\documentclass[11pt]{article}

\usepackage[utf8]{inputenc}
\usepackage{mathrsfs}
\usepackage{graphicx}
\usepackage{amsmath, mathtools}
\usepackage{amsthm}
\usepackage{amssymb,fourier}
\usepackage{microtype}
\usepackage{xparse}
\usepackage{xcolor}
\usepackage{enumitem}
\usepackage[colorlinks=true,citecolor=blue]{hyperref}
\usepackage[
    backend=biber,
    style=alphabetic,
    sortlocale=auto,
    natbib=true,
    url=false, 
    doi=true,
    eprint=true,
    safeinputenc=true,
    giveninits=true,
    maxbibnames=10
]{biblatex}
\AtBeginBibliography{\raggedright}
    \newtheorem{theorem}{Theorem}
    \newtheorem{lemma}[theorem]{Lemma}
    \newtheorem{proposition}[theorem]{Proposition}
    \newtheorem{corollary}[theorem]{Corollary}

    \newtheorem{assumption}[theorem]{Assumption}

\theoremstyle{definition} % For roman text in the body
    
    \newtheorem{definition}[theorem]{Definition}
    \newtheorem{remark}[theorem]{Remark}

\newcommand{\tr}{\operatorname{Tr}}

\renewcommand{\P}[1]{\mathds{P}\left[#1\right]}
\newcommand{\E}[1]{\mathds{E}\left[#1\right]}

\renewcommand{\i}{\mathbf{i}}

\newcommand{\1}{\mathbf{1}}
\newcommand{\A}{\mathcal{A}}

\renewcommand{\d}{\mathrm{d}}
\newcommand{\C}{\mathbb{C}}

\renewcommand{\E}{\mathbb{E}}
\newcommand{\D}{\mathbb{D}}
\newcommand{\F}{\mathscr{F}}

\newcommand{\N}{\mathbb{N}}
\renewcommand{\O}{\mathcal{O}}
\renewcommand{\P}{\mathbb{P}}
\newcommand{\R}{\mathbb{R}}

\renewcommand{\u}{\mathbf{U}}

\newcommand{\z}{\mathfrak{z}}
\newcommand{\Z}{\mathbb{Z}}

\newcommand{\Id}{\operatorname{Id}}

\numberwithin{equation}{section}
\numberwithin{theorem}{section}
\newcommand{\Exp}{\mathbb{E}}
\newcommand{\prob}{\mathbb{P}}

\renewcommand{\Pr}{\prob}
\DeclareDocumentCommand \one { o }
{%
\IfNoValueTF {#1}
{\mathbf{1}  }
{\mathbf{1}\left\{ {#1} \right\} }%
}

\newcommand{\lawequals}{\overset{\mathrm{law}}{=}}
\DeclareDocumentCommand{\Prto} {o} {
\IfNoValueTF {#1}
 {\overset{\Pr}{\longrightarrow}}
 { \xrightarrow[ #1 \to \infty]{\Pr }}
}
\DeclareDocumentCommand{\Asto} {o} {
\IfNoValueTF {#1}
 {\overset{\operatorname{a.s.}}{\longrightarrow}}
 {
 \xrightarrow[ #1 \to \infty]{\operatorname{a.s.} }
 }
}
\DeclareDocumentCommand{\Mgfto} {o} {
\IfNoValueTF {#1}
{\overset{\operatorname{mgf}}{\longrightarrow}}
{ \xrightarrow[ #1 \to \infty]{\operatorname{mgf} }}
}

\DeclareDocumentCommand{\Wkto} {o} {
\IfNoValueTF {#1}
 {\overset{\rm law}{\longrightarrow}}
 { \xrightarrow[ #1 \to \infty]{\operatorname{law}}}
}

\DeclareDocumentCommand \LPto { O{1} }
{\overset{\operatorname{\LP^{#1}}}{\longrightarrow}}

\newcommand{\Ai}{\operatorname{Ai}}
\newcommand{\Bi}{\operatorname{Bi}}
\newcommand{\SAi}{\operatorname{SAi}}

\newcommand{\Es}{\mathscr{U}}
\renewcommand{\i}{\mathbf{i}}

\renewcommand{\O}{\mathcal{O}}
\newcommand{\Q}{\mathbb{Q}}
\newcommand{\T}{\mathscr{T}}
\renewcommand{\u}{\mathrm{u}}
\renewcommand{\v}{\mathrm{v}}
\newcommand{\W}{\mathcal{W}}

\renewcommand\o{
  \mathchoice
    {{\scriptstyle\mathcal{O}}}% \displaystyle
    {{\scriptstyle\mathcal{O}}}% \textstyle
    {{\scriptscriptstyle\mathcal{O}}}% \scriptstyle
    {\scalebox{.7}{$\scriptscriptstyle\mathcal{O}$}}%\scriptscriptstyle
  }

\title{Strong approximation of Gaussian $\beta$ ensemble characteristic polynomials: the edge regime and the stochastic Airy function}
\date{\today}
\author{Gaultier Lambert\footnote{
KTH Royal Institute of Technology, Dept. of Mathematics, SE-100 44 Stockholm, Sweden.
\newline Email: \href{mailto:glambert@kth.se}{\nolinkurl{glambert@kth.se}}} 
~and Elliot Paquette\footnote{
 McGill University, 805 Sherbrooke Street West, Montreal, Quebec, Canada.
 \newline Email: \href{mailto:elliot.paquette@mcgill.ca}{\nolinkurl{elliot.paquette@mcgill.ca}}
 }}

\begin{document}

\maketitle

\begin{abstract}
For any  $\beta>0$, we investigate the properties of the stochastic Airy$_\beta$ operator, a random Schr\"odinger operator whose spectrum characterizes the soft-edge scaling limit of the $\beta$-ensemble. 
Our main results provide a spectral decomposition of this operator in terms of a random special function, which we call the \emph{stochastic Airy function} $\SAi$, and stability results for solutions of the stochastic Airy equation.
$\SAi$ can be viewed as the stochastic counterpart of the classical Airy function, and our stability results imply that the stochastic Airy function should describe the generic behavior of randomly perturbed second-order systems at a turning point.
This is illustrated by considering the characteristic polynomial of the Gaussian $\beta$-ensemble at the edge. 
Using the tridiagonal model, we build a coupling in which the stochastic Airy function is the almost sure scaling limit of the characteristic polynomial. In addition, the rescaling only involves Gaussian random variables.  
%  for general through its transfer matrix recurrence.
%We show that the rescaled characteristic polynomial converges to a random entire function in a neighborhood of the edge of the limiting spectrum.
%This random entire function, called the stochastic Airy function, is the unique (up to scaling) $L^2$ solution to the stochastic Airy equation, a family of second order stochastic differential equations. 
%Moreover, we obtain a coupling between $\varphi_N$ and a solution of the stochastic Airy equation which allows us to show that for any $\epsilon>0$, these two function are uniformly close by $N^{-1/6 + \epsilon}$ with overwhelming probability.
%These results build on the results of \cite{LambertPaquette02} in which the hyperbolic portion of the transfer matrix recurrence for $\varphi_N$ is analyzed.  %The results here show the parabolic portion is well--approximated by a solution of the stochastic Airy equation.  
\end{abstract}

%\paragraph{Keywords:} Stochastic Airy equation, Stochastic Airy equation, Gaussian $\beta$-ensembles, Airy$_\beta$ point process, characteristic polynomial. random transfer matrices recurrence, 

{
  \hypersetup{linkcolor=black}
\tableofcontents
}

% ===== begin sections/01-introduction.tex =====
\section{Introduction}

\subsection{Main results}

\paragraph{Stochastic Airy operator.}
The edge scaling limit of the  \emph{Gaussian $\beta$-ensemble}, or G$\beta$E, is given by a one-dimensional random Schr\"odinger operator depending on $\beta>0$ which is  called the \emph{stochastic Airy operator}, formally written as 
$\mathbf{H}_\beta=-\partial_{t}^{2} + t+ \tfrac{2}{\sqrt{\beta}}\d B_t $ on $\R_+$ where $\{B_t\}$ is a standard Brownian motion. 
 This operator was proposed in \cite{EdelmanSutton1} by viewing the tridiagonal random matrix model for G$\beta$E as a finite-difference scheme perturbed by independent noise, and the convergence was rigorously established in \cite{RRV09}. In particular, these works imply that the negative of the spectrum of the stochastic Airy operator, known as the Airy$_\beta$ point process, describes the local edge behavior of the G$\beta$E. Moreover, this limit is universal: it describes the soft-edge behavior of large classes of $\beta$-ensembles \cite{BEY14, KRV}, as well as that of generalized Wigner matrices in the orthogonal and unitary symmetry classes \cite{BEY14}.
 
 \medskip

The goal of this paper is to study the properties of solutions of the stochastic Airy equation and to relate its special solution to the characteristic polynomial of the G$\beta$E. 
Let $\{B(t) : t \in \R\}$ be a two-sided standard Brownian motion such that $B(0)=0$, so that
$\E[B^2(t)] =  |t|$ for all $t\in\R$.
The stochastic Airy equation is a system of stochastic differential equations, indexed by a parameter $\lambda\in\C$:
\begin{equation} \label{SA1}
\d\phi_\lambda'(t) = (t+\lambda) \phi_\lambda(t)\d t +\tfrac{2}{\sqrt{\beta}}\phi_\lambda \d B(t).
\end{equation} 
Given deterministic initial data $\left(c_{1}=\phi_\lambda(T), c_2=\phi_\lambda'(T) \right) \in \C^2$ for $T\in\R$, this equation can be posed as an It\^o differential equation for $\phi_\lambda' =\partial_t\phi_\lambda$ both for $\left\{ t > T \right\}$ and for $\left\{ t < T \right\}$,  in which case solutions are suitably adapted to the canonical filtration of $B$ pinned at $T$. 
Moreover, integrating by parts, it is possible to make sense of \eqref{SA1} for almost all realizations of $B$ and for random initial data. 
We define the (continuous) kernel
\begin{equation} \label{kernelU}
\Es_\lambda(t,u) = \frac{t^2-u^2}{2} +\tfrac{2}{\sqrt{\beta}}\big(B(t) - B(u)\big) + \lambda(t-u)  , \qquad \lambda \in \C, u, t \in\R.
\end{equation}
The equation \eqref{SA1}
%when viewed as a single-variable equation in $t$, 
can be written in terms of this kernel as 
$\d \phi_\lambda' = \phi_\lambda \d\Es_\lambda$ and if $\Phi_\lambda=\phi_\lambda'$ solves \eqref{SA1}, then it also satisfies the integral equation: 
\begin{equation} \label{SA2}
\begin{cases}
\Phi_\lambda(t) = c_2  { +} \int_T^t \Es_\lambda(t,u) \Phi_\lambda(u) \d u { +} c_1 \Es_\lambda(t,T) ,  \\
 \phi_\lambda(t)= c_1+ \int_T^t \Phi_\lambda(u) \d u , 
\end{cases}
 \qquad \lambda \in \C, t , T\in \R . 
\end{equation} 
We use this equation as the precise meaning of the \emph{stochastic Airy equation}.
This viewpoint is considered in \cite[Theorem 1]{Minami} to establish that the stochastic Airy operator $\mathbf{H}_\beta$ on $\R_+$, with Dirichlet boundary conditions, has a unique self-adjoint (Friedrichs) extension, the eigenvalue equation being interpreted as in \eqref{SA2} with $c_1=0$ at $T=0$.
In contrast to \eqref{SA1}, this approach does not require It\^o calculus, and solutions are deterministic functionals of the Brownian path.
In particular, one can take $c_1, c_2 \in\C$ to be random variables in \eqref{SA2}. 
Conversely, we show in Section \ref{sec:SAE} that any solution to \eqref{SA2} for deterministic $c_1$ and $c_2$ solves \eqref{SA1} in the It\^o sense with initial data $\left\{ \phi_\lambda(T)=c_1, \phi'_\lambda(T) = c_2 \right\}$. 

\medskip

From an analytic viewpoint, the equation \eqref{SA1} can be viewed as the usual Airy equation\footnote{The Airy equation typically shows up in critical wave phenomena at a turning point and it describes a transition between oscillatory and exponential behavior.
In particular, for $\beta=\infty$, \eqref{SA1} recovers the (classical) Airy equation $\phi_\lambda''(t) = (t+\lambda) \phi_\lambda(t)$ for $t\in\R$. This spectral equation has a unique (up to a multiplicative constant) solution in $L^{2}(\R_{+})$ given by $t\mapsto \Ai(t+\lambda)$ for $\lambda\in\C$ where $\Ai$ is the \emph{Airy function}. In addition, the asymptotics of the Airy function as $t\to\infty$ should be compared to \eqref{eq:SAIconvention} with $\beta=\infty, \lambda=0$: 
\[
\Ai(t) \simeq  t^{-1/4} \exp\big(-\tfrac{2}{3} t^{3/2} \big)/ \sqrt{4\pi} \qquad\text{as $t\to\infty$}. 
\]
} perturbed by a multiplicative white noise. It is established in  \cite{RRV09} that just like the Airy equation, \eqref{SA1} is \emph{limit-point} type at $+\infty$, meaning that $\Pr$-almost surely, there is a unique (up to a multiplicative constant) solution $\SAi_\lambda(t)$ which is in $L^2(\R_+)$. In fact, any other solution diverges like $e^{\frac 23 t^{3/2}}$ as $t \to \infty$.
This special solution is constructed in Section~\ref{sec:SAi} and we call it the \emph{stochastic Airy function}, denoted $\SAi$. 
We will show that it is an entire function of $\lambda\in\C$ which has the following almost sure asymptotics (see Proposition~\ref{prop:asympDN}) as $t\to\infty$: 
\begin{equation}
    \SAi_\lambda(t)
    \simeq
    t^{-\tfrac{1}{4}\bigl(1 + \tfrac 2 \beta\bigr)}
    \exp\biggl(
    -\tfrac{2}{3} (t + \lambda)^{3/2}
    -\tfrac{2}{\sqrt{\beta}}\int_0^t \mathfrak{X}(u)\,\d u + { \tfrac{2c_*}{\beta} }
    \biggr)/ \sqrt{4\pi},
%    \quad
%    \text{as }
%    t\to \infty\ \As
  \label{eq:SAIconvention}
\end{equation}
where $c_*\in\R$ is a (deterministic) constant and $\mathfrak{X}$ is the  continuous Gaussian process
\begin{equation}
  \mathfrak{X}(u) :=
  \int_0^u e^{\frac{4}{3}( t^{3/2} - u^{3/2})} \d B(t)
  \label{eq:SAIG}.
\end{equation}
In particular, the normalization is fixed such that (formally) $\lim_{\beta\to\infty}     \SAi_\lambda(t) =\ \Ai(t+\lambda)$ --- see Remark~\ref{rk:normalization}.
%We expect that $\SAi_\lambda$ is also an analytic function of the parameter $\gamma = \sqrt{2/\beta}$; see the discussion in Section~\ref{sec:discussion}.
Our goal in this paper is to develop the properties of the \emph{stochastic Airy function} and to show that it describes the edge-asymptotics of the Gaussian $\beta$-ensemble characteristic polynomial. 

\paragraph{Properties of the  stochastic Airy function.} First, in terms of the stochastic Airy function, we give an explicit diagonalization of the 
stochastic Airy operator:  
\begin{theorem}\label{thm:SAO}
Let $\mathcal{A} = \left\{ \lambda \in \C : \SAi_\lambda(0)=0 \right\} $ be the zero set of the stochastic Airy function. The spectrum of $\mathbf{H}_\beta$ is $-\mathcal{A}$, and for every $\lambda\in\mathcal{A}$, an associated eigenfunction with eigenvalue $-\lambda$ is $\SAi_\lambda$. In particular,  $\mathcal{A} \subset \R$ has the law of the Airy$_\beta$ point process.
\end{theorem}
\noindent This spectral theorem for the stochastic Airy operator is proved in Section~\ref{sec:SAi} (see Proposition~\ref{thm:spec}). We also provide several equivalent descriptions of the counting function for the random set $\mathcal{A}$ in Proposition~\ref{prop:cf}.
In contrast, in \cite{RRV09}, the spectrum of $\mathbf{H}_\beta$ is defined in terms of the quadratic form by the min-max principle, see \cite[Section 4.5]{AGZ}, and there is a  {\emph Sturm--Liouville} description of the eigenvalue counting function in terms of an It\^o process called the \emph{Riccati diffusion} (see Section~\ref{sec:SAR}). 

\medskip 

Second, as the underlying white noise $t\mapsto dB(t)$ enjoys a shift-invariance in law, it follows that the equation \eqref{SA1}  has a similar shift-invariance. In particular the stochastic Airy function has the following invariance property: 
Here and below, $\lawequals$ denotes equality in law.
\begin{proposition}\label{prop:SAIshift}
  For any $\varsigma \in \R,$
\begin{equation*}
  ( \SAi_\lambda(t) : \lambda \in \C , t \in \R)
  \lawequals
  ( \SAi_{\lambda-\varsigma}(t+\varsigma) : \lambda \in \C , t \in \R).
\end{equation*}
\end{proposition}

%{
%\noindent From this point, we can conclude from \eqref{eq:SAIconvention}
%that as the spectral parameter $\R \ni \lambda \to \infty$:
%\begin{equation} \label{SAasymp}
%  \SAi_\lambda(t)
%  \lawequals
%  (\lambda+t)^{-\tfrac{1}{4}\biggl(1+ \tfrac 2 \beta\biggr)}
%    \exp\biggl(
%    -\frac{2}{3} (t + \lambda)^{3/2}
%    -\int_0^{t+\lambda} \mathfrak{X}(u)\,\d u + { \frac{2c_*}{\beta}}  
%    \biggr) \bigg(\frac{1}{\sqrt{4\pi}}+\o_{\lambda,t,\beta}(1) \bigg),
%\end{equation}
%where the error is a random variable that converges uniformly on compact sets of $t \in\R$ to 0 in probability. 
%}

Third, we obtain high-probability estimates complementing the almost-sure asymptotics \eqref{eq:SAIconvention}.

\begin{theorem}  \label{thm:asymp}
Let $K\subset \C$ be a compact.
For any  $\epsilon \in (0, \tfrac 12],$  there is a constant $C = C_{K,\epsilon,\beta}$ so that for $\ell\in\{0,1\}$ and $T\ge 1$:
\[
\P \left( \sup_{\lambda \in K, t\ge T} \bigg|
\frac{\partial_t^\ell  \SAi_\lambda(t)\exp\big( \frac{2}{\sqrt{\beta}}\int_{0}^{t} \mathfrak{X}(u)\d u \big)}{\partial_t^\ell \Ai(t+\lambda)\E \exp\big( \frac{2}{\sqrt{\beta}}\int_{0}^{t} \mathfrak{X}(u)\d u \big)} -1 \bigg|  \ge  T^{\frac\epsilon3-\frac16} 
\right) \le  C  e^{- T^{1+\epsilon}}  
\]
and equivalently
\[
\P \left(    \sup_{\lambda \in K, t\ge T}  \bigg|    t^{\frac{1}{4}( (-1)^\ell + \frac 2 \beta)} \exp\biggl(
    \tfrac{2}{3} (t + \lambda)^{3/2} +
    \frac{2}{\sqrt{\beta}}\int_0^t \mathfrak{X}(u)\,\d u - { \frac{2c_*}{\beta} }
    \biggr) (-\partial_t)^\ell \SAi_\lambda(t) -\frac{1}{\sqrt{4\pi}}  \bigg| \geq  T^{-\epsilon} \right) \le  C  e^{- T^{3/2-3\epsilon}}  .
 \]
\end{theorem}

\paragraph{Stability  of the  stochastic Airy equation.}
In order to relate the stochastic Airy function to the characteristic polynomial of the Gaussian $\beta$-ensemble, we will need to establish that solutions of \eqref{SA2} are stable under perturbation of the initial data.

\begin{proposition}[Stability of the initial value problem] \label{prop:SAiwins}
  Fix $1<\nu<3/2$ and a compact interval $K\subset\R$, and set $\epsilon=\frac12-\frac\nu3$. There is a constant $c=c(\beta,\nu)>0$ such that the following holds for every $T\ge 1$. Define
  \[
  K' = \bigl\{\lambda\in\C : \Re\lambda\in K,\ |\Im\lambda|\le e^{-(1-1/\nu)T^\nu}\bigr\}.
  \]
  Let $\phi_\lambda$ be the solution of \eqref{SA2} with initial data
  \(
 \big( \phi_\lambda(T),  \phi'_\lambda(T)\big) = (c_{1,\lambda},  c_{2,\lambda} )
  \)
such that $|c_{1,\lambda}|+ |c_{2,\lambda}|/\sqrt{T}  \leq 1$ for $\lambda\in K'$.  
There is another solution $\chi_\lambda$ so that 
\[
    \phi_\lambda(t)
    =
    \Theta_\lambda \SAi_\lambda(t)
    +
    \chi_\lambda(t)
    \qquad
    \text{for all $t \in \R$}
\]
and, with probability at least $1 - e^{-cT^{\nu}}$,
\[
    \sup_{\lambda \in K' , t \in [-e^{T},T/2]}
      \bigl(
      |\chi_\lambda(t)|
      +|\chi'_\lambda(t)|
      \bigr)
 \le e^{-T^{3/2}/4} 
\qquad\text{and}\qquad 
  \Theta_\lambda
=  \frac{c_{1,\lambda}-c_{2,\lambda}/\sqrt{T} + \O(T^{-\epsilon})}{2\SAi_\lambda(T)\big(1+\O(T^{-\epsilon})\big)}
\]
where the $\O(T^{-\epsilon})$ terms are uniform in $\lambda\in K'$. 
Moreover, if the maps $\lambda\mapsto c_{i,\lambda}$ extend analytically to a neighborhood of $K'$ for $i=1,2$, then so do the random maps $\lambda\mapsto \Theta_\lambda$ and $\lambda\mapsto \chi_\lambda(t)$.
\end{proposition}

\noindent  We expect that the same conclusion holds uniformly on any fixed compact subset of $\C$, without the restriction on the imaginary part imposed in the definition of $K'$; that restriction is a technical artifact of our proof. The proof of Proposition~\ref{prop:SAiwins}, which is the technical heart of our results, occupies a substantial part of this paper (Sections~\ref{sec:SAK}--\ref{sec:Red5}). It relies on precise estimates for the kernel of the equation \eqref{SA2}  
(see Theorem~\ref{thm:Kest}). These kernel estimates are used to control both the stability of solutions under  perturbation of the initial data (Proposition~\ref{prop:SAiwins}) and under perturbation of the stochastic Airy equation (Proposition~\ref{prop:stability}).
They are derived by studying the fine asymptotics of the Riccati diffusion associated with the SDE \eqref{SA1}. 
We will provide a more extensive discussion of the properties of $\SAi$ in Section \ref{sec:discussion} and we now move on to its relationship to Gaussian $\beta$-ensembles.

\paragraph{G$\beta$E characteristic polynomial.}
The $N$--dimensional  \emph{Gaussian $\beta$-ensemble} is the $N$--point process on $\R$ with a joint density function
\begin{equation}\label{eq:GbE}
%d\mu_{G\beta E}(\lambda_1, \lambda_2, \dots,\lambda_N) =
{(\lambda_1, \lambda_2, \dots,\lambda_N) \mapsto}
\frac{1}{\mathcal{Z}_{N,\beta}} e^{- \sum_{i=1}^N \beta N \lambda_i^2} \prod_{i > j} \left|\lambda_i - \lambda_j\right|^\beta
\end{equation}
where $\mathcal{Z}_{N,\beta}>0$ is a normalizing constant.
When $\beta=1,2,4$, this is the law of the eigenvalues of the classical Gaussian invariant ensembles \cite[Chap.~3, 4]{AGZ}. 
In terms of these random points, we define the \emph{characteristic polynomial}
\[
\varphi_N(z) := {\textstyle \prod_{i=1}^N} ( z - \lambda_i)  , \qquad z\in\C. 
\]
Using the properties of the stochastic Airy equation developed in the previous paragraphs, our goal is to describe the \emph{edge scaling limit} of this characteristic polynomial. With the scaling \eqref{eq:GbE}, the empirical measure  converges to the semicircle law on $[-1,1]$ and we recall that the eigenvalue spacing at the edge is of order $N^{-2/3}$. We also introduce the \emph{weight}
\begin{equation}   \label{eq:wn}
  w_n(z)^{-1} \coloneqq (2\pi)^{1/4}e^{N z^2} 2^{-n} (Nz^2)^{-1/12}\sqrt{\tfrac{n!}{N^n}}  ,\qquad z\in\R\setminus0. 
\end{equation}
Our main result shows that the scaling limit of $\varphi_N$ around $z=1$ is given in terms of the stochastic Airy function by $\lambda \mapsto \SAi_\lambda(0)$. 

\begin{theorem}\label{thm:edge1}
Define $\Psi_N(\lambda) \coloneqq (w_N\varphi_N)\left(1+\frac{\lambda}{2 N^{2/3}}\right)$ for $\lambda\in\R$.
  There is a centered Gaussian variable $\mathfrak{G}_N$ with 
  \[
    \Exp \mathfrak{G}^2_N = \frac{2}{3\beta}\log N + \O_\beta(1)
    \quad\text{as}\quad N\to \infty
  \]
  such that as a random real-analytic function under the topology of locally uniform convergence of the function and all its derivatives:
  \[
    %{\exp( \sqrt{\tfrac{2}{\beta}}\mathfrak{G}_N)}
    \biggl(\Psi_N(\lambda)
    \tfrac
    {\Exp \exp(\mathfrak{G}_N)}
    {\exp(\mathfrak{G}_N)}
    ~:~\lambda \in \R\biggr)
    \Wkto[N]
    \bigl(\SAi_\lambda(0) : \lambda \in \R\bigr).
  \]
\end{theorem}

Hence, this paper provides a new characterization of the Airy$_\beta$ point process as the zero set 
of $\lambda\mapsto \SAi_{\lambda}(0)$, and  with the strong mode of convergence \eqref{cvgas} below, a new proof that it describes the (almost sure) limit of the G$\beta$E edge  eigenvalues.
The proof of  Theorem~\ref{thm:edge1}  relies on the tridiagonal representation of the G$\beta$E from  \cite{DumitriuEdelman} which provides a probabilistic coupling and a recursion for $(\varphi_N ;N\in\N)$ to study its asymptotics as $N\to\infty$.
This framework is reviewed in detail in Section~\ref{sect:TM}, and we give a refined version of Theorem~\ref{thm:edge1} which holds with overwhelming probability within this coupling; see Theorem~\ref{main:thm}.

As a direct consequence of Theorem~\ref{thm:edge1} (the independence follows from the more precise asymptotics given in Theorem~\ref{main:thm}), we obtain a central limit theorem (CLT) for the logarithm of the characteristic polynomial at the edge:
\begin{corollary}\label{cor:clt}
  For any $\lambda \in \R$,  
  \[
    \frac{\log|\Psi_N(\lambda)|
      +\frac{1}{3\beta}\log N}
      {\sqrt{\frac{2}{3\beta}\log N}}
    \Wkto[N]
    \mathcal{N}(0,1).
  \]
  The convergence holds jointly for finite-dimensional marginals of $\lambda \in \R$, and the limit is the same for all $\lambda\in\R.$  Moreover, the limit is independent of the stochastic Airy function in Theorem~\ref{thm:edge1}. 
\end{corollary}

There are several analogous CLTs for the log-characteristic polynomials of different random matrix ensembles in the bulk.
Let us mention \cite{KeatingSnaith, BHNY} for random unitary matrices and
the recent work \cite{Augeri} for Gaussian $\beta$-ensembles. 
There are also asymptotics for the determinant (at 0) of the Gaussian $\beta$--ensembles \cite{Duy17} and for general Wigner matrices $(\beta=1,2)$ in \cite{BM19}. 
Since the original version, an alternative proof of the edge CLT of Corollary~\ref{cor:clt} can be found in \cite{Johnstone}.
It has also been proved that this edge CLT is universal for one-cut regular $\beta$-ensembles in  \cite{BourgadeModyPain}. 
This work relies on optimal local laws for $\beta$-ensembles to establish the \emph{log-correlated structure} of the log-characteristic polynomial of $\beta$-ensembles. 

\subsection{\texorpdfstring{Gaussian $\beta$--ensembles}{Gaussian beta-ensembles}} \label{sect:TM}

\paragraph{Random tridiagonal matrices \& Transfer matrix recurrence.}
To study the asymptotics of the G$\beta$E characteristic polynomial, our starting point is the Dumitriu--Edelman tridiagonal matrix model for the Gaussian $\beta$--ensemble.  Recall that for any $\alpha>0$, a $\chi_{\alpha}$ random variable has density proportional to
$x^{\alpha -1} e^{- x^2 /2} \1_{x>0}$. 
% and $\chi_{\alpha}^2 \sim \Gamma(\frac\alpha2, 2)$ where $\Gamma(\frac\alpha2, 2)$ denotes a Gamma distribution with shape $\frac\alpha2$ and rate $\frac 12$.  
In terms of these variables, we define the semi-infinite tridiagonal matrix 
\begin{equation} \label{def:trimatrix}
  \mathbf{A} =
  \begingroup
  \setlength\arraycolsep{5pt}
  \begin{bmatrix}
    b_1 & a_1 & &\\[6pt]
    a_1 & b_2 & a_2 & \\[1pt]
    & a_2 & b_3 & \ddots  \\
    && \ddots & \ddots
  \end{bmatrix}
  \endgroup
\end{equation}
where $b_i \sim \mathcal{N}(0,2)$ and  $a_i \sim \chi_{\beta i}$ are independent random variables. 
By \cite{DumitriuEdelman}, the eigenvalues of the principal $N \times N$ minor of  the random matrix 
$ \mathbf{A}/\sqrt{4N\beta}$  have law \eqref{eq:GbE}.
In particular, we can realize the characteristic polynomial $\varphi_N(z) = \det[z-({4}{N\beta})^{-1/2}\mathbf{A}]_{N,N}$. 
%which is scaled so that its limiting density of states is given by the semicircle law on $[-1,1]$. 

\medskip

For each ambient scale $N$, we let
\begin{equation} \label{def:varphi}
\varphi_n^{(N)}(z)= \det[z-({4}{N\beta})^{-1/2}\mathbf{A}]_{n,n}
\qquad\text{ for $n\in\N$, $z\in\C$.}
\end{equation}
Since $N$ is fixed throughout the recurrence analysis, we suppress this dependence and write $\varphi_n=\varphi_n^{(N)}$.  In particular, $\varphi_N$ is the characteristic polynomial of the $N\times N$ matrix above.
By cofactor expanding the $n$--th column/row of this determinant, we are led to the following recurrence for any integer $n \geq 2$:
\begin{equation*}
  \begingroup
  \setlength\arraycolsep{5pt}
  \begin{bmatrix}
    \\[-0.85em]
     \varphi_{n}(z) \\[0.50em] \varphi_{n-1}(z)\\[0.35em]
  \end{bmatrix}
  =
  \begin{bmatrix}
      z-\frac{b_n}{{2}\sqrt{N\beta}} & -\frac{ a_{n-1}^2}{{4}N\beta}  \\
    1 & 0 
  \end{bmatrix}
  \begin{bmatrix}
    \\[-0.85em]
     \varphi_{n-1}(z) \\[0.50em] \varphi_{n-2}(z)\\[0.35em]
  \end{bmatrix}
  \eqqcolon T_{n-1}(z) 
  \begin{bmatrix}
    \\[-0.85em]
     \varphi_{n-1}(z) \\[0.50em] \varphi_{n-2}(z)\\[0.35em]
  \end{bmatrix},
  \endgroup
\end{equation*}
where by convention $\varphi_0 = 1$ and  
$ \varphi_1(z) = z- \frac{b_1}{ 2\sqrt{N\beta}}$. 
This shows that for any $n\ge 1$,
\begin{equation} \label{eq:recurrence}
  \begin{bmatrix}
     \varphi_{n}(z) \\[0.35em] \varphi_{n-1}(z)\\[0.35em]
  \end{bmatrix}
= T_{n-1}(z)  \cdots T_1(z)  \begin{bmatrix}   z- \frac{b_1}{ 2\sqrt{N\beta}}   \\ 1\end{bmatrix} . 
\end{equation}

In particular, the noise which naturally appears in the random matrices $ T_n(z)$ are independent and we denote for all $n \geq 1,$
\begin{equation} \label{def:XY}
X_n = \frac{b_{n+1}}{\sqrt{2}} 
\qquad\text{and}\qquad 
Y_n = \frac{a_{n}^2 - \beta n}{\sqrt{2\beta n}}. 
\end{equation}
These are independent, mean $0$ and variance $1$ random variables which,  after a mild truncation, are sub-Gaussian for $n\gg1$.
Our analysis applies to the characteristic polynomial of any random Jacobi with these properties and our results are \emph{universal}, meaning that they do not depend on the specific distributions of $(X_n,Y_n)_{n\ge 1}$ under the following assumptions: 

\begin{assumption} \label{def:noise}
We assume that the entries of the tridiagonal matrix model $\mathbf{A}$ are independent random variables and there are constant $ c, C>0$  such that \eqref{def:XY} satisfy the conditions: 
\begin{equation}\label{norm}
\E X_n =0 , \qquad \E X_n^2 =  1, \qquad 
\VERT X_n \VERT_1 \le C  ,\qquad  
\VERT X_n\1\{|X_n| < c \sqrt{n}\}\VERT_2 \le C, \qquad 
\text{for every $n\in\N$},  
\end{equation}
and similarly for $\{Y_n\}_{n\in\N}$. 
In this paper, $\VERT\cdot\VERT_q$ refers to the Orlicz norm defined in Section~\ref{sec:concentration}. 
In particular, these assumptions suffice to apply the results from \cite{LambertPaquette02}, and they are satisfied by the G$\beta$E tridiagonal models; see \cite[Appendix~B]{LambertPaquette02}. 
In the sequel, all constants are allowed to depend on the fixed parameters $\beta,c, C>0$.
We introduce $\sigma$-algebras $\F_0 \coloneqq \sigma(b_1)$ and
\(
\F_{n} \coloneqq {\sigma\bigl(b_1, X_k,Y_k:k\le n \bigr)}. 
\)
\end{assumption}

\paragraph{Approximate stochastic Airy equation.}
Let $w_n$ be as in \eqref{eq:wn} and define for $n\ge 1$,
\[
{\Psi}_n(\lambda) \coloneqq (w_n \varphi_n)( 1 + \tfrac{\lambda}{2N^{2/3}}) .
\]
We claim that at the edge $z= 1 + \frac{\lambda}{2N^{2/3}}$, in terms of the random variables \eqref{def:XY}, we obtain the approximation for $N\gg1$,
\begin{equation}\label{eq:recurrence0}
\Psi_n- 2\Psi_{n-1} + \Psi_{n-2} 
\approx
\biggl(\frac{N-n}{N}
+\frac{\lambda}{N^{2/3}}
-\frac{X_{n-1} \sqrt{2}}{\sqrt{\beta N}}
\biggr)
\Psi_{n-1}
-
\biggl(
\frac{Y_{n-1} \sqrt{2}}{\sqrt{\beta N}}
\biggr)
\Psi_{n-2}, 
\end{equation}
see \eqref{eq:recurrence2} for the exact statement.  
Then, for a suitable mesh of $t \in\R$ in \eqref{eq:recurrence0}, this is a natural discretization of \eqref{SA1} driven by the two-sided (scaled) random walk:
\begin{equation} \label{Bedge}
\tilde{B}_{t} =  -\sqrt{\frac{1}{2N^{1/3}}} \sum_{k= N- \lfloor t N^{1/3}\rfloor}^{N-1} (X_k +Y_k), 
  \quad
  t \in \R.
\end{equation}
% Here and below, a sum whose lower limit exceeds its upper limit is interpreted as an oriented sum, so that $\widetilde B_0=0$.
% One has $\tilde{B}_{t} \approx B_t$ where $\{B_t\}$ is a standard Brownian motion. For $n=N-t/\delta$, the recurrence formally reads
% \[
% % \frac{\Psi_n-2\Psi_{n-1}+\Psi_{n-2}}{\delta}
% \approx
% (t+\lambda)\Psi_{n-1}\delta
% +\tfrac{2}{\sqrt{\beta}}\Psi_{n-1}\bigl(\widetilde B_{t+\delta}-\widetilde B_t\bigr),
% 
% \]
% whose left-hand side is an increment of the discrete derivative; this is the discrete form of \eqref{SA1}.
To make this precise, we can build a probability space which supports the noise \eqref{def:XY} coupled to a Brownian motion $(B_t : t \in \R)$ with $\tilde{B}_t \approx B_t$ so that the G$\beta$E characteristic polynomial (rescaled at the edge) can be directly compared to the stochastic Airy function driven by $(B_t : t \in \R)$. 
A similar coupling has been introduced in our previous paper \cite{LambertPaquette02}.
Within this coupling, we can state a quantitative space-time refinement of Theorem~\ref{thm:edge1}. 

\paragraph{High probability approximation.}
Using this coupling, we now formulate a quantitative version of the stochastic Airy convergence at the edge. 

\begin{theorem}\label{main:thm} 
Let $K \subset \R$ be any compact set. 
There is a small $0<\epsilon < \frac16$ such that with $T(N) \coloneqq (\log N)^{1-\epsilon}$, 
it holds for  $t \in  [-e^{T},T/2] $ and $\lambda \in K$,
\[
  \Psi_{N- N^{1/3}t}(\lambda) = 
  \frac{\exp\big(\frac{2}{\sqrt{\beta}}\int_{0}^{T} \mathfrak{X}(u)\d u + \boldsymbol{g} + \varepsilon_{\lambda}^{(N)} \big)}
  {\E  \exp\big(\frac{2}{\sqrt{\beta}}\int_{0}^{T} \mathfrak{X}(u)\d u +  \boldsymbol{g} \big)}
  \biggl( \SAi_\lambda(t) + \Upsilon_{\lambda}^{(N)}(t)\biggr)
\]
where
\begin{itemize}[leftmargin=*]
\item $\mathfrak{X}$  is a Gaussian process as in  \eqref{eq:SAIG} driven by the same Brownian motion $B$ as $\SAi$, and 
$\boldsymbol{g}$ is a centered Gaussian which is independent of $B$ with variance $\frac2\beta \log\big(\frac{N^{1/3}}{2T^{1/2}}\big) +\o(1)$ as $N\to\infty$,
\item $\varepsilon_{\lambda}^{(N)}$ and  $\Upsilon_{\lambda}^{(N)}$ are random analytic functions. 
 These errors satisfy with probability at least $1- e^{-(\log N)^{1+\epsilon}},$ for any $\ell\in\N$ and uniformly for  $t \in  [-e^{T},T/2]$ and $\lambda \in K$, 
\[
\partial_\lambda^{\ell-1}\varepsilon_{\lambda}^{(N)} = \O\big((\log N)^{\epsilon-\frac16}\big)
\qquad\text{and}\qquad
\partial_\lambda^{\ell-1} \Upsilon_{\lambda}^{(N)}(t) =\O(N^{\ell\epsilon-1/6}) .
\]
\end{itemize} 
\end{theorem}

At step $N$, the spectral edges $z\in\{\pm1\}$ correspond to \emph{tunring points} for the recursion \eqref{eq:recurrence} where the behavior of the solutions transitions from an exponential behavior (\emph{hyperbolic regime}) to an oscillatory behavior (\emph{elliptic regime}). Our result shows that this transition is described by  $t\mapsto \SAi_\lambda(t)$ for $N\gg1$. 
We expect that this behavior is universal for randomly perturbed second-order systems at a turning point in the presence of critical noise.
We formulated Theorem~\ref{main:thm} for $K \subset \R$ for technical reasons. Using the results from \cite[Proposition 4.6]{Assiotis} (the entire function $\lambda\in\C \mapsto \SAi_\lambda(t)$ is in the Laguerre--P\'olya class for a fixed $t\in\R$), one can upgrade the convergence to uniform on compact sets $K\subset\C$.

 \medskip

The proof of Theorem~\ref{main:thm} is given in Section~\ref{sec:quant} and it relies on two main ingredients: 
\begin{itemize}[leftmargin=.5cm]
\item  The results from \cite{LambertPaquette02} to control the  \emph{hyperbolic part of the recursion} \eqref{eq:recurrence}, that is, the \emph{exponentially growing behavior of the characteristic polynomial} for steps $n\in[1, N- T(N)N^{1/3}]$. These results are reviewed in Theorem~\ref{thm:hypersimple} and are used to control the entrance behavior in equation \eqref{SA2} at time $T(N)$.
In particular, the Gaussian random variable $\boldsymbol{g}$ comes from the \emph{hyperbolic part}, which is why, with the appropriate coupling, it is independent of the \emph{parabolic objects} $(\SAi,\mathfrak{X})$. 
 \item Then, introducing a coupling between \eqref{Bedge} and a Brownian motion $\{B_t\}$, the rescaled characteristic polynomial $t \in \R \mapsto  \Psi_{N- N^{1/3}t}(\lambda)$ is an approximate (backward) solution of the equation \eqref{SA2} with a small random perturbation.  Using the stability result from Proposition~\ref{prop:SAiwins}, this allows us to compare the characteristic polynomial to the stochastic Airy function with an explicit prefactor and errors that converge to 0 with overwhelming probability as $N\to\infty$. 
\end{itemize}

Since $\E \big(\int_{0}^{T} \mathfrak{X}(u)\d u \big)^2 = \tfrac14( \log  T+ c_*) + \o(1) $ as $T\to\infty$,
Theorem~\ref{main:thm} with $t=0$ immediately implies Theorem~\ref{thm:edge1}, as well as the fact that the convergence holds with overwhelming probability:  for any compact set $K \subset \R \times \R$ and any $\ell \in \N_0$, as $N\to\infty$, 
\begin{equation} \label{cvgas}
  \max_{(t,\lambda) \in K}
  \biggl|
  \partial^\ell_{\lambda}\Big( \Psi_{N - tN^{1/3}}(\lambda)
  -\SAi_{\lambda}(t)
  \tfrac
  { \exp( \mathfrak{G}_N)}
  { \Exp \exp(\mathfrak{G}_N)}
  \Big)
  \biggr| \to 0.
\end{equation}
As for Corollary~\ref{cor:clt}, it follows from
\[
\E \biggl(\int_{0}^{T} \mathfrak{X}(u)\d u \biggr)^2 = \o(\log N),
\]
whereas
\[
\E \boldsymbol{g}^2
=\frac2\beta\log\biggl(\frac{N^{1/3}}{2T^{1/2}}\biggr)+\o(1)
=\frac{2}{3\beta}\log N+\O(\log\log N),
\]
and $\boldsymbol{g}$ is independent of $(\SAi,\mathfrak{X})$.

%but we expect that a similar result holds for $K \subset \C$ and may be proved using the method from this paper.

\begin{remark}
Let us observe that  $\varphi_n(-z)= (-1)^n \det([z-({4}{N\beta})^{-1/2}\widehat{\mathbf{A}}]_{n,n})$ for $n\in\N,$ where $\widehat{\mathbf{A}}$ is the random tridiagonal matrix like \eqref{def:trimatrix} with sequences $(a_k, -b_k)_{k\in\N}$. 
Hence, Theorem~\ref{main:thm} also describes the scaling limit of the characteristic polynomial near the other edge point at $-1$.  
Hence, we recover the fact that in the large $N$ limit, one observes two independent Airy$_\beta$ point processes around $\pm 1$. 
Moreover, it is apparent from our coupling that the two Brownian motions ${B}^{(\pm 1)}$ which drive the stochastic Airy functions $\SAi_{\lambda}^{(\pm1)}$ and Gaussian processes $\mathfrak{X}^{(\pm1)}$ near $\pm1$ are independent. It turns out that the Gaussian random variables $\boldsymbol{g}^{(\pm1)}$ also have small correlations; see Remark~\ref{rk:sym}. 
\end{remark}

\subsection{Parabolic regime} \label{sec:para}
The behavior of the recursion \eqref{eq:recurrence} 
can be explained by considering the expected characteristic polynomial. For any $\beta >0$ and $n\ge 0$,  $\widetilde\varphi_n(z)= \E\varphi_n(z)$ are monic Hermite polynomials orthogonal with respect to the measure $e^{-2N x^2}\d x$ on $\R$. This is a consequence of the fact that  
\begin{equation}\label{eq:hermite}
%\label{eq:Ttilde}
    \begin{bmatrix}
     \widetilde\varphi_{n}(z) \\[0.35em] \widetilde\varphi_{n-1}(z)\\[0.35em]
  \end{bmatrix}
  =
  \widetilde T_{n-1}
  \cdots
  \widetilde T_{1}
  \begin{bmatrix}
  z  \\[0.35em] 1\\[0.35em]
\end{bmatrix},
\qquad \text{where}
\quad \widetilde T_{k}  = \E T_k = 
\begin{bmatrix}
  z &- \frac{k}{{ 4}N} \\[0.35em]  1 & 0 \\[0.35em]
\end{bmatrix}
\quad 
\text{for all} 
\quad k \geq 1.
\end{equation}
Note that this also corresponds to the recurrence for $\beta=\infty$. 
Then, the deterministic normalization \eqref{eq:wn} is such that with $\widetilde\Psi_n = w_n\widetilde\varphi_n$, using the classical  \emph{Plancherel--Rotach} asymptotics  \cite{PR29}, as $N\to\infty$
\begin{equation}\label{eq:Airy}
\widetilde\Psi_{N-tN^{1/3}}(1+\tfrac{\lambda}{2N^{2/3}}) \to \Ai(\lambda + t)
\end{equation}
uniformly on compact sets of $\lambda \in \C$ and $t \in \R$, where $\Ai$ is the classical Airy function.
This is the deterministic equivalent of Theorem~\ref{main:thm}, where $\Ai(\lambda+\cdot) = \E \SAi_\lambda $ is the special solution of the ODE \eqref{SA1} with $\beta=\infty$. 
In fact, the asymptotics \eqref{eq:Airy} can be derived using ideas similar to those in Section~\ref{sec:fde}, by approximating the finite difference equations for $(\widetilde\Psi_n )$ coming from \eqref{eq:hermite} by the Airy ODE.
One of the main challenges of this paper is to control the effect of the \emph{white noise} and to develop the stability theory for solutions of the stochastic Airy equation \eqref{SA1}.

The recursions \eqref{eq:recurrence} and  \eqref{eq:hermite} have the same quantitative behavior with three different regimes (for a spectral parameter $z\in\R$):
\begin{itemize}[leftmargin=.5cm]
\item For $n \ll Nz^2$, the eigenvalues of $\widetilde T_{n}$  are real and separated, so that the solution has an \emph{exponentially growing behavior}. This is called the  \emph{hyperbolic regime}, and it has been studied in detail in our previous paper \cite{LambertPaquette02} using methods from \emph{perturbation theory}. 
\item For $n \gg Nz^2$, the eigenvalues of $\widetilde T_{n}$  are complex conjugate and the solution has an \emph{oscillatory behavior}. This is called the  \emph{elliptic regime}; it characterizes the bulk of the spectrum and is studied in our follow-up paper \cite{LambertPaquetteBulk}.
\item If $z \approx \pm\sqrt{n/N}$, there is a \emph{transition} where the transfer matrices are almost singular and the noise is most relevant. This is called the  \emph{parabolic regime} and it is the subject of this paper.
\end{itemize}

In the sequel, we denote by  $N_0 =N_0(z) \coloneqq \lfloor N z^2 \rfloor $ the \emph{turning point}  of the recursion and by $\mathfrak{L} =\mathfrak{L}(z) \coloneqq N_0^{1/3}$ the \emph{parabolic scale} in the transition window. 
 If $z\in[-1,1]$, in the spectrum, the turning point is directly relevant to study the characteristic polynomials $\big(\varphi_n(z); n\in[1,N]\big)$ and 
 $\varphi_N(z)$ is in the parabolic regime exactly at the edge $z=\pm1$. 
 In this case, Theorem~\ref{main:thm} shows that the \emph{transition from exponential to oscillatory  behavior} is described by the stochastic Airy function.
Moreover, in the bulk of the spectrum, for $z\in(-1,1)$, all three regimes are relevant to obtain the complete asymptotics of  $\varphi_N(z)$ as $N\to\infty$. 
In particular, the stochastic Airy function is also relevant to describe the complete bulk asymptotics as explained in \cite{LambertPaquetteBulk}. 

\medskip

To explain the \emph{parabolic scaling},
observe that according to \eqref{eq:wn}, by scaling, the Hermite polynomials satisfy for $z>0$ and $N, n\ge 1$:
\[
\widetilde\Psi_n(z) = w_n \widetilde\varphi_n(z) =  N_0^{1/12} \phi_n\big(2\sqrt{N}z\big)
=   N_0^{1/12}  \phi_n\big(2\sqrt{n}\sqrt{\tfrac{Nz^2}n}\big)
\]
where $(\phi_n; n\in\N_0)$ are the \emph{Hermite functions} (orthonormal basis of $L^2(\R)$ independent of $N\in\N$). 
Then, in the parabolic regime, one has the asymptotic approximation for $N_0 \gg 1$:
\[
\widetilde\Psi_{N_0-t \mathfrak{L}}\big(z\big(1 + \tfrac{\lambda}{2N_0^{2/3}}\big)\big)
=     N_0^{1/12}  \phi_n\big(2\sqrt{n}\sqrt{\tfrac{n+(t+\lambda) \mathfrak{L} + \lambda^2/4\mathfrak{L}}n}\big)
\approx   n^{1/12}   \phi_n\big( 2\sqrt{n}+\tfrac{t+\lambda}{n^{1/6}}\big) \approx \Ai(t+\lambda) .
\]

These asymptotics hold for any $z\in[-1,1]$, and the rescaling only depends on $N_0 = Nz^2$, instead of $N$. 
Similarly, for the randomly perturbed recursion \eqref{eq:recurrence}, at any point in the spectrum $z\in[-1,1]$ except for $z=0$,  the stochastic Airy function describes the behavior of the sequence of characteristic polynomials around the turning point and the approximation \eqref{eq:recurrence0}--\eqref{Bedge} holds replacing $N\leftarrow N_0(z)$ and the Brownian motion $\{B_t\} \leftarrow \{B_t^z\}$.
The precise approximation statement for $N_0\gg1$ is given in Theorem~\ref{main:approx} (Theorem~\ref{main:thm} follows from this result with $z=1$).

\subsection{Discussion and further developments}
\label{sec:discussion}

The first version of this work \cite{LambertPaquetteV1} included a list of open questions about the stochastic Airy function and its connection to G$\beta$E.  In the intervening years, several of these problems have been resolved.  We refer to \cite{LambertPaquetteV1} for the original formulations, and below we discuss the relevant developments alongside the surrounding context.
Recall that $\mathbf{H}_\beta=-\partial_{t}^{2} + t+ \tfrac{2}{\sqrt{\beta}}\d B_t $ denotes the stochastic Airy operator and we let $\{ \mathfrak{z}_n\}_{n\ge1}$ be the zero set of $\SAi_\lambda(0)$ ordered so that $\mathfrak{z}_1 > \mathfrak{z}_2 >\cdots$. In particular, this is a realization of the Airy$_\beta$ point process.

\paragraph{Related theory on the stochastic Airy$_\beta$ operator.}
We have constructed the stochastic Airy function, as the special solution to the stochastic Airy equation and established that  the zero set of $\lambda\mapsto \SAi_{\lambda}(0)$ gives a new characterization of the Airy$_\beta$ point process; see Proposition~\ref{prop:cf} below.
In particular, we obtain a new proof that  it describes the edge fluctuations of the G$\beta$E eigenvalues. 
Another point of view on the stochastic Airy operator is provided by Gorin and Shkolnikov \cite{GorinShkolnikov}, who found a pathwise representation for the \emph{stochastic Airy semigroup} $(e^{-t\mathbf{H}_\beta}; t\ge 0)$ using the Feynman--Kac formula.
The proof also relies on the tridiagonal representation of the G$\beta$E to study the Laplace functional of the edge eigenvalues.   
There are also similar representations  for the limiting point processes arising from \emph{spiked models} \cite{BloemendalViragII,LamarreShkolnikov, Lamarre}.
For a rank-1 spike, this point process can also be formulated as the negative of the spectrum of $\mathbf{H}_\beta$ but with Robin boundary conditions \cite{BloemendalViragI}.  In Section \ref{sec:SAi} we make a connection between the stochastic Airy function and the spectrum of $\mathbf{H}_\beta$  with general Robin type boundary conditions.  We do not know of a connection between the stochastic Airy semigroup and the stochastic Airy function, but believe it is an interesting research direction. 

The recent work of \cite{DumazLiValko} (see also \cite{AshburyBridgwood}) establishes that the inverse of $\mathbf{H}_\beta$ 
is a well-defined Hilbert--Schmidt operator and it is used to prove a hard-to-soft edge transition at the level of operators. 
Their analysis has similarities to our construction of the stochastic Airy kernel in Section \ref{sec:SAK}.

In \cite{GonzalezHolcomb}, the authors study the spectrum of the stochastic Airy operator $\mathbf{H}_\beta$ on $[t,\infty)$ with Dirichlet boundary conditions. Within our framework, the negative of this spectrum corresponds to the zero set $\{ \mathfrak{z}_n(t)\}_{n\ge1}$ of $\SAi_\lambda(t)$. They prove that  $\left\{ \mathfrak{z}_n(t) + t ; t\in\R \right\}_{n\ge1}$ is a stationary process of continuously-differentiable non-intersecting paths with the distribution of Airy$_\beta$. 
This is a consequence of Theorem \ref{thm:SAO} and the stationarity follows from the shift invariance in law from Proposition \ref{prop:SAIshift}.

For rational $\beta>0$, the Airy$_\beta$ point process can also be characterized as the unique solution of a hierarchy of loop equations \cite{bourgade2026loop}. These equations are derived from the explicit joint eigenvalue density of the Gaussian $\beta$-ensemble by repeated integration by parts, rather than from the Dumitriu--Edelman tridiagonal model. The main challenge is to prove that the hierarchy of equations for the Stieltjes transform of the Airy$_\beta$ point process has a unique solution.

\paragraph{Relations to KPZ theory and the Airy$_\beta$ line ensemble.}
There is an interesting connection between the Cole--Hopf solution of the KPZ equation (with special boundary conditions) and the Airy$_\beta$ point process (for $\beta=1,2$). This arises from the fundamental relationship between the Laplace functionals of the stochastic heat equation and the Airy$_\beta$ point processes \cite[Theorem 2.2]{BorodinGorin}.  
Several works rely on this connection to establish large deviations for the Cole--Hopf solution \cite{CorwinGhosal,Tsai}. %This has been extended in \cite{Zhong} to part of a large deviation principle for the Airy$_\beta$ point process. In particular, \cite{Zhong} also relies on a quantitative coupling between the extremal eigenvalues of G$\beta$E and the Airy$_\beta$ point process which has similarities with our method to study the G$\beta$E characteristic polynomial at the edge. 
%It is an interesting research direction whether the Laplace functional of the  Cole--Hopf solution admits another representation in terms of the stochastic Airy function.

The Airy$_\beta$ line ensemble is a random collection of curves that describes the trajectories of the largest eigenvalues in Dyson Brownian motion and the extreme particles in the G$\beta$E corners process. This process was recently introduced (for any $\beta>0$)
in \cite{GorinXuZhangAiryBeta} and characterized through its multitime joint moments, which are computed using a discretization involving multivariate Bessel generating functions.

\paragraph{Asymptotic properties of $\SAi$.}
Theorem~\ref{thm:asymp} gives the complete asymptotic expansion of $\SAi_\lambda(t)$ as $t\to\infty$. 
Using the invariance property from Proposition \ref{prop:SAIshift}, this yields a  type of distributional asymptotic approximation for $\SAi_\lambda(0)$ as the spectral parameter $\lambda \in\R \to \infty$.  
There is considerable room to develop the almost sure properties of the stochastic Airy function as the spectral parameter $\lambda\to\pm \infty$. 
In the oscillatory direction (as $t\to-\infty$ for $\lambda=0$), it has now been established in \cite[Theorem 3]{PainchaudPaquette} (see also \cite[Section 2]{LambertPaquetteBulk}) that any solution of the stochastic Airy equation \eqref{SA2} with $T=-1$ satisfies for $t \ge 1$,
\begin{equation} \label{oscasymp}
 \phi_0(-t) = {\rm c}  t^{-1/4 + 1/2\beta}\, e^{ \xi_\beta (t)} \cos \theta_\beta(t),
  \qquad
  \phi_0'(-t) = {\rm c}  t^{1/4 + 1/2\beta}\, e^{\xi_\beta(t)} \sin \theta_\beta(t),
\end{equation}
where $ {\rm c} = \sqrt{c_1^2+c_2^2}$ and $(\xi_\beta,\theta_\beta)$ are adapted (to the filtration $\{B_{-t}; t\ge 1\}$) real-valued processes.  The (random) phase $\theta_\beta$ \emph{oscillates fast} so that  $\theta_\beta(t) \mod 2\pi$ converges in law to a uniform random variable on $[0,2\pi]$ as $t\to\infty$, and  the log-amplitude correction $\xi_\beta (t)$ is approximately Gaussian with variance $\tfrac1{2\beta}\log t$. These asymptotics are in agreement with the high-probability estimates from our Lemma~\ref{lem:energy}. 
We also mention that the structure of the formula \eqref{oscasymp} is consistent with the large-$N$ pointwise asymptotics for the G$\beta$E characteristic polynomial in the bulk of the spectrum recently obtained in \cite{LambertPaquetteBulk}.

%these asymptotics do not reveal the almost sure behavior of $\SAi_\lambda(0)$ for large $\lambda\in\R$, and the question of these almost sure asymptotics is stated in \cite{LambertPaquetteV1} and remains open.
%\begin{question}\label{q:SAI_aslimit}
%  What are the almost sure asymptotics of $\SAi_\lambda$ as $\lambda \to \infty?$
%\end{question}
%In addition, in the $t \to \infty$ direction, while we have given the almost sure asymptotic, we have left open if the convergence holds in $L^p(\Pr)$ for some $p \in (0,\infty)$ or indeed potentially for all $p \in (0,\infty)$.
%\begin{question}\label{q:SAI_LPlimit}
%  For which $p \in (0,\infty)$ does \eqref{eq:SAIconvention} hold in $L^p(\Pr)?$
%\end{question}
 
%\begin{question}\label{q:SAI_oscillatory}
%What is the almost sure behavior of $\SAi_\lambda(t)$ as $\lambda \to -\infty$ or as $t \to -\infty?$
%\end{question}

\paragraph{Hadamard product formula.}
By the Weierstrass factorization theorem and some basic estimates on the growth of the counting function of the Airy$_\beta$ points, the stochastic Airy function admits the representation
\[
\SAi_\lambda(t) = e^{g_\lambda(t)} \prod_{n=1}^\infty \left( 1- \frac{\lambda}{\z_n(t)} \right)e^{\frac{\lambda}{\z_n(t)}},
\]
where $\lambda \mapsto g_\lambda(t)$ is a (random) affine function.
After this paper, the infinite-product formula for the \emph{normalized version} of $\SAi_\lambda(0)$ was obtained using different methods. First, it is established in  \cite{AshburyBridgwood} that such a ratio is a convergent product and that it coincides with the regularized Fredholm determinant of the resolvent of the stochastic Airy operator:
\begin{equation} \label{ratio1}
  \frac{\SAi_\lambda(0)}{\SAi_0(0)} =  \prod_{n=1}^\infty \left( 1- \frac{\lambda}{\z_n} \right)e^{\frac{\lambda}{\z_n}}
  = \det\,_2(\mathrm{I}+\lambda \mathbf{H}_\beta^{-1}) .
\end{equation}
In particular, $\mathbf{H}_\beta^{-1}$ is a Hilbert--Schmidt operator but is not trace class, which is why the regularized determinant $\det_2$ appears in \eqref{ratio1}.  The analysis of  \cite{AshburyBridgwood} relies on the variational characterization of the eigenvalues of the stochastic operator from \cite{RRV09}. 
An alternative approach has been developed by \cite{Assiotis} using the properties of the Laguerre--P\'olya class of entire functions  which can be applied to study the  scaling limit of the characteristic polynomial of many random matrix models. This includes the G$\beta$E at the edge; see \cite[Section 4.3]{Assiotis}, where the connection with the stochastic Airy function is established.
In particular,  \cite[Proposition 4.6]{Assiotis} upgrades the convergence from Theorem~\ref{thm:edge1} to uniform convergence on compact sets of $\C$. Other applications to unitarily invariant infinite Hermitian matrices and $\beta$-ensemble interlacing arrays are also given in \cite{Assiotis} establishing almost-sure product representations for the characteristic polynomials of these ensembles. 
Finally, a general framework for studying the convergence of random infinite products, including rescaled ratios of characteristic polynomials of random matrices, is formulated in \cite{NajnudelNikeghbali}. 
In the random matrix context, the conditions in this paper are related to \emph{eigenvalue rigidity}. 
As an application,  \cite[Section 6.4]{NajnudelNikeghbali} establishes that for G$\beta$E, in distribution as $N\to\infty$:
\begin{equation} \label{ratio2}
  \Bigg(\frac{\varphi_N(1+\tfrac{\lambda}{2 N^{2/3}})}{\varphi_N(1)} e^{-N^{1/3}\lambda}
  : \lambda \in \R\Bigg)
\to  \prod_{n=1}^\infty \left( 1- \frac{\lambda}{\z_n} \right)e^{\frac{\lambda}{\z_n}} .
\end{equation}
We note that the RHS is expressed in a different form in \cite{NajnudelNikeghbali} and that the almost sure limit \eqref{ratio2} is a direct consequence of Theorem~\ref{main:thm} and \eqref{ratio1}. In particular, the only relevant factor in this ratio comes from the  deterministic renormalization $\frac{w_N(1+\frac{\lambda}{2 N^{2/3}})}{w_N(1)} \simeq e^{-N^{1/3}\lambda}$ as $N\to\infty$; see
\eqref{eq:wn}.

\paragraph{The stochastic $\zeta$-function.} The counterpart of the stochastic Airy function associated with the Sine process was first constructed for $\beta=2$ in \cite{ChhaibiNajnudelNikeghbali} by considering the scaling limit of ratios of the characteristic polynomial of a  CUE matrix. In particular, the limiting object, known as the \emph{stochastic zeta function}, also has a convergent infinite-product representation:
\[
\xi_\infty^{(\beta)}(\lambda) = e^{\i \pi \lambda} ~\operatorname{pv}~\prod_{k\in\Z} (1-\lambda/\mathfrak{y}_k)  \qquad\text{where $(\mathfrak{y}_k)_{k\in\Z}$ is a realization of the Sine$_\beta$ point process}.
\]
This result follows in effect from the \emph{rigidity} of the Sine$_\beta$ point process, and it has been generalized to a family of $\beta=2$ random matrices in \cite{CHNNR} by relying on the determinantal structure. 
For general $\beta>0$, the stochastic zeta function has been introduced in \cite{VV20}. 
This object has several interesting descriptions: it has an infinite-product representation; it is the secular function (the counterpart of the characteristic polynomial) for the Sine$_\beta$ operator; and, like $\SAi$, it has an almost-sure representation built from Brownian motion. It is also the almost-sure scaling limit of the characteristic polynomial of the circular-$\beta$ ensemble (within the CMV coupling from \cite{KillipNenciu} --- the counterpart of the Dumitriu--Edelman tridiagonal matrix model for C$\beta$E).
The connection with the Riemann $\zeta$-function and number theory is also discussed in detail in \cite{NajnudelNikeghbali}.

This raises natural questions about the relation between $\SAi_\lambda$ and $\xi_\infty^{(\beta)}$ as $\lambda\to\infty$, or equivalently between the Airy$_\beta$ and Sine$_\beta$ point processes. 
To study this problem, it is actually simpler to consider the asymptotics of $\SAi_\lambda(t)$ as $t\to-\infty$, as discussed above; see \eqref{oscasymp}. In this regime, one can introduce a suitable \emph{phase function} that encodes the eigenfunctions and satisfies an appropriate SDE. In particular, the operator convergence (at the level of canonical systems) from  Airy$_\beta$ to Sine$_\beta$ has been recently established in \cite{PainchaudPaquette}. 
The stochastic zeta function is also expected to describe the microscopic scaling limit of the G$\beta$E characteristic polynomial inside the bulk of the spectrum (by universality, it  is expected for general $\beta$-ensembles). 
This has been established in our follow-up paper \cite{LambertPaquetteBulk}. Our main result is the counterpart of Theorem~\ref{main:thm} at any point in the bulk, $z\in [-1 +\tfrac{C_N}{N^{2/3}}, 1 -\tfrac{C_N}{N^{2/3}}]$, for $\Psi_N(\lambda;z) = (w_N\varphi_N)\left(z+\frac{\lambda}{\varrho(z)N}\right)$ for $\lambda\in\R$, where $\varrho(z)$ is the semicircle density and $C_N\gg1$. 
We obtain pointwise asymptotics for the rescaled characteristic polynomial which include both the Gaussian log-correlated fluctuations and the microscopic Sine$_\beta$ fluctuations. This complements the analysis of this paper by studying the \emph{elliptic part }of the G$\beta$E recursion (where the transfer matrices have complex conjugate eigenvalues for $n\gg Nz^2$). 
We emphasize that the results in this paper are also crucial to derive the bulk asymptotics and that Theorem \ref{thm:para} is the  main input of \cite{LambertPaquetteBulk}.

%\begin{question} \label{q:CUE}
%(1) For general $\beta>0$, show that $\xi_\infty^{(\beta)}(\lambda)$ for $\lambda\in\C$
%is the scaling limit of ratios of the G$\beta$E characteristic polynomials in the bulk. \\
%(2) Show that in a suitable scaling limit, $\SAi_\lambda(t)$ converges locally uniformly for $\lambda\in\C$ to $\xi_\infty^{(\beta)}(\lambda)$ as $t$ or $\lambda \to -\infty$.
%\end{question}

\paragraph{Joint moments.}
There is an extensive body of work on the computation of expectations of products and ratios of characteristic polynomials of random matrices, motivated by local universality considerations and the conjectural relationship to the Riemann $\zeta$-function.  We mention in particular \cite{BorodinStrahov} in which the expectation of ratios of characteristic polynomials of the GOE/GUE/GSE are computed in the bulk scaling limit (the edge behavior is still open to our knowledge) and also \cite{FyodorovStrahov} in which ratios of the characteristic polynomial of orthogonal polynomial ensembles $(\beta=2)$ are considered in more detail using determinantal formulae. These two works reduce the problem to classical asymptotics for orthogonal polynomials (Hermite polynomials in the Gaussian case). 
 Related ratio results for the characteristic polynomial of the classical compact groups have also been given \cite{CFKRS}. For G$\beta$E, we discovered that the analysis of \emph{products} of  characteristic polynomials was carried out in \cite{DesrosiersLiu} using the relationship to Jack polynomials (this work extends to the other classical $\beta$--Laguerre and $\beta$--Jacobi ensembles). At the edge, they identify these moments with a generalized Kontsevich integral, and obtain its large-$N$ asymptotics.  In view of Theorem~\ref{main:thm}, this gives explicit formulas for $\E\big[\prod_{j=1}^n\SAi_{\lambda_j}(0)\big]$ with $\lambda_j\in\R$ in terms of a multivariate Airy function.
We also mention that, based on the results in this paper, 
using the Feynman--Kac formula for solutions of the stochastic Airy SDE \eqref{SA1} and Girsanov's theorem, one can in principle compute $\E\big[\prod_{j=1}^n\SAi_{\lambda_j}(t_j)\prod_{j=1}^m\SAi_{\mu_j}'(s_j)]$ for $\lambda_j,\mu_j\in\C, t_j,s_j\in\R$ in terms of solutions of some ODE and deduce the counterpart for the G$\beta$E characteristic polynomial at the edge.
Such computations are carried out in the thesis \cite{Kirkiles}, leading to some explicit formulae. 
For instance, based on \cite[Theorem 1.2.4]{Kirkiles}, one obtains for a fixed $\lambda \in \R$, as $N\to\infty$, 
\begin{equation} \label{mom2}
  \Exp\bigl[\Psi_N^2(\lambda)\bigr]
  \sim
  \frac{N^{\frac2{3\beta}}}{4\pi^{3/2} \i}
  \int_{-\i\infty}^{\i\infty}
  e^{z^3/12 - \lambda z}\, z^{-2/\beta - 1/2}\, \d z,
\end{equation}
a formula which is originally due to \cite{su2010second}.
For the third and fourth moments, this leads to formulas analogous to \eqref{mom2} involving Whittaker functions.

Joint expectations of powers of characteristic polynomials of $\beta$-ensembles are also related to so-called \emph{Fisher--Hartwig asymptotics} and \emph{Gaussian multiplicative chaos} (GMC). These questions have been thoroughly investigated for $\beta=2$ using the determinantal structure and the corresponding Riemann--Hilbert problem for orthogonal polynomials. We refer to the works \cite{CharlierHankelFH,CharlierFahsWebbWongMultiCutFH} for asymptotics of powers of characteristic polynomials for general Hermitian in the one-cut and multi-cut respectively. 
In the recent work of  Keles \cite{KelesGMC}, this has been generalized to merging singularities for the stationary Hermitian matrix Ornstein--Uhlenbeck process (a stationary dynamical extension of the GUE). In particular, this implies convergence of renormalized powers of the characteristic polynomial to a GMC. The method is based on \cite{BourgadeFalconetLQG} which considers a similar problem for the circular unitary ensemble and establishes a first connection between random matrix dynamics and two-dimensional Liouville quantum gravity.

Finally, Assiotis and Najnudel \cite{AssiotisNajnudelCBEField} recently computed \emph{supercritical moments} for the partition function associated with the circular $\beta$-ensemble characteristic polynomial. These moments, and also the correlation functions of the Sine$_\beta$ point process, are expressed in terms of a variant of the stochastic $\zeta$- function.

\subsection{Organization} \label{sec:org}

Let us draw a roadmap for the proofs of the results presented above.
In Section~\ref{sect:prel}, we introduce the formalism for concentration of random variables that we will use in the remainder of this paper and the Riccati diffusion which is a crucial tool to obtain the asymptotic behavior of solutions of the stochastic Airy equation \eqref{SA1}.
Then, the paper is divided in two independent parts:
\begin{itemize}[leftmargin=*]
\item In Sections~\ref{sec:SAE}--\ref{sec:Red5}, we study the properties of the stochastic Airy equation, including its stability, the asymptotics of its solutions and its spectral theory.
\item In Section~\ref{sec:quant} and Appendix~\ref{sec:hyperuniform}, we consider the characteristic polynomial of the Gaussian $\beta$-ensemble in the \emph{parabolic regime}. These results rely crucially on our previous results from \cite{LambertPaquette02} on the \emph{hyperbolic regime} 
(Appendix~\ref{sec:hyperuniform}) and the stability properties of the stochastic Airy equation developed above. The goal is to give quantitative approximation results for the characteristic polynomial in terms of the stochastic Airy function. 
\end{itemize}

In Section~\ref{sec:SAE}, we show (global) existence and uniqueness of the solutions \eqref{SA1}, interpreted as a Volterra-type equation in the form \eqref{SA2},  and develop the basic properties of the solution map.
Section~\ref{sec:SAK} addresses the stability of solutions of  \eqref{SA2}.
Specifically, solutions can be represented in terms of a natural object that we call the \emph{stochastic Airy kernel}.
Our main result is an overwhelming probability estimate (cf.~Theorem~\ref{thm:Kest}) for this kernel which is instrumental in showing that the characteristic polynomial \eqref{eq:recurrence}, after the appropriate rescaling, approximately solves the stochastic Airy equation. In particular, we prove Proposition~\ref{prop:SAiwins} and another key stability result, Proposition~\ref{prop:stability}, at the end of this section. 

\medskip

In Section~\ref{sec:SAi}, we properly define the \emph{stochastic Airy function} $\SAi$.
This definition arises from the almost sure asymptotics of the Dirichlet and Neumann solutions of \eqref{SA2}.
In fact, we obtain the full asymptotics as $t\to\infty$ of any solution of \eqref{SA2} with overwhelming probability. This construction allows us to show that $\SAi_\lambda$ is the only solution of \eqref{SA2} lying in the Sobolev space $H^1(\R_+)$ and to prove Theorem~\ref{thm:asymp}.
This also makes a connection with the spectral theory of the \emph{stochastic Airy operator}.
In addition, we show that for any Robin-type boundary condition, this operator is \emph{diagonalized} by the stochastic Airy function, cf.~Theorem~\ref{thm:spec}.
In conjunction, we prove that the zero set of $\lambda \mapsto \SAi_\lambda(0)$ is a realization of the \emph{Airy$_\beta$ point process} and provide a description of its counting function in terms of the Riccati diffusion, cf.~Proposition~\ref{prop:cf}.
Hence, this provides an alternative approach to the results of \cite{RRV09} on the spectrum of the \emph{stochastic Airy operator}.
Our construction of the \emph{stochastic Airy function} is based on the fine asymptotics of the Riccati diffusion and the main steps are explained in Section~\ref{sec:asymp}. The technical details of the proofs are collected in 
Sections \ref{sec:Delta}--\ref{sec:X}. 

\medskip

In Sections~\ref{sec:SAestimate}--\ref{sec:Red5}, we give the proof of the kernel estimate Theorem~\ref{thm:Kest}.
Section~\ref{sec:SAestimate} provides a simple series of reductions after which
it suffices to estimate the growth of an arbitrary solution $\phi_\lambda(t)$ of \eqref{SA1} with $\lambda=0$.
In the oscillatory direction $\{t \le 1\}$, we obtain such bounds by a basic \emph{energy method} in Section~\ref{sec:Red3}.
Our bounds in the expanding direction $\{t \ge 1\}$ are presented in Section~\ref{sec:Red5}, see Proposition~\ref{prop:GE}.
Our approach relies again on the Riccati diffusion: in Sections~\ref{sec:bd}--\ref{sec:sc}, we give some technical moderate deviations estimates for the stability and the number of blow-downs of the Riccati diffusion.

\medskip

In Section~\ref{sec:quant}, we study the
characteristic polynomial of the Gaussian $\beta$-ensemble for $z\in [-1,1]$ in the parabolic regime. 
We introduce the relevant coupling to a Brownian motion (depending on $z$ and $N$) and we present two quantitative high-probability approximation results in terms of two different solutions of the stochastic Airy equation (Theorem~\ref{main:approx} in terms of the stochastic Airy function and Theorem~\ref{thm:para} in terms of an adapted solution). 
These results are proved in Section~\ref{sec:fde} by viewing the recurrence \eqref{eq:recurrence} as a discretization of the integral equation \eqref{SA2} with $T=(\log N_0)^{1-\epsilon}$.
Using Proposition~\ref{prop:stability}, this allows us to compare the characteristic polynomial to the solution of the stochastic Airy equation with the same initial conditions at time $T$. Since this equation is run backward in time with $T\gg1$, by Proposition~\ref{prop:SAiwins}, up to a negligible error, any such solution is a multiple of the stochastic Airy function and the multiplicative constant is determined by the initial data. Finally, we can compute this constant with high probability using our results from \cite{LambertPaquette02} on the hyperbolic part of the characteristic polynomial recursion and the asymptotics of the solutions to the stochastic Airy equation. At the edge of the spectrum ($z=1$), this result reduces to Theorem~\ref{main:thm}. 

\paragraph{Acknowledgements}

We would like to thank Diane Holcomb, conversations with whom helped launch this project.
We would like to thank Ovidiu Costin for helpful conversations.  
We would like to acknowledge support from the Park City Mathematics Institute 2017, at which this program was begun, and in particular acknowledge NSF grant DMS:1441467.
 GL acknowledges the support of the following grants: SNSF Ambizione grant S-71114-05-01, Swedish Research Council starting grant 2022-04882, as well as the Ragnar S\"oderbergs Foundation and G\"oran Gustafsson Foundation for their support. 
EP would like to acknowledge support from a Simons travel grant, 638152.

% ===== end sections/01-introduction.tex =====

% ===== begin sections/03-preliminaries.tex =====
\section{Preliminaries} \label{sect:prel}

Throughout this article, the parameter $\beta>0$ is fixed and we do not keep track of the $\beta$-dependency of various constants.
We make use of the symbols $\lesssim$ and $\gtrsim$ in the following form.  We write $f(x) \lesssim_{\alpha} g(x)$ if there is a finite function $C=C(\beta,\alpha)>0$ so that $|f(x)| \le Cg(x)$ for all relevant $x$. 
We alternatively use $f(x) =\O_{\alpha}\big(g(x) \big)$ or $f(x) \leq \O_\alpha\big(g(x)\big)$ to mean  $f(x) \lesssim_{\alpha} g(x)$.  If we omit the subscript $\alpha$ in either case, we mean the inequality holds with a constant $C_\beta$ which only depends on $\beta>0$.  We will sometimes use the notation $f(N) = N^{-\alpha_-}$ for a real number $\alpha > 0$ to mean a function such that for any $\epsilon > 0$ there is a $C = C(\beta,\epsilon)$ so that for all $N \gg 1$, $f(N) \leq CN^{\epsilon-\alpha}$.

% We also let  $\diag(M)$ denotes the diagonal matrix matching the diagonal of matrix $M.$
% We take the convention that for a sequence of matrices $\left\{ M_n \right\}$, 
% \[
%  {\textstyle \prod_{j=p}^n} M_j = M_n M_{n-1} \dots M_{p+1}M_p.
% \]
Throughout this article, $B$ is a standard Brownian motion and for any $s \in \R$, we define
\[
  \mathscr{G}_s = (\mathscr{G}_{s,t} : t \in \R)
  \quad\text{with}\quad
  \mathscr{G}_{s,t} = \sigma\big( B(u)-B(s) : u \in [s,t] \big),
\]
where we allow both $t > s$ and $t < s$  (this is a filtration for $t > s$ and a reversed filtration for $t < s$).
Moreover, we work on the event that $B$ is $\alpha$-H\"older for a fixed $0<\alpha<1/2$.

If $z\in \C$, we denote by $(z)_+ = \max\{ \Re z, 0\}$. 

\subsection{Concentration} \label{sec:concentration}

On top of the Gaussian coupling from the previous section, our analysis  makes crucial use of the theory of sub--Gaussian and sub--exponential random variables.  We review many standard concentration results in terms of this theory, following \cite[Chapter 2]{Vershynin}, where one may find the proofs of all the claims in this section. Classical estimates for Brownian motion are also formulated in terms of this theory in the Appendix~\ref{sec:bm}.

Define, for any $p \geq 1,$ and any complex valued random variable $X$, 
\[
  \VERT X \VERT_p = \inf\left\{ t \geq 0 : \Exp e^{|X|^p/t^p} \leq 2\right\}.
\]
For all $X$ for which $\VERT X \VERT_p < \infty,$ this defines a norm.  In the cases of $p=1$ and $p=2,$ these are the sub--exponential norm and the sub--Gaussian norm, respectively, and those are the only two cases we use in this paper.

By Markov's inequality, it follows that if $\VERT X \VERT_p < \infty,$ then for all $t \geq 0$
\begin{equation}\label{eq:Xtail}
  \Pr \left[ |X| \geq t \right] \leq 2\exp( -t^p/ \VERT X \VERT_p^p ),
\end{equation}
on observing the infimum in the definition of $\VERT \cdot \VERT_p$ is attained whenever it is finite.  Moreover, this concentration inequality is equivalent to the finiteness of $\VERT\cdot\VERT_p,$ in that if there exists $s \geq 0$ such that
\begin{equation}\label{eq:normup}
  \Pr \left[ |X| \geq t \right] \leq s\exp( -t^p ) \quad \text{for all } t \geq s,
\end{equation}
then $\VERT X \VERT_p \leq C_{p,s}$ for some absolute constant $C_{p,s} > 0.$ 

Finally we observe as a consequence of Young's inequality that for any $p,q \geq 1$ satisfying $1/p + 1/q = 1,$
there is an absolute constant $C_{p,q}$ so that
for any two random variables $X$ and $Y$, 
\begin{equation}\label{eq:VERT}
  \VERT XY \VERT_{1} \leq C_{p,q} \VERT X\VERT_p \VERT Y \VERT_q.
\end{equation}
In the particular case that $p=q=2,$ one can further take $C_{2,2}=1.$

We also make use of concentration inequalities for the local supremum of a Gaussian process, in the following form.  Suppose that $K \subset \C$ is a compact set and that $G : K \to \C$ is a centered Gaussian process which has the property that
\[
  \Exp |G(z) - G(x)|^2 \leq |z-x|^{\alpha}
\]
for some $\alpha \in (0,2].$  Then there is a constant $C=C(K,\alpha)$ so that
\begin{equation}\label{eq:Dudley1}
  \Big\VERT\sup_{z,x \in K} |G(z)-G(x)|\Big\VERT_2 \leq C.
\end{equation}
This follows as a corollary of standard chaining techniques (Talagrand majorization or Dudley's inequality), see \cite[Theorem 8.1.6]{Vershynin}.  If in addition we have that $\VERT G(y)\VERT_2 \leq C$ for some $y \in K$, then we may further conclude
\begin{equation}\label{eq:Dudley2}
  \VERT\sup_{z \in K} |G(z)|\VERT_2 \leq 2C.
\end{equation}

\subsection{Riccati diffusion}  \label{sec:SAR}

In this section, we present basic facts about solutions of the so-called \emph{Riccati SDE} \eqref{Ric0} and its connection to the stochastic Airy equation. 
Analyzing the behavior of this diffusion is our main tool to obtain the asymptotics as $t\to\infty$ of solutions of the stochastic Airy equation of Theorem~\ref{thm:asymp}, as well as to obtain moderate deviations bounds for the growth of any solution in Section~\ref{sec:Red5} 
(see Proposition~\ref{prop:GE}). 
These bounds are instrumental for obtaining the high probability asymptotics of the characteristic polynomial of the Gaussian $\beta$-ensemble around the spectral edges. 
%key estimates for the so-called stochastic Airy kernel (see Theorem~\ref{thm:Kest} which is proved in Section~\ref{sec:SAestimate}). 

\medskip

As we explained in the introduction, the equation \eqref{SA1} can formally be viewed as an eigenvalue problem for the stochastic Airy operator $\mathbf{H}_\beta$ with eigenvalue $-\lambda$, which turns out to be a generalized Sturm--Liouville problem.
With our conventions, $\{ t\le 0\}$ is the region where solutions of  \eqref{SA1} have    \emph{oscillatory behavior} while they generically have an exponentially growing behavior for  $\{t\ge 0\}$. 
It is well-known in Sturm--Liouville theory that it is possible to analyze these behaviors by performing a so-called \emph{Riccati substitution}: $\rho = \phi_\lambda'/\phi_\lambda$. 
Then, we verify that $\rho =\rho_{\lambda,s}$ solves the following equation\footnote{We refer to the proof of Lemma~\ref{lem:phipv}.}:
\begin{equation} \label{Ric0}
\d\rho(t) = (t+\lambda-\rho^2) \d t + \tfrac{2}{\sqrt{\beta}}\d B(t) , \qquad \rho(s) = \omega .
\end{equation}
where $s\in\R$ is fixed and $\lambda\in\C$ is a parameter.  
As this SDE is locally Lipschitz, \eqref{Ric0} has a (unique) strong solution
$\rho : [s, \mathfrak{z}_1) \to \R$ with $\mathfrak{z}_1 = \inf\big\{ t > s : \rho(t) = -\infty\big\}$.
Indeed, any solution is a continuous adapted process which can \emph{blow down} to $-\infty$ in finite time, which corresponds naturally to a zero of $\phi.$
In particular, $\mathfrak{z}_1$ is a stopping time with respect to the filtration 
 $\mathscr{G}_{s} = \sigma\big( B(u)-B(s) : u \ge s \big)$.
Moreover, before blowing down, two solutions with different initial conditions cannot cross. Hence, by a limiting procedure, it is possible to define a solution with initial data $ \omega = +\infty$ (for $s=0$, this corresponds to the Neumann solution studied in \cite{RRV09}). 
Hence, we can define a (unique) continuous solution $\rho : [s,\infty)\to \overline\R$ by induction as follows: For $k\in \N$, on the event $\{\mathfrak{z}_k <\infty\}$, the process $\rho : [\mathfrak{z}_k,\mathfrak{z}_{k+1}) \to  \overline\R $ is the solution of \eqref{Ric0} with initial data  $\rho(\mathfrak{z}_k ) = +\infty$ and $\mathfrak{z}_{k+1} = \inf\big\{ t >\mathfrak{z}_k : \rho(t) = -\infty\big\}$.

For a given $s\in\R$, these solutions $(\rho_{\lambda,s})_{\lambda\in\C, \omega \in (-\infty, +\infty]}$ are coupled together  and they \emph{interlace}\footnote{This means that $\rho_{\lambda,s}(t) \le \rho_{\mu,s}(t)$ for all $t \in [s, \mathfrak{z}_1)$ where $\mathfrak{z}_1 = \inf\big\{ t > s : \rho_{\lambda,s}(t) = -\infty\big\}$, then on the event $\{\mathfrak{z}_1<\infty\}$, we have $\rho_{\lambda,s}(t) \le \rho_{\mu,s}(t)$ for all $t \in [\mathfrak{z}_1 , \mathfrak{z}_1']$ where   $\mathfrak{z}_1' = \inf\big\{ t > s : \rho_{\mu,s}(t) = -\infty\big\}$, etc.}:
For $\lambda, \mu\in\R$, 
\begin{equation} \label{interlacing}
\rho_{\lambda,s} \prec \rho_{\mu,s} \qquad\text{ if  $\rho_{\lambda,s}(s) \le \rho_{\mu,s}(s)$ and $\lambda\le\mu$}.
\end{equation} 
For $\lambda\in\R$, let us denote for $t\in\R$, 
\[
\mathbf{N}_{[s,t)}(\rho_{\lambda,s})  = \#\big\{ u\in[s,t), \rho_{\lambda,s}(u) =\infty \big\} 
=  \#\big\{ u\in [s,t), \phi_\lambda(u) =0 \big\}. 
\]
These functions are c\`adl\`ag  non-decreasing and  \eqref{interlacing} implies that for all $t\ge s$, 
\[
\mathbf{N}_{[s,t)}(\rho_{\lambda,s}) \le \mathbf{N}_{[s,t)}(\rho_{\mu,s}) \le \mathbf{N}_{[s,t)}(\rho_{\lambda,s})+1 .
\]

Moreover, according to Proposition~\ref{prop:zerobounds} below, for any
$\lambda\in\R$ with $\lambda\ge 1$, 
\[
 \mathbf{N}_{[0,\infty)}(\rho_{\lambda,0}) = \lim_{t\to\infty} \mathbf{N}_{[0,t)}(\rho_{\lambda,0}) \text{ is finite almost surely.}
\]
Besides being used to count the number of zeros of a given solution of the stochastic Airy equation \eqref{SA1}, there is a precise relation between $\phi_\lambda$ and its Riccati transform which is summarized by the next lemma. 

\begin{lemma} \label{lem:phipv}
Almost surely, any non-trivial solution of \eqref{SA2} has only isolated simple zeros. 
Let $\phi_\lambda$ be the solution of \eqref{SA2}  with $\lambda\in\R$ and initial data  $\left\{ \phi_\lambda(s) = c_1, \phi_\lambda'(s)=c_2 \right\} \in \R^2$ with $c_1\neq 0$,
then the function $\rho_{\lambda} : \R \to \overline{\R}$ given by $\rho_{\lambda}(t)  = \phi_\lambda'(t) /\phi_\lambda(t)$ is continuous and it holds for all $t\in\R$, 
 \begin{equation} \label{Ric_representation}
\phi_\lambda(t) = c_1 \exp\biggl( \operatorname{pv} \int_{s}^{t}\rho_{\lambda}(u)\d u+  \i \pi \mathbf{N}_{[s,t)}(\rho_{\lambda}) \biggr) ,  
\end{equation}
 where $ \operatorname{pv}$ has the same meaning as in Lemma~\ref{lem:pv}. 
Moreover, if $c_1,c_2$ are fixed, $\rho_\lambda$ solves the SDE \eqref{Ric0} with initial condition $\omega =c_2/c_1$. 
\end{lemma}

\begin{proof}
%Like in the proof of Lemma~\ref{lem:Wronskian}, we verify that the Wronskian of $\psi$ with another (linearly independent) solution of \eqref{SA1} is constant for all $t\ge 0$. This shows that all zeros of $\phi$ are simple. 

The proof relies on the results of Section~\ref{sec:SAE} and  Lemma~\ref{lem:pv}. 
Namely, according to Lemma~\ref{lem:Wronskian}, (almost sure) existence and uniqueness of the solution map $\T$ implies that if any solution has a double-zero, then it is identically zero.  
In turn, as $\T$ takes values in $C^1(\R)$, if the zeros of a solution had an accumulation point $\mathfrak{z} \in \R$, then  $\mathfrak{z}$ would be a double-zero, which cannot happen. 
Then, the representation \eqref{Ric_representation} is a direct consequence of the deterministic  Lemma~\ref{lem:pv}.  
Finally, if $c_1,c_2$ are deterministic, $\phi_\lambda$ solves \eqref{SA1} in the  It\^o  sense and we verify that $\rho_{\lambda}$ is a solution of \eqref{Ric0} by applying It\^o formula. 
Namely, $\phi_\lambda'$ is a continuous semimartingale and $\d\phi_\lambda = \phi_\lambda' \d t$  since $\phi_\lambda\in C^1(\R)$, so that 
\[
\d\rho_{\lambda,s} = \frac{\d\phi_\lambda'}{\phi_\lambda} - \frac{(\phi_\lambda')^2}{\phi_\lambda^2}\d t =  (t+\lambda-\rho^2) \d t + \tfrac{2}{\sqrt{\beta}}\d B(t) . \qedhere
\]
\end{proof}

Lemma~\ref{lem:phipv} provides a representation for real-valued solutions of the stochastic Airy equation ($\lambda\in\R$ and $(c_1,c_2) \in\R^2$).  
Using the following observation, this can be extended to general solutions. 
%Since $\phi_{\overline{\lambda}} = \overline{\phi_\lambda}$

\begin{lemma} \label{lem:nozero}
  Let $s \in \R$ be fixed and suppose that $\eta = \Im\lambda \neq 0$.
Let $\phi_\lambda$ solve \eqref{SA2} with initial data  $(c_1,c_2) \in\C^2$. 
The function  $t \mapsto \phi_\lambda(t)$ has at most one zero.
Moreover, if  $(c_1,c_2) \in\R^2$ with $c_1\neq 0$, then $\phi_\lambda$ has no zeros. 
\end{lemma}

\begin{proof}
We consider the process $I(t) = \Im\big\{ \phi_\lambda'(t) \overline{\phi_\lambda(t)} \big\}$ for $t\in\R$.  
 $I$ is a continuous $\mathscr{G}_s$-adapted process and applying It\^o formula, 
 \[
\d I(t) = \eta |\phi_\lambda(t)|^2\d t.
\]
Hence, for any non-trivial solution, $\displaystyle I(t) = \Im\{c_2\overline{c_1} \} +\eta \int_s^t |\phi(u)|^2 \d u$ has at most 1 zero . 
In case $(c_1,c_2) \in\R^2$, this implies that neither $\phi_\lambda$ nor $\phi'_\lambda$  has any zero except possibly at $t=s$.
\end{proof}

This lemma shows that we can also perform the Riccati substitution for $\lambda \in\C$ in which case if $\phi_\lambda$ solves \eqref{SA2} for fixed $(c_1,c_2) \in\R^2$ with $c_1\neq 0$, then it holds for all $t\in\R$, 
 \begin{equation} \label{Ric_representation_C}
\phi_\lambda(t) = c_1 \exp\biggl( \int_{s}^{t}\rho_{\lambda,s}(u)\d u \biggr) 
\qquad\text{and}\qquad 
\mathbf{N}_{[s,\infty)}(\rho_{\lambda,s}) =0 . 
\end{equation}

% ===== end sections/03-preliminaries.tex =====

% ===== begin sections/04-sae-solution-operator.tex =====
\section{The solution operator of the stochastic Airy equation} \label{sec:SAE}

In this section, we show the equivalence of \eqref{SA1} with the integral equation \eqref{SA2} defined on the space of continuous functions, and we develop the properties of the solution map associated to this equation. 
%Recall that $B$ is a standard Brownian motion and that for any $s \in \R$,
%\[
%  \mathscr{G}_s = (\mathscr{G}_{s,t} : t \in \R)
%  \quad\text{with}\quad
%  \mathscr{G}_{s,t} = \sigma\big( B(u)-B(s) : u \in [s,t] \big),
%\]
%where we allow both $t > s$ and $t < s$  (this is a filtration for $t > s$ and a reversed filtration for $t < s$).

 \begin{proposition} \label{prop:Top}
For a fixed $s\in\R$, there exists a continuous function $\C \times \R^2 \ni (\lambda, u,t ) \mapsto  A_{\lambda}(t,u)$ which is entire in $\lambda$,  $\mathscr{G}_{s,t}$--measurable for all $u \in [s,t]$ such that for all $(c_1, c_2) \in \C^2$, the function
  \[
    \Phi_\lambda(t) = c_2 + c_1 \Es_\lambda(t,s) + \int_s^t A_{\lambda}(t,u)( c_2 + c_1 \Es_\lambda(u,s) )\,\d u, \qquad t\in \R , 
  \]
  is the unique solution to \eqref{SA2}.  Then for deterministic $(c_1,c_2)\in\C^2$,  the processes $\big(\Phi_\lambda(t)\big)_{t\ge s}$ and  $\big(\Phi_\lambda(t)\big)_{t\le s}$  are adapted to $\mathscr{G}_s$ and $    \phi_\lambda(t) = c_1 + \int_s^t     \Phi_\lambda(v) \d v $  is the unique solution of \eqref{SA1}. 
\end{proposition}

\begin{proof} 
  The integral equation \eqref{SA2} is of Volterra type, and so the kernel $A$ is easily constructed by Picard iteration.  Define inductively for $m \geq 2$ and $(t,s)\in\R^2$,  
  \begin{equation}\label{eq:Hms}
  \mathcal{K}^m_\lambda(t,s) := \int_s^t   \mathcal{K}^1_\lambda(t,u)   \mathcal{K}^{m-1}_\lambda(u,s)  \d u
  \quad
  \text{where}
  \quad
  \mathcal{K}^1_\lambda(t,u)
  : =
  \Es_\lambda(t,u).
\end{equation}
Each of these kernels is entire in $\lambda$ and continuous in its other parameters.  In addition, $    \mathcal{K}^m_\lambda(t,u)$ are $\mathscr{G}_{s,t}$--measurable for all $u \in [s,t]$.  For a fixed compact set $K\subset\C$, let
\[
I_{s,t}:=[s\wedge t,s\vee t],
\qquad
M_{K;s,t}:=
\sup_{\substack{\lambda\in K\\v,w\in I_{s,t}}}
|\Es_\lambda(v,w)|.
\]
Then induction gives, for every $m\geq1$,
\begin{equation}\label{eq:VoltEst}
sup_{\substack{\lambda\in K\\u\in I_{s,t}}}
\big|\mathcal K_\lambda^m(t,u)\big|
\le M_{K;s,t}^{\,m}\frac{|t-s|^{m-1}}{(m-1)!},
\qquad (s,t)\in\R^2.
\end{equation}
 Hence, the sum
$A_{\lambda,s}(t,u) \coloneqq \sum_{m=1}^\infty \mathcal{K}^m_\lambda(t,u)$
converges locally uniformly for $\lambda\in\C$, $t\in\R$ and $u\in[s,t]$. This defines a kernel which is entire in $\lambda\in\C$,  continuous in $(t,u)\in\R^2$ and has the desired adaptedness.
Moreover, it holds for any $(s,t)\in\R^2$, 
\begin{equation} \label{Dir'}
\int_s^t    \Es_\lambda(t,u)  A_{\lambda,s}(u,s)  \d u = A_{\lambda,s}(t,s) -  \Es_\lambda(t,s)  , 
 \end{equation}
interchanging the integral and sum which by virtue of \eqref{eq:VoltEst} is justified and using \eqref{eq:Hms}. 
This already shows that $\R\ni t \mapsto   A_{\lambda}(t,s)  $ is a  solution to \eqref{SA2} with \emph{initial data} $(c_1,c_2)=(1,0)$. 
Moreover by Fubini's Theorem, this implies that  with $\Phi_\lambda$ as in the statement of the proposition, we have
\[
  \begin{aligned}
  \int_s^t  \Es_\lambda(t,u) \Phi_\lambda(u)\,\d u
  =\int_s^t   A_\lambda(t,u)(c_2  + c_1 \Es_\lambda(u,s)) \,\d u
=  \Phi_\lambda(t) - \big( c_2 + c_1 \Es_\lambda(t,s) \big)
  \end{aligned}
\]
so we conclude that $\Phi_\lambda$  solves \eqref{SA2}.
The uniqueness of this solution follows from Banach's fixed point Theorem. 
Indeed, the difference of two solutions is a fixed point of the mapping 
\[
 f \in C(\R) \mapsto \Bigg( \R\ni t \mapsto  \int_s^t  \Es_\lambda(
 t,u) f(u)\,\d u \Bigg) .
\]
Using the almost sure H\"older--continuity of $B$ (and of the kernel $\Es_\lambda$, see \eqref{kernelU}), this mapping is almost surely (locally) a contraction, so its only  fixed point is $f=0$.

Finally, given deterministic initial data $c_1,c_2\in\C$, by construction the stochastic process $\big(\phi_\lambda(t) = c_1 + \int_s^t \Phi_\lambda(u) \d u \big)_{t\in\R}$ is $C^1$ and adapted to $\mathscr{G}_s$ in the above sense. Hence, the following It\^o integral is well-defined and by an integration by parts, it holds for $t\ge s$, 
\begin{equation} \label{Fubini}
\int_{s}^t \phi_\lambda(u) \big( (u +\lambda) \d u + \tfrac{2}{\sqrt{\beta}}\d B(u)\big)  =  c_1 \Es_\lambda(t,s) + \int_{s}^t \Phi_\lambda(v)\Es_\lambda(t,v) \d v .
\end{equation}
Here, we have used that by definition $\displaystyle \Es_\lambda(t,s)  =   \int_{s}^t \big((u+\lambda) \d u + \tfrac{2}{\sqrt{\beta}}\d B(u)\big)$  for  $t\ge s$ and $\phi_\lambda(s) =c_1$. 
Similarly, \eqref{Fubini} also holds for $t\le s$. 
This shows that $\phi_\lambda$ solves the It\^o stochastic equation \eqref{SA1} both for $t\le s$ and $t\ge s$.  
Conversely, any $C^1$ solution of \eqref{SA1} such that $\Phi_\lambda = \phi_\lambda'$ is adapted to   $\mathscr{G}_s$ gives rise to a solution of the integral equation \eqref{SA2}. 
 This completes the proof.
\end{proof}

\begin{remark} \label{rk:inverse}
Observe that if $\phi_\lambda$ solves the  It\^o stochastic equation \eqref{SA1} with fixed initial data, upon reversing time and changing the driving Brownian motion, 
$(\phi_\lambda(t) : t \leq s) \overset{\rm law}{=}(\phi_{-\lambda-s}(t) : t\ge 0)$. As \eqref{SA2} has a transparent meaning (almost surely) for any $t \in \R,$ we do not dwell on this.
\end{remark}

This allows us to define the \emph{solution map} $\T$ to the stochastic Airy equation \eqref{SA2}  by the formula:
\begin{equation}\label{eq:T}
  \T(c_1,c_2,\lambda,s,t)
  =
  c_2 + c_1 \Es_\lambda(t,s) + \int_s^t A_\lambda(t,u)( c_2 + c_1 \Es_\lambda(u,s) )\,\d u . 
\end{equation}
%From \eqref{eq:Hms}, it is clear that the kernels $H_m$ are continuous in $s,$ and from \eqref{eq:VoltEst}, the series $\mathscr{H}_{\lambda,s}(t,u)$ will also be continuous in $s$ as it converges absolutely.
%Thus we have shown
% Hence, in fact we conclude that $\T$ has all the desired properties in Theorem \ref{thm:T}. 

%As a corollary of Lemma~\ref{lem:sol}, we obtain that the initial value problem \eqref{SA2} has (almost surely) a strong solution $\phi_\lambda\in C^1(\R)$. 
%Our next goal is to show that we can define the \emph{solution map} $\T_{\lambda,s}$ consistently for all $s\in\R$ and $\lambda\in\C$. We begin by showing that the maps $(\T_{\lambda,s})_{s\in\R}$ from Lemma~\ref{lem:sol} have the same range.  

In the remainder of the section, we develop some basic properties of this solution map~$\T$.
While this map is entire for $\lambda\in\C$, jointly continuous for $(s,t) \in\R^2$ from the properties of the kernels $\mathscr{U}$ and $A$, there is a stronger connection between the map for different values of $s.$

\begin{lemma} \label{lem:Top}
  Almost surely, 
  for any $s_1,s_2 \in \R$ and any $c_{11},c_{12},\lambda \in \C$, 
if $\Phi_\lambda = \T(c_{11},c_{12},\lambda,s_1,\cdot)$ 
then 
$ \Phi_\lambda = \T(c_{21},c_{22},\lambda,s_2, \cdot)$ 
  where
  \begin{equation} \label{T:relation} 
    \begin{cases}
      c_{21}  = c_{11} + \int_{s_1}^{s_2} \Phi_\lambda(u) du,  \\
      c_{22} = c_{12} +  c_{11} \Es_\lambda({s_2},{s_1}) +  \int_{s_1}^{s_2} \Es_\lambda(s_2,u) \Phi_\lambda(u) du  = \Phi_\lambda(s_2).
    \end{cases}
  \end{equation}
\end{lemma}
\noindent This shows that the (random) space of functions $\big\{ (t,\lambda) \mapsto \T( c_1,c_2, \lambda,s,t), (c_1,c_2)\in\C^2 \big\}$  is a 2-dimensional vector space which does not depend on the choice of base point $s.$

\begin{proof}
By definition, we have for any $t\in\R$,
\[ \begin{aligned}
\Phi_\lambda(t) & = c_{12} + \int_{s_1}^t \Es_\lambda(t,u) \Phi_\lambda(u) du + c_{11} \Es_\lambda(t,{s_1}) \\
& =  c_{12} +  c_{11} \Es_\lambda({s_2},{s_1}) +   \int_{s_1}^{s_2} \Es_\lambda(t,u) \Phi_\lambda(u) du + \int_{s_2}^t \Es_\lambda(t,u) \Phi_\lambda(u) du   + c_{11} \Es_\lambda(t,{s_2})
\end{aligned}\] 
where we used that $ \Es_\lambda(t,{s_1}) =  \Es_\lambda(t,{s_2}) +  \Es_\lambda({s_2},{s_1})$.  Similarly, we have 
\[
\int_{s_1}^{s_2} \Es_\lambda(t,u) \Phi_\lambda(u) du  = \int_{s_1}^{s_2} \Es_\lambda(s_2,u) \Phi_\lambda(u) du +   \Es_\lambda(t,s_2) \int_{s_1}^{s_2} \Phi_\lambda(u) du .
\]
Hence, we conclude that $ \Phi_\lambda(t) = c_{22} + \int_{s_2}^t \Es_\lambda(t,u) \Phi_\lambda(u) du + c_{21} \Es_\lambda(t,{s_2}) $ where $c_{22}$ and $c_{21}$ are as in the statement of the lemma.  The uniqueness of solutions to \eqref{SA2}, established in Proposition~\ref{prop:Top}, concludes the proof.
\end{proof}

%\begin{proposition} \label{prop:Top}
%For a fixed $\lambda\in\C$, almost surely, there exists a map $\T_\lambda : \R\times \C^2 \mapsto C(\R)$ such that 
%for any $(s, c_1, c_2) \in\R\times\C$, $f=  \T_{\lambda,s}(c_1,c_2)$ is the unique solution of the integral equation \eqref{SA2}.
%\end{proposition}
%
%\begin{proof} 
%By Lemma~\ref{lem:sol}, almost surely, the family of maps $ \T_\lambda =  (\T_{\lambda,s})_{s\in\Q}$ is well-defined.
%Moreover, by Lemma~\ref{lem:Top},  for any $s_1, s_2 \in\Q$ and $(c_{11},c_{12}) \in \C^2$, there exists  $(c_{21},c_{22}) \in \C^2$ such that
%\[
%\T_{\lambda,s_1}(c_{11},c_{12}) =  \T_{\lambda,s_2}(c_{21},c_{22}) = f
%\]
%and the relationships \eqref{T:relation} hold. In particular, since both $f$ and $\Es_\lambda$ are continuous function, it holds  as $s_2 \to s_1$,
%\[
%c_{21}  \to c_{11} 
%\qquad\text{and}\qquad
% c_{22} \to c_{12} .
%\]
%This implies that we can extend the family $ \T_\lambda =  (\T_{\lambda,s})_{s\in\R}$  by continuity in such a way that almost surely, for any $s, c_1,c_2\in\R$, the continuous function $f=\T_{\lambda,s}(c_1,c_2)$  satisfies \eqref{SA2}. 
%%and its primitive satisfies the stochastic Airy equation \eqref{sol2}. 
%%Moreover, it is clear that the maps $\T_s$ all have the same image. 
%\end{proof}

%\begin{remark}
%The proof also shows that the mapping $s\in\R\mapsto \T_{\lambda,s} $ is continuous in the sense that if $c_1,c_2 :\R\to\C$ are continuous functions, then for any $T>0$, 
%\[
%\sup_{}
%\]
%\end{remark}

Since the equation \eqref{SA2} is linear, we can represent a general solution as a linear combination of two special solutions. 
We use the following basis which exist almost surely.

\begin{definition}[Dirichlet and Neumann solutions] \label{def:DNS}
Let us define for $\lambda\in\C$ and $s\in\R$, 
\[
\big( \mathrm{f}_{\lambda,s}(t) , \mathrm{g}_{\lambda,s}(t) \big) = \bigg( 1 + \int_{s}^t \T(1,0,\lambda,s,u) \d u , 
  \int_{s}^t \T(0,1,\lambda,s,u) \d u  \bigg)_{t\in\R} . 
\]
\end{definition}

Let us observe that we have seen in the proof of Proposition~ \ref{prop:Top}, c.f. formula \eqref{Dir'}, that
  \[
\mathrm{f}_{\lambda,u}'(t)
    =\T(1,0,\lambda,u,t)
= A_{\lambda}(t,u) . 
  \]

The following classical lemma shows that these form a fundamental system for solutions of the stochastic Airy equation. 

\begin{lemma} \label{lem:Wronskian}
Almost surely for all $\lambda \in \C$ and $s \in \R,$
the Wronskian
  $ \W = \mathrm{f}_{\lambda,s} \mathrm{g}_{\lambda,s}' - \mathrm{f}_{\lambda,s}' \mathrm{g}_{\lambda,s} = 1.$   
Moreover, if 
\[
  \phi_\lambda(t) = c_1 + \int_{\varsigma}^{t} \T(c_1,c_2,\lambda,\varsigma, u) \d u ,  \qquad t\in\R
\]
for $\varsigma \in\R, \lambda, c_1, c_2 \in\C$, then it holds for any $s\in\R$, 
\begin{equation} \label{lcs} 
  \phi_\lambda(t) =   \phi_\lambda(s) \mathrm{f}_{\lambda,s}(t) +  \phi_\lambda'(s) \mathrm{g}_{\lambda,s}(t)
  \qquad  t\in\R. 
\end{equation}
\end{lemma}

\begin{proof}
  For fixed $s \in \R$,  let $\big( \mathrm{f} , \mathrm{g}\big) = \big( \mathrm{f}_{\lambda,s} , \mathrm{g}_{\lambda,s} \big)$. 
By Proposition \ref{prop:Top}, both $\mathrm{f}$ and $\mathrm{g}$ are (adapted) solutions of the stochastic Airy equation \eqref{SA1} so that by It\^o formula for any $\lambda \in \C$ and any $t>s$, 
\[\begin{aligned}
\d \W & = \mathrm{f} \d\mathrm{g}' -  \mathrm{g}\d\mathrm{f}' \\
& =  \mathrm{f}(t) \mathrm{g}(t) (t \d t + \tfrac{2}{\sqrt{\beta}}\d B)   - \mathrm{g}(t)\mathrm{f}(t)(t \d t + \tfrac{2}{\sqrt{\beta}}\d B) =0 .
\end{aligned}
\]
The same computation for $t<s$ shows that for all $t\in\R$, $\W(t)= \W(s) = 1$ almost surely.  
As $\W$ is also a continuous function of $(s,\lambda),$ it follows that almost surely it is identically $1$ by applying the above argument to a countable dense set of $(s,\lambda)\in \R\times\C.$

This implies that $( \mathrm{f}, \mathrm{g})$ is a basis of solutions of the stochastic Airy equation. 
Moreover with $  \phi_\lambda$ as in the statement of this Lemma,  by Lemma \ref{lem:Top}, we have 
\[
  \phi_\lambda'(t) 
  =\T(c_1,c_2,{\lambda, \varsigma, t})
  =\T(   \phi_\lambda(s)  ,   \phi_\lambda'(s)  ,\lambda,s,t) . 
\]
By linearity of $\T$ in its first two variables, if we integrate this equation, we obtain \eqref{lcs} and this completes the proof. 
\end{proof}

We conclude by observing that the solution map has a type of invariance property in law.
\begin{proposition}\label{prop:shift}
  The operator $\T$ has an invariance in law: for any $\varsigma \in \R$
  \[
    \left( 
    \T(c_1,c_2,\lambda,t,s) : c_1,c_2,\lambda \in \C, t,s \in \R
    \right)
    \lawequals
    \left( 
    \T(c_1,c_2,\lambda-\varsigma,t+\varsigma,s+\varsigma) : c_1,c_2,\lambda \in \C, t,s \in \R
    \right).
  \]
\end{proposition}
\begin{proof}
  Let $\varsigma \in \R$ be fixed, let $\widehat{B}(t) = B(t+\varsigma) - B(\varsigma)$
  and
    \[
    \widehat{\Es_\lambda}(t,s) = 
    \tfrac{t^2-s^2}{2} + \tfrac{2}{\sqrt{\beta}}\big(\widehat{B}(t) - \widehat{B}(s)\big) + \lambda(t-s)  
    , \qquad \lambda \in \C, s, t \in\R.
  \]
 Then by \eqref{kernelU}, we verify that 
  \(
    \widehat{\Es_\lambda}(t,s)
    =\Es_{\lambda - \varsigma}( t+\varsigma,s+\varsigma).
  \)
Since the solution map $\widehat\T$ defined in terms of $\widehat{\Es_\lambda}$ has the same law as $\T$, we draw the desired conclusion.
\end{proof}

% ===== end sections/04-sae-solution-operator.tex =====

% ===== begin sections/05-stochastic-airy-kernel.tex =====
\section{The stochastic Airy kernel} \label{sec:SAK}

\paragraph{Definition and basic properties.}
In this section, for $\lambda\in\C$, we define a (random) function $\mathfrak{A}_\lambda \in C^1(\R\times\R \to \C)$ that we call the \emph{stochastic Airy kernel} and we can use to represent solutions of the \emph{stochastic Airy equation} with a \emph{forcing term}.

\begin{definition}[Stochastic Airy kernel] \label{def:SAK}
We define the kernel
\[
\mathfrak{A}_\lambda(t,u)   : = \mathrm{f}_{\lambda,0}(t) \mathrm{g}_{\lambda,0}(u)  - \mathrm{f}_{\lambda,0}(u) \mathrm{g}_{\lambda,0}(t)  , \qquad \lambda \in \C, u,t\in\R
\]
%We denote $\mathfrak{A}^\lambda_ t   = \frac{\d}{\d t} \mathfrak{A}^\lambda$, $\mathfrak{A}^\lambda_ s   = \frac{\d}{\d s} \mathfrak{A}^\lambda$ 
%and $\mathfrak{A}^\lambda_{t,u}   = \frac{\d^2}{\d t \d u} \mathfrak{A}^\lambda$. 
This kernel is (almost surely) well-defined and $C^1$ on $\C \times \R \times \R,$ and it is entire in $\lambda \in \C.$
 \end{definition}

%Note that the kernel $\mathfrak{A}_\lambda$ is anti-symmetric, we can assume it is defined for $t\ge u$.

It turns out that the {stochastic Airy kernel} is an intrinsic object which has several other representations, e.g. in terms of the Dirichlet and Neumann solutions from Definition~\ref{def:DNS}. 
This Lemma also reveals some interesting relationships between Dirichlet and Neumann solutions.

\begin{lemma} \label{lem:A}
For any $s\in\R$, we have $\mathfrak{A}_\lambda(t,u)   = \mathrm{f}_{\lambda,s}(t) \mathrm{g}_{\lambda,s}(u)  - \mathrm{f}_{\lambda,s}(u) \mathrm{g}_{\lambda,s}(t)$ for $u,t \in \R$ and $\lambda\in\C$. 
This implies that for all $u,t \in \R$,
\begin{equation} \label{ker1}
\mathfrak{A}_\lambda(t,u)= \mathrm{g}_{\lambda,t}(u) = - \mathrm{g}_{\lambda,u}(t) ,  \qquad
\partial_u\mathfrak{A}_\lambda(t,u)   =  \mathrm{f}_{\lambda,u}(t) 
\quad\text{and}\quad
\partial_t\mathfrak{A}_\lambda(t,u)   =  -\mathrm{f}_{\lambda,t}(u) . 
\end{equation}
Consequently, we also have for all  $u,t \in \R$,
\begin{equation} \label{ker2}
\partial_{tu}\mathfrak{A}_\lambda(t,u) = \mathrm{f}_{\lambda,u}'(t)   =  -\mathrm{f}_{\lambda,t}'(u) . 
\end{equation} 
\end{lemma}

\begin{proof}
By Lemma~\ref{lem:Wronskian}, we have for $u,t\in\R$
\[
\begin{pmatrix} 
\mathrm{f}_{\lambda,s}(t) & \mathrm{g}_{\lambda,s}(t) \\
 \mathrm{f}_{\lambda,s}(u) &  \mathrm{g}_{\lambda,s}(u)
\end{pmatrix}
= \begin{pmatrix} 
\mathrm{f}_{\lambda,0}(t) & \mathrm{g}_{\lambda,0}(t) \\
 \mathrm{f}_{\lambda,0}(u) &  \mathrm{g}_{\lambda,0}(u)
\end{pmatrix} 
\begin{pmatrix} 
\mathrm{f}_{\lambda,s}(0) & \mathrm{g}_{\lambda,s}(0) \\
 \mathrm{f}_{\lambda,s}'(0) &  \mathrm{g}_{\lambda,s}'(0)
\end{pmatrix} .
\]
By taking determinants,  this implies that 
\[
 \mathrm{f}_{\lambda,s}(t) \mathrm{g}_{\lambda,s}(u)  - \mathrm{f}_{\lambda,s}(u) \mathrm{g}_{\lambda,s}(t) = \W\, \mathfrak{A}_\lambda(t,u) . 
\]
Since  $\W =1$, this proves the first claim. 
By Definition~\ref{def:DNS}: $\mathrm{g}_{\lambda,u}(u) = \mathrm{f}_{\lambda,u}'(u) = 0 $ and $\mathrm{f}_{\lambda,u}(u)=\mathrm{g}_{\lambda,u}'(u)=1$  for any $u\in\R$ so that 
\[
\mathfrak{A}_\lambda(t,u)   = \big( \mathrm{f}_{\lambda,s}(t) \mathrm{g}_{\lambda,s}(u)  - \mathrm{f}_{\lambda,s}(u) \mathrm{g}_{\lambda,s}(t) \big)\big|_{s=u} = - \mathrm{g}_{\lambda,u}(t)  
\]
and 
\[
\partial_u\mathfrak{A}_\lambda(t,u)   = \big( \mathrm{f}_{\lambda,s}(t) \mathrm{g}_{\lambda,s}'(u)  - \mathrm{f}_{\lambda,s}'(u) \mathrm{g}_{\lambda,s}(t) \big)\big|_{s=u} =  \mathrm{f}_{\lambda,u}(t)  . 
\]
By a similar argument, we find that $\mathfrak{A}_\lambda(t,u)= \mathrm{g}_{\lambda,t}(u)$, $\partial_t\mathfrak{A}_\lambda(t,u)   =  -\mathrm{f}_{\lambda,t}(u)$ and $\partial_{tu}\mathfrak{A}_\lambda(t,u) = \mathrm{f}_{\lambda,u}'(t)  =  -\mathrm{f}_{\lambda,t}'(u)$. 
\end{proof}

\begin{remark} \label{rk:Wronskian}
The previous argument shows that we can write the kernel
$\mathfrak{A}_\lambda(t,u)  = \chi(t)  \phi(u)  - \chi(u) \phi(t) $
for any two solutions $ \phi,\chi$ of the  stochastic Airy equation \eqref{SA2} with fixed initial data at $s\in\R$ whose Wronskian is normalized to 1.
\end{remark}

\paragraph{Inhomogeneous stochastic Airy equation.}
Our next  proposition allows us to express the solution of the stochastic Airy equation $\d \phi' = \phi \d \Es_\lambda(t,s) + \d \zeta$ with forcing $\zeta$ (written as an integral equation \eqref{sol3}) in terms of the stochastic Airy kernel defined above. 
This is instrumental to obtain our coupling between the characteristic polynomial of the Gaussian $\beta$-ensemble and a solution of the stochastic Airy equation. 
The following proposition holds almost surely. 
 %(using the superposition principle together with Proposition~\ref{prop:Top}). 

\begin{proposition} \label{prop:sol}
Let $\lambda\in\C$ and $\zeta:\R\to\R$ be a c\`adl\`ag function with $\zeta(s)=0$ for a fixed $s\in\R$ and suppose that $h$  solves the equation 
\begin{equation} \label{sol3}
h(t) = \int_s^t \Es_\lambda(t,v) h(v)\d v  + \zeta(t)  , \qquad t\in \R . 
\end{equation}
Then, it holds for any $t\in\R$, 
\begin{equation}\label{sol5}
h(t)  =  \zeta(t)+  \int_{s}^t \partial_{tu}\mathfrak{A}_\lambda(t,u) \zeta(u) \d u . 
\end{equation}
Conversely, if $h$ is given by \eqref{sol5}, then it also solves the equation \eqref{sol3}.
\end{proposition}

\begin{proof}
Let $\displaystyle h(t)= - \int_{s}^t  \partial_t\mathfrak{A}_\lambda(t,u)  \d \zeta(u);$ this function is well-defined since $\zeta$ is a c\`adl\`ag process and for a fixed $t\in\R$, $ u\mapsto  \partial_t\mathfrak{A}_\lambda(t,u)$ is a $C^1(\R)$ function. Moreover, by an integration by parts,
\[
h(t) = \zeta(t) + \int_{s}^t \partial_{tu}\mathfrak{A}_\lambda(t,u)\zeta(u) \d u
\]
where we used that $\zeta(s)=0$ and $\partial_t\mathfrak{A}_\lambda(t,u) |_{u=t} = -\mathrm{f}_{\lambda,t}(t) =-1$ for any $t\in\R$; see \eqref{ker1}.  
This shows that the process $h$ defined above is given by  \eqref{sol5}.

Now, let us check that $h$ also solves the equation \eqref{sol3}. 
First note that by \eqref{ker1}, we also have  $\partial_t\mathfrak{A}_\lambda (t,s) = - \mathrm{g}_{\lambda,s}'(t)$ for $s,t\in\R$ where $\mathrm{g}_{\lambda,s}'$ is a solution of the equation \eqref{SA2} with $c_1=0$ and $c_2=1$. This implies that for all $u,t\in\R$, 
\[
  \int_{s}^t \partial_v\mathfrak{A}_\lambda(v,s)  \Es_\lambda(t,v)\d v
  = - \int_{s}^t \mathrm{g}_{\lambda,s}'(v)  \Es_\lambda(t,v)\d v
  =  1-\mathrm{g}_{\lambda,s}'(t)
  = \partial_t\mathfrak{A}_\lambda(t,s)  + 1 . 
\]
By Fubini's theorem, this shows that
\[
  \begin{aligned}
    \int_s^t \Es_\lambda(t,v) h(v) \d v 
    & = -  \int_{s}^t  \int_{u}^t \Es_\lambda(t,v) \partial_v\mathfrak{A}_\lambda(v,u) \d v \d \zeta(u) \\
    & =  -\int_{s}^t  \big(\partial_t\mathfrak{A}_\lambda(t,u) + 1 \big)  \d \zeta(u)  \\
    & = h(t)- \zeta(t). 
  \end{aligned}
\]
For the last step, we used again that $\zeta(s)=0$.
The converse statement follows from the fact that the equation \eqref{sol3} has a unique solution. 
Indeed, if $g$ is the difference of two solutions for \eqref{sol3}, then by linearity, $g$ satisfies
$\displaystyle g(t) = \int_s^t \Es_\lambda(t,v) g(v)\d v$.
This equation is of the form \eqref{SA2} with $c_1=c_2 = 0$. Hence, by Proposition~\ref{prop:Top}, we conclude that $g=0$ is the unique such solution. 
\end{proof}

\paragraph{Stability of solutions to the stochastic Airy equation}

As we shall see in Section \ref{sec:fde}, the 3-term recurrence \eqref{eq:recurrence} \emph{approximately} solves \eqref{SA2}.  From here we would like to conclude that $\left\{ \varphi_n \right\}$ is well approximated by a given solution of the stochastic Airy equation. 
Our main technical result which enables us to draw such a conclusion is a bound for the kernel $\partial_{tu}\mathfrak{A}_\lambda(t,u) =  \mathrm{f}_{\lambda,u}'(t) $ which holds with overwhelming probability. 
For $\lambda\in\C$ and $t,s\in\R$, let
\begin{equation}\label{eq:Et}
E_\lambda(t,s)=
    e^{\frac23\big((\Re \lambda + s \vee t)_+^{3/2}- (\Re \lambda + s \wedge t)_+^{3/2}\big)} .
%    \quad
%    \text{where}
%    \quad
%    t_+ = t\1_{t>0}.
\end{equation}

% \paragraph{An a priori parameter $T$.}
% The stability statements that follow, as well as the corresponding estimates of Sections~\ref{sec:SAi}, \ref{sec:SAestimate}, and \ref{sec:Red5}, are formulated in terms of an a priori large parameter $T \ge 1$ together with a fixed $\nu \in (1, 3/2)$, in place of the matrix size $N$.  When applied in Section~\ref{sec:dtoc} to the characteristic polynomial $\varphi_n$ of the Gaussian $\beta$-ensemble, we set
% \[
% \log N \coloneqq T^\nu ,
% \]
% equivalently $T = (\log N)^{1/\nu}$. In particular, polynomial-in-$N$ rates from \cite{LambertPaquette02} translate as $N^c = e^{c\, T^\nu}$, and $(\log N)^{1+\eta}$ becomes $T^{(1+\eta)\nu}$ — an exponent strictly greater than $1$ for $\eta > 0$, so the resulting probability tails beat any negative power of $N$.

\begin{theorem}\label{thm:Kest}
Fix $\nu \in (1, 3/2)$.  For $T \ge 1$, set
\[
J(T) \coloneqq [-e^{T} , T]
\quad\text{and}\quad
K_\nu(T) \coloneqq \big\{\lambda \in \C : \Re \lambda \in J(T) ,\, |\Im \lambda| \leq e^{-(1-1/\nu)\, T^\nu} \big\} .
\]
For any $\delta>0$, there exist constants $C=C(\beta, \delta, \nu)$ and $c=c(\beta, \delta, \nu)>0$ so that for  $T \gg 1$,
\[
\P\biggl[ \sup\big\{ | \partial_{tu}\mathfrak{A}_\lambda(t,u) | E_\lambda^{-1}(t,u) :
  t,u \in J(T) ,\, \lambda \in K_\nu(T) \big\}   > C\, e^{\delta\, T^\nu}\biggr] \le  e^{-c\, T^{(3\nu-1)/2}} .
\]
\end{theorem}
\noindent We give the proof in Section \ref{sec:SAestimate}.

\begin{remark}
To understand the bounds from Theorem~\ref{thm:Kest}, one can compare to the (deterministic) Airy kernel for the usual Airy equation which is given by 
\[
\widehat{\mathfrak{A}_\lambda}(t,u) = \pi\big( \Ai(t+\lambda)\Bi(u+\lambda)-  \Ai(u+\lambda)\Bi(t+\lambda) \big) . 
\]
In this case, using the well-known asymptotics of $\Ai'$ and $\Bi'$, see e.g.~\cite[Section 9.7]{DLMF}, we verify that for $t\ge u$, 
\[
\partial_{tu}\widehat{\mathfrak{A}_\lambda}(t,u) \simeq  -\pi  \Ai'(u+\lambda)\Bi'(t+\lambda)
\simeq (u+\lambda)^{1/4} (t+\lambda)^{1/4} E_\lambda(t,u)/2 . 
\]
\end{remark}

Using the previous bounds, we formulate a statement which shows that solutions of \eqref{SA2} are stable under perturbation, that is, suppose $h_\lambda$ solves the equation
\begin{equation}
 h_\lambda(t) = c_2 + \int_s^t (\Es_\lambda(t,u) + \mathcal{D}_{\lambda}^1(t,u)) h_\lambda(u) \d u + c_1 \Es_\lambda(t,s) + \mathcal{D}_{\lambda}^2(t),
  \label{eq:SAperturb}
\end{equation}
for all $t\in I$, a closed interval containing $s.$ 
If $\mathcal{D}^1$ and $\mathcal{D}^2$ are small, the following proposition allows us to conclude that $h_\lambda$ is close to the solution $\phi_\lambda$ of the same equation with $\mathcal{D}^1=\mathcal{D}^2 = 0$.

\begin{proposition} \label{prop:stability}
Let  $1<\nu <3/2$, $J(T)$ and $K_\nu(T)$ be as in Theorem~\ref{thm:Kest}. Assume that $h_\lambda(t)$ solves \eqref{eq:SAperturb} for $(\lambda,s) \in K_\nu  \times J $ with $|c_1|+|c_2|\le1$ and that for $\alpha>0$:
  \begin{equation}    \label{eq:D_stab}
 \sup_{u \in [s,t]}   |\mathcal{D}_\lambda^1(t,u)| \leq e^{-\alpha\, T^\nu}
    \quad
    \text{and}
    \quad
    |\mathcal{D}_\lambda^2(t)| \leq e^{-\alpha\, T^\nu} E_\lambda(t,s) ,
  \end{equation}
  for all $(\lambda,t) \in K_\nu \times J $.
Let $\Phi_\lambda$ solves the stochastic Airy equation \eqref{SA2}.  There are constants
$C=C(\beta, \alpha, \nu)$ and $c=c(\beta, \alpha, \nu)>0$ and $0<\alpha_- <\alpha$ so that for $T \gg 1$, with probability at least $1 - e^{-c\, T^{(3\nu-1)/2}}$: 
  \[
    \sup_{t \in J ,\, \lambda \in K_\nu}\biggl(|h_\lambda(t)-\Phi_\lambda(t)| E_\lambda^{-1}(t,s)\biggr) \leq C e^{-\alpha_-\, T^\nu} .
  \]
\end{proposition}
\begin{proof}
We work on the event of Theorem~\ref{thm:Kest}: for some small $\delta>0$ to be chosen and for all $t,u \in J(T)$ and $\lambda \in K_\nu(T)$,
  \begin{equation}    \label{eq:stab0}
    |\partial_{tu}\mathfrak{A}_\lambda(t,u)| \leq e^{\delta\, T^\nu} E_\lambda(t,u) .
  \end{equation}
By Lemma~\ref{lem:A}, we can write both $\mathrm{f}_{\lambda,s}'(t)$ and $\mathrm{g}_{\lambda,s}'(t)$ in terms of the kernel $\partial_{tu}\mathfrak{A}_\lambda(t,u)$. We obtain the crude bounds, for $t, s\in J$ and $\lambda\in K_\nu$:
\begin{equation} \label{fbound}
|\mathrm{f}_{\lambda,s}'(t)| \le e^{\delta\, T^\nu} E_\lambda(t,s) , \qquad
|\mathrm{f}_{\lambda,s}(t)| = |\mathrm{g}_{\lambda,s}'(t)| 
\le 1+ \int_{[t,s]}  |\partial_{ts}\mathfrak{A}_\lambda(u,s)| \d u \lesssim  e^{T + \delta\, T^\nu} E_\lambda(t,s) . 
\end{equation}
% and
% \[
% |\mathrm{g}_{\lambda,s}(t)| \le  \int_{[t,s]}  |\mathrm{f}_{\lambda,s}(t)| \d s \lesssim  e^{2T +\delta\, T^\nu} E_\lambda(t,s) .
% \]
We may represent any solution of \eqref{SA2} as 
\( \Phi_\lambda = c_1 \mathrm{f}_{\lambda,s}' + c_2 \mathrm{g}_{\lambda,s}'\).
If $|c_1|+|c_2|\lesssim 1$, this implies that for all $t\in J$ and $\lambda \in K_\nu$:
  \begin{equation}    \label{eq:stab1}
|\Phi_{\lambda}(t)| \le  e^{2\delta\, T^\nu} E_\lambda(t,s)
\end{equation}
if $T\gg1$ and $\nu >1$.

Now let us denote $\xi_\lambda := \Phi_\lambda - h_\lambda$.  By  \eqref{eq:SAperturb}, this function satisfies
  \begin{equation*}
    \xi_\lambda(t)
    =
    \int_s^t \Es_\lambda(t,u) \xi_\lambda(u) \d u
    + \zeta_\lambda(t)
    \qquad\text{where}\quad
	    \zeta_\lambda(t) :=\int_s^t \mathcal{D}_\lambda^1(t,u) (\underbrace{\xi_\lambda(u)-\Phi_\lambda(u)}_{-h_\lambda(u)}) \d u    -\mathcal{D}_\lambda^2(t).
  \end{equation*}
Hence, by Proposition~\ref{prop:sol}, we conclude
\begin{equation}\label{xirep}
    \xi_\lambda(t)
    =  \zeta_\lambda (t)+  \int_{s}^t \partial_{tu}\mathfrak{A}_\lambda(t,u) \zeta_\lambda(u) \d u .
\end{equation}

Let 
\(
  \mathcal{M} := \sup_{t \in J, \lambda \in K_\nu} \bigl(|\xi_\lambda(t)| E^{-1}_\lambda(t,s)\bigr)
\)
for the quantity that we want to estimate.
It follows from the conditions \eqref{eq:D_stab} and \eqref{eq:stab1} that
\[
|\zeta_\lambda(t)| \le e^{-\alpha\, T^\nu} (\mathcal{M}+ e^{2\delta\, T^\nu}) \int_{[s,t]} E_\lambda(u,s) \d u 
\le 2 e^{T-\alpha\, T^\nu} \big(\mathcal{M}+ e^{2\delta\, T^\nu}\big) E_\lambda(t,s) .
\]
Using the bound \eqref{eq:stab0} again in \eqref{xirep}, this implies that
with \(
  \mathcal{M}' := \sup_{t \in J, \lambda \in K_\nu} \bigl(|\zeta_\lambda(t)| E^{-1}_\lambda(t,s)\bigr)
\),
\[ \begin{aligned}
|\xi_\lambda(t)|  & \le  \mathcal{M}'\bigg(  E_\lambda(t,s)  + e^{\delta\, T^\nu} \int_{[s,t]} E_\lambda(t,u)  E_\lambda(u,s) \d u \bigg)  \\
& \le \mathcal{M}'  e^{2\delta\, T^\nu}   E_\lambda(t,s) \\
\mathcal{M} & \le e^{-(\alpha-3\delta)\, T^\nu} \big(\mathcal{M}+ e^{2\delta\, T^\nu}\big) = e^{- \alpha_- T^\nu} \big(1+\mathcal{M}e^{-2\delta\, T^\nu}\big)
\end{aligned}\]
with $\alpha_-=\alpha-5\delta>0$ ($\delta>0$ is arbitrary small).
This implies that on the event of Theorem~\ref{thm:Kest}, $  \mathcal{M} \lesssim e^{- \alpha_- T^\nu}$. This concludes the proof. 
\end{proof}

Let us also record the following consequence of Theorem~\ref{thm:Kest}. 

\begin{corollary} \label{cor:solbound}
Fix $1<\nu < 3/2$ and a small $\delta>0$.  
There exist constants $C=C(\beta, \nu, \delta)$ and $c=c(\beta, \nu, \delta)>0$ so that with probability at least $1 - e^{-c\, T^{(3\nu-1)/2}}$,
every solution $\phi_\lambda$ of the stochastic Airy equation \eqref{SA2} and $s \in J(T)$ with $|c_1|+|c_2|\le1$  satisfies  
\[
\sup_{t \in J ,\, \lambda \in K_\nu}\big|\phi_\lambda'(t)  E_\lambda^{-1}(t,s)\big| \leq C\, e^{ \delta T^\nu} 
\]
and a similar bound holds for $\phi_\lambda$. 
\end{corollary}

\begin{proof}
The first bound corresponds to \eqref{eq:stab1} --
$\phi_\lambda' =\Phi_{\lambda}$ and it holds on the event from  Theorem~\ref{thm:Kest}. A similar bound for $\phi_\lambda$ follows by integrating this equation and adjusting $\delta>0$ and the constants.
\end{proof}

%{\color{magenta}
%\paragraph{The stochastic Airy kernel and its relation to Sturm--Liouville operator.}
%
%\begin{remark}
%Fix $\lambda\in\C$ and let us define the operator $\displaystyle \mathscr{A}_\lambda : \zeta \in C_0(\R_-) \mapsto  -\int_{0}^t  \mathfrak{A}_\lambda(t,u) \zeta(u)  \d u $. 
%Using that for any $u\in\R_-$,  $-\mathfrak{A}^\lambda(\cdot,u) = \mathrm{g}_{\lambda,u} \in C_1(\R)$ with  $\mathrm{g}_{\lambda,u}(u) =0 , \mathrm{g}_{\lambda,u}'(u)=1$ and $\d \mathrm{g}_{\lambda,u}' = - \mathrm{g}_{\lambda,u}(t-\lambda)\d t - \mathrm{g}_{\lambda,u} \d B$,  we obtain the stochastic differential equation
%\[
%\d( \mathscr{A}_\lambda \zeta )' = \zeta  -(t-\lambda)(\mathscr{A}_\lambda\zeta) \d t - ( \mathscr{A}_\lambda \zeta)\d B
%\]
%This computation shows that formally if $\mathscr{L}  = -\frac{\d^2}{\d t^2} - t - \d B $, then $(\mathscr{L}-\lambda \mathrm{Id})\mathscr{A}_\lambda = - \mathrm{Id}$.
%
%\end{remark}
%
%
%Regarding the previous remark, in the deterministic case, we have 
%\[
%\mathfrak{A}_\lambda(t,u) =  \Ai(t-\lambda)\Bi(u-\lambda) -  \Ai(u-\lambda)\Bi(t-\lambda)
%\]
%which cannot be equal to the resolvent kernel
%\[
%\mathscr{R}_\lambda(t,u) = \sum_{k\in\N}  \frac{\Ai(t-\lambda_k)\Ai(u-\lambda_k)}{\lambda-\lambda_k}
%\]
%where $(\lambda_k)_{k\in\N}$ are the zeros of $\Ai$. \\
%}

% ===== end sections/05-stochastic-airy-kernel.tex =====

% ===== begin sections/06-asymptotics-and-airy-function.tex =====
\section{Asymptotics of solutions and construction of the stochastic Airy function}\label{sec:SAi}

In this section, we construct the stochastic Airy function and give its asymptotic behavior as $t \to \infty$.  We also discuss the stability of the initial value problem for the  stochastic Airy equation at large time by showing that the stochastic Airy function is the unique stable solution. 

We recall that $\displaystyle \mathfrak{X}(t) = \int_0^t e^{-\frac43(t^{3/2}-u^{3/2})}  \d B(u)$ where $B$ is a standard Brownian motion and we define for  $\lambda\in\C$, $t\in\R_+$, 
\begin{equation} \label{def:theta}
  \theta_\lambda(t) \coloneqq \sqrt{t} + \frac{\lambda}{2\sqrt{t}} - \frac{1}{4(t+1)} +\tfrac{2}{\sqrt{\beta}}\mathfrak{X} (t)
  - \frac{4}{\beta} 
  \int_{0}^t 
  \int_{0}^s 
  e^{-\frac43(t^{3/2}-u^{3/2}) -\frac43(s^{3/2}-u^{3/2})}  
  \d u \d s . 
\end{equation}
We use this function to describe the almost--sure and high--probability behavior of Dirichlet and Neumann solutions of \eqref{SA2}.
The last term in \eqref{def:theta} can be understood in terms of the covariance structure of the  (centered) Gaussian process $\mathfrak{X}:$ 
for $t\ge s \ge 0$, 
\begin{equation}\label{eq:covX1}
  \E[ \mathfrak{X}(t) \mathfrak{X}(s)] = \int_{0}^s e^{-\frac43(t^{3/2}-u^{3/2}) -\frac43(s^{3/2}-u^{3/2})}  \d u = \tfrac14 e^{-\frac43(t^{3/2}-s^{3/2})}\biggl(\frac{1}{\sqrt{s}} +\underset{s\to\infty}{\mathcal{O}}(s^{-3/2}) \biggr). 
\end{equation}
In particular, this implies that the covariance integral occurring in the last term on the RHS of \eqref{def:theta} equals
\begin{equation}\label{eq:covX2}
  \frac12 \partial_t \E\bigg[ \bigg( \int_0^t  \mathfrak{X}(s) \d s \bigg)^2\bigg] = \int_0^t  \E[ \mathfrak{X}(t) \mathfrak{X}(s)]  \d s 
  = \frac{1}{8t} +\underset{t\to\infty}{\mathcal{O}}(t^{-2}) .
\end{equation}

Our strategy to construct the stochastic Airy function is to take a locally uniform limit of the Dirichlet~$\mathrm{f}$ and/or Neumann $\mathrm{g}$ solutions to the stochastic Airy equation.

\begin{definition} \label{def:SAiprime}
  For all $\lambda\in\C$ and $s\in\R$
  \[
    \lim_ {t\to \infty} \mathrm{f}_{\lambda,s}(t)  \exp\biggl( -\int_{0}^{t}\theta_{\lambda}(u)\d u \biggr) 
    \eqqcolon - \sqrt{\pi}\SAi_\lambda'(s).
  \]
\end{definition}
If this limit exists locally uniformly for $\lambda\in\C$, the resulting function is an entire function in $\lambda$, being itself a locally uniform limit of entire functions.
It also turns out that $s\mapsto \SAi_\lambda'(s)$ is continuous on $\R$ and in $L^2(\R_+)$.  We can then define the stochastic Airy function and stochastic BAiry function by taking:

\begin{definition} \label{def:SAi}
  For $\lambda\in\C$ and $s\in\R$, we define the {\bf stochastic Airy function}:
\[
\SAi_\lambda(s) = - \int_s^\infty \SAi_\lambda'(u) \d u
=
\lim_ {t\to \infty} \mathrm{g}_{\lambda,s}(t) \exp\biggl( -\int_{0}^{t}\theta_{\lambda}(u)\d u \biggr) 
/\sqrt{\pi}. 
\]
The function $(t,\lambda) \mapsto \SAi_\lambda(t)$ is  $C^1(\R\times\C)$, entire for $\lambda\in\C$ and $\overline{\SAi_\lambda} = \SAi_{\overline\lambda}$. 
\end{definition}

%As defined, these functions also depend on the parameter $\beta>0$ and the pair $\left\{\SAi_\lambda, \SBi_\lambda \right\}$ are a fundamental set of solutions of the stochastic Airy equation with Wronskian 
%equals to $1/\pi$.

\begin{remark}\label{rk:SBi} (Stochastic BAiry)
  One can make a pair of fundamental solutions to the stochastic Airy equation by choosing another solution $f$ which has nonzero Wronskian with $\SAi_\lambda$. We note that essentially any solution we pick will have BAiry type asymptotics as $t \to \infty.$ 
  We do not know of a canonical choice, but it does not play an important role in our approach, in any case.

  For the (non-stochastic) Airy equation, there is a canonical second solution $\Bi,$ the BAiry solution. This solution is a real--valued (see \cite[(2.4)]{Vallee}) combination of the Airy equation solutions 
 \[
   \left\{ \Ai( e^{2\pi i/3} x), \Ai( e^{4\pi i/3}x) \right\}.
 \]
% or due to the Sturm--Liouville properties of the pair $\left( \Ai(-x),\Bi(-x) : x \leq 0 \right)$ (see \cite[2.2.5]{Vallee}).  
 Because of a lack of analytic structure in $t$, we do not know a way to generalize this solution to the stochastic Airy equation for positive $\beta >0$.  One other possible candidate for a second solution would generalize the Sturm--Liouville properties of the pair $\left( \Ai(-x),\Bi(-x) : x \leq 0 \right)$ (see \cite[2.2.5]{Vallee}).
\end{remark}

In fact we provide quantitative asymptotics for $\big\{ \mathrm{f}_{s}(t),\mathrm{g}_{s}(t)\big\}$ when $t \geq s \gg 0.$  A complete description of this behavior is complicated and we refer to Section \ref{sec:asymp} for the relevant statements.  Instead, we present the main corollary of these estimates when applied to the stochastic Airy function:

\begin{proposition}[Existence of the stochastic Airy function and its asymptotics] \label{prop:asympDN}
  Definitions \ref{def:SAiprime} and \ref{def:SAi} are well--posed, and the limits therein exist almost surely and locally uniformly in $\lambda \in \C.$ 
  Furthermore, for any compact $K \subset \C$ and any $\epsilon \in [0, \tfrac 12],$  there is a constant $C_{K,\epsilon}$ so that for $s\ge 1$, 
\begin{equation} \label{SAiasymp}
    \begin{aligned}
   & \Pr\biggl(
    \sup_{\lambda \in K}
    \sup_{t \geq s}
    \biggl|2\sqrt{\pi}\exp\biggl(\int_{0}^{t}\theta_{\lambda}(u)\d u \biggr)\SAi'_{\lambda}(t) + 1\biggr| \geq s^{-\epsilon} 
    \biggr)
    \leq C_{K,\epsilon,\beta} \exp(- s^{3/2-3\epsilon}) ,\\
    &\Pr\biggl(
    \sup_{\lambda \in K}
    \sup_{t \geq s}
    \biggl|2\sqrt{\pi t}\exp\biggl(\int_{0}^{t}\theta_{\lambda}(u)\d u \biggr)\SAi_{\lambda}(t)- 1\biggr| \geq s^{-\epsilon}
    \biggr)
    \leq C_{K,\epsilon,\beta} \exp(-s^{3/2-3\epsilon}) .
    \end{aligned}
  \end{equation}
In particular, we have the almost sure asymptotics
  \begin{equation} \label{asympAB}
    \begin{aligned}
      &\lim_ {s\to \infty} \biggl\{ \SAi_\lambda(s) \sqrt{s},\SAi_\lambda'(s) \biggr\}\exp\biggl(\int_{0}^{s}\theta_{\lambda}(u)\d u \biggr)  = \biggl\{\frac{1}{2\sqrt{\pi}},-\frac{1}{2\sqrt{\pi}}\biggr\}, 
%      &\lim_{t\to\infty}
%      \biggl\{
%	\SBi_\lambda(t),
%	\SBi'_\lambda(t)/\sqrt{t}
%      \biggr\}
%      \exp\biggl( -\int_{0}^{t}\theta_{\lambda}(u)\d u \biggr)
%      =
%      \left\{ \frac{1}{\sqrt{\pi}}, \frac{1}{\sqrt{\pi}} \right\},
    \end{aligned} 
  \end{equation}
  from which it follows that, for each fixed $\lambda\in\C$, the space of solutions of the stochastic Airy equation \eqref{SA2} that belong to $L^p(\R_+)$ for every $p\in[1,\infty]$ is one-dimensional and spanned by $\SAi_\lambda$. The asymptotics \eqref{asympAB} fix the displayed scalar normalization of $\SAi_\lambda$.
\end{proposition}

%We do not give a tail bound for $\left\{ \SBi_\lambda(t), \SBi'_\lambda(t) \right\}.$ 
%This would be useful in the representation of solutions to the stochastic Airy equation with initial data at large positive time.  However, it is possible to circumvent the need for this tail bound, in showing that essentially any fixed initial data at large positive time produces a solution which is a multiple of stochastic Airy, at least for nearly real $\lambda.$   

The other main result of this section is Proposition~\ref{prop:SAiwins} (stability of the initial value problem), stated in the introduction. As anticipated there, this result follows from the large-$t$ asymptotics of Proposition~\ref{prop:asympDN} together with the stochastic Airy kernel estimates of Proposition~\ref{prop:stability}. The main steps of the proofs of Propositions~\ref{prop:asympDN} and \ref{prop:SAiwins} are presented in Section~\ref{sec:proofasymp}. We fill in the technical details in Sections~\ref{sec:entrance}--\ref{sec:X2}.

\begin{remark} \label{rk:normalization}
In Proposition~\ref{prop:asympDN}, we made a convention to define $\SAi'$ that we now explain by comparing the asymptotics \eqref{asympAB} with their counterpart for the (deterministic) Airy equation. Since the Wronskian of the Airy and BAiry functions, $\mathcal{W}(\Ai,\Bi)=1/\pi$, the Dirichlet solutions for the Airy equation are given by
\[
\widehat{\mathrm{f}_{\lambda,s}}(t) =\pi \big( \Bi'(\lambda+s)\Ai(\lambda+t) - \Ai'(\lambda+s)\Bi(\lambda+t) \big) , \qquad s,t\in\R , \lambda\in\C . 
\]
%According to \eqref{eq:covX2}, 
Let us define
\[
\widehat{\theta_{\lambda}}(t)  :=  \partial_t \log \E\bigg[ \exp\biggl( \int_{0}^{t}\theta_{\lambda}(u)\d u \biggr) \bigg] =   \sqrt{t} + \frac{\lambda}{2\sqrt{t}} - \frac{1}{4(t+1)}   =  \theta_{\lambda}\big|_{\beta=\infty}(t) , 
\]
cf. \eqref{eq:covX2}.
Using the well-known asymptotics of $\Ai$ and $\Bi$, see e.g.~\cite[Section 9.7]{DLMF}, we verify that
\[
\lim_{t\to\infty}\Ai(\lambda+t)  \exp\biggl( -\int_{0}^{t}\widehat{\theta_{\lambda}}(u)\d u \biggr) = 0
\quad\text{and}\quad
\lim_{t\to\infty}\Bi(\lambda+t)  \exp\biggl( -\int_{0}^{t}\widehat{\theta_{\lambda}}(u)\d u \biggr) = \frac{1}{\sqrt{\pi}} . 
\]
Hence, this shows that for $s\in\R$ and $\lambda\in\C$, 
\[
\widehat{\SAi_\lambda}'(s) : =  -\tfrac{1}{  \sqrt{\pi}}
\lim_{t\to\infty} \widehat{\mathrm{f}_{\lambda,s}}(t) \exp\biggl( -\int_{0}^{t}\widehat{\theta_{\lambda}}(u)\d u \biggr) =  \Ai'(\lambda+s). 
\]
We further verify that this is consistent with \eqref{asympAB} in the sense that
\[
\lim_{s\to\infty} \widehat{\SAi_\lambda}'(s) \exp\biggl(\int_{0}^{s}\widehat{\theta_{\lambda}}(u)\d u \biggr) 
=  \lim_{s\to\infty}\Ai'(\lambda+s)  \exp\biggl(\int_{0}^{s}\widehat{\theta_{\lambda}}(u)\d u \biggr)=  -\frac{1}{2\sqrt{\pi}} .
\]
Note that $ \widehat{\SAi_\lambda}$ is also the unique solution of the (deterministic) Airy equation with such asymptotics. 
In fact, we claim that $\lim_{\beta\to\infty}\SAi_\lambda(t)$ exists and is a solution of the (deterministic) Airy equation which has the same asymptotics as $\widehat{\SAi_\lambda}(t)$ when $t\to\infty$ (we do not justify these claims here). 
Then our convention implies that $\lim_{\beta\to\infty}\SAi_\lambda(t) =  \Ai(\lambda+t)$. 
In fact, one can rewrite the asymptotics from Proposition~\ref{prop:asympDN} as $t\to\infty$,  
\begin{equation} \label{asympDN2}
\begin{aligned}
 \SAi'_\lambda(t)  
   & \simeq  \Ai'(\lambda+t) \exp\biggl( -\frac{2}{\sqrt{\beta}}\int_{0}^{t} \mathfrak{X}(u)\d u \biggr)\E  \exp\biggl(\frac{2}{\sqrt{\beta}}\int_{0}^{t} \mathfrak{X}(u)\d u \biggr)  \\
      \SAi_\lambda(t)  
   & \simeq  \Ai(\lambda+t)  \exp\biggl( -\frac{2}{\sqrt{\beta}}\int_{0}^{t} \mathfrak{X}(u)\d u \biggr) \E  \exp\biggl(\frac{2}{\sqrt{\beta}}\int_{0}^{t} \mathfrak{X}(u)\d u \biggr) 
\end{aligned}  
\end{equation}
locally uniformly  in $\lambda\in \C$. 
\end{remark}

Let us define
\begin{equation} \label{def:cstar}
c_* := \lim_{T\to\infty}\left(
2\int_0^T\int_0^t\int_0^s
e^{-\frac43(t^{3/2}-u^{3/2})-\frac43(s^{3/2}-u^{3/2})}
\,\d u\,\d s\,\d t-\frac14\log T
\right).
\end{equation}
This limit exists by \eqref{eq:covX1}--\eqref{eq:covX2} and according to \eqref{def:theta}, we obtain the asymptotics
\[ \begin{aligned}
\int_{0}^{t}\theta_{\lambda}(u)\d u 
& = \int_{0}^{t} \bigg(  \sqrt{u} + \tfrac{\lambda}{2\sqrt{u}} +  \tfrac{2}{\sqrt{\beta}}\mathfrak{X} (u)  - \tfrac{1/4}{u+1} + \tfrac{1/2\beta}{u+1} \bigg) \d u
+ \frac{2c_*}{\beta} +\underset{t\to\infty}{\mathcal{O}_\beta}(t^{-1}) \\
& = \tfrac{2}{3} (t + \lambda)^{3/2} +\frac{2}{\sqrt{\beta}}\int_0^t \mathfrak{X}(u)\,\d u  + \big( \tfrac{1}{2\beta}   -\tfrac{1}{4}\bigg) \log(t) +\frac{2c_*}{\beta} +\underset{t\to\infty}{\mathcal{O}_{\beta, \lambda}}(t^{-1/2})  
\end{aligned}\]
where the errors are deterministic and locally uniform for $\lambda\in\C$.  
By \eqref{asympAB}, this implies the asymptotics~\eqref{eq:SAIconvention}.

\subsection{Connecting the stochastic Airy function to the stochastic Airy operator}

It turns out that Theorem~\ref{thm:SAO}, the spectral theorem for the stochastic Airy operator, is another consequence of Proposition~\ref{prop:asympDN}.
In the following, we present a result for general boundary conditions.
Let $\mathbf{H}_{\beta, \omega} = -\partial_t^2 + t + \tfrac{2}{\sqrt{\beta}}\d B$ be the operator acting formally on $L^2(\R_+)$ with boundary condition:
$f'(0) = \omega f(0) $ for $\omega \in (-\infty, +\infty]$ and $f(+\infty)=0$. 
Note that in the case of the stochastic Airy operator $\mathbf{H}_\beta$, the Dirichlet boundary condition corresponds naturally to the limiting case $\omega=+\infty$. 
In our setting, the operators $\mathbf{H}_{\beta, \omega}$ are coupled together (with the same driving Brownian motion) and their spectra are understood as in \eqref{SA2}. 
Namely, $\phi_{\lambda,\omega}$ is an eigenfunction for $\mathbf{H}_{\beta, \omega}$ with eigenvalue $-\lambda$ if it satisfies the equation \eqref{SA2} with $s=0$, $c_2 =\omega c_1$ and $\phi_{\lambda,\omega}\in H^1(\R_+)$ (equivalently, $\phi_{\lambda,\omega},\Phi_{\lambda,\omega}\in L^2(\R_+)$). 

\begin{theorem} \label{thm:spec}
Let $\mathcal{A}_\omega =\big\{\lambda\in\C : \SAi_\lambda'(0) =  \omega  \SAi_\lambda(0)  \big\} $. 
The spectrum of the operator $\mathbf{H}_{\beta,\omega}$ is given by 
$(-\lambda, \SAi_\lambda)_{\lambda\in\mathcal{A}_\omega}$ and in particular, 
$\mathcal{A}_\omega \subset \R$ is a countable set, bounded from above with no accumulation point. 
\end{theorem}
 
\begin{proof}
Let us recall from \cite{BloemendalViragI} that the operators $\mathbf{H}_{\beta, \omega}$ are symmetric, bounded from below and have pure point spectrum for any $\omega \in (-\infty, +\infty]$ (these facts are established by using the standard variational characterization for the spectrum). 
By Proposition \ref{prop:asympDN}, we have characterized all solutions of \eqref{SA2} which lie in  $H^1(\R_+)$.
By matching the boundary conditions at 0, this yields an explicit description of the spectrum of $\mathbf{H}_{\beta, \omega}$.
\end{proof}

According to \cite[Proposition 3.4]{RRV09}, it is also possible to describe the spectrum of the stochastic Airy operator $\mathbf{H}_\beta$, and more generally of $\mathbf{H}_{\beta,\omega}$, in terms of the Riccati diffusion \eqref{Ric0}.  For simplicity we restrict to the case of $\omega=+\infty.$ 
Namely, for $E\in\R$, we define the conventional eigenvalue counting function of the negated Airy point process $-\A=\operatorname{spec}(\mathbf H_\beta)$ by
\[
\begin{aligned}
\mathbf{N}^{\mathbf{H}}_{E}
&= \#\big\{\mu\in\operatorname{spec}(\mathbf H_\beta):\mu\le E\big\} \\
&= \#\big\{ t\ge 0 : \rho_{-E}^*(t) = -\infty \big\}, \\
\rho_{-E}^*&\text{ solves \eqref{Ric0} with parameter $-E$ and }
\rho_{-E}^*(0) =+\infty.
\end{aligned}
%\mathrm{p}_\lambda = \operatorname{Riccati}(\mathrm{g}_{\lambda,0}). 
\]
The function $\R\ni E\mapsto \mathbf{N}^{\mathbf{H}}_{E}$ is non-decreasing and c\`adl\`ag, it jumps only by $+1$, and $\lim_{E\to-\infty} \mathbf{N}^{\mathbf{H}}_{E} =0$.
We also denote the stochastic Airy zero counting function in the spectral coordinate $E$ by
\[
\mathbf{N}^{\SAi}_{E} = \#\big\{ t \ge 0 : \SAi_{-E}(t) =0 \big\},
\]
and the conventional counting function of the negated set $-\A$, where $\A= \big\{ \lambda\in\R : \SAi_\lambda(0) =0 \big\}$, by
\[
\mathbf{N}^{-\A}_{E} = \#\big\{\nu\in-\A:\nu\le E\big\}
= \#\big\{\lambda\ge-E:\SAi_\lambda(0)=0\big\}.
\]

It turns out that these three counting functions coincide, which gives a direct proof of the fact that the zero set of the stochastic Airy function $\lambda \mapsto \SAi_\lambda(0) =0 $ is the Airy$_\beta$ point process.
The following proposition is a consequence of the fact that the stochastic Airy operator $\mathbf{H}_\beta$ is a Sturm--Liouville operator with a random potential. 

\begin{proposition} \label{prop:cf}
Almost surely,  
$\mathbf{N}^{\mathbf{H}}_{E} = \mathbf{N}^{\SAi}_{E} =\mathbf{N}^{-\A}_{E}$ for all $E\in\R$.
\end{proposition}

\begin{proof}
In this proof, we let  $\mathrm{q}_\lambda = (\log\SAi_\lambda)'$ be the Riccati transform of the solution $\SAi_\lambda$.
We also abbreviate $\mathrm{g}_\lambda = \mathrm{g}_{\lambda,0}$ and  
$\rho_\lambda^* = (\log\mathrm{g}_\lambda)'$.
Then using the notation of Section~\ref{sec:SAR},  $\mathbf{N}^{\mathbf{H}}_{-\lambda} = \mathbf{N}_{[0,\infty)}(\rho_{\lambda}^*)$,
$\mathbf{N}^{\SAi}_{-\lambda} =  \mathbf{N}_{[0,\infty)}( \mathrm{q}_\lambda)$
and by \eqref{interlacing}, $\mathrm{q}_\lambda \prec \rho_\lambda^*$, so that it holds for any $\lambda\in\R$, 
\[
\mathbf{N}^{\mathbf{H}}_{-\lambda} \le  \mathbf{N}^{\SAi}_{-\lambda}  \le \mathbf{N}^{\mathbf{H}}_{-\lambda}+1.
\] 
In particular, this guarantees that the counting function $\mathbf{N}^{\SAi}_{E}$ is well-defined. 
We rely on the following characterization of the stochastic Airy function which follows from Proposition~\ref{prop:asympDN}:
if $\phi_\lambda$ is any solution of \eqref{SA2}, then 
$\big\{{\displaystyle \lim_{t\to\infty}} (\log\phi_\lambda(t))' \to -\infty\big\} =
\big\{ (\log\phi_\lambda)' = \mathrm{q}_\lambda  \big\} = \big\{ \phi_\lambda \propto \SAi_\lambda\big\}$
up to a null event.  This is a single null event off of which the events coincide for all $\lambda \in \R,$ and we shall suppress this going forward. 

Then the key observation is that since $\mathrm{q}_\lambda \prec \rho_\lambda^*$ and $\lim_{t\to\infty} \mathrm{q}_\lambda(t) =-\infty $, it holds for any $\lambda\in\R\setminus\A$ and  $k\in\N$:  $\big\{ \mathbf{N}^{\SAi}_{-\lambda} =k \big\} \subseteq \big\{ \mathbf{N}^{\mathbf{H}}_{-\lambda} =k\big\}$. 
Indeed, for  $ \mathbf{N}^{\SAi}_{-\lambda} > \mathbf{N}^{\mathbf{H}}_{-\lambda}$ it must be that the Riccati $\mathrm{q}_\lambda$ has blown down one more time, and so is larger than the Riccati $\rho_\lambda^*$ for all time sufficiently large.  Hence as $\lim_{t\to\infty} \mathrm{q}_\lambda(t) =-\infty,$
$\big\{  \mathbf{N}^{\SAi}_{-\lambda} =k,  \mathbf{N}^{\mathbf{H}}_{-\lambda} =k-1\big\} 
\subset \big\{{\displaystyle \lim_{t\to\infty}}  \rho_\lambda^*(t) =-\infty \big\} 
= \big\{ \mathrm{q}_\lambda =  \rho_\lambda^* \big\}.$
 Thus we have for any $\lambda\in\R$ and $k\in\N$,
\[
\big\{\mathbf{N}^{\SAi}_{-\lambda} =k ,\mathbf{N}^{\mathbf{H}}_{-\lambda} =k-1\big\}
\subseteq \big\{\SAi_\lambda \propto \mathrm{g}_\lambda\big\} \subseteq \{ \lambda \in \A\}.
\] 
These facts imply that 
 \[
 \big\{ \mathbf{N}^{\SAi}_{-\lambda} =k \big\}= \big\{ \mathbf{N}^{\mathbf{H}}_{-\lambda} =k\big\} , \qquad\text{for all } \lambda\in\R\setminus\A \text{ and } k\in\N_0.
\]

Now recall that the function $(\lambda , t) \in\R^2 \mapsto  \SAi_\lambda(t)$ is $C^1$, analytic in $\lambda$, 
and  $\big\{ t\in \R_+ :   \SAi_\lambda(t) =   \SAi_\lambda'(t) =0  \big\} = \varnothing$. This implies that $-\mathcal{A}$ is a countable set with no accumulation point, that the counting function $E\mapsto\mathbf{N}^{\SAi}_E$ is c\`adl\`ag and non-decreasing, and that it jumps by $+1$ if and only if $E\in-\A$ (if a zero of $\SAi_{-E}(t)$ were to appear or disappear in the bulk as $E$ increases, it would have to be a double zero).
 By convention, the counting functions $\mathbf{N}^{\mathbf{H}}$ and $\mathbf{N}^{-\A}$ are also c\`adl\`ag.
 Hence, we conclude that $\mathbf{N}^{\mathbf{H}}= \mathbf{N}^{\SAi} $  and in particular 
  $\lim_{E\to-\infty} \mathbf{N}^{\SAi}_{E} =0$.
 This implies that $-\mathcal{A}$ is bounded from below, or equivalently that $\mathcal A$ is bounded from above, and that
 $\mathbf{N}^{\SAi} =\mathbf{N}^{-\A}$ (as both have the same jumps and converge to $0$ as $E\to-\infty$).
\end{proof}

\subsection{Locally uniform asymptotics for the Riccati diffusions}%: Proofs of Propositions \ref{prop:asympD} and \ref{prop:asympN}}
\label{sec:asymp}
Our main tool to prove Proposition~\ref{prop:asympDN} is a precise study of the asymptotic behavior of the Riccati diffusions introduced in Section~\ref{sec:SAR}.
Throughout this section, we fix a compact set $K \subset \C$. 
We are interested in asymptotics for the Dirichlet and Neumann solutions (see Definition~\ref{def:DNS}) which are uniform for $\lambda=  (\mu +\i\eta) \in K$.
Let us recall from Lemma~\ref{lem:phipv} that for any $n\ge 1$ (which need not be an integer) and $s_n>n$, on the event $\big\{ \mathbf{N}_{[s_n,\infty)}(\rho_{\lambda,n}) =0\big\} $, 
it holds for any $t\ge s_n$, 
\begin{equation} \label{Ric2}
\phi(t) = \phi(s_n) \exp\biggl( \int_{s_n}^{t}\rho_{\lambda,n}(u)\d u \biggr) 
\qquad\text{for } \phi \in \big\{ \mathrm{f}_{\lambda,n}, \mathrm{g}_{\lambda,n} , \lambda\in K \big\} . 
\end{equation}

Let $\rho=\rho_{\lambda,n}$ and decompose $\rho:= \sigma + \i \varpi$. Then, the Riccati equation \eqref{Ric0} becomes a coupled system:
\begin{equation}
  \begin{aligned}
  \d\sigma(t) &= -\sigma^2(t) + \varpi^2(t) + \mu + t + \tfrac{2}{\sqrt{\beta}}\,\d B(t) ,  \quad \sigma(n) = 0/\infty, \\
  \d\varpi(t) &= -\big(2\sigma(t)\varpi(t) - \eta\big) \d t ,
 \hspace{2.1cm} \varpi(n) = 0 ,
  \end{aligned}
  \label{eq:complexRic}
\end{equation}
where the initial condition $ \sigma(n) = 0/\infty$ depends on whether we consider the Dirichlet or Neumann solution. 
The equation for $\varpi$ being a first order ODE is solved by
\begin{equation}   \label{eq:omega}
  \varpi_{\lambda,n}(t) =  \eta\int_{n}^t e^{-2\int_u^t \sigma_{\lambda,n}(v)\,\d v}\,\d u . \\
\end{equation}
Then, we deduce from \eqref{eq:complexRic} an autonomous equation for the diffusion $\sigma_{\lambda,n}$ (up to the first blow-down time). 
We expect that as $t\to\infty$, this system approaches equilibrium in the sense that $ \sigma_{\lambda,n}(t)$ fluctuates around the stable parabola $\sqrt{t+\mu}$ and $\varpi_{\lambda,n}(t)$ fluctuates around $\frac{\eta}{2\sqrt{t+\mu}} $.  
Moreover this equilibrium is reached quickly if the starting time $n$ is large and we can explicitly describe the fluctuations of the Riccati diffusion $\rho_{\lambda,n}$.
 We now overview these facts. %and the proof of Proposition \ref{prop:asympDN}.

\subsubsection*{Step 1: Entrance behavior.}

If the initial time $n>0$ is large, the diffusions  $\sigma_{\lambda,n}$ reach the stable parabola in a very specific way both in the Dirichlet $(\sigma_{\lambda}(n)=0)$  and Neumann  $(\sigma_{\lambda}(n)=\infty)$ cases. In the Dirichlet case, one can approximate \eqref{eq:complexRic} on a short time interval $[n,s_n]$ by the ODE:
\begin{equation} \label{tanheq}
\psi' = -\psi^2 + n, \qquad \psi(n)=0. 
\end{equation}
The solution being given explicitly by $\psi_n(t) \coloneqq \tanh((t-n)\sqrt{n})\sqrt{n}$.  In the Neumann case, instead it is convenient to work with the inverse Riccati diffusion.  But again, we show the diffusion is well approximated by an ODE.  We refer to \eqref{eq:com-inversericatti} and \eqref{tanheq2} below for the precise approximation. We give a summary result which applies to both Dirichlet and Neumann cases.

\begin{proposition} \label{prop:entrance}
% Let $(\sigma_{\lambda}(t) : t \geq s,\lambda \in K)$ be the real-part of the Riccati diffusions either in the Dirichlet or Neumann cases. 
Let $s_n = n+ \frac{\log n}{\sqrt n}$ and $\kappa_n = n+ \frac{\kappa}{\sqrt n}$ where $\kappa$ is any fixed constant.
 For any $\epsilon \in (0, \tfrac 34)$, there are constants  $c=c_{\beta,K,\epsilon}$ and  $C=C_{\beta,K,\epsilon, \kappa}$ so that for $n>0$ large enough (depending only on $K$), it holds with probability $1-e^{-\frac{n^{3/2-2\epsilon}}{C\log(\log n)}}$,
\[    
|  \varpi_{\lambda,n}(s_n)|  \le c n^{-1/2} 
 \quad\text{and}\quad
 \big| \sigma_{\lambda,n}(s_n) - \sqrt{s_n}\big|
 \le c n^{1/2-\epsilon} .
\]
Moreover,   for $\ell\in\{0,1\}$,  it holds  uniformly for  $t \in [n,s_n]$ and $\lambda\in K$, 
\[
  \begin{aligned}
 \mathrm{f}_{\lambda,n}(t) =  \cosh((t-n)\sqrt{n})e^{\O(n^{-\epsilon}\log n)}
  \quad
  \text{and}
  \quad
  \partial_t^\ell \mathrm{g}_{\lambda,n}(t)=  \partial_t^\ell\Big(\tfrac{\sinh((t-n)\sqrt{n})}{\sqrt{n}}\Big)e^{\O(n^{-\epsilon}\log n)} . 
\end{aligned}
\]
For the Dirichlet solutions, we have instead  uniformly for  $t \in [\kappa_n,s_n]$ and $\lambda\in K$, 
\[
\mathrm{f}_{\lambda,n}'(t) 
= \partial_t\big(  \cosh((t-n)\sqrt{n}) \big) e^{\O(n^{-\epsilon}\log n)} .
\]
  \end{proposition}
\noindent The proof is given in Section~\ref{sec:entrance}. 
 Let us observe that in the Dirichlet case, for very small $t\ge 0$, 
 $\mathrm{f}_{\lambda,n}'(n+t) \approx \sigma_{\lambda,n}(n+t) = \tfrac{2}{\sqrt{\beta}}\bigl(B(n+t)-B(n)\bigr) + \O(t^{1-\epsilon})$ so that this function cannot be approximated by $\sqrt{n} \sinh(t\sqrt{n}) \simeq n t$ up to a vanishing multiplicative error  for $t\ll \kappa/\sqrt{n}$. 
 
\subsubsection*{Step 2: the linearized process.}
Once in a neighborhood of the parabola, it is possible to give a good linear diffusion approximation for $\sigma_{\lambda,s}(t)-\sqrt{t+\mu}.$  Specifically, we consider the following diffusion: for $\mu\in\R$ with $\mu<n<t$, 
\begin{equation}
  \d \mathfrak{h}_{\mu,n}(t) 
  =
  - 2\Big( \sqrt{t+\mu}\,\mathfrak{h}_{\mu,n}(t) + \frac{1/4}{\sqrt{t+\mu}} \Big) \d t + \tfrac{2}{\sqrt{\beta}}\,\d B(t),
  \quad
  \mathfrak{h}_{\mu,n}(n) = 0.
  \label{eq:linearized}
\end{equation}
This diffusion is explicitly solvable.  In preparation to give its solution, we introduce the integrating factor for $\mu\in\R$, 
\begin{equation} \label{def:y}
y_\mu(t) =  \exp\big(\tfrac43 (t+\mu)^{3/2}  \big)  
\quad\text{so that for }t \geq s \geq \mu,\quad
y_\mu(t) y_\mu^{-1}(s)  = \exp\bigg(2\int_s^t  \sqrt{u+\mu}  \d u \bigg) . 
\end{equation}
Throughout the remainder of this paper, we make use of the following basic estimates:  for any $\gamma \in \R$,  
\begin{equation}\label{eq:yeezy}
  \begin{aligned}
    y_\mu^{-1}(t) \int_c^t y_\mu(u) u^{\gamma}\d u = \tfrac{1}{2} t^{\gamma-1/2} + \underset{t\to\infty}{\O( t^{\gamma-2})} \quad \text{and} \quad
    y_\mu(t) \int_t^\infty y_\mu^{-1}(u) u^{\gamma}\d u  =  \tfrac{1}{2} t^{\gamma-1/2} + \underset{t\to\infty}{\O( t^{\gamma-2})} 
  \end{aligned}
\end{equation}
where both errors are uniform for $\mu\in [-c,c]$. 

The solution to the stochastic differential equation \eqref{eq:linearized} has the explicit representation
\begin{equation*}
\mathfrak{h}_{\mu,n}(t) =   y_\mu^{-1}(t) \Bigg(- \int_{n}^t \frac{y_\mu(u)}{{2}\sqrt{u+\mu}} \d u
+ \tfrac{2}{\sqrt{\beta}}\int_{n}^t y_\mu(u)  \d B(u) \Bigg) .
\end{equation*}
This is therefore an uncentered Gaussian process, and we can use the asymptotics \eqref{eq:yeezy} to express
 \begin{align}  \label{def:mathfrakh} 
\mathfrak{h}_{\mu,n}(t)  
= -\frac{1}{4(t+1)} +  \tfrac{2}{\sqrt{\beta}}\mathfrak{X}_{\mu, n}(t)  + \O_K(t^{-2}),
\qquad\text{where}\quad
\mathfrak{X}_{\mu, n}(t) \coloneqq  y_\mu^{-1}(t) \int_{n}^t y_\mu(u)  \d B(u) 
\end{align}
and the error is deterministic and uniform for $\mu\in K$ if $n$ is large enough. 
We need estimates for the process $\mathfrak{X}_{\mu,n}$ to ensure that it does not perform a large excursion since this would result in a loss of control on the difference between $\sigma_{\lambda,s}(t) - \sqrt{t + \mu}$. 
In Section~\ref{sec:X}, we obtain the following tail bounds:

\begin{proposition} \label{prop:Xest}
  For any $\epsilon \in [-\tfrac 14, \tfrac 34)$,  and $n \geq \kappa \ge 0$,
    \[
%      \Big\VERT \sup_{|\mu|, |\nu| \leq \kappa, t\ge s} \bigg| (1+t)^{1/2+\epsilon} \tfrac{\mathfrak{X}_{\mu, s}(t)-\mathfrak{X}_{\nu, s}(t)}{\mu-\nu}\bigg| \Big\VERT_2 
%      +
      \Big\VERT \sup_{\mu\in[-\kappa,\kappa], t\ge n} | (1+t)^{-1/2+\epsilon} \mathfrak{X}_{\mu, n}(t)| \Big\VERT_2 
      \lesssim_{ \kappa}   (1+n)^{-3/4+\epsilon}\log(2+n).
    \]
%and 
%\[
%\bigg\VERT \sup_{|\mu| \le \kappa , t\ge s} \Big|  (1+t)^{-\epsilon} \int_s^t \mathfrak{X}_{\mu, s}(u) \d u \Big| \bigg\VERT_2  \lesssim_{\kappa} \beta^{-1}  (1+s)^{- \epsilon}  .
%\]
\end{proposition}
\noindent We will be interested in $\epsilon>0$ small. It follows by simple computations for the mean of $\mathfrak{h}$ that both of the conclusions of the proposition also hold for $\mathfrak{h}$ in place of $\mathfrak{X}.$

\subsubsection*{Step 3: control of the linearization error.}

We define the \emph{linearization error} as the difference:
\begin{equation} \label{delta0}
\Delta_{\lambda,n}(t)   
\coloneqq \sigma_{\lambda,n}(t)- \sqrt{t+\mu} -\mathfrak{h}_{\mu,n}(t),
%+ \int_s^t  \tfrac{y_\mu(u)}{y_\mu(t)} \mathfrak{h}_{\mu,s}^2(u)  \d u  
\qquad t\ge n.  
\end{equation}
We define a (random) integrating factor for $s\ge n$, 
\begin{equation}\label{eq:if}
  v_{\lambda,s}(t) \coloneqq 
  \frac{y_\mu(t)}{ y_\mu(s)}\exp\biggl( \int_s^t (\Delta_{\lambda,n}(u) + 2\mathfrak{h}_{\mu,n}(u))\d u \biggr)
  =\exp\biggl( \int_s^t \big(\sqrt{u+\mu} + \sigma_{\lambda,n}(u) + \mathfrak{h}_{\mu,n}(u) \big)\d u \biggr).
\end{equation}
One can rewrite the first equation in \eqref{eq:complexRic} in terms of $v$ and $\Delta$ as 
\begin{equation*} 
\d (v_{\lambda,s}(t) \Delta_{\lambda,n}(t)) = v_{\lambda,s}(t)\big( \varpi_{\lambda,n}(t)^2 - \mathfrak{h}_{\mu,n}(t)^2 \big) \d t. 
\end{equation*}
Note that this equation has no stochastic integral term,  hence its solution is continuously differentiable.   Moreover, this yields for $t\ge s\ge n$, 
\begin{equation} \label{delta2}
\Delta_{\lambda,n}(t)
= v_{\lambda,s}^{-1}(t)\Delta_{\lambda,n}(s)  + \int_s^t \frac{v_{\lambda,s}(u)}{v_{\lambda,s}(t)}\big( \varpi_{\lambda,s}^2(u) - \mathfrak{h}_{\mu,n}^2(u) \big) \d u. 
\end{equation}

We introduce a typical event on which we have control of the Riccati diffusion after a short time and of the linearization error: let  $s_n = n+ \frac{\log n}{\sqrt{n}}$ for $n\ge 2$,  $c\ge 1$, a small $\epsilon>0$ and define
    \begin{equation} \label{Eevent}
    \mathcal{E}_{n} \coloneqq \left\{   \big| \varpi_{\lambda,n}(s_n)\big|  \le c n^{-1/2}  , \big| \Delta_{\lambda,n}(s_n)\big|  \le c n^{1/2-\epsilon}   :  \lambda\in K \right\} \cap \left\{  
   \big|\mathfrak{X}_{\mu,n} (t)| \le c t^{1/2-\epsilon} :   t \geq n, \mu \in \Re(K) \right\} .
    \end{equation}
By Propositions~\ref{prop:entrance} and~\ref{prop:Xest},  there is a constant $C=C_{K,\beta, \epsilon,c}$ so that for any $\epsilon \in (0, \tfrac 12]$,
 \begin{equation} \label{PEevent}
 \P\big( \mathcal{E}_{n}^c \big) \le C \exp\big(-n^{3/2-5\epsilon/2}\big). 
\end{equation}

Working on this event, we can show that deterministically, the linearization remains negligible with respect to the principal $\sqrt{t}$ term:

\begin{proposition}\label{lem:sig}
For $n\ge n_K$ large enough and $\epsilon \in (0, \tfrac 34)$, there is a constant $C=C_{K,\beta, \epsilon,c}$ so that on the event \eqref{Eevent},  it holds for all $t \geq s_n$,
  \[
    | \varpi_{\lambda,n}(t) | \le     Ct^{-1/2}
    \quad \text{and}\quad 
    | \Delta_{\lambda,n}(t) | \le     Ct^{1/2-\epsilon} . 
%    \max\left\{ 
%      \tfrac{|\sigma_{\lambda,n}(t)-\sigma_{\nu,n}(t)|}{|\lambda-\nu|},
%      \tfrac{|\varpi_{\lambda,n}(t)-\varpi_{\nu,n}(t)|}{|\lambda-\nu|}
%    \right\}
%    \leq Ct^{-1/2}.
  \]
Moreover, for all $t \geq n+\frac{4\log n}{\sqrt n}$,
\[
\bigl(\Delta_{\lambda,n}(t)\bigr)_+ \leq Ct^{-3/2}.
\]
\end{proposition}

%    \begin{equation}\label{eq:sig-event}
%    \begin{aligned}
%    \mathcal{E}_{n}(\epsilon) &= 
%    \biggl\{ \forall \lambda,\nu\in K: 
%      \max\left\{
%	\varpi_{\lambda,n}(s_n),
%	\tfrac{|\sigma_{\lambda,n}(s_n)-\sigma_{\nu,n}(s_n)|}{|\lambda-\nu|},
%	\tfrac{|\varpi_{\lambda,n}(s_n)-\varpi_{\nu,n}(s_n)|}{|\lambda-\nu|}
%      \right\}
%      \leq c s_n^{-1/2},
%     % |\varpi_{\lambda,n}(s_n)| \le \frac{c}{\sqrt{s_n}}, 
%      |\Delta_{\lambda,n}(s_n)| \le c{s_n^{1/2-\epsilon}}\biggr\}  \\ %\right\} \\
%     &\qquad  %\left\
%     \cap
%     \biggl\{ 
%      \forall \mu \in \Re K, t \geq n :  
%      |\mathfrak{h}_{\mu,n} (t)| \le {c}{t^{1/2-\epsilon}}
%      %,  %   	\bigg|  \int_{n}^t \mathfrak{h}_{\mu, n}(u) \d u \bigg| \le c t^\epsilon
%    \biggr\} , 
%    \end{aligned}
%  \end{equation}

\noindent See Section \ref{sec:Delta} for the proof.
In particular, on the event $\mathcal{E}_{n}$,  it holds for all $\lambda\in K$ and $t\ge s_n$ 
\begin{equation} \label{Ric3}
|\rho_{\lambda,n}(t)- \sqrt{t}| \le C t^{1/2-\epsilon} . 
\end{equation}

\subsubsection*{Step 4: control on the integrated errors.}

Proposition~\ref{lem:sig} implies that  $ \mathbf{N}_{[s_n,\infty)}(\rho_{\lambda,n}) =0$ and formula \eqref{Ric2} makes sense. Moreover we can rewrite
\begin{equation} \label{rhoint1}
\int_{s_n}^{t}\rho_{\lambda,n}(u)\d u
= \int_{s_n}^{t}  \bigg( \sqrt{u+\mu} + \i \frac{\eta}{2\sqrt{u}} +\mathfrak{h}_{\mu,n}(u) \bigg) \d u  + \int_{s_n}^{t}\bigg(\Delta_{\lambda,n}(u) + \i \bigg( \varpi_{\lambda,n}(u) - \frac{\eta}{2\sqrt{u}}\bigg) \bigg) \d u . 
\end{equation}
The first integrand accounts for the first three terms in the definition of $\theta$, see \eqref{def:theta}, 
and we now want to estimate the contribution from the integrated error that is, the second term on the RHS of \eqref{rhoint1}. 
It turns out that the mean of the squared--Gaussian process $\mathfrak{h}^2_{\mu,n}$ in the expression \eqref{eq:if} for $\Delta$  when averaged in $t$ yields a non-trivial contribution. 
Namely, we show in Section \ref{sec:X2} that: 
\begin{proposition} \label{prop:2021}
For $n\ge n_K$ large enough and $\epsilon \in (0, \tfrac 34)$, there is a constant $C=C_{K,\beta, \epsilon,c}$  so that
  \[
    \Pr\biggl(
    \mathcal{E}_{n}
    \cap
    \biggl\{
    \sup_{\lambda \in K}
    \sup_{T \geq s_n} 
    \bigg| \int_{s_n}^T
    \biggl(
    \Delta_{\lambda,n}(t) + \frac{4}{\beta}\int_0^t \Exp[\mathfrak{X}(t)\mathfrak{X}(u)] \d u\biggr)  \d t \bigg|  \geq Cn^{-\epsilon}
  \biggr\}
    \biggr)
    \leq \exp\bigg(-\frac{n^{3/2-2\epsilon}}{(\log n)^2}\big) .
  \]
\end{proposition}

In Section \ref{sec:X2}, we also show that the integrated error coming from the imaginary part of the Riccati diffusion is small for large enough $n$. 

\begin{proposition} \label{prop:2022}
  Under the same setup as Proposition~\ref{prop:2021}, 
  \[
    \Pr\biggl(
    \mathcal{E}_{n}
    \cap
    \biggl\{
    \sup_{\lambda \in K}
    \sup_{T \geq s_n} 
    \bigg| \int_{s_n}^T \bigg( \varpi_{\lambda,n}(t) -\frac{\eta}{2\sqrt{t}}\bigg)  \d t \bigg|  \geq C n^{-1/2}
  \biggr\}
    \biggr)
    \leq \exp\bigg(-\frac{n^{3/2-2\epsilon}}{(\log n)^2}\big) .
  \]
\end{proposition}

The a priori bounds from Proposition~\ref{lem:sig} are crucial to prove both Propositions~\ref{prop:2021} and~\ref{prop:2022} as it allows us 
to approximate the integrated factor \eqref{eq:if} by $v_{\lambda,s}(t)\approx y_{\mu}(t)y_{\mu}^{-1}(s)$ in an averaged sense. 
Then, in terms of \eqref{delta2} and \eqref{def:mathfrakh}, we obtain
$\Delta_{\lambda,n}(t) \approx - \frac{4}{\beta}\int_s^t \frac{y_{\mu}(u)}{y_{\mu}(t)} \mathfrak{X}_{\mu,n}^2(u) \d u$ for $t$ large enough. 
In addition, we will require concentration estimates for such integrals of squared--Gaussian processes; see Lemma~\ref{lem:S}.

From this point, it is a simple extension to produce the asymptotics for the integral of the Riccati diffusion and to deduce that of the Dirichlet and Neumann solutions. 
Recalling that by definition \eqref{def:theta}--\eqref{eq:covX2}:
\[
\theta_\lambda(t) = \sqrt{t} + \frac{\lambda}{2\sqrt{t}} - \frac{1}{4(t+1)} +\frac{2}{\sqrt{\beta}}\mathfrak{X} (t) - \frac{4}{\beta}\int_0^t  \E[ \mathfrak{X}(t) \mathfrak{X}(s)]  \d s , 
  \] 
we can conclude that:

\begin{proposition} \label{prop:X2int}
For $n\ge n_K$ large enough and $\epsilon \in (0, \tfrac 34)$, there is a constant $C=C_{K,\beta, \epsilon,c}$  so that
  \[
    \Pr\biggl(
    \mathcal{E}_{n}
    \cap
    \biggl\{
    \sup_{\lambda \in K}
    \sup_{T \geq s_n} 
    \bigg| \int_{s_n}^T
    \bigl(
    \rho_{\lambda,n}(t) - \theta_{\lambda}(t)
    \bigr)
    \d t \bigg|  \geq C n^{-\epsilon}
  \biggr\}
    \biggr)
    \leq \exp\bigg(-\frac{n^{3/2-2\epsilon}}{(\log n)^2}\big) .
  \]
\end{proposition}
\noindent See Section \ref{sec:X2} for the proof.

%\paragraph{Step 5: Existence of the stochastic Airy function limit and high probability estimates.}
\subsection{Existence of the stochastic Airy function and high probability estimates: Proofs of Propositions~\ref{prop:asympDN} and \ref{prop:SAiwins}}
\label{sec:proofasymp}

\paragraph{Typical events.}
Recall that $s_n = n + \frac{\log n}{\sqrt n}$ and $\kappa_n = n + \frac{\kappa}{\sqrt n}$ for a fixed $\kappa\ge 1$ and $n\in \R_+$ with $n \ge e^\kappa$. For $\epsilon\in (0,\tfrac 12]$ and $C>0$, we define the events
\[
\mathcal{B}_n \coloneqq \left\{
\bigg|\frac{\rho_{\lambda,n}(t)}{\sqrt t} -1\bigg| \le C t^{-\epsilon} , 
  \bigg| \int_{s_n}^t
  \bigl(
  \rho_{\lambda,n}(u) - \theta_{\lambda}(u)
  \bigr)
  \d u \bigg|  \leq C n^{-\epsilon}  : t\ge s_n , \lambda\in K
\right\}  
\]
where $\rho_{\lambda,n}$ can be the Riccati transform of either $\mathrm{f}_{\lambda,n}$ or $\mathrm{g}_{\lambda,n}.$
 By Proposition~\ref{prop:X2int}, \eqref{Ric3}  and \eqref{PEevent}, in both cases,
\[
 \P\big( \mathcal{B}_{n}^c \big) \le  \exp\big(-n^{3/2-2\epsilon}\big) +  \P\big( \mathcal{E}_{n}^c \big)  \lesssim \exp\big(-n^{3/2-5\epsilon/2}\big) . 
\]

We also work on the entrance event from  Proposition~\ref{prop:entrance}, 
    \begin{align}\notag
\mathcal{C}_n & \coloneqq \left\{
\begin{gathered}
 \partial_t^\ell \mathrm{f}_{\lambda,n}(t)
 = \partial_t^\ell \cosh((t-n)\sqrt{n})e^{\O_\ell(n^{-\epsilon})},\\
  \partial_t^\ell \mathrm{g}_{\lambda,n}(t)
  = \partial_t^\ell\tfrac{\sinh((t-n)\sqrt{n})}{\sqrt{n}}e^{\O_\ell(n^{-\epsilon})}  :\,
  t\in [\kappa_n,s_n] , \lambda\in K , \ell \in \{0,1\}
\end{gathered}
\right\} \\
& \label{Cevent}\qquad \cap 
\left\{ \big| \theta_\lambda(t)  - \sqrt{t} \big|\le Cn^{3/4-5\epsilon/4}t^{-1/4} : t\ge n , \lambda\in K \right\}
\end{align}
To control the last event, by Proposition~\ref{prop:Xest}, we verify that for $n\ge 1$, 
\[
\P\left( 
\sup_{t\ge n}|\mathfrak{X}(t)|t^{1/4} > n^{3/4-5\epsilon/4} \right) \lesssim  \exp\big(-n^{3/2-5\epsilon/2}\big) .
\]
By definition of the process $\theta$, this implies that
\[
 \P\big( \mathcal{C}_{n}^c \big)  \lesssim \exp\big(-n^{3/2-5\epsilon/2}\big) . 
\]

Hence there is a constant $C_{\epsilon,K,\beta}$ so that for any $n\ge 1$,
\begin{equation} \label{gtb}
 \sum_{k \in (\N)^{2/3} : k\ge n} \P( \mathcal{B}_{k}^c  \cup  \mathcal{C}_{k}^c ) \le 
 C_{\epsilon,K,\beta} \exp\big(-n^{3/2-3\epsilon}\big) .
\end{equation}
We are going to use that by the Borel--Cantelli lemma, these events hold almost surely for all sufficiently large $n \in (\N)^{2/3}$ (the relevance of which will be explained shortly). On these typical events, we can proceed to prove the results announced in Section~\ref{sec:SAi}. 
First, let us make a few extra basic observations.

\medskip

On the event \eqref{Cevent}, for any $\ell \in \N_0$,
\[
\bigg(\partial_t^\ell\frac{\sinh((t-n)\sqrt{n})}{\sqrt{n}}\bigg)  \exp\biggl( -\int_{n}^{t}\theta_{\lambda}(u)\d u \biggr) \bigg|_{t=s_n} = \frac{n^{(\ell-1)/2}}{2} e^{ \O_\ell(n^{-\epsilon})}. 
\]
This shows that on $\mathcal{C}_n$, it holds uniformly for $\lambda\in K,$
\begin{equation} \label{scontrol}
\begin{aligned}
\mathrm{f}_{\lambda,n}(t)e^{-\int_{n}^{t}\theta_{\lambda}(u)\d u}
 \bigg|_{t=s_n}  
    =  \frac{e^{ \O(n^{-\epsilon})} }{2} 
\qquad\text{and}\qquad
\mathrm{g}_{\lambda,n}(t) e^{-\int_{n}^{t}\theta_{\lambda}(u)\d u}
  \bigg|_{t=s_n}  
    =  \frac{e^{ \O(n^{-\epsilon})} }{2\sqrt{n}} .
\end{aligned}
\end{equation}

%    \exp\biggl( -\int_{n}^{t}\theta_{\lambda}(u)\d u \biggr)
We now consider the events 
    \begin{equation} \label{Fevent}
    \begin{aligned}
    \mathcal{F}_n & \coloneqq
\left\{ \mathrm{f}_{\lambda,n}(t) e^{-\int_{n}^{t}\theta_{\lambda}(u)\d u}  = 
    \frac{e^{\O(n^{-\epsilon})}}{2} \text{ and }     \mathrm{g}_{\lambda,n}(t) e^{-\int_{n}^{t}\theta_{\lambda}(u)\d u}
 =  \frac{e^{\O(n^{-\epsilon})}}{2 \sqrt{n}}   \text{ uniformly for }t\ge s_n , \lambda\in K \right\} \\
     \mathcal{G}_n & \coloneqq
\left\{ \mathrm{f}_{\lambda,n}'(t) \frac{e^{-\int_{n}^{t}\theta_{\lambda}(u)\d u}}{\sqrt t}  = 
    \frac{e^{\O(n^{-\epsilon})}}{2} \text{ and }  \mathrm{g}_{\lambda,n}'(t) \frac{e^{-\int_{n}^{t}\theta_{\lambda}(u)\d u}}{\sqrt t}
 =  \frac{e^{\O(n^{-\epsilon})}}{2 \sqrt{n}}   \text{ uniformly for }t\ge s_n , \lambda\in K \right\}, 
\end{aligned}\end{equation}
where the implied constant in the $\O(n^{-\epsilon})$ term is any fixed number.

On $\mathcal{B}_n$, we have $\mathbf{N}_{[s_n,\infty)}(\rho_{\lambda,n}) =0$ and by formula \eqref{Ric2}, for $ \phi \in \big\{ \mathrm{f}, \mathrm{g}\big\}$ and for all $t\ge s_n$, 
\begin{equation*} 
\begin{aligned}
\phi_{\lambda,n} (t) e^{-\int_{n}^{t}\theta_{\lambda}(u)\d u}
&= \phi_{\lambda,n} (s_n) e^{-\int_{n}^{s_n}\theta_{\lambda}(u)\d u}
\exp\biggl( \int_{s_n}^{t}\big(\rho_{\lambda,n}(u) - \theta_\lambda(u) \big)\d u \biggr)\\
&= \phi_{\lambda,n} (s_n) e^{-\int_{n}^{s_n}\theta_{\lambda}(u)\d u + \O(n^{-\epsilon})}  .
\end{aligned}
\end{equation*}
Similarly, by differentiating formula \eqref{Ric2}, it holds on the event $\mathcal{B}_n$,  for $ \phi \in \big\{ \mathrm{f}, \mathrm{g}\big\}$ and for all $t\ge s_n$, 
\begin{equation*} 
\begin{aligned}
\phi_{\lambda,n} '(t) \frac{e^{-\int_{n}^{t}\theta_{\lambda}(u)\d u}}{\sqrt t}
&= \phi_{\lambda,n} (s_n) e^{-\int_{n}^{s_n}\theta_{\lambda}(u)\d u}
\frac{\rho_{\lambda,n}(t)}{\sqrt t}
\exp\biggl( \int_{s_n}^{t}\big(\rho_{\lambda,n}(u) - \theta_\lambda(u) \big)\d u \biggr)\\
&= \phi_{\lambda,n} (s_n) e^{-\int_{n}^{s_n}\theta_{\lambda}(u)\d u + \O(n^{-\epsilon})}  .
\end{aligned}
\end{equation*}

By \eqref{scontrol}, this shows that $\mathcal{B}_{n}  \cap  \mathcal{C}_{n} \subset \mathcal{F}_n \cap  \mathcal{G}_n$ for $n\ge 1$. 
Hence,  we can work on the events $ \mathcal{F}_{n}$ and $ \mathcal{G}_{n}$ which control the asymptotics of the Dirichlet and Neumann solutions and their derivative as $t\to\infty$.
In particular these events occur for all sufficiently large $n \in (\N)^{2/3}$.

\begin{proof}[Proof of Proposition~\ref{prop:asympDN}.]
We introduce the notation for $n\in\R$ and $\lambda\in \C$, 
\begin{equation} \label{def:a}
\mathrm{a}_{\lambda,n}^{\mathrm{f}} \coloneqq \lim_{t\to\infty}\bigg\{ \mathrm{f}_{\lambda,n}(t) \exp\biggl( -\int_{n}^{t}\theta_{\lambda}(u)\d u \biggr)  \bigg\}
  \quad\text{and}\quad
  \mathrm{a}_{\lambda,n}^{\mathrm{g}} \coloneqq \lim_{t\to\infty}\bigg\{ \mathrm{g}_{\lambda,n}(t) \exp\biggl( -\int_{n}^{t}\theta_{\lambda}(u)\d u \biggr)  \bigg\}.
\end{equation}

\paragraph{Existence of $\SAi_\lambda'$.}
We begin by showing that the limits \eqref{def:a} exist almost surely uniformly in $\lambda \in K$ (where $K\subset\C$ is an arbitrary compact set). In particular, the resulting functions of $\lambda\in\C \mapsto \mathrm{a}_{\lambda,n}^{\mathrm{f}}, \mathrm{a}_{\lambda,n}^{\mathrm{g}}$ are entire, as $\lambda\in \C \mapsto\theta_\lambda(t) , \mathrm{f}_{\lambda,n}(t),\mathrm{g}_{\lambda,n}(t)$ are all entire functions.

\medskip

Recall that for any time $n,k \in\R$, we can always decompose a solution $\phi_{\lambda,n} $  in the Dirichlet--Neumann basis,
\begin{equation}\label{eq:aD1}
\phi_{\lambda,n}(t) = 
\phi_{\lambda,n}(k) \mathrm{f}_{\lambda,k}(t) + {\phi_{\lambda,n}'(k)}\mathrm{g}_{\lambda,k}(t)  . 
\end{equation}
On the event $   \mathcal{F}_n \cap \mathcal{G}_n$, for $\phi \in \big\{ \mathrm{f}, \mathrm{g}\big\}$, both
\[
  \sup_{k \geq n , \lambda\in K} |\phi_{\lambda,n}(k)|\exp\biggl( -\int_{n}^{k}\theta_{\lambda}(u)\d u \biggr)
  < \infty
  \quad\text{and}\quad
  \sup_{k \geq n, \lambda\in K} \tfrac{|\phi'_{\lambda,n}(k)|}{\sqrt{k}}\exp\biggl( -\int_{n}^{k}\theta_{\lambda}(u)\d u \biggr)
  < \infty.
\]
Moreover, on $\bigcap_{k \in  (\N)^{2/3} : k\ge n}   \mathcal{F}_k $, we have 
\begin{equation}\label{eq:aD2}
  \lim_{ (\N)^{2/3}\ni k \to \infty}
  \sup_{t \geq s_k, \lambda\in K}
  \biggl\{
  \biggl|
  2\mathrm{f}_{\lambda,k}(t)\exp\biggl( -\int_{k}^{t}\theta_{\lambda}(u)\d u \biggr)-1\biggr|
  +
  \biggl|
  2\sqrt{k}\mathrm{g}_{\lambda,k}(t)\exp\biggl( -\int_{k}^{t}\theta_{\lambda}(u)\d u \biggr)
-1\biggr|\biggr\}
=0 .
\end{equation}
Then, it follows from \eqref{eq:aD1} that $t \mapsto \phi_{\lambda,n}(t)\exp\biggl( -\int_{n}^{t}\theta_{\lambda}(u)\d u \biggr)$ is a Cauchy sequence.  This implies that almost surely, $\mathrm{a}_{\lambda,n}^{\mathrm{f}}$ and $\mathrm{a}_{\lambda,n}^{\mathrm{g}}$ exist uniformly for $\lambda\in K$ for all $n \in (\N)^{2/3}$ sufficiently large.

\medskip

Using \eqref{eq:aD1} once more choosing $k \in (\N)^{2/3}$ large enough, we conclude that on the co-null event  
\[
  \widehat{\mathcal{F}} \coloneqq\bigcup_{n \in (\N)^{2/3}}\bigcap_{k \in  (\N)^{2/3} : k\ge n}   \mathcal{F}_k 
\]
both of  $\mathrm{a}_{\lambda,n}^{\mathrm{f}}, \mathrm{a}_{\lambda,n}^{\mathrm{g}}$ actually exist for all $n \in \R$. 
Let us recall that according to Definition~\ref{def:SAiprime}, 
\begin{equation} \label{SAi'lim}
  \SAi'_{\lambda}(n) \coloneqq 
  -\frac{1}{\sqrt{\pi}}
  \lim_{t\to\infty}\bigg\{ \mathrm{f}_{\lambda,n}(t) \exp\biggl( -\int_{0}^{t}\theta_{\lambda}(u)\d u \biggr)  \bigg\}
  =  -\frac{1}{\sqrt{\pi}}
  \mathrm{a}_{\lambda,n}^{\mathrm{f}}
  \exp\biggl( -\int_{0}^{n}\theta_{\lambda}(u)\d u \biggr). 
  %\quad\text{and}\quad
  %\SAi_{\lambda}(s) \coloneqq 
  %\frac{1}{\sqrt{\pi}}
  %\mathrm{a}_{\lambda,s}^{\mathrm{f}}
  %\exp\biggl( \int_{0}^{s}\theta_{\lambda}(u)\d u \biggr)
\end{equation}
Thus, we have proved that the function $\SAi'$ exists almost surely as a locally uniform limit in $\lambda \in \C$, and hence as an entire function. 
Next, we show that the expression $\SAi_{\lambda}(t) = -\int_t^\infty \SAi'_{\lambda}(s)\d s$ makes sense and we connect it to $\mathrm{a}_{\lambda,n}^{\mathrm{g}}.$

\paragraph{Asymptotics of $\SAi_\lambda'$.}
From \eqref{SAi'lim} determining the asymptotics $\SAi_\lambda'$ as $t \to \infty$ is equivalent to determining those of $\mathrm{a}_{\lambda,n}^{\mathrm{f}}$ as $n \to \infty.$
As we have the bound,
\begin{equation} \label{eq:aD0prime}
  |\mathrm{a}_{\lambda,n}^{\mathrm{f}} - \frac{1}{2}|
  \leq
  \limsup_{t \to \infty}
  \biggl|
	  \mathrm{f}_{\lambda,n}(t) e^{-\int_{n}^{t}\theta_{\lambda}(u)\d u}
  - \frac{1}{2}
  \biggr|
\quad\text{and}\quad
  |\sqrt{n}\mathrm{a}_{\lambda,n}^{\mathrm{g}} - \frac{1}{2}|
  \leq
  \limsup_{t \to \infty}
  \biggl|
  \mathrm{g}_{\lambda,n}(t) e^{-\int_{n}^{t}\theta_{\lambda}(u)\d u}
  - \frac{1}{2}
  \biggr|
\end{equation}
it follows that on the event $\widehat{\mathcal{F}}$ we have that restricting to $k \in (\N)^{2/3}$
\begin{equation} \label{eq:aD0}
  \mathrm{a}_{\lambda,k}^{\mathrm{f}} 
  =\frac{1}{2} + \O(k^{-\epsilon})
\quad\text{and}\quad
  \sqrt{k}\mathrm{a}_{\lambda,k}^{\mathrm{g}} 
  =\frac{1}{2} + \O(k^{-\epsilon})
\end{equation}
uniformly for $\lambda\in K$.  Moreover, for a given $n > 0$ on the event $\cap_{ k \in (\N)^{2/3} :  k \geq n } \mathcal{F}_k$ we have that \eqref{eq:aD0} holds for all $k \geq n$ which are in $(\N)^{2/3}$.
%, and the implied constants in the $\O(k^{-\epsilon})$ (essentially) are those used to define \eqref{Fevent}.

We proceed by identifying an event of overwhelming probability on which we can extend \eqref{eq:aD0} to hold for all $k \in \R$ with $k \geq n.$
By \eqref{SAi'lim}, this translates to a statement about the stochastic Airy function.  
We must show that in between values of $(\N)^{2/3},$ these processes do not oscillate too much. 
We again use the Dirichlet--Neumann representation \eqref{eq:aD1}, although by Lemma~\ref{lem:A}, we can  represent the solutions solely in terms of $\mathrm{f}_{\lambda,k}$ and $\mathrm{g}_{\lambda,k},$ which produces for $t, n \in \R_+$ and choosing $k \in (\N)^{2/3}$
\begin{equation}\label{eq:aD3}
 \mathrm{f}_{\lambda,n}(t) = 
 \mathrm{g}'_{\lambda,k}(n)
 \mathrm{f}_{\lambda,k}(t)  
 -{\mathrm{f}_{\lambda,k}'(n)}
 \mathrm{g}_{\lambda,k}(t)  .
\end{equation}

We now choose the constant $ \kappa>0$ small enough so that for all $m \in \N$,  
\[
 \kappa  \sqrt{ m^{-2/3}} \le  (m+1)^{2/3} - m^{2/3} \le  \kappa^{-1}\sqrt{ m^{-2/3}} .   
\]
This is the basic reason for working with the mesh $(\N)^{2/3}$. Then, on the event \eqref{Cevent}, we control uniformly both  $ \mathrm{f}'_{\lambda,k}(n),  \mathrm{g}'_{\lambda,k}(n)$ for $n\in[m^{2/3},(m+1)^{2/3}]$ by choosing $k =(m-1)^{2/3}$. 
Therefore, on the event $\mathcal{C}_k \cap\mathcal{F}_k$, it holds for  $n\in[m^{2/3},(m+1)^{2/3}]$  and $t\ge s_n$,
\[ \begin{aligned}
\mathrm{f}_{\lambda,n}(t) e^{-\int_{n}^{t}\theta_{\lambda}(u)\d u } 
 & =    \frac{e^{\O(k^{-\epsilon})}}{2} \mathrm{g}'_{\lambda,k}(n) e^{\int_{k}^{n}\theta_{\lambda}(u)\d u}
 -   \frac{e^{\O(k^{-\epsilon})}}{2 \sqrt{k}}  {\mathrm{f}_{\lambda,k}'(n)} e^{\int_{k}^{n}\theta_{\lambda}(u)\d u} \\
 &= \frac12 \bigg( \cosh((n-k)\sqrt{k}) e^{(n-k)\sqrt{k}+\O(n^{-\epsilon})}
 - \sinh((n-k)\sqrt{k})e^{(n-k)\sqrt{k}+\O(n^{-\epsilon})} \bigg)  . 
%\exp\biggl((s-n)\sqrt{n} + \O_K(n^{-\epsilon}\log n) \biggr)
\end{aligned}\]
where we used that  $\theta_\lambda(t)  = \sqrt{t} + \O(k^{1/2-\epsilon})$ uniformly for $t\ge k$ at the second step and that $(n-k)  \le 2 \kappa^{-1}/\sqrt{k}$ to control the errors. 
This shows that for $m\in \N$, on the event $\bigcap_{k\in (\N)^{2/3} : k\ge (m-1)^{2/3}} \mathcal{B}_k \cap\mathcal{C}_k$, it holds uniformly for $\lambda\in K$ and $n\in[m^{2/3},(m+1)^{2/3}]$, 
\begin{equation} \label{Dasymp}
\sup_{t\ge s_n}\bigg| \mathrm{f}_{\lambda,n}(t) e^{-\int_{n}^{t}\theta_{\lambda}(u)\d u }  -  \frac 12  \bigg|
=\O(n^{-\epsilon}) .
\end{equation}
In particular, taking the limit as  $t\to\infty$, we obtain that $\mathrm{a}_{\lambda,n}^{\mathrm{f}} 
= \tfrac 12 + \O(n^{-\epsilon}) $ almost surely as $\R\ni n\to\infty$.
Similarly, using the representation for  $t, n \in \R_+$, 
\begin{equation}\label{eq:aDg}
 \mathrm{g}_{\lambda,n}(t) = 
- \mathrm{g}_{\lambda,k}(n)
 \mathrm{f}_{\lambda,k}(t)  
 +{\mathrm{f}_{\lambda,k}(n)}
 \mathrm{g}_{\lambda,k}(t)  .
\end{equation}
and proceeding exactly in the same way, we obtain that on the same event, it holds uniformly for $\lambda\in K$ and $n\in[m^{2/3},(m+1)^{2/3}]$, 
\begin{equation} \label{Nasymp}
\sup_{t\ge s_n}\bigg| \sqrt{n}\mathrm{g}_{\lambda,n}(t) e^{-\int_{n}^{t}\theta_{\lambda}(u)\d u }  -  \frac 12  \bigg|
=\O(n^{-\epsilon})
\end{equation}
and that  $\sqrt{n}\mathrm{a}_{\lambda,s}^{\mathrm{g}} = \tfrac 12 + \O(n^{-\epsilon})$ almost surely as $\R\ni n\to\infty$.
By \eqref{gtb}, these estimates imply that with probability at least $1-C_{\epsilon,K,\beta} \exp\big(-n^{3/2-3\epsilon}\big)$, 
\begin{equation} \label{lima}
\sup_{s\in \R : s\ge n} \big| \mathrm{a}_{\lambda,s}^{\mathrm{f}} - \tfrac12\big| = \O(n^{-\epsilon})   
\quad \text{and}\quad 
\ \sup_{s\in \R : s\ge n} \big|  \sqrt{s}\,\mathrm{a}_{\lambda,s}^{\mathrm{g}} -  \tfrac12 \big| = \O(n^{-\epsilon})  .
\end{equation}
By formula \eqref{SAi'lim}, this establishes the first tail-bound \eqref{SAiasymp} and it also shows that almost surely, it holds uniformly for $\lambda \in K$ as $s\to\infty$, 
\[
 \SAi_\lambda'(s) \exp\biggl( \int_{0}^{s}\theta_{\lambda}(u)\d u \biggr) 
  =-\frac{1+\O(s^{-\epsilon}) }{2\sqrt{\pi}} . 
\]

In particular, given the growth properties of the process $\theta$, this establishes that almost surely $\SAi_\lambda' \in L^p(\R_+)$ for any $p\in [1,\infty]$.  
Indeed, according to \eqref{def:cstar}, as $t\to\infty$, 
\[
\int_{0}^{t}\theta_{\lambda}(u)\d u 
 = \tfrac{2}{3} (t + \lambda)^{3/2} +\frac{2}{\sqrt{\beta}}\int_0^t \mathfrak{X}(u)\,\d u  + \big( \tfrac{1}{2\beta}   -\tfrac{1}{4}\bigg) \log(t) +\frac{2c_*}{\beta} +\underset{t\to\infty}{\mathcal{O}_{\beta, \lambda}}(t^{-1/2})  
\]
and we have seen that $\mathfrak{X}(t) =\O(t^{-1/4})$ almost surely as $t\to\infty$. 

\begin{remark} \label{rk:derv}
By definition for $\phi \in \big\{ \mathrm{f}, \mathrm{g} \big\}$, $\phi_{\lambda,n}'(t) = \rho_{\lambda,n}(t)\phi_{\lambda,n}(t)$ for $n, t\in\R$. Since on the event $\mathcal{B}_n$, $\sup_{t\ge n}\big|\frac{\rho_{\lambda,n}(t)}{\sqrt t} -1\big| \le C n^{-\epsilon}$, \eqref{Dasymp} and \eqref{Nasymp} imply that with probability at least $1-C_{\epsilon,K,\beta} \exp\big(-n^{3/2-3\epsilon}\big)$,  it holds for $s\ge n$, 
\[
\sup_{t\ge s}\bigg| \mathrm{f}_{\lambda,s}'(t)\frac{e^{-\int_{s}^{t}\theta_{\lambda}(u)\d u }}{\sqrt t}  -  \frac 12  \bigg|
=\O(s^{-\epsilon}) 
\quad \text{and}\quad
\sup_{t\ge s}\bigg| \mathrm{g}_{\lambda,s}'(t)\frac{e^{-\int_{s}^{t}\theta_{\lambda}(u)\d u }}{\sqrt t}  -  \frac{1}{2\sqrt{s}}  \bigg|
=\O(s^{-1/2-\epsilon}) . 
\]
In particular, this shows that almost surely, uniformly for $\lambda\in K$
\[ \lim_{t\to\infty}\phi_{\lambda,n}'(t) \tfrac{e^{-\int_{n}^{t}\theta_{\lambda}(u)\d u}}{\sqrt t} = \mathrm{a}_{\lambda,n}^{\phi} 
\qquad \text{for } \phi \in \big\{ \mathrm{f}, \mathrm{g}\big\}.
\]

%$\bullet$ Formula \eqref{eq:aD1} also shows that on the event  $\bigcap_{k \in  (\N)^{2/3} : k\ge n} \mathcal{F}_k $,  if $n$ is large enough (depending on a fixed $C\ge 2$), then
%\begin{equation} \label{bdas}
%  \sup_{t\ge n}\Bigl| \phi_{\lambda,n}(t)  e^{-\int_{n}^{t}\theta_{\lambda}(u)\d u } \Bigr|
%  \le C   \Big(  \Bigl| \phi_{\lambda,n}(k)  e^{-\int_{n}^{k}\theta_{\lambda}(u)\d u } \Bigr|
% +\Bigl| k^{-1/2} \phi_{\lambda,n}'(k)  e^{-\int_{n}^{k}\theta_{\lambda}(u)\d u } \Bigr| \Big)
%\end{equation}
%for $k \in  (\N)^{2/3}$ with $k\ge n$ and $\lambda\in K$. 
%Thus, the existence of both  limits on the RHS of \eqref{bdas} establishes that almost surely,
%\begin{equation}\label{eq:aD6}
%\limsup_{n\to\infty}  \sup_{t \geq n, \lambda \in K} \Bigl| \phi_{\lambda,n}(t)  e^{-\int_{n}^{t}\theta_{\lambda}(u)\d u } \Bigr| <\infty.
%\end{equation}
\end{remark}

\paragraph{Existence of $\SAi_\lambda$, its representation in terms of Neumann solutions and tail bounds.}

%According to \eqref{def:cstar}, as $t\to\infty$, 
%\[
%\int_{0}^{t}\theta_{\lambda}(u)\d u 
% = \tfrac{2}{3} (t + \lambda)^{3/2} +\int_0^t \mathfrak{X}(u)\,\d u  + \big( \tfrac{1}{2\beta}   -\tfrac{1}{4}\bigg) \log(t) +\frac{2c_*}{\beta} +\underset{t\to\infty}{\mathcal{O}_{\beta, \lambda}}(t^{-1/2})  
%\]
%with a deterministic error. Moreover, making an integration by parts, the random term is given by 
%\[
%\int_0^t \mathfrak{X}(u)\,\d u = V(u) \mathfrak{X}(u) + \int_0^t V(u) \d B(u)
%\]
%where $V(u) = y(u) \int_u^\infty y^{-1}(s)\d s$.
%By \eqref{eq:yeezy}, $V(u) = 1/\sqrt{u} + \O(u^{-2})$ as $u\to\infty$ and 
%$ \int_0^t V(u) \d B(u) \overset{\rm law}{=} B( \int_0^t V(u)^2  \d u)$ by the martingale representation Theorem.
%So, by Proposition~\ref{prop:Xest} and the law of iterated logarithm, we have almost surely  
%\[
%\lim_{u\to\infty}  V(u) \mathfrak{X}(u) =0
%  \quad
%  \text{and}
%  \quad
%  \limsup_{t \to \infty} \frac{|\int_0^t V(u) \d B(u)|}{\sqrt{\tfrac{2}{\beta} \log t \log\log\log t}} = 1.
%  \]

Hence we can define, as in Definition~\ref{def:SAi},
  \[
    \SAi_\lambda(t) 
    \coloneqq
    -\int_t^\infty \SAi_\lambda'(s)\d s
    =
    \int_t^\infty
    \mathrm{a}_{\lambda,s}^{\mathrm{f}}\exp\biggl( -\int_{0}^{s}\theta_{\lambda}(u)\d u \biggr)/\sqrt{\pi}\d s
    .
  \]
We would like to establish its connection to $\mathrm{a}^{g}_{\lambda,t},$ which will allow us to make statements about its asymptotics.
Recall that by  Lemma~\ref{lem:A}, 
\[
 \mathrm{g}_{\lambda,s}(t)  = \int_s^t  \mathrm{f}_{\lambda,v}(t) \d v, \qquad
  s,t \in\R \text{ and } \lambda\in\C.
\]
This implies that 
\[
  \mathrm{a}_{\lambda,s}^{\mathrm{g}}
  \exp\biggl(-\int_{0}^{s}\theta_{\lambda}(u)\d u \biggr)
  = \lim_{t\to\infty}  \left( \int_s^t  \mathrm{f}_{\lambda,v}(t) \exp\biggl( -\int_{0}^{t}\theta_{\lambda}(u)\d u \biggr)  \d v  \right). 
\]
By \eqref{Dasymp}, observe that $\limsup_{n\to\infty}  \sup_{t \geq n, \lambda \in K} \Bigl| \mathrm{f}_{\lambda,n}(t)  e^{-\int_{n}^{t}\theta_{\lambda}(u)\d u } \Bigr| =1/2$ almost surely
and the function $v \in \R_+\mapsto \exp\biggl( -\int_{0}^{v}\theta_{\lambda}(u)\d u \biggr)$ is integrable.
Then, by  \eqref{SAi'lim}  and the dominated convergence theorem, it holds for $s\in\R$ and $\lambda\in \C$, 
\begin{equation}\label{eq:aD7}
  \mathrm{a}_{\lambda,s}^{\mathrm{g}}
  \exp\biggl(-\int_{0}^{s}\theta_{\lambda}(u)\d u \biggr)
  = - \sqrt{\pi}\int_s^\infty \SAi_\lambda'(v) \d v =  \sqrt{\pi} \SAi_\lambda(s) .
\end{equation}

In addition, it follows from \eqref{lima} that  almost surely as $s\to\infty$, 
\[
  \sqrt{s} \SAi_\lambda(s) \exp\biggl( \int_{0}^{s}\theta_{\lambda}(u)\d u \biggr) 
  = \frac{1+\O(s^{-\epsilon})}{2\sqrt{\pi}}
\]
 uniformly for $\lambda \in K$. 
This also establishes the second tail-bound \eqref{SAiasymp} as well as the fact that $\SAi_\lambda \in L^p(\R_+)$ for any $p\in [1,\infty]$.  
Finally, going back to formula \eqref{eq:aD1}, for $s, t\in\R$
 \begin{equation*}
  \SAi_\lambda(s) 
  = \frac{1}{\sqrt \pi} \lim_{n \to\infty} \bigg\{ \mathrm{g}_{\lambda,s}(n) e^{-\int_{0}^{n}\theta_{\lambda}(u)\d u } \bigg\}
= \frac{1}{\sqrt \pi} \lim_{n\to\infty} \left\{
\mathrm{g}_{\lambda,s}(t) \mathrm{f}_{\lambda,t}(n) e^{-\int_{0}^{n}\theta_{\lambda}(u)\d u } + {\mathrm{g}_{\lambda,s}'(t)}\mathrm{g}_{\lambda,t}(n)  e^{-\int_{0}^{n}\theta_{\lambda}(u)\d u }\right\} . 
\end{equation*}
By \eqref{SAi'lim}, this shows that 
 \begin{equation*}
  \SAi_\lambda(s) 
=
- \mathrm{g}_{\lambda,s}(t)   \SAi_\lambda'(t) + {\mathrm{g}_{\lambda,s}'(t)}   \SAi_\lambda(t) 
\end{equation*}
By Lemma~\ref{lem:A}, $\mathrm{g}_{\lambda,s}(t)   =- \mathrm{g}_{\lambda,t}(s)$ and $\mathrm{g}_{\lambda,s}'(t) = \mathrm{f}_{\lambda,t}(s)$ so we conclude that 
 \begin{equation*}
  \SAi_\lambda(s) 
=\SAi_\lambda(t)   \mathrm{f}_{\lambda,t}(s) +  \SAi_\lambda'(t)   \mathrm{g}_{\lambda,t}(s)
\end{equation*} 
is indeed an entire solution of the stochastic Airy equation. 
\end{proof}

\begin{proof}[Proof of Proposition~\ref{prop:SAiwins}.]
In the context of this proposition, $K \subset \R$ is a fixed compact and $K'= \{z\in\C : \Re z\in K , |\Im z| \le e^{-\epsilon\, T^\nu}  \}$ for $\nu>1$ and $T\gg 1$. We decompose $\phi_{\lambda}$ in the basis $\left\{ \SAi_\lambda, \mathrm{f}_{\lambda,s} \right\}$ to  in terms of the Dirichlet solution at times $s=T/2$. The Wronskian of these solutions (evaluated at $s$) is $-\SAi_\lambda'(T/2).$  Then
\[
\phi_\lambda
=
\Theta_\lambda
\SAi_\lambda
+b_\lambda \mathrm{f}_{\lambda,s}
\quad
\text{where}
\quad
\begin{pmatrix}
  \Theta_\lambda \\
  b_\lambda
\end{pmatrix}
=
\frac{1}{-\SAi_\lambda'(T/2)}
\begin{pmatrix}
  \mathrm{f}'_{\lambda,s}(T) & -\mathrm{f}_{\lambda,s}(T) \\
  -\SAi'_\lambda(T) & \SAi_\lambda(T) \\
\end{pmatrix}
\begin{pmatrix}
  c_{1,\lambda} \\
 c_{2,\lambda} \\
\end{pmatrix}.
\]
In particular,  if $|c_{1,\lambda}|+ |c_{2,\lambda}|/\sqrt{T}  \leq 1$, then
\[
|b_\lambda| \le \frac{\max\big\{ |\SAi_\lambda'(T)|,|\sqrt{T} \SAi_\lambda(T)| \big\}}{|\SAi_\lambda'(T/2)|} .
\]
Using the tail bounds \eqref{SAiasymp} with $\epsilon=0$,  with probability at least $1- \exp(-c(T/2)^{3/2})$, for any  $\lambda\in K$, 
\[
|b_\lambda| \lesssim \exp\bigg( - \int_{T/2}^T \theta_\lambda(u) \d u \bigg) \le \exp(-T^{3/2}/2)
\]
where we  used the event $\mathcal{C}_{T/2}$ given by \eqref{Cevent} to uniformly control $\theta_\lambda$ and that $T\gg 1$. This event has probability at least  $1- C e^{-T^{\nu}}$ if $\nu<3/2$ (choosing $\epsilon>0$ sufficiently small). 

By Corollary~\ref{cor:solbound}, with probability at least $1- e^{-c T^{\nu}}$, we also control the Dirichlet solution: 
\[
|\mathrm{f}_{\lambda,s}(t)| , |\mathrm{f}_{\lambda,s}'(t)| \lesssim
e^{ \delta T^\nu} E_\lambda(t,s)
\qquad\text{for all $t\in [-e^{T},s]$ and   $\lambda\in K_\nu$.}
\]
Using that $E_\lambda(t,s) \lesssim e^{\frac23 (T/2)^{3/2}}= e^{T^{3/2}/3 \sqrt{2}} $, this shows that with probability at least $1- e^{-c T^{\nu}}$, for $\ell=0,1,$
\[
\sup_{t\in [-e^{T},s] , \lambda\in K_\nu } \big|b_\lambda \partial_t^\ell \mathrm{f}_{\lambda,s}(t) \big|  \le     \exp(-T^{3/2}/4) .
\]
We also have the asymptotics: with probability at least $1- e^{-c T^{\nu}}$, there is $\epsilon = \frac12- \frac\nu3>0$ so that uniformly for $t\ge s \gg1$ and $\lambda\in K_\nu$: 
\[
 \mathrm{f}_{\lambda,s}'(t)\frac{e^{-\int_{s}^{t}\theta_{\lambda}(u)\d u }}{\sqrt t}  =  \frac 12+\O(s^{-\epsilon}) , \qquad 
 \mathrm{f}_{\lambda,s}(t) e^{-\int_{s}^{t}\theta_{\lambda}(u)\d u } =  \frac 12+\O(s^{-\epsilon}) ,  
\]
see e.g.~\eqref{Dasymp} and Remark~\ref{rk:derv}. 
By \eqref{SAiasymp}, we also have similar asymptotics for $\SAi_\lambda'(T/2)$ and $\SAi_\lambda(T)$.
This implies that with probability at least $1- e^{-c T^{\nu}}$:
\[
  \Theta_\lambda
  = \frac{c_{2,\lambda}\mathrm{f}_{\lambda,s}(T)- c_{1,\lambda}\mathrm{f}_{\lambda,s}'(T) + \O(T^{-\epsilon})}{\SAi_\lambda'(T/2)}
  = \sqrt{T} \frac{c_{2,\lambda}/\sqrt{T}-c_{1,\lambda} + \O(T^{-\epsilon})}{2\SAi_\lambda'(T/2)e^{-\int_{T/2}^{T}\theta_{\lambda}(u)\d u } }
=  \frac{c_{1,\lambda}-c_{2,\lambda}/\sqrt{T} + \O(T^{-\epsilon})}{2\SAi_\lambda(T)\big(1+\O(T^{-\epsilon})\big)} . 
\]
This concludes the proof. 
\end{proof}

\subsection{Control of the entrance behavior: Proof of Proposition~\ref{prop:entrance}} \label{sec:entrance}

We separately give proofs for the Dirichlet and Neumann cases, which while similar, retain important differences.

\begin{lemma}[Dirichlet case]\label{lem:s-dirichlet-entrance}
  Let $(\sigma_{\lambda,s}(t) : \lambda \in K)$ be the real-part of the Riccati diffusions for $(\mathrm{f}_{\lambda,s} : \lambda \in K)$ with $s\gg 1$. 
For any $\epsilon \in (0, \tfrac 34)$ and $\gamma \in (1, s^{\epsilon}]$,
there is a constant $C=C_{\beta,K,\epsilon}$ so that it holds with probability at least $1-e^{-\frac{s^{3/2-2\epsilon}}{C\log\gamma}}$,
\vspace{-.3cm}
\[
  \big| \sigma_{\lambda,s}(t) - \psi_s(t)\big|
  \leq 
  s^{1/2-\epsilon},
  \quad
  |\varpi_{\lambda,s}(t)| \leq C|\eta| t^{-1/2}
  \quad
  \text{for all $t \in\big[s, s+ \gamma/\sqrt{s}\big]$ and $\lambda \in K.$}
%  \quad
%  \max\left\{ 
%    \frac{|\sigma_{\lambda,s}(t)-\sigma_{\nu,s}(t)|}{|\lambda-\nu|},
%    \frac{|\varpi_{\lambda,s}(t)-\varpi_{\nu,s}(t)|}{|\lambda-\nu|}
%  \right\}
%  \leq Ct^{-1/2}.
\]
\end{lemma}

\begin{proof}
  Let $\kappa = \max_{\lambda \in K}|\Re \lambda|$. 
By \eqref{eq:complexRic}, we note that all Riccati diffusions $\big(\sigma_{\lambda,s}(t): t \geq s,\lambda \in K\big)$ remain above the (strong) solution to the SDE:
  \begin{equation}\label{eq:su-1-s}
    \d\rho_\kappa(t) = -\rho_\kappa^2(t) -\kappa+ t + \tfrac{2}{\sqrt{\beta}}\,\d B(t), \quad \rho_\kappa(s) = 0,
  \end{equation}
  which we define to be absorbed at $-\infty.$
Hence, to obtain a lower-bound, it suffices to estimate the probability that $\rho_\kappa$ remains above $\psi_s(t) -s^{1/2-\epsilon} $ for the short time-interval that we consider. 
Let 
  \[
\tau :=   \inf\Big\{ t \in \big[s,  s+ \tfrac{\gamma}{\sqrt{s}}\big]\, : 
  | \rho_\kappa(t)  - \psi_s(t) | > s^{1/2-\epsilon}   \Big\}.
  \]

By \eqref{eq:omega}, it holds for any $t \in[s, \tau]$, 
  \[ \begin{aligned}
  \max_{\lambda\in K} 
  | \varpi_{\lambda,s}(t) |
  & \le |\eta| \int_s^{t} e^{- 2\int_u^t \rho_\kappa(v)\d v}\,\d u \\
  &  \leq  |\eta| \int_s^{t}  e^{2\gamma s^{-\epsilon}-2\int_u^t \psi_s(v)\,\d v }\,\d u \\
 &   \lesssim_{\epsilon} |\eta| \int_0^{t-s} \biggl(\frac{\cosh( u\sqrt{s})}{\cosh( (t-s)\sqrt{s})} \biggr)^2\,\d u
   \lesssim_{\epsilon} |\eta|  s^{-1/2}.
  \end{aligned}\]
 Hence, this quantity is bounded uniformly for $\lambda \in K$ and by \eqref{eq:complexRic},  all Riccati diffusions $\big(\sigma_{\lambda,s}(t): t \geq s,\lambda \in K\big)$ also remain below the solution of
  \[
    \d\rho_{\kappa'}(t) = -\rho_{\kappa'}^2(t) +\kappa' + t + \tfrac{2}{\sqrt{\beta}}\,\d B(t), \quad \rho_{\kappa'}(s) = 0,
  \]
  for another constant $\kappa'>\kappa$ depending on $K,\epsilon$.  Note that since we look at time $t \leq \tau,$ the process $\rho_\kappa$ has not blown down, and so $\sigma_{\lambda,s}$ for $\lambda \in K$ cannot have jumped above $\rho_{\kappa'}.$
This SDE is of the same type as \eqref{eq:su-1-s}, therefore it suffices to bound the probability of the event $\big\{  \tau <  s+ \tfrac{\gamma}{\sqrt{s}}\big\}$ uniformly over the parameter $\kappa \in [-C_{\epsilon,K} ,C_{\epsilon,K}]$.  
The idea is to compare the diffusion  \eqref{eq:su-1-s} with the solution $\psi_s$ of the ODE \eqref{tanheq}. 
Let us denote $A(t) = \int_s^t (\rho_\kappa(u) + \psi_s(u))\, \d u$.
Using \eqref{tanheq},  we can write \eqref{eq:su-1-s} as
  \begin{equation*}
   \d\big(\rho_\kappa(t) - \psi_s(t) \big) = - \big(\rho_\kappa(t) - \psi_s(t) \big)\d A(t) + 
   \kappa + (t-s) + \tfrac{2}{\sqrt{\beta}}\,\d B(t) . 
  \end{equation*}
Hence, we can write the solution of  \eqref{eq:su-1-s} as
    \begin{equation}
      \rho_\kappa(t) - \psi_s(t) = \int_s^t e^{-A(t)+A(u) } (u-s + \kappa) \d u
    +  \tfrac{2}{\sqrt{\beta}}e^{-A(t)}\int_s^t e^{A(u)}\d B(u).
    \label{eq:su-2}
  \end{equation}
  For all $t\in[s,\tau]$, we have the estimate
    \begin{equation*}
      |A(t) - \hat{A}(t)| \leq (t-s)s^{1/2-\epsilon},
      \quad
      \text{where}
      \quad
      \hat{A}(t) = 2 \int_s^t \psi_s(u) \d u =2\log(\cosh((t-s)\sqrt{s})).
  \end{equation*}
  
This shows that when $\epsilon\in(0,\frac34)$,  $ |A(t) - \hat{A}(t)| \lesssim_\epsilon 1$ for all $t\in[s,\tau]$. In particular, we can bound the mean in \eqref{eq:su-2} by   
   \begin{equation}
    \sup_{ t \in [s,  s+\gamma/\sqrt{s}], |\kappa| \le C_{\epsilon,K} }  \bigg| \int_s^t e^{-A(t)+A(u) } (u-s + \kappa) \d u \bigg| 
    \lesssim_{\epsilon, K }  \int_0^{t-s} e^{- \hat{A}(t)+  \hat{A}(u+s)}(1+u) \d u      \lesssim  \frac{1}{\sqrt{s}} + \frac{\gamma}{s}. 
    \label{eq:su-3}
  \end{equation}

As for the martingale part in \eqref{eq:su-2}, let $M(t) = \tfrac{2}{\sqrt{\beta}}\int_s^{t \wedge \tau} e^{A(u)}\d B(u).$  
The quadratic variation of this martingale is bounded above and below (up to a multiplicative constant) by
  \begin{equation*} 
    \langle M(t)\rangle
    =
    \frac{4}{\beta}\int_s^t e^{2A(u)}\d u
    \asymp_{\beta,\epsilon}
    \int_s^t e^{2 \hat{A}(u)}\d u
    = \int_0^{t-s} \cosh^4( u\sqrt{s})\,\d u . 
  \end{equation*}

On the one-hand, for short times, this immediately implies the sub-Gaussian estimate
\[ \begin{aligned}
 \Big\VERT \sup_{ t \in \big[s,  s+ \frac{1}{\sqrt{s}}\big]} e^{-A(t)} M(t)
  \Big\VERT_2
\lesssim    \sqrt{\big\langle M(s+1/\sqrt{s})\big\rangle} 
 \lesssim_{\beta,\epsilon} s^{-1/4}.
\end{aligned}\] 

On the other-hand, for $t\ge s+ 1/\sqrt{s}$,
\[
    \langle M(t)\rangle \asymp_{\beta,\epsilon}  \int_0^{t-s} \cosh^4( u\sqrt{s})\,\d u \asymp \frac{e^{4(t-s)\sqrt{s}}}{4\sqrt{s}} 
    \quad\text{and}\quad
    e^{-\hat{A}(t)}= \cosh((t-s)\sqrt{s})^{-2} \asymp e^{-2(t-s)\sqrt{s}} . 
\]
By the martingale representation theorem
 $M(t) \overset{\rm law}{=} B_{ \langle M(t)\rangle}$ and the previous estimates imply that
 \[ \begin{aligned}
  \Big\VERT \sup_{ t \in \big[s+ \frac{1}{\sqrt{s}}, s+ \frac{\gamma}{\sqrt{s}}\big]} e^{- A(t)} M(t) \Big\VERT_2
 & \lesssim_\epsilon
    \Big\VERT \sup_{ t \in \big[s+ \frac{1}{\sqrt{s}}, s+ \frac{\gamma}{\sqrt{s}}\big]} e^{- \hat{A}(t)} M(t) \Big\VERT_2 \\
   & \asymp 
       \Big\VERT \sup_{ u\in [1, \gamma]} e^{-2 u } B(e^{4u}/c\sqrt{s})\Big\VERT_2 
 \end{aligned}\]
for a small constant $c = c_{\beta,\epsilon}>0$. 
This comparison with an Ornstein--Uhlenbeck process (using the scaling property of $B$) shows that 
 \[
   \Big\VERT \sup_{ t \in \big[s+ \frac{1}{\sqrt{s}}, s+ \frac{\gamma}{\sqrt{s}}\big]} e^{- A(t)} M(t) \Big\VERT_2  \lesssim_{\beta,\epsilon} s^{-1/4} \sqrt{\log\gamma} . 
 \]
 Combining the previous estimates, we conclude that 
 \[
  \Big\VERT \sup_{ t \in \big[s,  s+ \frac{\gamma}{\sqrt{s}}\big]} \Big|\tfrac{2}{\sqrt{\beta}}e^{-A(t)} \int_s^t e^{A(u)}\d B(u)\Big|
  \Big\VERT  \lesssim_{\beta,\epsilon} s^{-1/4} \sqrt{\log\gamma} . 
 \]
 Using the representation \eqref{eq:su-2} together with the estimate  \eqref{eq:su-3} for the mean, 
 this yields a sub-Gaussian bound for $\sup_{ t \in \big[s,  s+ \frac{2\log s}{\sqrt{s}}\big]}\big|\rho_\kappa(t) - \psi(t)\big| $ which is locally uniform in the parameter $\kappa\in\R$.
 Thus, $\P\big[\tau\le s+ \frac{\gamma}{\sqrt{s}}\big] \le e^{-\frac{s^{3/2-2\epsilon}}{C_{K,\epsilon,\beta}\log\gamma}}$ and as claimed above.  This completes the proof.  
\end{proof}

In the Neumann case, as $\sigma_{\lambda,s}(s)=\infty$, it is convenient to work instead with its inverse $x_\lambda=1/\sigma_\lambda$. 
By \eqref{eq:complexRic}--\eqref{eq:omega} and It\^o formula, we verify that the diffusions 
$(x_{\lambda}(t) : t \geq s,\lambda \in K)$  satisfy the following SDEs:
\begin{equation}
  \begin{aligned}
  \d x_\lambda(t) &= 
  \big(1-(\varpi_\lambda^2(t) + t + \mu)x_\lambda^2(t) + \tfrac{4}{\beta}x_\lambda^3(t)\big)\d t - \tfrac{2}{\sqrt{\beta}}x_\lambda^2(t) \d B(t),  \qquad   x_\lambda(s)=0 ,  \\
\varpi_\lambda(t) &= \eta\int_s^t e^{-\int_u^t 2x_\lambda^{-1}(v)\,\d v}\,\d u ,
  \end{aligned}
  \label{eq:com-inversericatti}
\end{equation}
up to the first time $t>s$ that $x_\lambda(t)$ hits either $0$ or $\infty$. 
Observe that in this interval, $x_\lambda(t)>0$ so that this expression for $\varpi_\lambda$ makes sense and $0\le \varpi_{\lambda,s}(t)\le \eta(t-s)$ for $\eta \ge 0$. 

For $s\gg 1$, one can approximate \eqref{eq:com-inversericatti} on a short time interval by the ODE:
\begin{equation} \label{tanheq2}
\mathfrak{m}' = 1 -s\mathfrak{m}^2 , \qquad \mathfrak{m}(s)=0. 
\end{equation}
The solution being given explicitly by $\mathfrak{m}_s(t):= \tanh((t-s)\sqrt{s})/\sqrt{s}=\psi_s(t)/s $, we obtain the following behavior.

\begin{lemma}[Neumann case]\label{lem:s-neumann-entrance}
For any $\epsilon \in (0, \tfrac 34)$ and $\gamma \in (1, s^{\epsilon}]$,
there is a constant $C=C_{\beta,K,\epsilon}$ so that it holds with probability at least $1-e^{-\frac{s^{3/2-2\epsilon}}{C\log\gamma}}$,
  \[
     \big| x_{\lambda,s}(t) - \mathfrak{m}_s(t)\big|
    \leq \mathfrak{m}_s(t) s^{-\epsilon}
    \qquad
    \text{for all } t \in\big[s, s+ \gamma/\sqrt{s}\big] \text{ and }\lambda \in K.
  \]
\end{lemma}

\begin{proof} 
Fix $\epsilon>0$, let $\mathfrak{m}=\mathfrak{m}_s$ and
  \[
\tau :=   \inf\Big\{ t \in \big[s,  s+ \gamma/\sqrt{s}\big]\, : 
 \exists \lambda\in K,  | x_\lambda(t)  - \mathfrak{m}(t) | > \mathfrak{m}(t) s^{-\epsilon} \Big\}.
  \]
For $t\in[s,\tau]$, the definition of $\tau$ gives $x_{\lambda,s}^{-1}(t)\ge (1+s^{-\epsilon})^{-1}\mathfrak{m}^{-1}(t)$.  Moreover, for $s<u<t$,
\[
\int_u^t\mathfrak{m}^{-1}(r)\,\d r
=\log\!\left(\frac{\sinh((t-s)\sqrt{s})}{\sinh((u-s)\sqrt{s})}\right).
\]
Consequently, writing $v=(t-s)\sqrt{s}$ and using \eqref{eq:com-inversericatti},
\[
\begin{aligned}
|\varpi_{\lambda,s}(t)|
&\le |\eta|\int_s^t
\exp\!\left(-\frac{2}{1+s^{-\epsilon}}\int_u^t\mathfrak{m}^{-1}(r)\,\d r\right)\d u\\
&=\frac{|\eta|}{\sqrt{s}}\int_0^v
\left(\frac{\sinh w}{\sinh v}\right)^{2/(1+s^{-\epsilon})}\d w
\le \frac{|\eta|}{\sqrt{s}}\int_0^v\frac{\sinh w}{\sinh v}\,\d w
\le \frac{|\eta|}{\sqrt{s}}.
\end{aligned}
\]
Thus $|\varpi_{\lambda,s}(t)|\lesssim_K s^{-1/2}$ uniformly for $t\in[s,\tau]$ and $\lambda\in K$.  In particular, $\mu+\varpi_{\lambda,s}^2(t)$ remains in a fixed compact interval.  This sandwiches $(x_{\lambda,s}(t):t\in[s,\tau],\lambda\in K)$ between two solutions of
  \begin{equation}\label{eq:newsandwich}
  \d x_\kappa(t) = 
  \big(1-(t + \kappa )x_\kappa^2(t) + \tfrac{4}{\beta}x_\kappa^3(t)\big)\d t - \tfrac{2}{\sqrt{\beta}}x_\kappa^2(t) \d B(t), \quad x_\kappa(s) =0 , 
\end{equation}
  for $\kappa \in \{\pm C_K\}$. 
 Hence, it suffices to control these two fixed comparison diffusions. 
   Set $A(t) = \int_s^t (x_\kappa(u) + \mathfrak{m}(u))\, \d u$ and observe that
  \[
\d\big(x_\kappa(t)  - \mathfrak{m}(t) \big) =  - s \big(x_\kappa(t)  - \mathfrak{m}(t) \big)\d A(t)
- (t-s+\kappa) x_\kappa^2(t) \d t+ \tfrac{4}{\beta}x_\kappa^3(t)\d t - \tfrac{2}{\sqrt{\beta}}x_\kappa^2(t) \d B(t)
  \]
Hence, by integrating this equation, we obtain
  \begin{equation}\label{eq:x-v}
    \begin{aligned}
x_\kappa(t)  - \mathfrak{m}(t)
    =
    -\int_s^t e^{-s(A(t)-A(u))}\big((u-s + \kappa)x_\kappa^2(u) - \tfrac{4}{\beta}x_\kappa^3(u)\big) \d u
    +
    e^{-sA(t)} M_\kappa(t)
  \end{aligned}
  \end{equation}
where $\displaystyle  M(t) \coloneqq  \tfrac{2}{\sqrt{\beta}}\int_s^t     e^{sA(u)}    x_\kappa^2(u)\,\d B(u)$ is a martingale. 
Moreover, for all $t\in[s,\tau]$, we have the estimate
    \begin{equation*}
       s |A(t) - \hat{A}(t)| \leq s^{1-\epsilon} \hat{A}(t)
      \quad
      \text{where}
      \quad
      \hat{A}(t) = 2 \int_s^t \mathfrak{m}(u) \d u =2\log(\cosh((t-s)\sqrt{s}))/s.
  \end{equation*}
This shows that 
\[
  s |A(t) - \hat{A}(t)| \lesssim \gamma s^{-\epsilon} \lesssim_\epsilon 1
\]
 and we verify that the mean in \eqref{eq:x-v} is controlled by
\[ \begin{aligned}
\sup_{t\in [s,\tau]} \bigg| \int_s^t e^{-s(A(t)-A(u))}\big((u-s + \kappa)x_\kappa^2(u) - \tfrac{4}{\beta}x_\kappa^3(u)\big) \d u \bigg| 
 & \lesssim_{\epsilon,K,\beta} \int_s^t e^{-s(\hat{A}(t)-\hat{A}(u))} \mathfrak{m}(u)^2 \d u \\
 & \le \mathfrak{m}(t)/s , 
\end{aligned}\]
where we used that the function $\mathfrak{m}$ is increasing. 
We now turn to controlling the martingale term, its quadratic variation being given by 
\[
\langle M(t) \rangle = \tfrac{4}{\beta}  \int_s^t     e^{2sA(u)}    x_\kappa^4(u)\,\d u
\asymp_{\epsilon,\beta} \int_s^t     e^{2s\hat{A}(u)} \mathfrak{m}^4(u)\,\d u .
 \]

This shows that $\langle M(t) \rangle  \lesssim_{\epsilon,\beta} \mathfrak{m}(t)^2/s^{3/2}$ for all $ t \in \big[s,  s+ \frac{1}{\sqrt{s}}\big]$ where we used the uniform bound $ \mathfrak{m}(t)\le 1/\sqrt{s}$. 
Hence, for short times, we obtain the sub-Gaussian estimate
\[ \begin{aligned}
 \Big\VERT \sup_{ t \in \big[s,  s+ \frac{1}{\sqrt{s}}\big]} \big| \mathfrak{m}(t)^{-1} e^{- sA(t)} M(t) \big|  \Big\VERT_2
 \lesssim_{\epsilon} 
  \Big\VERT \sup_{ t \in \big[s,  s+ \frac{1}{\sqrt{s}}\big]} \big| \mathfrak{m}(t)^{-1} B_{ \mathfrak{m}(t)^2/s^{3/2}} \big|  \Big\VERT_2
  \lesssim_{\epsilon,\beta}  s^{-3/4}.
\end{aligned}\] 

On the other-hand, for larger time, we have 
$\langle M(t) \rangle \asymp_{\epsilon,\beta}  s^{-5/2} e^{2s\hat{A}(t)}$ and $e^{s\hat{A}(t)} \asymp e^{2(t-s)\sqrt{s}}$. 
Then, by the martingale representation theorem
 $M(t) \overset{\rm law}{=} B_{ \langle M(t)\rangle}$ and we obtain 
  \[ \begin{aligned}
   \Big\VERT \sup_{ t \in \big[s+ \frac{1}{\sqrt{s}}, s+ \frac{\gamma}{\sqrt{s}}\big]} \big| e^{-s A(t)} M(t) \big|  \Big\VERT_2 
& \asymp_\epsilon    \Big\VERT \sup_{u\in[1,\gamma]} e^{-2u} B(e^{4u}/cs^{5/2}) \Big\VERT_2 . 
 \end{aligned}\]
for a small constant $c = c_{\beta,\epsilon}>0$. 
This comparison with an Ornstein--Uhlenbeck process shows that 
 \[
   \Big\VERT \sup_{ t \in \big[s+ \frac{1}{\sqrt{s}}, s+ \frac{\gamma}{\sqrt{s}}\big]} \big| \mathfrak{m}(t)^{-1}e^{-s A(t)} M(t) \big| \Big\VERT_2  \lesssim_{\beta,\epsilon} s^{-3/4} \sqrt{\log \gamma} . 
 \]
 where we used that $\mathfrak{m}(t) \ge c/\sqrt{s}$ for all $t\ge s+1/\sqrt{s}$.

 Consequently, using the representation \eqref{eq:x-v}, these estimates yield control of the sub-Gaussian norm of 
 $\sup_{ t \in \big[s,  s+ \frac{\gamma}{\sqrt{s}}\big]} \big|x_\kappa(t)\mathfrak{m}(t)^{-1} -1\big|$  by $C_{\beta,\epsilon, K} s^{-3/4} \sqrt{\log \gamma}$. 
 Thus, we conclude that 
 $\P[\tau \le s+\gamma/\sqrt{s}] \lesssim e^{-\frac{s^{3/2-2\epsilon}}{C_{\beta,\epsilon, K}\log \gamma}}$ and this proves the claim.
 \end{proof}
 
Our next result upgrades the estimate from Lemma~\ref{lem:s-neumann-entrance} for short times.
 
 \begin{proposition} \label{prop:Entrance_inverse_ric}
 Let $\epsilon\in(0,\tfrac34)$, $\kappa > 0$ be a constant and $n\gg 1$. 
Let us consider the event $\mathsf{Ent}_n^{\mathrm N} \coloneqq \big\{  \big| x_{\lambda,n}(t) - \mathfrak{m}_n(t)\big|
    \leq \mathfrak{m}_n(t) n^{-\epsilon} :  t \in\big[n, n+ \frac{\kappa}{\sqrt{n}}\big], \lambda\in K \big\}.$ 
 On  $\mathsf{Ent}_n^{\mathrm N}$ it holds   $|\varpi_{\lambda,n}(t)| \le \frac{\kappa|\eta|}{\sqrt{n}}$ for $t\in [n,n+ \frac{\kappa}{\sqrt{n}}]$  and  there is a constant $C = C_{K,\epsilon,\kappa, \beta}$ so that
    \[
    \P\left( \mathsf{Ent}_n^{\mathrm N} \cap \Big\{ \sqrt{n}\big| x_{\lambda,n}(n+\tfrac{v}{\sqrt n}) - \mathfrak{m}_{n}(n+\tfrac{v}{\sqrt n}) \big| \le C n^{-\epsilon} v^{5/2-\epsilon} : v\in [0,\kappa], \lambda\in K \Big\}  \right) \ge 1- e^{-n^{3/2-2\epsilon}/ C}
    \]
 \end{proposition}

\begin{proof}
Recall the equations \eqref{eq:com-inversericatti} and that $x_{\lambda,n}=\sigma_{\lambda,n}^{-1}$. Let us also denote $\kappa_n=  n+ \frac{\kappa}{\sqrt{n}}$. 
Without loss of generality, we assume that $\eta\ge 0$ so that the imaginary part of the Riccati diffusion $\varpi_{\lambda,n}\ge 0$. 
On $\mathsf{Ent}_n^{\mathrm N}$, we verify that $\big|x_{\lambda,n}^{-1}(t)-\mathfrak{m}_n^{-1}(t) \big| \le 2 \mathfrak{m}_n^{-1}(t) n^{-\epsilon} $ for  $t \in\big(n, \kappa_n\big]$. 
In particular $x_\lambda(t)>0$ on this interval and  by the discussion following formula \eqref{eq:com-inversericatti}, $|\varpi_{\lambda,n}(t)| \le \frac{\kappa \eta}{\sqrt{n}}$ for $t\in [n,\kappa_n]$ as claimed.

Then, using \eqref{tanheq2}, we can rewrite  \eqref{eq:com-inversericatti} as
\begin{equation} \label{invric1}
\d\big( x_{\lambda,n}(t) - \mathfrak{m}_n(t) \big)
 = n \big(\mathfrak{m}_n^2(t)- x_{\lambda,n}^2(t) \big)\d t   -\big(\varpi_{\lambda,n}^2(t) + (t-n) + \mu- \tfrac{4}{\beta}x_{\lambda,n}(t)\big)x_{\lambda,n}^2(t)\d t - \tfrac{2}{\sqrt{\beta}}x_{\lambda,n}^2(t) \d B(t) . 
\end{equation}
On $\mathsf{Ent}_n^{\mathrm N}$, the bound on $\varpi_{\lambda,n}$ and \eqref{eq:newsandwich} sandwich every $x_{\lambda,n}$ between the two fixed comparison diffusions $x_{a,n}$ with $a\in\{-C_K,C_K\}$.  For either such $a$, let $\vartheta_{a,n}$ be the first $v\in[0,\kappa]$ for which
\[
\big|x_{a,n}(n+v/\sqrt n)-\mathfrak{m}_n(n+v/\sqrt n)\big|>\mathfrak{m}_n(n+v/\sqrt n)n^{-\epsilon},
\]
with $\vartheta_{a,n}=\kappa$ if there is no such $v$, and define the stopped martingale
\[
M_{a,n}(v)=\frac{2\sqrt n}{\sqrt\beta}\int_n^{n+(v\wedge\vartheta_{a,n})/\sqrt n}x_{a,n}^2(u)\,\d B(u).
\]
Integrating the comparison equation \eqref{eq:newsandwich} up to this stopping time gives, for $t\le n+\vartheta_{a,n}/\sqrt n$,
\begin{equation} \label{invric2}
\big|x_{a,n}(t)-\mathfrak{m}_n(t)\big|
\le 2n^{1-\epsilon}\int_n^t\mathfrak{m}_n^2(u)\,\d u
+\frac{2}{\sqrt\beta}\bigg|\int_n^t x_{a,n}^2(u)\,\d B(u)\bigg|.
\end{equation}
The quadratic variation of $M_{a,n}$ is bounded above, up to fixed constants, by
\[
q_n(v)=n\int_n^{n+v/\sqrt n}\mathfrak{m}_n^4(u)\,\d u
=\frac1{n^{3/2}}\int_0^v\tanh(u)^4\,\d u.
\]
For $v\in[0,\kappa]$, $q_n(v)\asymp_\kappa n^{-3/2}v^5$.  The martingale representation theorem and Lemma~\ref{lem:BMestimate} therefore give
\begin{equation} \label{invric3}
\Big\VERT\max_{a\in\{-C_K,C_K\}}\sup_{v\in[0,\kappa]}
\frac{|M_{a,n}(v)|}{v^{5/2-\epsilon}}\Big\VERT_2
\lesssim_{\beta,\kappa,\epsilon,K}n^{-3/4}.
\end{equation}
After rescaling \eqref{invric2}, its deterministic term is bounded by $C_{\kappa,\epsilon}n^{-\epsilon}v^{5/2-\epsilon}$, since
\[
n^{3/2}\int_n^{n+v/\sqrt n}\mathfrak{m}_n^2(u)\,\d u
=\int_0^v\tanh(u)^2\,\d u\le v^3/3
\]
and $v^3\le C_{\kappa,\epsilon}v^{5/2-\epsilon}$ on $[0,\kappa]$.  The estimate \eqref{invric3} controls the two stopped martingales with failure probability at most $e^{-n^{3/2-2\epsilon}/C}$.  Lemma~\ref{lem:s-neumann-entrance}, applied also to the two fixed comparison parameters, shows that neither stopping time occurs early with the same probability scale.  Comparison then transfers the estimate to every $\lambda\in K$; intersecting with $\mathsf{Ent}_n^{\mathrm N}$ proves the proposition.
\end{proof} 

We are now ready to complete our proof of Proposition~\ref{prop:entrance}.
  
  \begin{proof}[Proof of Proposition~\ref{prop:entrance}]
  We let $s_n = n + \frac{\log n}{\sqrt{n}}$ and  $\kappa_n = n+ \kappa/\sqrt{n}$ for a constant $\kappa>0$. 
  Since $\sigma_{\lambda,s}(s)=0/\infty$ for the Dirichlet/Neumann diffusion, we treat the two cases separately. We use Lemma~\ref{lem:s-dirichlet-entrance} and~\ref{lem:s-neumann-entrance} respectively with $\gamma= \log n$.  

\paragraph{Dirichlet case.} On the event of Lemma~\ref{lem:s-dirichlet-entrance}, 
the Riccati diffusion $\rho_{\lambda,n}(t)= \psi_n(t) + \O(n^{1/2-\epsilon}) $ uniformly for $t\in [n, s_n]$ and $\lambda\in K$. In particular, the first claim follows immediately from the fact that
 $\psi_n(s_n) -\sqrt{n}= \big(\tanh(\log n)-1\big)\sqrt{n} =\O(n^{-3/2})$.
Then, on this event, by formula \eqref{Ric2} and since $\mathrm{f}_{\lambda,n}(n)=1$, 
\[
\mathrm{f}_{\lambda,n}(t) = \exp\biggl( \int_{n}^{t} \big(  \psi_n(u) + \O(n^{1/2-\epsilon}) \big) \d u \biggr) 
= \cosh((t-n)\sqrt{n})  e^{\O(n^{-\epsilon}\log n)} 
\]
using that $\partial_t\log\cosh((t-n)\sqrt{n}) = \psi_n(t)$. 
This proves the second claim. For the third claim we use that $\mathrm{f}_{\lambda,n}'(t) =  \rho_{\lambda,n}(t)\mathrm{f}_{\lambda,n}(t) $ and that $\rho_{\lambda,n}(t)= \psi_n(t) e^{\O(n^{-\epsilon})}$ 
for $t\in [\kappa_n, s_n]$. This yields
\[
\mathrm{f}_{\lambda,n}'(t)  = \psi_n(t) \cosh((t-n)\sqrt{n})  e^{\O(n^{-\epsilon}\log n)} 
= \partial_t\big(  \cosh((t-n)\sqrt{n}) \big) e^{\O(n^{-\epsilon}\log n)} .
\]
This  completes the proof in the Dirichlet case with the required uniformity.

\paragraph{Neumann case.}
On the event of Lemma~\ref{lem:s-neumann-entrance}, we verify that
\begin{equation} \label{invric5}
\big|x_{\lambda,n}^{-1}(t)-\mathfrak{m}_n^{-1}(t) \big| \le 2 \mathfrak{m}_n^{-1}(t) n^{-\epsilon} , \quad
\sqrt{n} \le \mathfrak{m}_n^{-1}(t)\le C\sqrt n
\quad\text{for } t \in\big[\kappa_n, s_n\big] 
\quad\text{and}\quad |\varpi_{\lambda,n}(\kappa_n)| \le \frac{\kappa |\eta|}{\sqrt n}  ; 
\end{equation}
see the first claim of Proposition~\ref{prop:Entrance_inverse_ric}. 
Thus, by formula \eqref{eq:com-inversericatti}, we obtain the bounds for $t\in [\kappa_n, s_n]$, 
\[
|\varpi_{\lambda,n}(t)| \le |\varpi_{\lambda,n}(\kappa_n)| + |\eta| \int_{\kappa_n}^t e^{- 2\int_u^t x_{\lambda,n}^{-1}(v)\,\d v}\,\d u 
\le \frac{\kappa|\eta|}{\sqrt n} +  |\eta| \int_{\kappa_n}^t e^{- c\sqrt{n}(t-u)}\,\d u  .
\]

This shows that on the event of Lemma~\ref{lem:s-neumann-entrance}, $\varpi_{\lambda,n}(t)  \lesssim \frac{|\eta|}{\sqrt n} $  for all  $t\in [n, s_n]$ and $\lambda\in K$. 
Since  $\mathfrak{m}_n^{-1}(s_n) -\sqrt{n}= \big(\tanh(\log n)^{-1}-1\big)\sqrt{n} =\O(n^{-3/2})$,  this proves the first claim of Proposition~\ref{prop:entrance} in the Neumann case with the required probability. 

\medskip

Then, by formula \eqref{Ric2}, it holds on the event of Lemma~\ref{lem:s-neumann-entrance} for $n<s<t \le s_n$, 
\[
\mathrm{g}_{\lambda,n}(t) = \mathrm{g}_{\lambda,n}(s)  \exp\biggl( \int_{s}^{t} \big( x_{\lambda,n}^{-1}(u) + \O(n^{-1/2}) \big)\d u \biggr) 
\]
Using that explicitly, $\exp(\int_s^t \mathfrak{m}_n(u)^{-1} \d u )  = \frac{\sinh((t-n)\sqrt{n})}{\sinh((s-n)\sqrt{n})}$,
this gives 
\[
\mathrm{g}_{\lambda,n}(t) = \frac{\mathrm{g}_{\lambda,n}(s)}{\sinh((s-n)\sqrt{n})}  \sinh((t-n)\sqrt{n})\exp\biggl( \int_{s}^{t} \big( x_{\lambda,n}^{-1}(u)  - \mathfrak{m}_n(u)^{-1}  \big)\d u + \O(n^{-\epsilon} \log n) \biggr) 
\] 

By Proposition~\ref{prop:Entrance_inverse_ric}, it holds on the event of Lemma~\ref{lem:s-neumann-entrance} with probability at least $1- e^{-n^{3/2-2\epsilon}/ C}$, for $v\in [0,\kappa]$, 
\begin{equation} \label{invric4}
\frac{1}{\sqrt{n}}\Big| x_{\lambda,n}^{-1}(n+\tfrac{v}{\sqrt n}) - \mathfrak{m}_{n}^{-1}(n+\tfrac{v}{\sqrt n}) \Big| \le
2 \frac{\big| x_{\lambda,n}(n+\tfrac{v}{\sqrt n}) - \mathfrak{m}_{n}(n+\tfrac{v}{\sqrt n}) \big| }{\sqrt{n}\mathfrak{m}_{n}^2(n+\tfrac{v}{\sqrt n})} \le 
\frac{2 n^{-\epsilon} v^{5/2-\epsilon}}{\tanh(v)^2}  \lesssim n^{-\epsilon} v^{1/2-\epsilon}
\end{equation}
Hence, we can take the limit as $s\to n$  since $\mathrm{g}_{\lambda,n}$ is $C^1$ with $\mathrm{g}_{\lambda,n}'(n)=1$, this implies that for $t\in[n,s_n]$
\[
\mathrm{g}_{\lambda,n}(t) = \frac{\sinh((t-n)\sqrt{n})}{\sqrt n}  \exp\biggl( \int_{\kappa_n}^{t} \big( x_{\lambda,n}^{-1}(u)  - \mathfrak{m}_n^{-1}(u)  \big)\d u + \O(n^{-\epsilon}) \biggr) .
\] 
Note that we have restricted the integral to $u\in[\kappa_n, t]$ by using the bound \eqref{invric4}. 
Using \eqref{invric5} again and that 
\[
\int_{\kappa_n}^{s_n}  \mathfrak{m}_n^{-1}(u)  \d u = \int_{\kappa}^{\log n} \frac{\d u}{\tanh(u)} = \O(\log n)
\]
we conclude that on the event of Lemma~\ref{lem:s-neumann-entrance}, with probability at least $1- e^{-n^{3/2-2\epsilon}/ C}$,  it holds  uniformly for $t\in[n,s_n]$ and $\lambda\in K$, 
\[
\mathrm{g}_{\lambda,n}(t) = \frac{\sinh((t-n)\sqrt{n})}{\sqrt n}  e^{\O(n^{-\epsilon} \log n)} 
\]

For the third claim we use that $\mathrm{g}_{\lambda,n}'(t) =  \rho_{\lambda,n}(t)\mathrm{g}_{\lambda,n}(t) $ and that $\rho_{\lambda,n}(t)= \mathfrak{m}_n^{-1}(t) e^{\O(n^{-\epsilon})}$ 
for $t\in [n, s_n]$ so that 
\[
\mathrm{g}_{\lambda,n}'(t) = \mathfrak{m}_n^{-1}(t) \frac{\sinh((t-n)\sqrt{n})}{\sqrt n}  e^{\O(n^{-\epsilon} \log n)} 
= \partial_t\bigg(\frac{\sinh((t-n)\sqrt{n})}{\sqrt n} \bigg)  e^{\O(n^{-\epsilon} \log n)} . 
\]
This completes the proof in the Neumann case with the required uniformity.
 \end{proof}

  \subsection{Estimates for the Gaussian \texorpdfstring{process}{process} \texorpdfstring{$\mathfrak{X}$}{X}: Proof of Proposition~\ref{prop:Xest}} \label{sec:X}

Let us recall that for $\mu\in\R$ and  $s\ge -\mu$, we have defined
\[
\mathfrak{X}_{\mu, s}(t) = y_\mu^{-1}(t) \int_{s}^t y_\mu(u)  \d B(u) 
\quad\text{where}\quad y_\mu(t) =  \exp\big(\tfrac43 (t+\mu)^{3/2}  \big). 
\]  
Many of our estimates will be based on $t\mapsto y_\mu(t) \mathfrak{X}_{\mu,s}(t)$ being a Gaussian martingale, and for this reason we define for real $\mu,t$,
\begin{equation} \label{eq:G}
  G_\mu(t) \coloneqq 
  \int_0^t y_\mu^2(u) \d u
  =\int_\mu^{t+\mu} y_0^2(u) \d u,
\end{equation}
and we record some basic estimates for it.
\begin{lemma}\label{lem:G}
\begin{enumerate}[label=(G-\roman*), ref=(G-\roman*), leftmargin=*]
  \item \label{eq:Gupper}
    For all $t,\mu \geq 0,$
    \(
    \,\,G_\mu(t) \leq \frac{ y_\mu^{2}(t) }{\max\{\sqrt{t+\mu},1\}},
    \)
    and hence $\log(G_\mu(t)) \leq \tfrac{8}{3}(t+\mu)^{3/2}.$
  \item \label{eq:GlogGlower}
    For all $t,\mu \geq 0,$
    \(
    G_\mu(t)\bigl(\log(1+G_\mu(t))\bigr)^{1/3} \leq 2{y_\mu^{2}(t)}.
    \)
  \item \label{eq:Glower}
    For all $t,\mu \geq 0$ with $t \geq \min\{1,\mu^{-1/2}\},$
    \(
    G_\mu(t)
    \geq
    \frac{y_\mu^2(t)e^{-4}}{\sqrt{t+\mu}},
    \)
    and if in addition $t+\mu \geq 4,$ then
    \(
    \log G_\mu(t) \geq (t+\mu)^{3/2}.
    \)
\end{enumerate}
\end{lemma}
\begin{proof}

We begin by verifying \ref{eq:Gupper}.
For $\sqrt{t+\mu} \leq 1,$ we have $t \in [0,1],$ and  
\begin{equation}\label{eq:Arnault2}
 G_\mu(t) = \int_0^t y_\mu^{2}(u) \d u
 \leq y_\mu^{2}(t)
 = \frac{ y_\mu^{2}(t) }{\max\{\sqrt{t+\mu},1\}}.
\end{equation}
Thus it remains to consider the case that $t+\mu \geq 1.$
For such $t$ we verify that
\[
  \frac{\d}{\d t}
  \biggl(G_\mu(t) - \frac{ y_\mu^{2}(t) }{\sqrt{t+\mu}}\biggr)
  =
  \biggl(1-4+ \frac{1}{2(t+\mu)^{3/2}}\biggr)y_\mu^{2}(t)
  <0.
\]
Hence, it suffices to verify that $G_\mu(t_0) \leq \frac{ y_\mu^{2}(t_0)}{\sqrt{t_0+\mu}},$ for  $t_0 = (1-\mu)_+$.
In case $t_0 = 0$, we trivially have $G_\mu(t_0) =0 \leq \frac{ y_\mu^{2}(0)}{\sqrt{\mu}}$ since $\mu\ge 1$.  
In the other case, we may apply the bound \eqref{eq:Arnault2} at time $t_0$ to get the desired conclusion. This completes the proof of \ref{eq:Gupper}.

For \ref{eq:GlogGlower}, we consider first the case that $G_\mu(t) \le 1.$  In this case, $t\le 1$ and we have using \ref{eq:Arnault2}
\[
  G_\mu(t)\bigl(\log(1+G_\mu(t))\bigr)^{1/3} 
  \leq G_\mu(t)\bigl(\log 2\bigr)^{1/3} 
  \leq {y_\mu^{2}(t)}.
\]
If $G_\mu(t) > 1$ then multiplying \ref{eq:Gupper} by $(\log G_\mu(t))^{1/3}$ we arrive at
\[
  \begin{aligned}
  G_\mu(t)
  \bigl(\log(1+G_\mu(t))\bigr)^{1/3}
  &\leq
  G_\mu(t)
  \bigl(\log 2 + \log(G_\mu(t))\bigr)^{1/3} \\
  &\leq
  G_\mu(t)
  \bigl( (\log 2)^{1/3} + (\log(G_\mu(t)))^{1/3}\bigr) \\
  &\leq
  y_\mu^2(t)
  \bigl( (\log 2)^{1/3} + (\tfrac 83)^{1/3}\bigr) \\
  \end{aligned}
\]
where we have used subadditivity of the map $x \mapsto x^{1/3}$ and \ref{eq:Gupper}.  This numerical constant is less than $2.$

Finally for \ref{eq:Glower}, let us observe that for any $t_1 \in [0,t)$,
\[
G_\mu(t) \geq G_\mu(t_1) + (t-t_1) y_\mu^2(t_1) . 
\]
Under the condition $t \geq \min\{1,\mu^{-1/2}\},$ 
we may take $t_1 = t-(t+\mu)^{-1/2} \ge 0$.
Dropping $G_\mu(t_1) \ge 0$, this implies that
\[
G_\mu(t) \geq (t+\mu)^{-1/2} y_\mu^2(t)  \exp\biggl(- 4 \int_{t_1}^t \sqrt{u+\mu} \, \d u
  \biggr)  \geq (t+\mu)^{-1/2} y_\mu^2(t)  \exp\big( -4(t-t_1) \sqrt{t+\mu} \big) 
\]
Given our choice of $t_1$, this proves the first bound.  The final numerical bound follows from the inequality $x^{-1}e^{5/3 x^3 - 4} \geq 1$ for all $x \geq 2.$
\end{proof}

\begin{lemma} \label{lem:XGN}
For any $\epsilon \in (0,1]$ and $s\ge 0$, 
\[
\Big\VERT \sup_{t\ge s} (1+t)^{1/4-\epsilon}| \mathfrak{X}_{0, s}(t)| \Big\VERT_2
\lesssim (1+s)^{-\epsilon} \log(2+ s) . 
\]
\end{lemma}

\begin{proof}
Fix $\gamma <1/4$ and let $y=y_0$. By the martingale representation theorem, 
\[
\big( (1+t)^\gamma\mathfrak{X}_{0, s}(t) \big)_{t\ge s} \overset{\rm law}{=} 
\big( (1+t+s)^{\gamma} y^{-1}(t+s) B(G_s(t)) \big)_{t\ge 0} 
\]
where $G_s$ is as in \eqref{eq:G}. 
In particular, for any $\alpha>0$ 
\begin{equation}\label{eq:Arnault0} 
  \begin{aligned}
\Big\VERT \sup_{t\ge s} (1+t)^\gamma| \mathfrak{X}_{0, s}(t)| \Big\VERT_2
 & \le  \max_{x\ge s}\big( (1+x)^\gamma  y^{-1}(x)\big)
 \Big\VERT \sup_{t\in [0,\alpha]} |B(G_s(t))| \Big\VERT_2\\
&\quad
+ \Big\VERT \sup_{t\ge \alpha}| (1+t+s)^\gamma y_s^{-1}(t) B(G_s(t))|  \Big\VERT_2.
%& \lesssim \beta^{-1}  (s+1)^\gamma e^{-4 s^{3/2}/3} \sqrt{G_s(\alpha)} 
%+ \Big\VERT \sup_{t\ge \alpha}| (1+t+s)^\gamma 
%\tfrac{B(G_s(t))}{
%  G_s(t)^{1/2}\bigl(\log G_s(t)\bigr)_+^{1/6}}  \Big\VERT_2 \\
%\Big\VERT \sup_{u\ge G(\alpha)} |(1+G^{-1}(u)+s)^\gamma y_s^{-1}(G^{-1}(u)) B(u)|  \Big\VERT_2 . 
\end{aligned}\end{equation}
Provided $t \geq \min\{1,s^{-1/2}\}$ and that $t+s \geq 4$, we have by \ref{eq:Glower}  
\[
  (1+t+s)^\gamma \leq (5/4)^\gamma (t+s)^\gamma 
  \leq 2 \bigl(\log(1+ G_s(t))\bigr)^{2\gamma/3}.
\]
Similarly, from \ref{eq:GlogGlower}, we have
\[
  y_s^{-1}(t) \leq 2 G_s^{-1/2}\bigl(\log(1+G_s(t))\bigr)^{-1/6},
\]
and therefore from \eqref{eq:Arnault0},
\begin{equation*}\label{eq:Arnault1}
  \begin{aligned}
\Big\VERT \sup_{t\ge s} (1+t)^\gamma| \mathfrak{X}_{0, s}(t)| \Big\VERT_2 
\le 
\max_{x\ge s}\big( (1+x)^\gamma  y^{-1}(x)\big) \Big\VERT \sup_{t\in [0,\alpha]} |B(G_s(t))| \Big\VERT_2 + 4\Big\VERT \sup_{t\ge \alpha} \biggl| \frac{B(G_s(t)) / \sqrt{G_s(t)}} {  \bigl( \log(1+ G_s(t))\bigr)^{1/6-2\gamma/3}} \biggr|  \Big\VERT_2
\end{aligned}\end{equation*}
Applying Lemma \ref{lem:BMestimate}, for any $\gamma \in (0,\tfrac14)$ and $\alpha \geq 4$ so that $G_s(\alpha)\ge 3$ and we can use the previous estimates, we obtain
\[
\Big\VERT \sup_{t\ge s} (1+t)^\gamma| \mathfrak{X}_{0, s}(t)| \Big\VERT_2 
\lesssim_{\gamma,\beta}
(s+1)^\gamma e^{-4 s^{3/2}/3}\sqrt{G_s(\alpha)}
+
\bigl(\log G_s(\alpha)\bigr)^{\tfrac23(\gamma-\tfrac 14)} \log(\log G_s(\alpha)).
\]
Taking for instance $\alpha$ so that $G_s(\alpha) = c e^{8 s^{3/2}/3}/ \sqrt{s+1}$ for a constant $c>0$ so that $\alpha\ge 4$ (which can be done using \ref{eq:Gupper} and \ref{eq:Glower}), we conclude the bound claimed in the lemma (by letting $\gamma=1/4-\epsilon$). 
\end{proof}

%We now extend the bound from Lemma~\ref{lem:XGN} so that it holds locally uniformly for  the parameter $\mu\in\R$. 
We now control the oscillations of the process in the parameter $\mu.$ By definition, the family of processes $(\mathfrak{X}_{\mu, n})_{\mu\in[-\kappa,\kappa]}$ are solutions of the coupled SDEs:
\[
\d\mathfrak{X}_{\mu, n}(t)= - 2 \sqrt{t+\mu}\, \mathfrak{X}_{\mu, n}(t)  \d t + \d B(t), \qquad \mathfrak{X}_{\mu, n}(n) =0 . 
\]

Thus the difference $\mathfrak{X}_{\mu, n}(t)-\mathfrak{X}_{0, n}(t)$ satisfies an ODE, which when solved, implies that for any $\mu\in[-\kappa,\kappa]$, 
\begin{equation} \label{XSDE_rep}
\mathfrak{X}_{\mu, n}(t) = \mathfrak{X}_{0,n}(t) +  \mathcal{R}_{\mu,n}(t) , \qquad t\ge n,
\end{equation}
where 
\[
\mathcal{R}_{\mu,n}(t) \coloneqq  y_\mu^{-1}(t)  \int_{n}^t y_\mu(u)  \mathfrak{X}_{0, n}(u) \big(\sqrt{u+\mu} -\sqrt{u}\big)  \d u . 
\]

Moreover, we can use Lemma~\ref{lem:XGN} to obtain uniform tail bounds for the processes $(\mathcal{R}_{\mu,n})_{\mu\in[-\kappa,\kappa]}$. 

\begin{lemma} \label{lem:Ups}
For any $\epsilon \in (0,1]$, $\kappa\ge 0$ and $n\ge \kappa$, 
\[
\Big\VERT \sup_{t\ge n, \mu \in [-\kappa,\kappa]} (1+t)^{5/4-\epsilon} |\mathcal{R}_{\mu,n}(t)|  \Big\VERT_2
\lesssim_{\kappa} (1+n)^{-\epsilon} \log(2+ n) . 
\]
\end{lemma}

\begin{proof}
 For any $\mu \in [-\kappa,\kappa]$ and $t\ge n \ge \kappa$ and $\gamma <1/4$ 
\[ \begin{aligned}
|\mathcal{R}_{\mu,n}(t)| & \lesssim_{\kappa,\gamma}
\sup_{u\ge n}\Big\{ \big| (1+u)^\gamma |\mathfrak{X}_{0, n}(u) \big| \Big\}
y_\mu^{-1}(t)  \int_{n}^t y_\mu(u)  (1+u)^{-\gamma-1/2} \d u . 
\end{aligned}\]
By \eqref{eq:yeezy} and  Lemma~\ref{lem:XGN}, taking norms, this implies that 
\[
\Big\VERT \sup_{t\ge n} (1+t)^{\gamma+1} |\mathcal{R}_{\mu,n}(t)|  \Big\VERT_2
 \lesssim_{\kappa,\gamma}(1+n)^{\gamma-1/4} \log(2+ n) .  \qedhere
\]
\end{proof}

Using the representation \eqref{XSDE_rep} and the previous Lemma, we obtain the uniform tail bound, for $\epsilon\in (0,1]$
\[
\Big\VERT \sup_{t\ge n, \mu \in [-\kappa,\kappa]} (1+t)^{1/4-\epsilon} |\mathfrak{X}_{\mu,n}(t)|  \Big\VERT_2
\lesssim_{\kappa} (1+n)^{-\epsilon} \log(2+ n) . 
\]
This completes the proof of Proposition~\ref{prop:Xest}. 

\subsection{Deterministic control of the processes \texorpdfstring{$\Delta$}{Delta} and \texorpdfstring{$\varpi$}{varpi}: Proof of Proposition~\ref{lem:sig}} \label{sec:Delta}

In this section, we fix $\epsilon \in (0,\tfrac34)$ and work on the event \eqref{Eevent}
\[
    \mathcal{E}_{n}= \left\{   \big| \varpi_{\lambda,n}(s_n)\big|  \le c n^{-1/2}  , \big| \Delta_{\lambda,n}(s_n)\big|  \le c n^{1/2-\epsilon}   :  \lambda\in K \right\} \cap \left\{  
   \big|\mathfrak{X}_{\mu,n} (t)| \le c t^{1/2-\epsilon} :   t \geq n, \mu \in \Re(K) \right\} .
\]
where $s_n = n + \tfrac{\log n}{\sqrt{n}}$ and $n\gg 1$.
We define a stopping time, 
  \begin{equation*}
    \tau :=
    \inf\left\{ 
      t \geq s_n : \exists\lambda\in K\text{ such that }
      | \varpi_{\lambda,n}(t) | >     C t^{-1/2}
      \text{ or }
      | \Delta_{\lambda,n}(t) | >     C t^{1/2-\epsilon}
    \right\}.
  \end{equation*}
It suffices to show that we can choose a $C > c$ sufficiently large depending on $c, \epsilon, \beta$ and $K$ so that $\tau = \infty.$ We suppress the $\lambda,n$ subscripts and write $\varpi = \varpi_{\lambda,n}, s=s_n$ etc.
Let us also observe that according to the definitions \eqref{delta0}, we can assume that for $t \in[s,\tau]$, 
\[
|\mathfrak{h}(t)| \le 2 c t^{1/2-\epsilon}
 \qquad\text{and}\qquad
 \big|\sigma(t)- \sqrt{t}\big| \le 2C t^{1/2-\epsilon}
\]
so that $\sigma(t) \ge  3\sqrt{t}/4$   for $n\gg 1$.

  \paragraph{The imaginary part.}
By  \eqref{eq:omega}, we have on $\mathcal{E}_{n}$ for $t \in[s,\tau]$, 
    \[
    | \varpi(t)| \le     | \varpi(s_n)| e^{-2\int_{s}^t \sigma(v)\,\d v}  +|\eta| \int_{s}^t e^{-2\int_u^t \sigma(v)\,\d v}\,\d u
    \le \frac{c}{\sqrt{n}} e^{- t^{3/2} + s^{3/2}} + |\eta| \int_{s}^t e^{-t^{3/2}+u^{3/2}}\,\d u
 \] 
Hence, choosing $n, C$ large enough, we can ensure that $| \varpi(t)| \le  \frac{C}{2\sqrt{t}}$. 
  
    \paragraph{The linearization error.}
  According to \eqref{eq:if},  for $t \in[s,\tau]$, 
  \[
   v(t)
  =\exp\biggl( \int_s^t \big(\sqrt{u+\mu} + \sigma_{\lambda,n}(u) + \mathfrak{h}_{\mu,n}(u) \big)\d u \biggr)
  \ge  e^{t^{3/2}-s^{3/2}}
  \]
 by using that $\sqrt{u+\mu} + \mathfrak{h}_{\mu,n}(u)  \ge 3\sqrt{u}/4$  for $u\ge n$  if  $n\gg 1$. 
  By formula \eqref{delta2}, this implies that for $t \in[s,\tau]$, 
  \[
  \big| \Delta(t)- v(t)^{-1}\Delta(s) \big| \le \int_{s_n}^t \tfrac{v(u)}{v(t)} \max\big\{ \varpi(u)^2 ,  \mathfrak{h}(u)^2 \big\} \d u
  \le \int_{s_n}^t  e^{u^{3/2}-t^{3/2}} \max\big\{ \tfrac{C^2}{u} , 4c^2 u^{1-2\epsilon} \big\} \d u
  \]
  Again, if $n$ is sufficiently large, we can ensure that the RHS is bounded by $t^{1/2-\epsilon}$  for $t \in[s,\tau]$ --- we even have a gain.  
On the other hand, on $\mathcal{E}_{n}$, $v(t)^{-1}|\Delta(s)| \le c n^{1/2-\epsilon} e^{- t^{3/2} + s^{3/2}} $ so we can ensure that $|\Delta(t)| \le C t^{1/2-\epsilon}/ 2$ for all  $t \in[s,\tau]$.

This argument shows that on $\mathcal{E}_{n}$, $| \varpi(t)|\le C t^{-1/2}/2$ and $|\Delta(t)| \le C t^{1/2-\epsilon}/ 2$ for all  $t \in[s,\tau]$. Hence, by continuity of  the processes $\varpi, \Delta$, we conclude that $\tau=\infty$.

\paragraph{The positive part.}
Let $s'=n+\frac{4\log n}{\sqrt n}$. Applying \eqref{delta2} with $s=s_n$ and discarding the nonpositive contribution from $-\mathfrak{h}(u)^2$, we obtain, for $t\ge s'$,
\[
\bigl(\Delta(t)\bigr)_+
\leq
v(t)^{-1}\bigl(\Delta(s_n)\bigr)_+
+\int_{s_n}^t \frac{v(u)}{v(t)}\varpi(u)^2\,\d u
\lesssim
n^{1/2-\epsilon}e^{-t^{3/2}+s_n^{3/2}}
+\int_{s_n}^t e^{u^{3/2}-t^{3/2}}\frac{\d u}{u}.
\]
Since $s'-s_n=3\log n/\sqrt n$, the first term is $\O(t^{-3/2})$ for $t\ge s'$. The elementary Laplace estimate
\[
\int_{s_n}^t e^{u^{3/2}-t^{3/2}}\frac{\d u}{u}
\lesssim t^{-3/2}
\]
gives the same bound for the second term and proves the final assertion of Proposition~\ref{lem:sig}.
\qed

\subsection{Control of the integrated errors: Proof of Proposition~\ref{prop:X2int}}\label{sec:X2}

Our goal to prove Proposition~\ref{prop:asympDN} is to obtain an approximation for the integral of the Riccati diffusions $\int_{s_n}^{t}\rho_{\lambda,n}(u)\d u$ uniformly for $t\ge s_n$ and $\lambda\in K$, up to a vanishing as $n\to\infty$. 
We further want to achieve this control with overwhelmingly good control of the failure probability. 
By formula \eqref{rhoint1}, this amounts to control the integrals $ \int_{s_n}^{t}\Delta_{\lambda,n}(u) \d u$ of the linearization error and of the imaginary part of the Riccati diffusion $ \int_{s_n}^{t} \big( \varpi_{\lambda,n}(u) - \frac{\eta}{2\sqrt{u}}\big) \d u$.  
As explained in Section~\ref{sec:asymp},  the issue is to obtain concentration bounds for integrals of the squared Gaussian processes  $\mathfrak{X}_{\mu,n}^2$. 
The required bounds are formulated as Lemma~\ref{lem:X2} and~\ref{lem:X3}. 
After proving these results, we proceed to the proofs of Proposition~\ref{prop:2021}, Proposition~\ref{prop:2022} and finally Proposition~\ref{prop:X2int}.

\begin{lemma}\label{lem:X2}
Recall that $\mathfrak{X}_{n}(t) = y_n^{-1}(t) \int_{0}^t y_n(u)  \d B(u)$ where $y_n(t) = \exp\big(\frac 43(t+n)^{3/2} \big)$ for $t,n \ge 0$. 
For any $\epsilon\in(0,1]$, there is a constant $C_\epsilon>0$ so that for any $n\ge 1$, 
\[
    \Pr\biggl(
  \sup_{T \geq 0} 
\bigg| \int_0^{T}  \mathfrak{X}_{n}^2(t) \frac{\d t}{2\sqrt{t+n}} - \int_{n}^{T+n}\int_0^t \Exp[\mathfrak{X}(t)\mathfrak{X}(u)] \d u  \d t \bigg|  \geq C_\epsilon n^{-\epsilon}
  \biggr)
  \leq
\exp(-  n^{3/2-2\epsilon}) 
\]
where $ \mathfrak{X} =  \mathfrak{X}_0$. 
\end{lemma}

\begin{proof}
We recall $G_\mu(t) \coloneqq \int_0^t y_\mu^2(u) \d u$  so that for $t \ge  u\ge 0$
\begin{equation} \label{mean2}
  \Exp[\mathfrak{X}_\mu(t)\mathfrak{X}_\mu(u)]  = \frac{G_\mu(u)}{y_\mu(t)y_\mu(u)} =  \frac{y_\mu(u)}{y_\mu(t)}    \E\mathfrak{X}_{\mu}^2(u) .
\end{equation}
Throughout the proof, we are going to use that $\E\mathfrak{X}_{\mu}^2(t) = \O(1)$ uniformly for all $t,\mu\ge0$ and the bounds for $k=1,2,$
\[
\int_\mu^\infty y_0^{-k}(t) \d t  \le \frac{ y_0^{-k}(\mu)}{\sqrt{\mu}}  \, , \qquad
\int_0^T  y_0^k(t) \d t   \le C \frac{ y_0^k(T)}{\sqrt{T}} . 
\]

The proof is divided in three parts, we first control the mean and then show that the fluctuations are at most of size $\mu^{-\epsilon}$ with the required probability using basic concentration inequalities for  martingales.

\paragraph{Step 1: Contribution from the mean.}
Our first goal is to show that  there is a  constant $C$ depending only on $\beta$ so that 
\begin{equation} \label{mean0}
 \sup_{T \geq 0} 
\bigg| \int_0^{T} \E[\mathfrak{X}_{\mu}^2(t)] \frac{\d t}{2\sqrt{t+\mu}} - \int_\mu^{T+\mu}\int_0^t \Exp[\mathfrak{X}(t)\mathfrak{X}(u)] \d u  \d t \bigg|   \le \frac{C}{\mu} 
\end{equation}
First observe that, since $\E[\mathfrak{X}^2(t)]= \O(1)$, we have
\[
\int_\mu^T\int_0^\mu  \Exp[\mathfrak{X}(t)\mathfrak{X}(u)] \d u \d t
 \lesssim   \int_\mu^\infty y_0^{-1}(t) \d t \int_0^\mu y_0(u) \d u  \lesssim \mu^{-1} . 
 \]
 Hence we can rewrite
\begin{equation} \label{mean1}
\int_{\mu}^{T+\mu}\int_0^t \Exp[\mathfrak{X}(t)\mathfrak{X}(u)] \d u  \d t 
= \int_0^{T}\int_0^{t} \Exp[\mathfrak{X}(t+\mu)\mathfrak{X}(u+\mu)] \d u  \d t  + \O(\mu^{-1}) 
\end{equation}
 Then, according to \eqref{mean2}, it holds 
  \[
\E[\mathfrak{X}(t+\mu)^2]  - \E[\mathfrak{X}_\mu(t)^2]  
= y_\mu^{-2}(t)  \big(G_0(t+\mu)-G_\mu(t) \big)   
= y_\mu^{-2}(t)  \int_0^\mu y_0^2(u) \d u 
= y_\mu^{-2}(t) y_0^{2}(\mu)  \O(\mu^{-1/2}).
\]
so that 
\[
\Exp[\mathfrak{X}(t+\mu)\mathfrak{X}(u+\mu)] 
=  \E[\mathfrak{X}_\mu(t)^2]   \frac{y_\mu(u)}{y_\mu(t)}  
+  \frac{y_0^{-1}(\mu+t)  y_0^{-1}(\mu+u)}{y_0^{-2}(\mu)  }   \O(\mu^{-1/2}) .
\]
 
By Fubini's Theorem and \eqref{mean1}, this implies that 
 \[\int_{\mu}^{T+\mu}\int_0^t \Exp[\mathfrak{X}(t)\mathfrak{X}(u)] \d u  \d t 
 = \int_0^{T}   \E[\mathfrak{X}_\mu(u)^2]  \int_u^{T}   \frac{y_\mu(u)}{y_\mu(t)} \d t \d u + \O(\mu^{-1}).
\]
Using that 
$\displaystyle \int_u^T y_\mu^{-1}(t) \d t  = \frac{y_\mu^{-1}(u)}{2\sqrt{u+\mu}} -   \frac{y_\mu^{-1}(T)}{2\sqrt{T+\mu}} +   \O\bigg(\frac{y_\mu^{-1}(u)}{(u+\mu)^2} \bigg)$ and $ \E[\mathfrak{X}_\mu(t)^2]  = \O(1)$ uniformly for $t\ge 0$, this shows that
\[
\int_{\mu}^{T+\mu}\int_0^t \Exp[\mathfrak{X}(t)\mathfrak{X}(u)] \d u  \d t =  \int_0^T \frac{ \E[\mathfrak{X}_{\mu}^2(u)]}{2\sqrt{u+\mu}}  \d u  -  \frac{y_\mu^{-1}(T)}{\sqrt{T+\mu}} 
\O(1)\int_0^{T}  y_\mu(u) \d u + \O(\mu^{-1}) . 
\]
Hence, we conclude that 
\[
\int_{\mu}^{T+\mu}\int_0^t \Exp[\mathfrak{X}(t)\mathfrak{X}(u)] \d u  \d t =  \int_0^T \frac{ \E[\mathfrak{X}_{\mu}^2(u)]}{2\sqrt{u+\mu}}  \d u + \O(\mu^{-1}) .
\] 
This completes the proof of \eqref{mean0}.

\paragraph{Step 2: Reduction to a martingale concentration problem.}
Therefore it suffices to prove that if $\mu$ is large enough,  for any $\epsilon\in(0,1]$, there is a constant $c_\epsilon$ (depending on $\epsilon,\beta$) so that
\begin{equation*} 
    \Pr\biggl(
  \sup_{T \geq 0} 
\bigg| \int_0^{T} \Big(  \mathfrak{X}_{\mu}^2(t) - \E \mathfrak{X}_{\mu}^2(t) \Big) \frac{\d t}{2\sqrt{t+\mu}}  \bigg|  \geq C_\epsilon  \mu^{-\epsilon}
  \biggr) \le \exp(-\mu^{3/2-2\epsilon}) 
 \end{equation*}
 As processes $\big(\mathfrak{X}_\mu(t) \big)_{t\ge 0} \overset{\rm law}{=}   \big( y_\mu^{-1}(t) B(G_\mu(t)) \big)$.  To obtain the display above, it suffices to show that 
  \begin{equation} \label{Mbound}
    \Pr\Biggl(
  \sup_{T \geq 0} 
\bigg| \int_0^{T} M(G_\mu(t)) y_\mu^{-2}(t)  \frac{\d t}{2\sqrt{t+\mu}}  \bigg|  \geq C_\epsilon \mu^{-\epsilon}
  \Biggr) \le \exp(- \mu^{3/2-2\epsilon}) 
 \end{equation}
 where $M(t) \coloneqq B(t)^2-t$ is a continuous martingale. 
 Let us denote $ \displaystyle w_\mu(t) \coloneqq \int_t^\infty y_\mu^{-2}(u) \frac{\d u}{\sqrt{u+\mu}}$. 
By performing an integration by parts,
\begin{equation} \label{Mtrick}
\int_0^{T} M(G_\mu(t)) y_\mu^{-2}(t)  \frac{\d t}{2\sqrt{t+\mu}}   =  \frac12\bigg(-w_\mu(T) M(G_\mu(T))  + \int_0^T w_\mu(t) \d M(G_\mu(t)) \bigg)
\end{equation}
where we used that $G_\mu(0) =M(0) = 0$. 
According to Lemma~\ref{lem:BMestimate} and \eqref{eq:VERT}, it holds for any $\epsilon\in(0,1]$ and $t\ge 0$, 
\begin{equation} \label{MA2}
\Big\VERT \sup_{t\in[0,2]} |M(t)|
\Big\VERT_1 \lesssim_\beta 1
\qquad\text{and}\qquad
\Big\VERT  \sup_{t\ge 2} \tfrac{|M(t)|}{t (\log t )^\epsilon} \Big\VERT_1 
  \lesssim_{\beta,\epsilon} 1. 
\end{equation}
We need the following bound
\begin{equation} \label{wbound}
  w_\mu(t) 
  =\int_t^\infty y_\mu^{-2}(u) \frac{\d u}{\sqrt{u+\mu}}
  \leq\int_t^\infty y_\mu^{-2}(u) \frac{\sqrt{u+\mu}}{{t+\mu}}\d u
  \leq  \frac{y_\mu^{-2}(t)}{4({t+\mu})}.
\end{equation}
Using the first bound of \eqref{MA2} it holds for any $\kappa\in\N$,
\begin{equation}\label{eq:Arnault5}
\Big\VERT \sup_{T\ge0 : G_\mu(T) \le 2 } \Big| w_\mu(T) M(G_\mu(T)) \Big| \Big\VERT_1\lesssim_\beta w_\mu(0)   \lesssim_{\beta,\kappa} \mu^{-\kappa}.
\end{equation}
As for larger $T$, we use the bounds \eqref{wbound} and \ref{eq:Gupper}, which shows that
\begin{equation} \label{wBbound}
 w_\mu(t) 
 \le ({1/4})({t+\mu})^{-1} y_\mu^{-2}(t)
 \le \frac{1}{4G_\mu(t)({t+\mu})^{3/2}}
 \le \frac{2}{3G_\mu(t)\log G_\mu(t)}.
% \qquad\text{and}\qquad G_\mu(t) \le \frac{(1-\mu^{-3/2})^{-1}}{4\sqrt{t+\mu}} y_\mu^{2}(t)
\end{equation}
Then from the second bound in \eqref{MA2}, it holds for any $\epsilon \in [0,\tfrac 32],$
\[
\Big\VERT \sup_{T\ge0 : G_\mu(T) \ge 2 } \Big| w_\mu(T) M(G_\mu(T)) \Big| \Big\VERT_1
%\lesssim_{\beta, \epsilon} \sup_{T\ge0 : G_\mu(T) \ge 2 } \Big|  w_\mu(T) G_\mu(T) \big( \log G_\mu(T) \big)^{2\epsilon/3} 
%\Big| 
\lesssim_{\beta, \epsilon}  \mu^{-3/2+\epsilon}. 
\]
Combining this bound with \eqref{eq:Arnault5}
\begin{equation}  \label{M1}
\Big\VERT \sup_{T\ge 0} w_\mu(T)|M(G_\mu(T))|  \Big\VERT_1  \lesssim_{\beta, \epsilon}  \mu^{-3/2+\epsilon}  . 
\end{equation}
This reduces the problem to a tail bound for the martingale term on the RHS of \eqref{MA2}. 
 
   \paragraph{Step 3: Control of the martingale part.} 
 For the martingale term, it suffices to control its total quadratic variation:
\begin{equation} \label{QV2}
\begin{aligned}
 \left[  \int_0^\infty w_\mu(t) \d M(G_\mu(t))  \right] 
&= \int_0^\infty w_\mu^2(t) \d [M(G_\mu(t))] \\
&= 4\bigg( \int_0^{\mu^{-1/2}} \hspace{-.3cm} w_\mu^2(t) B(G_\mu(t))^2 G_\mu'(t) \d t 
 + \int_{T_\mu}^\infty w_\mu^2(G_\mu^{(-1)}(u))  B(u)^2 \d u \bigg)
\end{aligned}
\end{equation}
where we used that $\d [M(t)] = 4 B(t)^2 \d t$,  $T_\mu \coloneqq G_\mu(\mu^{-1/2}) $ and $G_\mu^{(-1)}$ denotes the inverse function of $G_\mu(t)=\int_0^t y_\mu^2(u) \d u$.
Using  Lemma~\ref{lem:BMestimate} again, for any $\epsilon>0$, 
\[
\Big\VERT \sup_{t \in\R_+} \tfrac{B(t)^2}{t(\log(1+ t))^\epsilon} \Big\VERT_1 \lesssim_{\beta, \epsilon}1 . 
\]
Hence, using that $G_\mu'(t) =  y_\mu^{2}(t)$, \eqref{wBbound} and $\log(1+ G_\mu(t)) \lesssim (\mu+t)^{3/2}$, this implies that
\[ \begin{aligned}
\Big\VERT \int_0^{\mu^{-1/2}} \hspace{-.3cm} w_\mu^2(t) B(G_\mu(t))^2 G_\mu'(t) \d t  \Big\VERT_1
&\lesssim_{\beta, \epsilon} \int_0^{\mu^{-1/2}} \hspace{-.3cm}   y_\mu^{-2}(t) G_\mu(t)  \frac{\d t}{(t+\mu)^2} \\
& \lesssim \mu^{-3+\epsilon} . 
\end{aligned}\]
This controls the first term on the RHS of \eqref{QV2}. For the second term, observe that $T_\mu \simeq \frac{1/4}{\sqrt{\mu}}  e^{\frac 83 \mu^{3/2}}$ as $\mu\to\infty$ (given the bounds \ref{eq:Gupper} and \ref{eq:Glower}) and that by \eqref{wBbound},
$w_\mu^2(G_\mu^{(-1)}(u))  \le  \frac{2}{3(u \log u)^2}$.
Hence, by applying Lemma~\ref{lem:S} with $L(u) =1/(\log u)$ and $T = T_\mu $, we have $m_T = L(T_\mu) \simeq c\mu^{-3/2}$ and this yields
\[
\P\bigg( \int_{T_\mu}^\infty w_\mu^2(G_\mu^{(-1)}(u))  B(u)^2 \d u \ge C \mu^{-3/2} \bigg) \le \exp(- c \mu^{3/2}) . 
\] 
for constants $C,c>0$ depending only on $\beta$. 
Going back to \eqref{QV2}, these estimates imply that  for any $\epsilon>0$, 
\[
\P\bigg(  \left[  \int_0^\infty w_\mu(t) \d M(G_\mu(t))  \right]  \ge C_\epsilon \mu^{-3/2} \bigg)
\le \exp(- \mu^{3/2-\epsilon}) . 
\]

To complete the argument, we work on the event that the bracket is small and conclude that for all $\mu$ sufficiently large:
 \begin{align} \notag
& \P\bigg(  \sup_{T\ge 0} \bigg| \int_0^T w_\mu(t) \d M(G_\mu(t)) \bigg|  \ge C \mu^{-\epsilon} \bigg) \\
&\notag\qquad\qquad\le \P\Bigg(    \sup_{T\ge 0}  \bigg| \int_0^T w_\mu(t) \d M(G_\mu(t)) \bigg|  \ge C \mu^{-\epsilon} ,   \left[  \int_0^\infty w_\mu(t) \d M(G_\mu(t))  \right]  \le  \mu^{-3/2} \Bigg)
+ \exp(-  \mu^{3/2-\epsilon})  \\
&\label{M2} 
\qquad\qquad \le  \exp(- \mu^{3/2-2\epsilon})  .
\end{align}
By combining formula \eqref{Mtrick} with the bounds \eqref{M1} \eqref{M2}, we conclude that for any $\epsilon>0$, 
\[
\P\Bigg(  \sup_{T\ge 0} \bigg| 
 \int_0^T M(G_\mu(t)) y_\mu^{-2}(t)  \frac{\d t}{\sqrt{t+\mu}} \bigg|
 \ge C_\epsilon \mu^{-\epsilon} \Bigg)  \le  \exp(- \mu^{3/2-2\epsilon})   . 
\]
This completes the proof of \eqref{Mbound} and the overall argument. 
\end{proof}

%\begin{remark}
%Let  $ \mathbf{S}= \int_0^\infty L(u)^2 B(u)^2 \frac{\d u}{u^{2\alpha}} $ where  $L\ge 0$ is a slowly varying function and $\alpha>1$.
%A variation of the proof shows that in this case,
%\[
%\P[\mathbf{S} \ge 2 m] 
%\]
% 
%
%
% the mean of $ \int_0^\infty L(u)^2 B(u)^2 \frac{\d u}{u^{2\alpha}}$ is given by
%\[
%m= \tr \mathcal{K} =  \int_0^\infty  \frac{ L(u)^2}{u^{2\alpha-1}} \d u 
%\]
%then, using the bound  $ \tr [\mathcal{K}^{\ell}] \le  (\tr \mathcal{K})^\ell$ valid for all $\ell\in\N$, 
%for $0<\lambda \le \frac{1}{2 m }$
%\[
% \Exp\big[ e^{\lambda( \mathbf{S}- m)/2}\big] \le \exp\bigg(\tfrac{\lambda^2 m^2}{2} \bigg)
% \]
% Then applying Markov's inequality with $\lambda =  \frac{1}{2 m }$, we obtain \eqref{}.
%\end{remark}

\begin{lemma}\label{lem:X3}
For any $\epsilon\in(0,1/2]$, there exists a constant $c_\epsilon$ so that for $n\ge 1$, 
\[
    \Pr\biggl( \int_0^\infty    \frac{ \mathfrak{X}_{n}^2(u) }{(u+n)^{1/2+\epsilon}} \d u
 \geq C_\epsilon n^{-\epsilon}
  \biggr)
  \leq
\exp(-n^{3/2-2\epsilon}) .
\]
\end{lemma}

\begin{proof}
We follow the steps from the proof of Lemma~\ref{lem:X2} except that the estimates are simpler.
First, we show that the mean is appropriately small. 
By \eqref{mean2} and the subsequent estimates, $\E\mathfrak{X}_{n}^2(u) \lesssim (u+n)^{-1/2}$
so that 
\[
 \int_0^\infty    \frac{\E \mathfrak{X}_{n}^2(u) }{(u+n)^{1/2+\epsilon}} \d u
 \lesssim n^{-\epsilon} . 
\]
So we can again reduce the problem to martingale tail bounds.
\[
 \int_0^\infty    \frac{\mathfrak{X}_{n}^2(u)- \E \mathfrak{X}_{n}^2(u)}{(u+n)^{1/2+\epsilon}} \d u
\,  \overset{\rm law}{=}\, 
\int_0^\infty \frac{M(G_n(u)) y_n^{-2}(u)}{(u+n)^{1/2+\epsilon}} \d u
\]
 where $M(t) \coloneqq B(t)^2-t$ is a continuous martingale. By an integration by parts, the RHS satisfies 
\[
\int_0^\infty \frac{M(G_n(u)) y_n^{-2}(u)}{(u+n)^{1/2+\epsilon}} \d u =  \int_0^\infty v_n(t) \d M(G_n(t)) 
\]
where $ \displaystyle v_n(t) \coloneqq \int_t^\infty \frac{y_n^{-2}(u)}{(u+n)^{1/2+\epsilon}} \d u \le w_n(t) $ for $t\ge 0$ as in \eqref{Mtrick} and we used that  $v_n(T) M(G_n(T))\to 0$ almost surely  as $T\to\infty$. In particular, 
the (total) quadratic variation of this martingale is stochastically dominated by \eqref{QV2} and we conclude by Step 3 of the proof of Lemma~\ref{lem:X2} that there is a constant $C_\epsilon>0$ so that
\[
\P\Bigg( \bigg| \int_0^\infty \frac{M(G_n(u)) y_n^{-2}(u)}{(u+n)^{1/2+\epsilon}} \d u \bigg|  \ge C_\epsilon n^{-\epsilon} \Bigg) \le \exp(- n^{3/2-2\epsilon})   . \qedhere
\]
\end{proof}

We are finally in a position to prove Proposition~\ref{prop:2021}. Fix $\epsilon\in (0,\frac12]$ and $C>0$ large enough, we can work on the event for $n\gg 1$, 
\begin{equation} \label{goodevent}
\mathsf{Ctrl}_n \coloneqq \Bigg\{
\begin{aligned}
  |\varpi_{\lambda,n}(t)|&\le Ct^{-1/2},
  &|\Delta_{\lambda,n}(t)|&\le Ct^{1/2-\epsilon},\\
  |\mathfrak{X}_{\mu,n}(t)|&\le Ct^{1/2-\epsilon},
  &|\mathcal{R}_{\mu,n}(t)|&\le Ct^{-1/2-\epsilon}
\end{aligned}
\ :\ t\ge s_n,\ \lambda\in K\Bigg\} .
\end{equation}
Indeed, by Proposition~\ref{lem:sig},
\[
\mathcal{E}_{n} \subset \Big\{    | \varpi_{\lambda,n}(t) | \le     Ct^{-1/2}, 
    | \Delta_{\lambda,n}(t) | \le     Ct^{1/2-\epsilon}  \text{  for all }t \geq s_n ,  \lambda\in K \Big\}
 \]
%\[
%\P\Big[\Big\{    | \varpi_{\lambda,n}(t) | \le     Ct^{-1/2}, 
%    | \Delta_{\lambda,n}(t) | \le     Ct^{1/2-\epsilon}  \text{  for all }t \geq s_n ,  \lambda\in K
%\Big\}^c \Big]     \le  \P\big( \mathcal{E}_{n}^c \big) \le C \exp\big(-n^{3/2-5\epsilon/2}\big) 
%\]
and, by Proposition~\ref{prop:Xest} and Lemma~\ref{lem:Ups}, if $n\gg 1$, 
\[
\P\Big[\Big\{ 
  |\mathfrak{X}_{\mu,n} (t)| \le C t^{1/2-\epsilon} , |\mathcal{R}_{\mu,n}(t)|  \le C t^{-\epsilon}  \text{  for all }t \geq s_n ,  \lambda\in K
\Big\}^{c}\Big]    
\le \exp\bigg(-\frac{n^{3/2-2\epsilon}}{(\log n)^2}\big) . 
\]

Then, our first goal is to control the contribution from the real part of the Riccati diffusion that is,  to show that
  \[
    \Pr\Bigg( \mathsf{Ctrl}_n \cap
    \biggl\{
    \sup_{\lambda \in K}
    \sup_{T \geq s_n} 
    \bigg| \int_{s_n}^T
    \biggl(
    \Delta_{\lambda,n}(t) + \frac{4}{\beta}\int_0^t \Exp[\mathfrak{X}(t)\mathfrak{X}(u)] \d u\biggr)  \d t \bigg|  \geq C_\epsilon n^{-\epsilon}
    \biggr\} \Bigg)
    \le
    \exp\left( - n^{3/2-2\epsilon} \right).
  \]

\begin{proof}[Proof of Proposition~\ref{prop:2021}]
   We start by recalling \eqref{delta2} that with $s=s_n$, 
  \begin{equation*} 
    \Delta_{\lambda,n}(t)
    = v_{\lambda,s}^{-1}(t)\Delta_{\lambda,n}(s)  + \int_{s}^t \frac{v_{\lambda,s}(u)}{v_{\lambda,s}(t)}\big( \varpi_{\lambda,n}^2(u) - \mathfrak{h}_{\mu,n}^2(u) \big) \d u. 
  \end{equation*}
  where 
  \[
   v_{\lambda,s}(t) \coloneqq 
  \frac{y_\mu(t)}{y_\mu(s)}\exp\biggl( \int_s^t (\Delta_{\lambda,n}(u) + 2\mathfrak{h}_{\mu,n}(u))\d u \biggr)
  \]
 and by \eqref{def:mathfrakh}, $\mathfrak{h}_{\mu,n}(t)  = \E \mathfrak{h}_{\mu,n}(t) + \tfrac{2}{\sqrt{\beta}}\mathfrak{X}_{\mu, n}(t) $ where $ \E \mathfrak{h}_{\mu,n}(t)  = \O(t^{-1})$ for $t\ge n$. 
 In particular, on the event \eqref{goodevent},   $|\mathfrak{h}_{\mu,n} (t)| \lesssim  t^{1/2-\epsilon} $ for all $t  \ge s$ and $\lambda\in K$ and we can replace
\begin{equation} \label{hX2}
 \mathfrak{h}_{\mu,n}^2(u) =  \frac{4}{\beta}\mathfrak{X}_{\mu,n}^2(u) +  \O(u^{-1/2-\epsilon}) . 
 \end{equation}

Let us denote $\displaystyle V_{\lambda,T}(u) \coloneqq \int_u^T  \frac{v_{\lambda,s}(u)}{v_{\lambda,s}(t)} \d t$. 
 In terms of this quantity, by Fubini's Theorem, 
\begin{equation} \label{Dint1}
 \int_s^T  \Delta_{\lambda,n}(t) dt = \Delta_{\lambda,n}(s)   V_{\lambda,T}(s) 
 +  \int_s^T   V_{\lambda,T}(u) \big( \varpi_{\lambda,n}^2(u) - \mathfrak{h}_{\mu,n}^2(u) \big)  \d u  
 \end{equation}
 where we used that $v_{\lambda,s}(s) =1$.
 Moreover, by an integration by parts,  
\[
 V_{\lambda,T}(u)= - v_{\lambda,s}(u)\Bigg\{ \bigg[  \frac{v_{\lambda,s}^{-1}(t)}{2\sqrt{t+\mu}} \bigg]_{t=u}^T + \int_u^T v_{\lambda,s}^{-1}(t) \bigg( \frac{\Delta_{\lambda,n}(t) + 2\mathfrak{h}_{\mu,n}(t)}{2\sqrt{t+\mu}} + \frac{1/4}{(t+\mu)^{3/2}} \bigg) \d t \Bigg\}
\]
This implies that on the event \eqref{goodevent},
\[
 V_{\lambda,T}(u) =  \frac{1/2}{\sqrt{u+\mu}}  - \frac{1/2}{\sqrt{T+\mu}}  \frac{v_{\lambda,s}(u)}{v_{\lambda,s}(T)} + \O_K(u^{-\epsilon})  V_{\lambda,T}(u) .
\]
This shows that on our event, provided that $s$ is sufficiently large, it holds (deterministically) uniformly for all $T\ge u\ge s$ and $\lambda\in K$,
\begin{equation} \label{Vest}
 V_{\lambda,T}(u)  =   \frac{1/2}{\sqrt u}  + \O(u^{-1/2-\epsilon})  -\O(T^{-1/2})   \frac{v_{\lambda,s}(u)}{v_{\lambda,s}(T)}. 
\end{equation}
In particular, $| V_{\lambda,T}(u)| \le C/\sqrt{u}$ so that $\Delta_{\lambda,n}(s)   V_{\lambda,T}(s) = \O(s^{-\epsilon})$. By \eqref{Dint1}, this implies that on $\mathsf{Ctrl}_n$,  
\begin{equation} \label{Dint2}
 \int_s^T  \Delta_{\lambda,n}(t) dt = \O(s^{-\epsilon})
 + \mathrm{I} - \mathrm{II} - \mathrm{III} 
 \end{equation}
  where 
\begin{equation} \label{Derrors} \begin{aligned}
\mathrm{I} &= \O(1) \int_s^T   \varpi_{\lambda,n}^2(u) \frac{\d u}{\sqrt u} \, , \\
 \mathrm{II} &= \O(T^{-1/2})  \int_s^T   \frac{v_{\lambda,s}(u)}{v_{\lambda,s}(T)} \big( \varpi_{\lambda,n}^2(u) - \mathfrak{h}_{\mu,n}^2(u) \big) \d u \, ,  \\ 
 \mathrm{III} &=  \int_s^T   \mathfrak{h}_{\mu,n}^2(u)  \Big( \frac{1/2}{\sqrt u}  + \O(u^{-1/2-\epsilon}) \Big) \d u \, . \\
% \mathrm{IV} &=\int_s^T \frac{\mathfrak{h}_{\mu,n}^2(u) - \E\mathfrak{h}_{\mu,n}^2(u)}{2\sqrt{u}} \d u \, . 
\end{aligned} 
\end{equation}

On  $\mathsf{Ctrl}_n$,  we have  the  crude bounds
\[
\mathrm{I}  \lesssim \int_s^\infty u^{-3/2} \d u \lesssim s^{-1/2} 
\qquad\text{and}\qquad
| \mathrm{II}| \lesssim T^{1/2-2\epsilon}   \int_s^T   \frac{v_{\lambda,s}(u)}{v_{\lambda,s}(T)}  \d u
\]
Observe that making another integration by parts, setting $w(t) \coloneqq  y_\mu^{-1}(t)  \int_s^t y_\mu(u) \d u$ for $t\ge s$, we obtain
\[
  \int_s^T   \frac{v_{\lambda,s}(u)}{v_{\lambda,s}(T)}  \d u
=   w(T)-  \int_s^T  \frac{v_{\lambda,s}(u)}{v_{\lambda,s}(T)}  w(u)    (\Delta_{\lambda,n}(u) + 2\mathfrak{h}_{\mu,n}(u))\d u
\]
Using that $w(t) \le C/\sqrt{t}$, this implies that on $\mathsf{Ctrl}_n$, 
\[
  \int_s^T   \frac{v_{\lambda,s}(u)}{v_{\lambda,s}(T)}  \d u \le w(T) + 
\O(s^{-\epsilon})   \int_s^T   \frac{v_{\lambda,s}(u)}{v_{\lambda,s}(T)}  \d u
\]
so that $\displaystyle  \int_s^T   \frac{v_{\lambda,s}(u)}{v_{\lambda,s}(T)}  \d u \lesssim T^{-1/2} $. 
We conclude that $| \mathrm{II}| \lesssim T^{-2\epsilon} $. 
According to \eqref{Dint2}--\eqref{Derrors}, this shows that on $\mathsf{Ctrl}_n$,  
\begin{equation} \label{Dint3}
 \int_s^T  \Delta_{\lambda,n}(t) dt = \O(s^{-\epsilon}) - \mathrm{III} 
 \end{equation}
where by \eqref{hX2},  
$\displaystyle \mathrm{III} =  \frac{4}{\beta}\int_s^T   \mathfrak{X}_{\mu,n}^2(u)  \Big( \frac{1/2}{\sqrt u}  + \O(u^{-1/2-\epsilon}) \Big) \d u +  \O(s^{-\epsilon})  $ only depends on the Gaussian process   $\mathfrak{X}$. 

\smallskip

Now, we want to replace $\mathfrak{X}_{\mu,n}^2$ by $\mathfrak{X}_{0,n}^2$ using \eqref{XSDE_rep}, 
\[
\mathfrak{X}_{\mu, n}^2(t) =  \mathfrak{X}_{0,n}^2(t) +2  \mathcal{R}_{\mu,n}(t)  \mathfrak{X}_{0,n}(t) +  \mathcal{R}_{\mu,n}^2(t) , \qquad 
\]
 Let us denote $\displaystyle \mathrm{IV}  \coloneqq  \int_n^\infty   \mathfrak{X}_{0,n}^2(u)   \frac{\d u}{u^{1/2+\epsilon}}$. 
Then by Cauchy--Schwarz inequality, 
\[
 \mathrm{III} = \frac{4}{\beta}\int_n^T   \mathfrak{X}_{0,n}^2(u)   \frac{\d u}{2\sqrt{u}} 
 +\O(\mathrm{IV})+  \O\bigg( \sqrt{ \bigl(\mathrm{IV}\bigr)\cdot \int_n^T \mathcal{R}_{\mu,n}^2(u)  \frac{\d u}{u^{1/2-\epsilon}}  } \bigg) + \O\bigg( \int_n^T \mathcal{R}_{\mu,n}^2(u)  \frac{\d u}{2\sqrt{u}}\bigg)
 + \O(n^{-\epsilon}) . 
\]

On the event  $\mathsf{Ctrl}_n$ and by Lemma~\ref{lem:X3}, it holds for any $\epsilon \in (0,1/2]$, 
\begin{equation} \label{PIV}
 \int_n^\infty \mathcal{R}_{\mu,n}^2(u)  \frac{\d u}{u^{1/2-\epsilon}} \lesssim n^{-1/2-\epsilon}  
\qquad\text{and}\qquad
\P\bigg(  \int_n^\infty   \mathfrak{X}_{0,n}^2(u)   \frac{\d u}{u^{1/2+\epsilon}} \ge C_\epsilon n^{-\epsilon} \bigg) \le e^{- n^{3/2-2\epsilon}}  . 
\end{equation}
This shows that with probability at least $1-e^{- n^{3/2-2\epsilon}}$ and the required uniformity, 
\[
 \mathrm{III} = \frac{4}{\beta}\int_n^T   \mathfrak{X}_{0,n}^2(u)   \frac{\d u}{2\sqrt{u}} 
 + \O(n^{-\epsilon}) . 
\]

In order to obtain \eqref{PIV} and compute the leading order of $\mathrm{III}$, 
we use that for a given $n\ge 1$, we have the equality of processes
\[
\left( \mathfrak{X}_{0,n}(t+n) =  y_0^{-1}(t+n) \int_n^{t+n} y_0(u) \d B(u)  \right)_{t\ge 0}
  \overset{\rm law}{=}  \, 
  \left( \mathfrak{X}_n(t)  =   y_n^{-1}(t) \int_0^{t} y_n(u) \d B(u)  \right)_{t\ge 0}
 \]
since white noise is invariant by translation. Then for $T\ge n$, 
\[
 \mathrm{IV}  \,  \overset{\rm law}{=}\,  \int_0^\infty   \mathfrak{X}_{n}^2(u)   \frac{\d u}{(u+n)^{1/2+\epsilon}} 
 \qquad\text{and}\qquad 
 \int_n^T   \mathfrak{X}_{0,n}^2(u)   \frac{\d u}{2\sqrt{u}}  \,  \overset{\rm law}{=}\, 
  \int_0^{T-n}   \mathfrak{X}_{n}^2(u)   \frac{\d u}{2\sqrt{u+n}} .
\]
Using Lemma~\ref{lem:X2}, we conclude 
that uniformly for $T\ge s_n$ and $\lambda \in K$ with probability at least $1-2e^{- n^{3/2-2\epsilon}}$, 
\[
 \mathrm{III} =  \frac{4}{\beta}\int_n^{T}\int_0^t \Exp[\mathfrak{X}(t)\mathfrak{X}(u)] \d u  \d t
 + \O(n^{-\epsilon}) . 
\]
This shows that the contribution from $ \mathrm{III} $ is deterministic but positive. 
By formula \eqref{mean2} and using again that $\E\mathfrak{X}^2(t) =\O(1)$ uniformly for all $t\ge 0$, 
\[
\int_n^{s}\int_0^t \Exp[\mathfrak{X}(t)\mathfrak{X}(u)] \d u  \d t  =  \O(1) \int_n^{s} \frac{\d u}{\sqrt {u}}
=  \O((\log n)/n) . 
\]
By \eqref{Dint3}, this completes the proof. 
  \end{proof}

We now turn to control the imaginary part of the Riccati diffusion. In this case, our goal is to show that
  \[
    \Pr\Bigg( \mathsf{Ctrl}_n \cap
    \biggl\{
    \sup_{\lambda \in K}
    \sup_{T \geq s_n} 
    \bigg| \int_s^T \bigg( \varpi_{\lambda,n}(t) -  \frac{\eta}{2\sqrt{t}} \bigg) \d t  \bigg|  \geq C_\epsilon n^{-1/2}
    \biggr\} \Bigg)
    \le
    \exp\left( - n^{3/2-2\epsilon} \right).
  \]
The method is similar to the Proof of Proposition~\ref{prop:2021} and we rely on some of the estimates in this proof.

\begin{proof}[Proof of Proposition~\ref{prop:2022}]
According to \eqref{eq:complexRic}, for $t\ge n$, 
\[
\d   \varpi_{\lambda,n} (t) = - 2\sqrt{t}  \varpi_{\lambda,n}(t) \d t +\big(\eta -2\big(  \sigma_{\lambda,n}(t)-\sqrt{t} \big)\varpi_{\lambda,n}(t) \big) \d t .
\]
Hence, with $y=y_0$ and $s= s_n = n+\frac{\log n}{\sqrt{n}}$, it holds for $t\ge s$, 
\[\begin{aligned}
 \varpi_{\lambda,n} (t) 
&=  \varpi_{\lambda,n} (s) \frac{y(s)}{y(t)}  +   \int_s^t \frac{y(u)}{y(t)} \big(\eta -2\big(  \sigma_{\lambda,n}(u)-\sqrt{u} \big)\varpi_{\lambda,n}(u) \big)  \d u \\
&= \frac{\eta}{2\sqrt{t}} + \varpi_{\lambda,n}(s) \frac{y(s)}{y(t)}  + \O_{K}(t^{-3/2}) -  2 \int_s^t \frac{y(u)}{y(t)} \big( \sigma_{\lambda,n}(u)-\sqrt{u} \big)\varpi_{\lambda,n}(u)  \d u
\end{aligned}\]
where the error is deterministic. 
Now we also recall that with a deterministic error, 
\[
\sigma_{\lambda,n}(u)-\sqrt{u}=  \Delta_{\lambda,n}(u) + \O(u^{-1/2}) + \tfrac{2}{\sqrt{\beta}}\mathfrak{X}_{\mu, n}(u)  . 
\]
Using that on the event \eqref{goodevent}, 
$\displaystyle\int_s^t \frac{y(u)}{y(t)}\varpi_{\lambda,n}(u)  \frac{\d u}{\sqrt u} = \O(t^{-3/2})$
uniformly for all $t\ge s$, this implies that
\[
\varpi_{\lambda,n}(t) -  \frac{\eta}{2\sqrt{t}}  =  \varpi_{\lambda,n}(s) \frac{y(s)}{y(t)} -  2 \int_s^t \frac{y(u)}{y(t)} \big(\Delta_{\lambda,n}(u) + \tfrac{2}{\sqrt{\beta}}\mathfrak{X}_{\mu, n}(u)  \big)  \varpi_{\lambda,n}(u)  \d u + \O(t^{-3/2})  . 
\]
Since the error is integrable, by Fubini's Theorem, 
\[
\int_s^T \bigg( \varpi_{\lambda,n}(t) -  \frac{\eta}{2\sqrt{t}} \bigg) \d t = \O\bigg( \frac{1+|\varpi_{\lambda,n}(s)|}{\sqrt{s}} \bigg) -  2 \int_s^T  \big(\Delta_{\lambda,n}(u) + \tfrac{2}{\sqrt{\beta}}\mathfrak{X}_{\mu, n}(u)  \big)  \bigg(  \frac{1}{\sqrt{u}} + \O(u^{-3/2}) \bigg) \varpi_{\lambda,n}(u)  \d u
\]
where both implied constants are deterministic and we used that $\displaystyle \int_u^T  \frac{y(u)}{y(t)} \d t = \frac{1}{2\sqrt{u}} + \O(u^{-3/2}) $ uniformly   for all  $T\ge u \ge1$.
We conclude that  on $\mathsf{Ctrl}_n$, 
\begin{equation}  \label{Dint5}
\int_s^T \bigg( \varpi_{\lambda,n}(t) -  \frac{\eta}{2\sqrt{t}} \bigg) \d t =\O(s^{-1/2}) -   \int_s^T  \Delta_{\lambda,n}(u) \varpi_{\lambda,n}(u)  \frac{\d u}{\sqrt{u}}
-   \tfrac{2}{\sqrt{\beta}}\int_s^T \mathfrak{X}_{\mu, n}(u)\varpi_{\lambda,n}(u)  \frac{\d u}{\sqrt{u}} . 
\end{equation}
 
Now using representation \eqref{delta2} and proceeding like in the proof of Proposition~\ref{prop:2021} (see formula \eqref{Dint1})
\begin{equation}  \label{Dint4}
\int_s^{T}   \Delta_{\lambda,n}(u) \varpi_{\lambda,n}(u)  \frac{\d u}{\sqrt{u}}
= \Delta_{\lambda,n}(s)   W_{\lambda,T}(s) 
 +  \int_s^T   W_{\lambda,T}(u) \big( \varpi_{\lambda,n}^2(u) - \mathfrak{h}_{\mu,n}^2(u) \big)  \d u 
 \end{equation}
where $\displaystyle W_{\lambda,T}(u) \coloneqq \int_u^T  \frac{v_{\lambda,s}(u)}{v_{\lambda,s}(t)} \varpi_{\lambda,n}(t)  \frac{\d t}{\sqrt{t}}$.
We can control $ W_{\lambda,T}$, by an integration by parts,
\[
 W_{\lambda,T}(u)= - v_{\lambda,s}(u)\Bigg\{ \bigg[  \frac{v_{\lambda,s}^{-1}(t)  \varpi_{\lambda,n}(t)}{2\sqrt{t(t+\mu)}} \bigg]_{t=u}^T + \int_u^T \bigg( \big( \Delta_{\lambda,n}(t) + 2\mathfrak{h}_{\mu,n}(t) + \O(t^{-3/2}) \big) \varpi_{\lambda,n}(t) -\varpi_{\lambda,n}'(t) \bigg) \frac{v_{\lambda,s}^{-1}(t) \d t}{{2\sqrt{t(t+\mu)}}} \Bigg\}
\]
where the error is deterministic and we used that the process  $\varpi$ is $C^1$. 
Since  $\varpi'(t)= \eta-2\sigma(t)\varpi(t) = \O(1)$ 
and $ \big( \Delta_{\lambda,n}(t) + 2\mathfrak{h}_{\mu,n}(t) + \O_K(t^{-3/2}) \big) \varpi_{\lambda,n}(t) = \O(1)$ on the event \eqref{goodevent}, this implies that uniformly for all $\lambda\in K$ and $u\in [s,T]$,
\[
 W_{\lambda,T}(u) = \frac{\varpi_{\lambda,n}(u)}{2\sqrt{u(u+\mu)}} -\O(T^{-3/2}) \frac{v_{\lambda,s}(u)}{v_{\lambda,s}(T)}  - \O(u^{-1})  V_{\lambda,T}(u)
\]
where $V$ is as in formula \eqref{Dint1}.
In particular, as $V_{\lambda,T}(u) =\O(u^{-1/2})$,  this shows that $W_{\lambda,T}(u)=\O(u^{-3/2})$ for $u\in [s,T]$, so that by \eqref{Dint4}, 
\[
\int_s^{T}  \Delta_{\lambda,n}(u) \varpi_{\lambda,n}(u)  \frac{\d u}{\sqrt{u}}
= \O(s^{-1}) - \O(1)  \int_s^T  \mathfrak{h}_{\mu,n}^2(u)  \frac{\d u}{u^{3/2}}
\]

Hence, we have reduced the problem to bound an integral depending only on the Gaussian process~$ \mathfrak{h}$. 
Like in the proof of Proposition~\ref{prop:2021}, on  the event $\mathsf{Ctrl}_n$,  $ \mathfrak{h}_{\mu,n}^2(u) =  \frac{4}{\beta}\mathfrak{X}_{0,n}^2(u)  + \O(n^{-2\epsilon})$ uniformly for $\lambda\in K$ and $u\ge n$ so that 
\[
\int_s^T  \mathfrak{h}_{\mu,n}^2(u)  \frac{\d u}{u^{3/2}}
\le s^{-1/2+\epsilon}  \int_s^\infty  \mathfrak{X}_{0,n}^2(u)  \frac{\d u}{u^{1/2+\epsilon}} + \O(n^{-1/2-2\epsilon}) . 
\]
Using the tail bound \eqref{PIV}, we conclude that 
\begin{equation}  \label{Dint6}
\P\Bigg(  \mathsf{Ctrl}_n  \cap\bigg\{ \sup_{\lambda\in K, T\ge s} \bigg| \int_s^T  \Delta_{\lambda,n}(u) \varpi_{\lambda,n}(u)  \frac{\d u}{\sqrt{u}} \bigg| \ge C_\epsilon n^{-1/2} \bigg\} \bigg) \le e^{- n^{3/2-2\epsilon}}  . 
\end{equation}

This provides control of the second term on the RHS of \eqref{Dint5}. 
For the third term on the RHS of \eqref{Dint5}, using again the representation \eqref{XSDE_rep}, we obtain on $\mathsf{Ctrl}_n$,
\[
\int_s^T \mathfrak{X}_{\mu, n}(u)\varpi_{\lambda,n}(u)  \frac{\d u}{\sqrt{u}}
= \int_s^T \mathfrak{X}_{0, n}(u)\varpi_{\lambda,n}(u)  \frac{\d u}{\sqrt{u}} + \O(n^{-1/2-\epsilon})
\]
Now, by the Cauchy--Schwarz inequality, for any $\epsilon>0$
\[
\bigg| \int_s^T \mathfrak{X}_{0, n}(u)\varpi_n(u)  \frac{\d u}{\sqrt{u}}  \bigg| \le 
\sqrt{ \int_s^\infty \frac{\mathfrak{X}_{0, n}^2(u)}{u^{1/2+\epsilon}} \d u \int_s^\infty \frac{\varpi_{\lambda,n}^2(u)}{u^{1/2-\epsilon}} \d u  } .
\]
On $ \mathsf{Ctrl}_n$, the second integral is  $\O(n^{-1/2+\epsilon})$ and the tail bound \eqref{PIV} controls the first integral. We conclude that 
\begin{equation} \label{Dint7}
\P\Bigg(  \mathsf{Ctrl}_n  \cap\bigg\{ \sup_{\lambda\in K, T\ge s}  \bigg| \int_s^T \mathfrak{X}_{\mu, n}(u) \varpi_{\lambda,n}(u)  \frac{\d u}{\sqrt{u}} \bigg| \ge C_\epsilon n^{-1/2} \bigg\} \bigg) \le e^{- n^{3/2-2\epsilon}}  . 
\end{equation} 
Combining the tail bounds \eqref{Dint6}, \eqref{Dint7} with the estimate \eqref{Dint5}, this completes the proof. 
\end{proof}

We can now perform the last step in the proof of Proposition~\ref{prop:asympDN}: namely, to show that for $n$ large enough, we can replace the Riccati diffusion by the simple Gaussian process $\theta_\lambda(t)$ up to an error which is uniformly small in an integrated sense and with overwhelming probability. 

\begin{proof}[Proof of Proposition~\ref{prop:X2int}]
Recall that according to \eqref{def:mathfrakh}--\eqref{delta0}, for $\lambda\in K$ and $t\ge n$,
\[
\rho_{\lambda,n}(t) =     \sigma_{\lambda,n}(t) + \i  \varpi_{\lambda,n}(t)
\qquad\text{and}\qquad 
\Delta_{\lambda,n}(t)   =\sigma_{\lambda,n}(t)- \sqrt{t} - \frac{\mu}{2\sqrt{t}}  +\frac{1}{4(t+1)} -  \tfrac{2}{\sqrt{\beta}}\mathfrak{X}_{\mu, n}(t)  + \O_K(t^{-3/2}) 
\]
with a deterministic error. By definition of $\theta_\lambda$, this implies that 
\[
\rho_{\lambda,n}(t) -  \theta_\lambda(t)
=  \Delta_{\lambda,n}(t) + \frac{4}{\beta}\int_0^t \Exp[\mathfrak{X}(t)\mathfrak{X}(u)]\d u
+ \i \bigg( \varpi_{\lambda,n}(t)- \frac{\eta}{2\sqrt t} \bigg) + \frac{2}{\sqrt{\beta}}\bigl(\mathfrak{X}_{\mu, n}(t)-\mathfrak{X}(t)\bigr)+ \O_K(t^{-3/2})  .
\]
Combining Propositions~\ref{prop:2021} and \ref{prop:2022}, this implies that with probability  at least $1-3e^{- \frac{n^{3/2-2\epsilon}}{(\log n)^2}}$, 
\[
 \sup_{\lambda \in K}
    \sup_{T \geq s_n} 
    \bigg| \int_{s_n}^T
    \bigl(\rho_{\lambda,n}(t) -  \theta_\lambda(t)\bigr)  \d t \bigg|  \leq C_\epsilon n^{-\epsilon}
   +  \sup_{\lambda \in K}
    \sup_{T \geq s_n} 
    \bigg| \int_{s_n}^T
    \tfrac{2}{\sqrt{\beta}}\bigl(\mathfrak{X}_{\mu,n}(t) -  \mathfrak{X}(t)\bigr)  \d t \bigg| . 
\]
where $s_n = n+ \frac{\log n}{\sqrt{n}}$. 
Hence, to obtain Proposition~\ref{prop:X2int}, it suffices e.g. to show that 
\begin{equation} \label{Xvar1}
\P\Bigg( \sup_{\lambda \in K, T \geq s_n} 
    \bigg| \int_{s_n}^T
    \bigl( \mathfrak{X}_{\mu,n}(t)-  \mathfrak{X}_{0, n}(t)\bigr)  \d t \bigg|  \ge  C n^{-1/4} 
 \text{ or }
    \int_{s_n}^\infty
    \bigl| \mathfrak{X}(t)-  \mathfrak{X}_{0, n}(t)\bigr|  \d t \ge C n^{-2} \Bigg) \le C e^{- n^{3/2-2\epsilon}}. 
\end{equation}

First using the representation \eqref{XSDE_rep}, it holds
\[\begin{aligned} 
\int_{s_n}^T
    \bigl(   \mathfrak{X}_{\mu, n}(t)-\mathfrak{X}_{0,n}(t)\bigr)  \d t  
&=   \int_{s_n}^T  \int_{n}^t \frac{y_\mu(u)}{y_\mu(t)}  \mathfrak{X}_{0,n}(u) \big(\sqrt{u+\mu} -\sqrt{u}\big)  \d u    \d t \\
& = \int_n^T \mathfrak{X}_{0,n}(u) Y_{\mu,n}(u) \frac{\d u}{\sqrt u}
  \end{aligned}\]
where $\displaystyle Y_{\mu,n}(u) \coloneqq  \big(\sqrt{u(u+\mu)} - u\big) \int_{u \vee s_n}^T \frac{y_\mu(u)}{y_\mu(t)}  \d t  = \O(u^{-1/2}) $ for $u\ge n$ uniformly in $\lambda\in K$.
By Cauchy--Schwarz's inequality, this implies that  for any $\epsilon>0$
\[
\bigg|  \int_{s_n}^T   \bigl(   \mathfrak{X}_{\mu, n}(t)-\mathfrak{X}_{0,n}(t)\bigr)  \d t  \bigg| \lesssim
\sqrt{ \int_n^\infty \frac{\mathfrak{X}_{0, n}^2(u)}{u^{1/2+\epsilon}} \d u \int_{n}^\infty \frac{ Y_{\mu,n}^2(u)}{u^{1/2-\epsilon}} \d u  } .
\]
Using the tail bound \eqref{PIV}, this proves the first part of \eqref{Xvar1}. For the second part, we use that by definition, for all $t\ge n$, 
\[
 \mathfrak{X}(t)-  \mathfrak{X}_{0, n}(t) =  \frac{y_0(n)}{y_0(t)} \mathfrak{X}(n) . 
\]
This implies that 
\[
\int_{s_n}^\infty   \bigl| \mathfrak{X}(t)-  \mathfrak{X}_{0, n}(t)\bigr|  \d t  \lesssim  \frac{y_0(n)}{y_0(s_n)}
\frac{|\mathfrak{X}(n) |}{\sqrt{n}} \lesssim \frac{|\mathfrak{X}(n) |}{n^{5/2}}
\]
using that  $s_n = n+ \frac{\log n}{\sqrt{n}}$ and $y_0(s) =  e^{\tfrac43 s^{3/2}}$. 
Using Lemma~\ref{lem:XGN} again, this completes the proof of \eqref{Xvar1}. 
\end{proof}

{}

% ===== end sections/06-asymptotics-and-airy-function.tex =====

% ===== begin sections/07-airy-kernel-bound.tex =====
\section{Bound for the stochastic Airy kernel: Proof of Theorem~\ref{thm:Kest}}\label{sec:SAestimate}

Let  $J(T) = [ -e^{ T} , T]$  and
 recall that $K_\nu(T) = \big\{\lambda \in \C : \Re\lambda \in J(T), |\Im \lambda| \leq e^{-(1-1/\nu)\, T^\nu} \big\}$ for $\nu\in(1,3/2)$.
In this section, we prove  the estimates from Theorem \ref{thm:Kest} on the stochastic Airy kernel in a series of reductions.
Our first reduction is to a local estimate for $|\lambda| <e^{-(1-1/\nu)\, T^\nu}$,  which is to say that we reduce Theorem \ref{thm:Kest} to the following lemma.

\begin{lemma}\label{lem:K1}
  For any $\delta, \nu$ with $\nu \in (1, 3/2)$ and $0 < \delta < 1 - 1/\nu$,
there exist constants $C=C(\beta, \delta,\nu)$ and $c=c(\beta, \delta,\nu)$ so that
 \[
\P\biggl[ \sup\big\{ | \partial_{tu}\mathfrak{A}_\lambda(t,u) | E_\lambda^{-1}(t,u) : u,t \in J(T) , |\lambda| < e^{-(1-1/\nu)\, T^\nu} \big\}   > Ce^{\delta\, T^\nu}\biggr] \leq C e^{-cT^{(3\nu-1)/2}} .
\]
\end{lemma}

\noindent This is done by using invariance properties of the solution map $\T$ (see Proposition~\ref{prop:shift}) and then by a net argument; see Section~\ref{sec:Red1} for further details. 
% choosing a suitable net of $\lambda$ which has polynomial $N^{\O(1)}$ cardinality.

\medskip 

The second reduction is that by continuity of the stochastic Airy kernel in the parameter $\lambda\in\C$, it suffices to have control for $ \partial_{tu}\mathfrak{A}_\lambda(t,u) = \mathrm{f}_{\lambda,u}'(t)$  at $\lambda = 0$, i.e.\ we reduce Lemma \ref{lem:K1} to:

\begin{lemma}\label{lem:K2}
  For any $\nu \in (1, 3/2)$ and $\delta > 0$, there exist constants $C=C(\beta, \delta,\nu)$ and $c=c(\beta, \delta,\nu)$ so that
  \[
    \begin{aligned}
      \Pr
      \biggl[
	\sup
	\big\{  |\mathrm{f}_{0,u}'(t) | E_0^{-1}(t,u)
	  :  u, t\in J(T)
	\big\}
	>  C e^{\delta\, T^\nu}
      \biggr]
      \leq C e^{-cT^{(3\nu-1)/2}}.
    \end{aligned}
  \]
\end{lemma}
\noindent  
The proof of Lemma~\ref{lem:K2} is analogous to that of Proposition~\ref{prop:stability}, and it relies on the fact that one can express $\mathrm{f}_{\lambda}'$ in terms of $\mathrm{f}_{0}'$ as a Volterra or Fredholm equation; see \eqref{Volterra1} in Section~\ref{sec:Red2} for further details. 

\medskip
 
The proof of Lemma~\ref{lem:K2} now amounts to controlling (with overwhelming probability) 
the real-valued Dirichlet solutions $\mathrm{f}_0'$ of the stochastic Airy equation  \eqref{SA1}. 
In the oscillatory direction $t\le 0$, this can be achieved by a simple \emph{energy estimate} which is presented in Section~\ref{sec:Red3}.
%As a consequence, we obtain the following bounds. 

\begin{lemma}\label{lem:K4}
  For any $\nu \in (1, 3/2)$ and $\delta > 0$, there exist constants $C=C(\beta, \delta,\nu)$ and $c=c(\beta, \delta,\nu)$ so that  any solution $\phi = \phi_0$ of \eqref{SA1} with $\lambda=0$ and initial data $\phi^2(1) + \phi'(1)^2 \leq 1$
  satisfies,
  \[
    \begin{aligned}
      \Pr
      \biggl[
	\sup
	\big\{
	  |\phi(t)| , |\phi'(t)|  %\frac{E(0, u)}{E(0, t)}
	  :
	  -e^{ T}
	  \leq t \leq 1
	\big\}
	\ge C e^{\delta\, T^\nu}
      \biggr]
     % & \\
      \leq  C e^{-cT^{2\nu - 1}}.
    \end{aligned}
  \]
\end{lemma}

We can obtain bounds for real-valued solutions of the stochastic Airy equation in the expanding direction by using Lemma~\ref{lem:phipv} and controlling the behavior of the Riccati diffusion for positive times.  
The arguments are somewhat similar to those of Section~\ref{sec:SAi}, except that we require quantitative bounds which hold with overwhelming probability. The following lemma summarizes these bounds.

\begin{lemma}\label{lem:K3}
For any $\nu \in (1, 3/2)$ and $\delta > 0$, there exist constants $C=C(\beta, \delta,\nu)$ and $c=c(\beta, \delta,\nu)$ so that for any $s\in[1, T]$, it holds for  any solution $\phi = \phi_0$ of \eqref{SA1} with $\lambda=0$ and initial data $\phi^2(s) + \phi'(s)^2 \leq 1$,
\[
\P\bigg[ \sup\big\{ |\phi(t)|  E_0^{-1}(t,s) , |\phi'(t)|  E_0^{-1}(t,s)
% | \mathrm{f}_{0,s}(t)| ,  | \mathrm{g}_{0,s}(t)| ,  | \mathrm{f}'_{0,s}(t)| ,  | \mathrm{g}'_{0,s}(t)|
: t\in [s, T] \big\}  \ge C e^{\delta\, T^\nu} \bigg] \le C e^{-cT^{(3\nu-1)/2}} .
\]
\end{lemma}

\noindent The proof of Lemma~\ref{lem:K3} is rather involved and is given in Section \ref{sec:Red5} --- see in particular Proposition~\ref{prop:GE} and Section~\ref{sec:K3}. 

\medskip

To complete the proof of Theorem~\ref{thm:Kest}, it remains to show that 
Lemma~\ref{lem:K4} and Lemma~\ref{lem:K3} imply Lemma~\ref{lem:K2}. 
Let us first observe that by combining these two Lemmas, we immediately obtain that for any solution $\phi = \phi_0$ of \eqref{SA1} with $\lambda=0$ and initial data $\phi^2(1) + \phi'(1)^2 \leq 1$,
\begin{equation} \label{UBsol}
\P\bigg[ \sup\big\{ |\phi(t)|  E_0^{-1}(t,1)  , |\phi'(t)|  E_0^{-1}(t,1) 
: t\in J(T) \big\}  \ge C e^{\delta\, T^\nu} \bigg] \le C e^{-cT^{(3\nu-1)/2}} . 
\end{equation}
This provides the required bound for $\mathrm{f}_{0,u}'$ when $u=1$. We are now going to show that this bound holds for all $u\in J(T)$ with a similar probability by a net argument.

\begin{proof}[Proof of Lemma~\ref{lem:K2}]
Let us recall from Lemma \ref{lem:A} that $\mathrm{f}_{0,u}'(t) =   \partial_{tu}\mathfrak{A}_0(t,u)$, so by anti-symmetry of the kernel $\mathfrak{A}_0$, it suffices to show that 
\[
  \Pr
      \biggl[
	\sup
	\big\{  |\mathrm{f}_{0,u}'(t) | E_0^{-1}(t,u) 
	  :  u, t\in J(T) , u <t
	\big\}
	\ge C e^{2\delta\, T^\nu}
      \biggr]
      \leq C e^{-cT^{(3\nu-1)/2}}.
\]
Note that the small parameter $\delta>0$ which controls the error is only involved in the constants $C,c$.
By Lemma \ref{lem:A}, we can represent for $u,t\in\R$,
\[
  \mathrm{f}_{0,u}'(t)= 
  \mathrm{f}'_{0,1}(t)
  \mathrm{g}'_{0,1}(u)
  -\mathrm{f}'_{0,1}(u)
  \mathrm{g}'_{0,1}(t) , 
\]
so that by \eqref{UBsol}, we have for $u,t \in J(T)$ with $u\le 1 \wedge t$,
\begin{equation} \label{f'UB1}
  |\mathrm{f}_{0,u}'(t)| \lesssim e^{\delta\, T^\nu} \big( |  \mathrm{f}'_{0,1}(t)| +   |\mathrm{g}'_{0,1}(t)|  \big)\lesssim e^{2\delta\, T^\nu} E_0(t,u)
\end{equation}
with probability at least $1- C e^{-cT^{(3\nu-1)/2}} $. Here we have used that 
$E_0(1,u) \le E_0(t,u)$ and
$E_0(1,u) \lesssim 1$ for $u\le 1 \wedge t$. 

In the case $1\le u \le t$, we cannot use the same argument and we rely on the continuity of $u\mapsto \mathrm{f}'_{0,u}$ and Lemma~\ref{lem:K3} instead.
Let  $\left\{ u_k \right\}$ be a uniform net of the interval $[1, T^{(1-\delta)\nu}]$ with cardinality at most $e^{4T}$.
By Lemma~\ref{lem:K3}, we have
\[
 \sup_{k}\big\{ | \mathrm{g}'_{0,u_k}(t)|  E_0^{-1}(t,u_k) , | \mathrm{f}'_{0,u_k}(t)|  E_0^{-1}(t,u_k)
: t\in [u_k, T] \big\}  \le C e^{\delta\, T^\nu}
\]
with probability $1- C e^{-cT^{(3\nu-1)/2}} $ (up to adapting the constants $C,c$).
By Lemma \ref{lem:A} again, we have
\[
  \mathrm{f}'_{0,u}(t) =   \mathrm{f}'_{0,u_k}(t)
  \mathrm{g}'_{0,u_k}(u)
  -\mathrm{f}'_{0,u_k}(u)
  \mathrm{g}'_{0,u_k}(t) , \qquad u \in [u_k,u_{k+1}] , \quad t\ge u.
\]
By \eqref{eq:Et}, $E_0(u_k,u_{k+1}) \le 2$, 
so that for all $t\ge u\ge 1$, 
\begin{equation} \label{f'UB2}
  |\mathrm{f}_{0,u}'(t)| \lesssim e^{\delta\, T^\nu} \big( |  \mathrm{f}'_{0,u_k}(t)| +   |\mathrm{g}'_{0,u_k}(t)|  \big)\lesssim e^{2\delta\, T^\nu} E_0(t,u) .
\end{equation}
 
 By combining the bounds \eqref{f'UB1} and  \eqref{f'UB2}, this gives the required (uniform) control for $  |\mathrm{f}_{0,u}'(t)|$. 
\end{proof}

\subsection{Reduction 1: from a global estimate in \texorpdfstring{$\lambda$}{lambda} to a local estimate.}
\label{sec:Red1}

In this section, we assume Lemma \ref{lem:K1} and give the proof of Theorem \ref{thm:Kest}.
The key insight here is that the kernel $\partial_{tu}\mathfrak{A}_\lambda(t,u)$ has an invariance in law: for any fixed $\alpha \in \R$, 
\begin{equation}
  (\mathfrak{A}_\lambda(t,u) : \lambda \in \C, t,u \in \R)
  \lawequals
  (\mathfrak{A}_{\lambda-\alpha}(t+\alpha,u+\alpha) : \lambda \in \C, t,u \in \R).
  \label{eq:SAK_invariance}
\end{equation}

Let $(\alpha_k)$ be a uniform $2 e^{-(1-1/\nu)\, T^\nu}$-net of the interval $J(T)$, with cardinality at most $e^{T + (1-1/\nu)\, T^\nu}$. If $T$ is sufficiently large, $K_\nu(T) \subset \bigcup_{k}\big\{ \lambda\in\C : |\lambda-\alpha_k| \le 2 e^{-(1-1/\nu)\, T^\nu} \big\}$, so that by a union bound,
\[ \begin{aligned}
& \P\biggl[ \sup\big\{ | \partial_{tu}\mathfrak{A}_\lambda(t,u) | E_\lambda^{-1}(t,u) :
u,t\in J(T) , \lambda \in K_\nu(T) \big\}   > e^{\delta\, T^\nu}\biggr]  \\
&\qquad\qquad \le \textstyle{\sum_{k}}
\P\biggl[ \sup\big\{ | \partial_{tu}\mathfrak{A}_\lambda(t,u) | E_\lambda^{-1}(t,u) :
u,t\in J(T) ,  |\lambda-\alpha_k| \le 2 e^{-(1-1/\nu)\, T^\nu}  \big\}   > e^{\delta\, T^\nu}\biggr] \\
&\qquad\qquad \le e^{T + (1-1/\nu)\, T^\nu} \,
 \P\biggl[ \sup\big\{ | \partial_{tu}\mathfrak{A}_\lambda(t,u) | E_\lambda^{-1}(t,u) :
u,t\in 2J(T) ,  |\lambda| \le 2 e^{-(1-1/\nu)\, T^\nu}  \big\}   > e^{\delta\, T^\nu}\biggr]
%C e^{-T^{(1+\epsilon)\nu}} .
\end{aligned}\]
where we used \eqref{eq:SAK_invariance} to bound the probabilities at the second step. By
Lemma \ref{lem:K1}, this probability is at most  $C e^{-T^{(1+\epsilon)\nu}}$ for some $
0< \epsilon< 1-1/\nu$ and a constant $C=C(\beta, \delta,\nu, \epsilon)$.
To complete the proof of Theorem~\ref{thm:Kest}, it remains only to justify \eqref{eq:SAK_invariance}.
Recall that from Definition \ref{def:SAK}, 
\[
  \mathfrak{A}_\lambda(t,u)
  =
  \mathrm{f}_{\lambda,0}(t)
  \mathrm{g}_{\lambda,0}(u)
  -
  \mathrm{f}_{\lambda,0}(u)
  \mathrm{g}_{\lambda,0}(t)
  \quad
  \text{for all}
  \quad \lambda \in \C, t,u \in \R
\]
From Proposition \ref{prop:shift}, for any fixed $\alpha \in \R$,  this entire process therefore has the same law as
\[
  \mathrm{f}_{\lambda-\alpha,\alpha}(t+\alpha)
  \mathrm{g}_{\lambda-\alpha,\alpha}(u+\alpha)
  -
  \mathrm{f}_{\lambda-\alpha,\alpha}(u+\alpha)
  \mathrm{g}_{\lambda-\alpha,\alpha}(t+\alpha)
  \quad
  \text{for all}
  \quad \lambda \in \C, t,u \in \R.
\]
By Lemma \ref{lem:A}, this kernel equals $\mathfrak{A}_{\lambda-\alpha}(t+\alpha,u+\alpha)$.

\subsection{Reduction 2: \texorpdfstring{from a local estimate in $|\lambda|\le e^{-\delta\, T^\nu}$ to an estimate at $\lambda=0$}{from a local estimate to an estimate at zero}.}
\label{sec:Red2}

We show that (deterministically), on the event
\begin{equation} \label{event2}
\big\{ |\mathrm{f}'_{0,u}(t)| \le  e^{\delta\, T^\nu} E_0(t,u) : \forall\, t, u\in J(T) \big\} ,
\end{equation}
a similar bound holds for $\mathrm{f}_\lambda$ uniformly  for all $|\lambda| \le e^{-(1-1/\nu)\, T^\nu}$  provided that $1-1/\nu>\delta$ and $T$ is sufficiently large.
This establishes that Lemma \ref{lem:K2} implies Lemma \ref{lem:K1}. The proof relies on the Volterra-type structure of the stochastic Airy equation \eqref{SA2}. 
Let us denote by 
$\mathrm{h}_{\lambda,s} := \mathrm{f}_{\lambda,s}' - \mathrm{f}_{0,s}'$ for $s\in\R$, $|\lambda|<e^{-\delta\, T^\nu}$
and $\displaystyle K(t,u) =  \int_u^t  \mathrm{f}_{0,t}(v) \d v  $. 
First observe that according to \eqref{eq:Et}, 
\[
\sup\big\{ E_\lambda(t,u)  E_0^{-1}(t,u)  : t, u\in J(T) , |\lambda| \le e^{-(1-1/\nu)\, T^\nu} \big\}
= \exp(\o_\nu(1))
\]
as $T\to\infty$, so it suffices to show that with $E=E_0$,
 \[
\big\{ | \mathrm{f}_{\lambda,u}'(t) | \le e^{\delta\, T^\nu} E(t,u): \forall u,t \in J(T) , |\lambda| < e^{-\delta\, T^\nu} \big\}  .
\]

Second, observe that on the event \eqref{event2}, 
the bound \eqref{fbound} holds with $\lambda=0$ for all $t, s \in J(T)$ and we can bound 
\begin{equation} \label{boundK}
|K(t,u)|  
\le E(t,u)   \int_u^t |\mathrm{f}_{0,t}(v)|  E(t,v)^{-1}  \d v 
\le  (1+ |J(T)|)^2 e^{\delta\, T^\nu}  E(t,u)  . 
\end{equation}

Now, since both $\mathrm{f}_{\lambda,s}' , \mathrm{f}_{0,s}'$ solve \eqref{SA2} with $c_1=1, c_2=0$ and $\Es_\lambda(t,s) = \Es_0(t,s) + \lambda(t-s)$, the function $\mathrm{h}_{\lambda,s}$ solves the equation:
\[
  \mathrm{h}_{\lambda,s}(t)
  =\int_s^t \Es_0(t,v)
  \mathrm{h}_{\lambda,s}(v)
  \,\d v + \lambda \zeta(t) ,
   \quad\text{where }
   \quad
	   \zeta(t) =   \int_s^t \bigl(1+ (t-v)\mathrm{f}'_{\lambda,s}(v)\bigr)\d v  = \int_s^t \mathrm{f}_{\lambda,s}(v)\d v . 
  \]
Since $\zeta\in C^1(\R)$, applying Proposition~\ref{prop:sol} with $ \mathfrak{A}= \mathfrak{A}_0$, we obtain
after integrating by parts twice
\[
  \mathrm{h}_{\lambda,s}(t)
  = -\lambda \int_s^t \partial_t\mathfrak{A}(t,u) \d \zeta(u)
  = -\lambda \int_s^t   \partial_t\mathfrak{A}(t,u) \mathrm{f}_{\lambda,s}(u) \,\d u
  = \lambda K(t,s) + \lambda \int_s^t K(t,u)\mathrm{f}'_{\lambda,s}(u) \d u,
\]
where we used that $\partial_t\mathfrak{A}(t,u) = - \mathrm{f}_{0,t}(u) = \partial_u K(t,u)$
(cf.~Lemma~\ref{lem:A}).
In summary, we arrive at a Volterra-type equation for  $\mathrm{h}_{\lambda,s}$ where the parameter $\lambda$ is small: 
\begin{equation} \label{Volterra1}
  \mathrm{h}_{\lambda,s}(t) = 
  \lambda \int_s^t   \mathrm{f}_{0,t}(u)  \mathrm{f}_{0,s}(u) \,\d u 
  + \lambda \int_s^t K(t,u)\mathrm{h}_{\lambda,s}(u) \,\d u.
\end{equation}

We can use this equation to deduce a (uniform) bound for $\mathrm{h}_{\lambda,s}$.  
On the event \eqref{event2}, 
\[
\bigg|  \int_s^t   \mathrm{f}_{0,t}(u)  \mathrm{f}_{0,s}(u) \,\d u  \bigg|
\le  E(t,s) \int_s^t  | \mathrm{f}_{0,t}(u)|  E(t,u)^{-1}  |\mathrm{f}_{0,s}(u)|  E(u,s)^{-1} \,\d u  
\le |J(T)| e^{2\delta\, T^\nu}  E(t,s) . 
\]
Let us denote
$ \mathcal{M} := \sup_{t \in J(T), |\lambda| \le e^{-\delta\, T^\nu}} \bigl(| \mathrm{h}_{\lambda,s}(t)|E(t,s)^{-1}\bigr).$
The previous bound and \eqref{boundK} imply that 
\[
|  \mathrm{h}_{\lambda,s}(t) | \le    |\lambda| |J(T)| e^{2\delta\, T^\nu}  E(t,s) +  |\lambda| \mathcal{M}  (1+ |J(T)|)^2 e^{\delta\, T^\nu} 
\int_s^t  E(t,u)E(s,u) \d u    .
\]
Hence, if $|\lambda| \le  e^{-\delta\, T^\nu}$ and $\delta< 1-1/\nu$, we obtain when $T$ is sufficiently large (depending only on the parameters $\delta,\nu$),
\[
\mathcal{M} \le (e^{\delta\, T^\nu} +  \mathcal{M})/2 .
\]
By rearranging, we conclude that $\mathcal{M} \le e^{\delta\, T^\nu}$ and
$|\mathrm{f}_{\lambda,s}(t)'| \le \bigl(\mathcal{M} + | \mathrm{f}_{0,s}'(t)|E(t,s)^{-1}\bigr)E(t,s) \le 2 e^{\delta\, T^\nu}E(t,s)$,
 this bound being uniform over all $s,t \in J(T)$ and $|\lambda| \le  e^{-\delta\, T^\nu}$.

\subsection{Control of solutions in the oscillatory direction} \label{sec:Red3}

We give the proof of Lemma~\ref{lem:K4} by formulating a tail bound for
$\sup_{t\in[-\tau,s]} \big\{|\phi(t)| + \sqrt{1+t_-}|\phi'(t)| \big\}$ where $ \tau, s>0$ and we assume that $\phi^2(s) + \phi'(s)^2 \leq 1$.
Since we are interested in the behavior of a single solution for large $\tau$, it will be convenient to reverse time  in the stochastic Airy equation \eqref{SA1}. Hence, we suppose that $\phi$ is a strong solution of the diffusion: for $t\ge -s$, 
\begin{equation} \label{DN1}
  \d \phi'(t) = \bigl(-t\d t + \tfrac{2}{\sqrt{\beta}}\,\d B(t)\bigr)\phi(t) , \qquad \text{ where } \quad \begin{cases} \phi(-s)=  \alpha, \\  \phi'(-s)=  \pm \sqrt{1-\alpha^2}, \end{cases}
\end{equation}
with $\alpha\in[-1,1]$ and $B$ is a standard Brownian motion (this Brownian motion must be a reversal of the one that appears in \eqref{SA1}, but only the law of $\phi$ will be important here). 

Let us choose a function $U : \R \to (0,1]$ such that $U \in C^1$, $U' \le 0$, $U(t) = 1/t$ for $t\ge 2$ and $U(t)=1$ for $t\le 0$. We define the \emph{Lyapunov function}:
\begin{equation} \label{eq:Potential}
  \mathcal{H}(t)
  \coloneqq 
  U(t) \phi'(t)^2 + \phi^2(t)
%  \quad\text{where}\quad
%  \mathcal{U}(t) = 
%  \begin{cases} \frac{1}{t}, &\text{if } t \geq 1, \\
%    2-t, &\text{if } -1 \leq t \leq 1. \\
%  \end{cases}
\end{equation}
and prove the following tail-bound.

\begin{lemma}\label{lem:energy}
  Let $\phi$ be any solution of \eqref{DN1}, then for any $\tau \ge 1$ and $\Lambda>0$,
  \[
    \P\big[  \sup_{-s \leq t \leq \tau} \mathcal{H}(t) \geq e^{ \frac{(s+2)^2}{2}+\Lambda} \big]
    \leq \exp\bigg( -\frac{\beta \Lambda^2}{8(2+s+ \log \tau)} \bigg).
  \]
\end{lemma}
Before proceeding to the proof, we observe that control on $\mathcal{H}$ implies control on both $\phi$ and $\phi',$ since $| \phi(t)| \le  \sqrt{\mathcal{H}(t)}$ and $| \phi'(t)| \le  \sqrt{(1+t_+)\mathcal{H}(t)}$ --- here we used that by construction $U(t) \ge \frac{1}{1+t_+}$ for all $t\in\R$. Hence taking
$\tau = e^T$ and $\Lambda = \delta T^\nu - \tfrac12 T$, we obtain for $s\in [0, T]$ and $T$ sufficiently large,
\begin{equation} \label{energybound}
 \Pr
      \biggl[
	\sup
	\big\{
	  |\phi'(t)| e^{-\frac{(s+2)^2}{2}}
	  : t\in [-s, e^T]
	\big\}
	> e^{\delta\, T^\nu}
      \biggr]
      \le    \exp\bigg( -\frac{c\beta  T^{2\nu} }{1+T} \bigg) ,
      %\leq e^{-CT^{(1+\epsilon)\nu}}.
\end{equation}
where $c = c(\delta)>0$.
Note that the bound \eqref{energybound} displays the wrong growth for large $s$, namely a factor $e^{s^2/2}$ instead of $e^{\frac 23 s^{3/2}}$, besides it remains accurate for all  $s\in [0, \delta' T^{\nu/2}]$ with $\delta'>0$ sufficiently small. 
Choosing $s=1$, the bound \eqref{energybound} gives $|\phi'(t)| \le Ce^{\delta T^\nu}$ for $t\in[-1, e^T]$ with probability at least $1 - Ce^{-cT^{2\nu-1}}$. On the same event, the energy lemma furnishes $\mathcal{H}(t) \le e^{9/2 + \delta T^\nu - T/2}$, so $|\phi(t)| \le \sqrt{\mathcal{H}(t)} \le e^{\delta T^\nu}$ for $T$ sufficiently large. Reverting to the original time convention $t \mapsto -t$ maps $[-1, e^T]$ to $[-e^T, 1]$, concluding the proof of Lemma~\ref{lem:K4}.

\begin{proof}
Note that $  \mathcal{H}>0$ as non-trivial solutions of the stochastic Airy equation have (almost surely) no double zero. 
Let us denote by $V(t) = \frac{2U(t)\phi(t)\phi'(t)}{  \mathcal{H}(t)} $
and $M(t) = \tfrac{2}{\sqrt{\beta}}\int_{-s}^t V(u) \d B(u)$. 
Since $U\in C^1$ and $U>0$, by applying It\^ o's formula and \eqref{DN1}, we verify that
\[
     \d \log \mathcal{H}(t) =  \big( \big(U(t)^{-1}-t\big) V(t) +  U'(t)\mathcal{H}(t)^{-1} \phi'(t)^2- \tfrac2\beta V(t)^2 \big) \d t + \d M(t). 
\]
Thus, using that $U'<0$ and $U(t) = 1/t$ for $t\ge 2$, this implies that for $t\ge -s$, 
\[
     \d \log \mathcal{H}(t) \le (2-t)_+ \d t + \d M(t) .
 \]
As $ \mathcal{H}(-s) \le 1 $, this allows to control for $t\ge 0$, 
\[
\log \mathcal{H}(t) \le \tfrac{(s+2)^2}{2} + M(t). 
\]

Moreover $V(t)^2 \le U(t)$ and $U\in(0,1]$,  the quadratic variation of the martingale $M$ satisfies for  $t\ge1$, 
  \[
    \langle M \rangle (t)= \frac4\beta \int_{-s}^t  V(u)^2\d u
\le  \frac4\beta \int_{-s}^t  U(u) \d u \le \tfrac4\beta \big(2+s+ \log(t) \big) .
    \] 
Hence by Freedman's inequality, for any $\Lambda>0$ and $\tau\ge 1$
\[
 \P\big[ \sup_{-s \leq t \leq \tau} M(t)  \ge \Lambda \big] \le \exp\big(- \tfrac{\beta \Lambda^2}{8(2+s+\log \tau)}\big) .
\]
This provides the required bound for $\mathcal{H}$.
\end{proof}
\section{Control of solutions in the expanding direction} \label{sec:Red5}

The goal of this section is to prove Lemma~\ref{lem:K3}. That is, we would like to control the growth\footnote{Indeed by \eqref{SA1}, $\phi_0'$ can be represented in terms of  $\phi_\lambda$ as an It\^o integral, and it suffices to control the growth of $\phi_0$.  See Section~\ref{sec:K3} for further details.} of a solution $\phi_0$ of the stochastic Airy equation \eqref{SA1} with parameter $\lambda=0$ in the expanding direction $t\ge s\ge 1$.
In addition, using the invariance in law of the stochastic Airy equation (see Proposition~\ref{prop:shift}):
$\big(\phi_0(t)\big)_{t\ge \lambda}  \overset{\rm law}{=}   \big(\phi_\lambda(t)\big)_{t\ge 0} $ for any $\lambda\ge 1$.
Hence, it suffices to prove the following estimates.

\begin{proposition} \label{prop:GE}
Fix $\nu \in (1,2)$ and $\lambda\in[1,T]$ for $T \ge 1$ sufficiently large.
Let  $\phi_\lambda$ solve \eqref{SA1} with initial data $\left\{ \phi_\lambda(0)=c_1 > 0, \phi_\lambda'(0)=c_2 \in \R \right\}.$
%c_1>0$ and $c_2\in\R$.
 There exists two constants $C_\nu, C_{\beta,\nu}>0$ such that for any $R\in[ C_\nu, T^{(2-\nu)/2}]$,  the following estimate holds with probability at least  $1-C_{\beta,\nu}e^{-c_\beta R\, T^\nu}$,
\begin{equation} \label{phiest}
\sup_{t\in[0,T]} \big(|\phi_\lambda(t)| E_\lambda^{-1}(t)  \big) \le  (|c_1|+2|c_2|)  \exp\big(\sqrt{ R}  T^{(3\nu+1)/4} \big) ,
\end{equation}
where 
\begin{equation} \label{def:E}
E_\lambda(t) \coloneqq E_\lambda(t,0) =  \exp\bigg(\int_0^t  \sqrt{v+\lambda}  \d v \bigg)
=  \exp\bigg(\frac23 (t+\lambda)^{3/2} -\frac23 \lambda^{3/2}  \bigg)  , \qquad t\ge 0. 
\end{equation}
\end{proposition}

The proof of Proposition~\ref{prop:GE} consists in expressing $\phi_\lambda$ in terms of its Riccati transform $\rho$ and then in analyzing the long-time behavior of $\rho$. 
Namely, by Lemma~\ref{lem:phipv}, we have 
\begin{equation}\label{phipv}
|\phi_\lambda(t)| = |c_1| \exp\biggl( \operatorname{pv} \int_{0}^{t}\rho(u)\d u\biggr).  
\end{equation}
Note that the estimate \eqref{phiest} essentially comes from the fact $\rho$ solves \eqref{Ric0} with $s=0$ and the drift forces the solutions  to become stationary: $\rho(t) \sim \sqrt{t+\lambda}$ (since the branch $-\sqrt{t+\lambda}$ is unstable).

Our proof strategy is to quantitatively control the integral on the RHS of \eqref{phipv} by cutting short regions around the zeros $\left\{ \mathfrak{z}_k \right\}$ of $\phi_\lambda$ where $\rho$ blows down, for which we provide quantitative estimates by taking advantage of the cancellation implicit in the principal value.  Away from the zeros, we use comparisons between the Riccati diffusion $\rho$ and a linearized process, which in comparison to the regions near zeros, is much simpler.

For $\nu\in(1,2)$, we use $T^{(\nu+1)/4}$ as the Riccati cutoff. Let $ \mathbf{N}_t  = \mathbf{N}_{[0,t)}(\rho) = \#\big\{ k :   \mathfrak{z}_k \in[0, t)\big\}$ be the counting function of the zeros of $\phi_\lambda$.
%as in Lemma \ref{lem:phipv}.
Define the stopping times $0 < \tau_1 <\sigma_1 < \tau_2 < \sigma_2< \cdots $ as follows: 
\begin{equation} \label{stoppingtimes}
\tau_k= \inf\big\{ t > \sigma_{k-1} : \rho(t) = -T^{(\nu+1)/4} \big\} 
\qquad\text{and}\qquad 
\sigma_k= \inf\big\{ t > \mathfrak{z}_k : \rho(t) = T^{(\nu+1)/4} \big\}  . 
\end{equation}
\emph{A priori} it could be that the process returns from $\tau_k$ to $\sigma_k$ without blowing down, but we shall show this is a low-probability event that can be discarded and it is a consequence of our estimates that with overwhelming probability each interval $(\tau_k,\sigma_k)$ contains exactly one zero $ \mathfrak{z}_k$.

Setting $\sigma_0=0$, it holds for any $t\ge 0$, 
\[
\operatorname{pv} \int_{0}^{t}\rho(u)\d u
\le\sum_{k=0}^{\mathbf{N}_t} \int_{\sigma_k}^{\tau_{k+1}\wedge t} \hspace{-.5cm} \rho(s)ds
+  \sum_{k=1}^{\mathbf{N}_{t}} \mathrm{pv} \int_{\tau_k}^{\sigma_k}\rho(s)ds   
\]
so that by \eqref{phipv}, we obtain
\begin{equation*}
|\phi_\lambda(t)|  \le c_1 \exp\bigg( \sum_{k=0}^{\mathbf{N}_t} \int_{\sigma_k}^{\tau_{k+1}\wedge t} \hspace{-.5cm} \rho(s)ds \bigg) \prod_{k=1}^{\mathbf{N}_{t}} \bigg| \frac{\phi_\lambda(\sigma_k)}{\phi_\lambda(\tau_k)} \bigg| . 
\end{equation*}

Suppose that we construct a continuous diffusion $x$ such that $x(t) \ge \rho(t)- \sqrt{t+\lambda}$ for all $t\in \bigcup_{k\ge 0} [\sigma_k, \mathfrak{z}_{k+1}) $. Then by \eqref{def:E}, the previous bound implies that 
\begin{equation*}
|\phi_\lambda(t)| E_\lambda^{-1}(t)   \le c_1 \exp\bigg( \sum_{k=0}^{\mathbf{N}_t} \int_{\sigma_k}^{\tau_{k+1}\wedge t} \hspace{-.5cm} x(s)ds \bigg) \prod_{k=1}^{\mathbf{N}_{t}} \bigg| \frac{\phi_\lambda(\sigma_k)}{\phi_\lambda(\tau_k)} \bigg| . 
\end{equation*}

By linearity, it suffices to prove the estimate for the two normalized solutions $\mathrm f_{\lambda,0}$ and $\mathrm f_{\lambda,0}+\mathrm g_{\lambda,0}$, whose initial data are respectively $(1,0)$ and $(1,1)$.  We fix either solution, so that $c_1=1$ and $c_2\in\{0,1\}$. We denote for $T\ge 0$, 
\[ 
\mathrm{I}_T  =   \sup_{t\in[0,T]} \bigg( \sum_{k=0}^{\mathbf{N}_t}\int_{\sigma_k}^{\tau_{k+1} \wedge t} \hspace{-.3cm} x(u)\d u \bigg)
\qquad\text{and}\qquad
\mathrm{II}_T   =   \prod_{k=1}^{\mathbf{N}_{T}} \bigg| \frac{\phi_\lambda(\sigma_k)}{\phi_\lambda(\tau_k)} \bigg| . 
\]
In the sequel, we refer to $\mathrm{I}$ as the \emph{non-singular contribution} and $\mathrm{II}$ as the \emph{singular contribution}.  
In terms of these quantities, we obtain the following uniform estimate for any $T
\ge 0$,  
\begin{equation} \label{rb1}
\sup_{t\in[0,T]} \big(|\phi_\lambda(t)| E_\lambda^{-1}(t)  \big) \le  \mathrm{II}_T  \exp(\mathrm{I}_T) .
\end{equation}

%This shows that for any $T>0$,
%\[
%\sup_{t\in[0,T]}\Bigg( \operatorname{pv} \int_{0}^{t}\rho(u)\d u - \frac{2}{3}
%\big((t+\lambda)^{3/2} - \lambda^{3/2} \big) \Bigg)  \le \mathrm{I}_T + \mathrm{I}_T .
%\]

In Section~\ref{sec:bd}, we begin our analysis by providing bounds for the counting function $\mathbf{N}_{t}$ of blow-downs when the parameter $\lambda \ge 1$. 
Notice that as $\lambda$ (or equivalently time) increases, these blow-downs become progressively more expensive, and we exploit this fact to show that with overwhelming probability, there are at most $\sqrt{T^\nu}$ blow-downs; see Proposition~\ref{prop:zerobounds}. 

In Section~\ref{sec:nsc}, we prove the following estimate for the non-singular contribution. 

\begin{proposition} \label{prop:nsc}
Let $\nu\in(1,2)$, $\lambda\ge1$, and $T\ge1$ be sufficiently large. It holds for any $R\in[C_\nu,\sqrt T]$,
\begin{equation*}
\P\big[   \mathrm{I}_T  \ge  \sqrt{ R}  T^{(3\nu+1)/4} \big]
\le C_{\beta,\nu}e^{-c_\beta R\, T^\nu} .
\end{equation*}
\end{proposition}

The proof is based on constructing a suitable diffusion $x$ such that  $x(t) \ge \rho(t)- \sqrt{t+\lambda}$ away from the blow-downs $\{ \mathfrak{z}_k\}$ and our previous control of the counting function $\mathbf{N}_{T}$.

To control the singular contribution $\mathrm{II}$, let us denote by 
$\v_k = \mathrm{f}_{\lambda,\tau_k}$ and  $\u_k= \mathrm{g}_{\lambda,\tau_k}$ the  Dirichlet, respectively Neumann, solutions of the  stochastic Airy equation \eqref{SA1} at the stopping time $\tau_k$ for $k\ge 1$. 
These two processes are adapted and by linearity of the equation \eqref{SA1}, it holds that for all $t\ge \tau_k$, 
\[
\phi_\lambda(t) = \phi_\lambda(\tau_k) \v_k(t) + \phi_\lambda'(\tau_k) \u_k(t) . 
\]
In particular,  this implies that 
$\frac{\phi_\lambda(\sigma_k)}{\phi_\lambda(\tau_k)} = \v_k(\sigma_k) -T^{(\nu+1)/4} \u_k(\sigma_k)$ for any $k\ge 1$.
The intuition is that with overwhelming probability, $(\sigma_k-\tau_k) \lesssim T^{-(\nu+1)/4}$ and both $ \v_k$ and $\u_k$ are approximately linear in this short interval.
This allows us to show that $\big|\frac{\phi_\lambda(\sigma_k)}{\phi_\lambda(\tau_k)}\big| \lesssim 1$
with overwhelming probability; this being made precise by Proposition~\ref{prop:singular} below. 
In particular, we prove the following proposition in Section~\ref{sec:sc}. 

\begin{proposition} \label{prop:sc}
For any $\nu\in(1,2)$, fix $\lambda\in[1,T]$ and $T$ sufficiently large that $T\le c_\beta T^{(\nu+1)/2}$, where $c_\beta>0$ is sufficiently small.
There exist constants $c_\nu, C_\beta>0$ such that
\[
\P\big[  \mathrm{II}_T  \le  e^{c_\nu\, T^{3\nu/4}}\big] \ge 1- C_\beta e^{- c_\beta  T^{(\nu+2)/2}} . \qedhere
\]
\end{proposition}

By combining Propositions \ref{prop:nsc} and \ref{prop:sc} with the upper-bound \eqref{rb1}, we conclude that for any $R\in[C_\nu,T^{(2-\nu)/2}]$, both normalized solutions satisfy, with probability at least  $1-C_{\beta,\nu}e^{-c_\beta R\, T^\nu}$,
\[
\sup_{t\in[0,T]} \big(|\phi_\lambda(t)| y_\lambda^{-1}(t)  \big) \le \exp\big(2 \sqrt{ R}  T^{(3\nu+1)/4} \big) .
\]
For arbitrary initial data,
\[
\phi_\lambda=(c_1-c_2)\mathrm f_{\lambda,0}+c_2(\mathrm f_{\lambda,0}+\mathrm g_{\lambda,0}),
\]
and $|c_1-c_2|+|c_2|\le |c_1|+2|c_2|$.  A union bound for the two normalized estimates therefore gives the asserted initial-data prefactor.  This completes the initial-data reduction in Proposition~\ref{prop:GE}.

\subsection{Proof of Lemma \ref{lem:K3}} \label{sec:K3}

Having Proposition \ref{prop:GE}, we deviate briefly to show that this implies the desired estimate in Lemma \ref{lem:K3}.
Fix $\nu\in(1,3/2)$ and $\delta>0$.
  First observe that by applying Proposition \ref{prop:GE} to the initial data $\left\{ \psi(0) = 1, \psi'(0) = 0  \right\}$ and 
$\left\{ \psi(0) = 1, \psi'(0) = 1 \right\}$ with $R=\frac{\delta^2}{16} T^{(\nu-1)/2}$.  Since $\nu<3/2$, we have $(\nu-1)/2<(2-\nu)/2$, so this choice belongs to $[C_\nu,T^{(2-\nu)/2}]$ for all sufficiently large $T$ and Proposition~\ref{prop:GE} applies.  We conclude that if $\lambda\in[1,T]$,
\[
  |\mathrm{f}_{\lambda,0}(t)| \leq E_\lambda(t)e^{\delta\, T^\nu/2}
  \quad
  \text{and}
  \quad
  |\mathrm{f}_{\lambda,0}(t) + \mathrm{g}_{\lambda,0}(t)| \leq 2 E_\lambda(t)e^{\delta\, T^\nu/2}
  \quad
  \text{for all}
  \quad
  0 \leq t \leq  T,
\]
with  probability $1 -  Ce^{-c \delta^2\, T^{(3\nu-1)/2}}.$  Taking linear combinations, we obtain a similar bound for any solution of \eqref{SA1} with $\phi_\lambda^2(0) + \phi_\lambda'(0)^2 \leq 1$, 
\[
  |\mathrm{\phi}_{\lambda}(t)|  \leq 3 E_\lambda(t)e^{\delta\, T^\nu/2}
  \qquad
  \text{for all}
  \quad
  0 \leq t \leq T,
\]

%To get a similar control for $\mathrm{\phi}_{\lambda}'$, we can use the SDE \eqref{SA1}, 
%\[
%\mathrm{\phi}_{\lambda}'(t) = \mathrm{\phi}_{\lambda}'(0)+ \int_0^t u \mathrm{\phi}_{\lambda}(u) \,\d u
%  + M_\lambda(t) 
%  \qquad\text{where}\quad
% M_\lambda(t)  =  \int_\lambda^t \mathrm{\phi}_{\lambda}(u)\d B(u).
%\]

Moreover, using the invariance in law from Proposition \ref{prop:shift}, 
this implies that for any fixed $\lambda\in[1,T]$, 
\[
  |\mathrm{f}_{0,\lambda}(t)| \leq  E_\lambda(t-\lambda)e^{\delta\, T^\nu/2}
  \quad
  \text{and}
  \quad
  |\mathrm{g}_{0,\lambda}(t)| \leq 3 E_\lambda(t-\lambda)e^{\delta\, T^\nu/2}
  \quad
  \text{for all}
  \quad
  \lambda \leq t \leq  T,
\]
with the same probability.
To get a similar control for $\mathrm{f}_{0,\lambda}', \mathrm{g}_{0,\lambda}'$, we use the SDE  \eqref{SA1}, for $t\ge \lambda$, 
\begin{equation} \label{f'int}
  \mathrm{f}'_{0,\lambda}(t) = \int_\lambda^t u \mathrm{f}_{0,\lambda}(u) \,\d u
  + M_\lambda(t) 
  \qquad\text{where}\quad
 M_\lambda(t)  =  \tfrac{2}{\sqrt{\beta}}\int_\lambda^t \mathrm{f}_{0,\lambda}(u) \,\d B(u).
\end{equation}
Let $M_\lambda^\tau$ denote $M_\lambda$ stopped the first time the preceding envelope for $\mathrm f_{0,\lambda}$ fails.  For $0\le s\le T-\lambda$,
\[
 \big\langle M_\lambda^\tau\big\rangle(\lambda+s)
 \le \frac4\beta e^{\delta T^\nu}\int_0^sE_\lambda(u)^2\,\d u.
\]
By \eqref{def:E}, the same calculation as in Lemma~\ref{lem:G}\ref{eq:GlogGlower} gives
\[
 \bigg(\int_0^sE_\lambda(u)^2\,\d u\bigg)
 \bigg(\log\bigg(1+\int_0^sE_\lambda(u)^2\,\d u\bigg)\bigg)^{1/3}
 \lesssim E_\lambda(s)^2.
\]
Brownian scaling and the same martingale-representation argument as in the proof of Lemma~\ref{lem:XGN} (see \eqref{eq:Arnault0} and Lemma~\ref{lem:BMestimate}) therefore yield
\[
 \bigg\VERT\sup_{0\le s\le T-\lambda}
 \frac{|M_\lambda^\tau(\lambda+s)|}{E_\lambda(s)}\bigg\VERT_2
 \lesssim_\beta e^{\delta T^\nu/2}.
\]
Its probability of exceeding $e^{\delta T^\nu}$ is consequently at most $C_\beta\exp(-c_\beta e^{\delta T^\nu})$.
On the event supplied by Proposition~\ref{prop:GE}, $\tau_\lambda=T$ and hence $M_\lambda^\tau=M_\lambda$.  Combining the two exceptional probabilities and returning to \eqref{f'int}, we conclude that, with probability at least $1-Ce^{-cT^{(3\nu-1)/2}}$, for every $t\in[\lambda,T]$,
\[
 \begin{aligned}
 |\mathrm f'_{0,\lambda}(t)|E_\lambda^{-1}(t-\lambda)
 &\le e^{\delta T^\nu/2}E_\lambda^{-1}(t-\lambda)
 \int_0^{t-\lambda}(u+\lambda)E_\lambda(u)\,\d u
 +e^{\delta T^\nu}\\
 &\le e^{\delta T^\nu/2}\sqrt t+e^{\delta T^\nu}
 \lesssim e^{\delta T^\nu}.
 \end{aligned}
\]
The same stopped-martingale argument gives the analogous bound for the Neumann solution $\mathrm g'_{0,\lambda}$.  Since $E_\lambda(t-\lambda)=E(t,\lambda)$ according to \eqref{eq:Et}--\eqref{def:E}, this completes the proof. \qed

\subsection{Estimates for the number of blow-downs} \label{sec:bd}

Let us first estimate the probability that the Riccati diffusion hits 0 when started from a generic point. 

\begin{proposition} \label{prop:transitlaplace}
Let $\zeta = \inf\big\{ t\ge 0  : \rho(t) = 0  \} $ where $\rho$ is the Riccati diffusion \eqref{Ric0} with $s=0$ and $\lambda \ge 1$.
If $ \omega \in \big( 0, \sqrt{\lambda/8}]$, then it holds  for all $u \ge \frac{\beta \lambda^2}{8} \vee \frac{8}{\beta}$,
\[
  \Exp \exp(-u \zeta)
  \leq \exp\biggl(- \frac{\omega}{2\sqrt{2}}\bigg( \sqrt{\beta u}+ \frac{\lambda\beta}{2\sqrt{2}} \bigg)\biggr) . 
\]
\end{proposition}

\begin{proof}
The idea is to compare the Riccati diffusion to a simpler process $x$ which is driven by the same Brownian motion. 
  Let us define the following stopping times:  $\vartheta_0 = 0$ and for $k\ge 0$, 
\[
\zeta_{k+1} = \inf\big\{ t\ge \vartheta_{k} :  x(t) = 2\omega \text{ or } x(t) = 0  \big\}
\qquad
\text{and}
\qquad
\vartheta_{k+1} = \inf\big\{ t\ge \zeta_{k+1} : \rho(t) = \omega  \big\} ,
\]
where the process  $x$ is defined (piecewise) by
\[
\begin{cases}
\d x = \lambda/2 \d t + \tfrac{2}{\sqrt{\beta}}\,\d B ,  & t\in [\vartheta_k, \zeta_{k+1}),  \qquad x(\vartheta_k)= \omega  \\
x(t)= \omega  , & t\in [\zeta_{k+1}, \vartheta_{k+1}]
\end{cases} .
\]
We easily verify that $x(t) \le \rho(t)$ for all $t \le  \tau = \inf\big\{ t\ge 0  : x(t) = 0  \} $, so that $\tau\le \zeta$ almost surely.  
Indeed, since $\omega \le \sqrt{\lambda/8}$, we have by construction $x\le 2\omega$ and
\[
\begin{cases}
\d\rho \ge \lambda/2 \d t + \tfrac{2}{\sqrt{\beta}}\,\d B , & t\in  [\vartheta_k, \zeta_{k+1}) \cap \{ \rho \le 2\omega\},   \\
\rho(t) \ge \omega  , & t\in [\zeta_{k+1}, \vartheta_{k+1})
\end{cases} .
\]

Since $\tau\le \zeta$, the Laplace transform of $\zeta$ satisfies for any $u\ge 0$,
\begin{equation} \label{Laplace1}
  \Exp e^{-u \zeta} \le   \Exp e^{-u \tau} . 
\end{equation}
We can use that $\big(x(t)\big)_{t\le \zeta_1}$ is a stopped Brownian diffusion with variance parameter $4/\beta$ and drift $\lambda/2$ to compute the RHS of~\eqref{Laplace1}. 
We have the domination
$ \tau \ge \sum_{k\in\N} \1_{\tau = \zeta_k} \sum_{j =0}^{k-1} (\zeta_{j+1}-\vartheta_{j})$ 
and by the (strong) Markov property, $(\zeta_{j+1}-\vartheta_{j})_{j\in\N_0}$ are i.i.d. with law $\zeta_1$. 
We further claim that the event $\big\{ x(\zeta_1) = 0 \big\}$ is independent from the stopping time $\zeta_1$ which has Laplace transform:
\begin{equation} \label{Laplace2}
 \Exp[e^{- \frac{2u^2}{\beta} \zeta_1}] = \frac{\cosh\big(\tfrac{\lambda\beta}{8}\omega\big)}{\cosh\big(\omega \sqrt{u^2+\lambda^2\beta^2/64}\big)}  , \qquad u\ge 0 . 
\end{equation}
 
If $\wp : = \P[ x(\zeta_1) = 2\omega ]$, by \eqref{Laplace1}, this shows that 
\[
  \Exp e^{-u \zeta} \le  \sum_{k\in\N}  (1-\wp) \wp^{k-1} \big(  \Exp e^{-u \zeta_1}  \big)^k 
  = \frac{ (1-\wp)\E[e^{- u \zeta_1}]}{1-\wp \E[e^{- u \zeta_1}]}  . 
\]
We claim that $\wp= \frac{\exp\big(\tfrac{\lambda\beta}{8}\omega\big)}{2\cosh\big(\tfrac{\lambda\beta}{8}\omega\big)} $, so that we verify that if
$u \ge (\beta\lambda/4) \vee 2$, 
%$\sqrt{u^2+\lambda^2\beta^2/64} \ge \tfrac{\lambda\beta}{8}+\frac{\log 2}{\omega} $, 
 then
\[
 \frac{ (1-\wp) \Exp[e^{- \frac{2u^2}{\beta} \zeta_1}] }{1-\wp  \Exp[e^{- \frac{2u^2}{\beta} \zeta_1}] }  
 =  \frac{\exp\big(-\tfrac{\lambda\beta}{8}\omega\big)}{2\cosh\big(\omega \sqrt{u^2+\lambda^2\beta^2/64}\big) - \exp\big(\tfrac{\lambda\beta}{8}\omega\big)}
  \le  e^{- \frac{\omega}{2}\big(u+ \frac{\lambda\beta}{4} \big)} . 
\]
Note that up to factor $1/2$, the previous bound is sharp for large $u$. 
Hence, we conclude that for all $u \ge (\beta\lambda/4) \vee 2$, 
\[
 \Exp\big[e^{- \frac{2u^2}{\beta} \zeta}\big]   \le   e^{- \frac{\omega}{2}\big(u+ \frac{\lambda\beta}{4} \big)} ,
 \]
which is the claimed bound after rescaling $u\ge 0$.

To complete the proof, it remains to verify our formula for the probability $\wp$ and  \eqref{Laplace2}. 
This is a classical argument based on Girsanov's theorem. 
Recall that $\E[B(t)^2] =  t$ for $t\ge 0$ and let us make the change of measure given by
\[
\frac{\d\Q}{\d\P}\bigg|_t =  
\exp\big( - \tfrac{\lambda\sqrt{\beta}}{4} B(t)-  \tfrac{\lambda^2\beta}{32}t  \big) =
\exp\big(-  \tfrac{\lambda\beta}{8} (x(t)-\omega) + \tfrac{\lambda^2\beta}{32}t   \big) ,
\]
where the second identity holds for $t\le \zeta_1$. 
By Girsanov's theorem, under $\Q$, the process $\big(x(t)\big)_{t\le \zeta_1}$ is $\omega+\tfrac{2}{\sqrt\beta}B^{\Q}(t)$, where $B^{\Q}$ is a standard Brownian motion under $\Q$. If  $\mathrm{f}_\Q$ denotes the probability density of the stopping time $\zeta_1$ (under $\Q$), then it holds for any $u\ge 0$, 
\begin{equation} \label{Laplace3}
\Q[e^{- \frac{2u^2}{\beta} \zeta_1}] = \int_0^{+\infty} e^{-\frac{2u^2}{\beta} t}\, \mathrm{f}_\Q(t)\d t = \frac{1}{\cosh(\omega u)} . 
\end{equation}
By definitions, observe that
\[
\Q\big[ \{ x(\zeta_1) =0 , \zeta_1=t\} \big] = \P\bigg[\{\zeta_1 = t  ,  x(t) = 0 \}  \frac{\d\Q}{\d\P}\bigg|_t \bigg]  =  \P\big[\{ x(\zeta_1) = 0, \zeta_1 = t  \} \big] \exp\big(\tfrac{\lambda\beta}{8}\omega + \tfrac{\lambda^2\beta}{32}t   \big) . 
\]
By symmetry $\Q\big[ \{ x(\zeta_1) =0 , \zeta_1=t\} \big] =\mathrm{f}_\Q(t)/2$, so that for any $t\ge 0$,
\[
 \P\big[\{\zeta_1 = t  ,  x(\zeta_1) = 0 \} \big]  =  \frac{\mathrm{f}_\Q(t)}{2} \exp\big(-\tfrac{\lambda\beta}{8}\omega - \tfrac{\lambda^2\beta}{32}t   \big) . 
\]
Similarly, we also verify that
$ \P\big[\{\zeta_1 = t  ,  x(\zeta_1) = 2\omega \} \big]  =  \frac{\mathrm{f}_\Q(t)}{2} \exp\big(\tfrac{\lambda\beta}{8}\omega - \tfrac{\lambda^2\beta}{32}t   \big)$. 
These formulae show that under $\P$ , the event $\{ x(\zeta_1) = 2\omega \}$ is independent of the stopping time $\zeta_1$ with probability $\wp= \frac{\exp\big(\tfrac{\lambda\beta}{8}\omega\big)}{2\cosh\big(\tfrac{\lambda\beta}{8}\omega\big)} $ as claimed and that $\zeta_1$ has probability density function
\[
\mathrm{f}_\P(t) = \mathrm{f}_\Q(t)  \exp\big(- \tfrac{\lambda^2\beta}{32}t   \big)  \cosh\big(\tfrac{\lambda\beta}{8}\omega\big) . 
\]

Using the explicit formula \eqref{Laplace3}, this allows to compute the Laplace transform $ \Exp e^{-u \zeta_1} $, hence verifying \eqref{Laplace2} and completing the proof. 
\end{proof}

Recall that $\{ \mathfrak{z}_k\}$ denotes the blow-down times of the Riccati diffusion $\rho$ such that $0< \mathfrak{z}_1 < \mathfrak{z}_2$ etc. and $ \mathbf{N}_t  =  \#\big\{ k :   \mathfrak{z}_k <t\big\}$ is the associated counting function.  
Proposition \ref{prop:transitlaplace} provides us with a lower-bound for the blow-down time of the Riccati diffusion started from $+\infty$ which allows us to deduce overwhelming probability upper-bounds for the process $\mathbf{N}_{t}$.  
Namely,  by the (strong) Markov property, it holds for any $\lambda \ge 1$, $t \ge 0$ and $k\in \N$, 
\begin{equation} \label{eq:genericnbound}
 \P\big[  \mathbf{N}_t   \ge k + 1\big]   \leq \P[Z_1+ \cdots + Z_k \le t ],
\end{equation}
where $\left\{ Z_j \right\}$ are i.i.d. random variables which satisfy for   $u \ge \frac{\beta \lambda^2}{8} \vee \frac{8}{\beta}$,
\begin{equation} \label{Laplace4}
  \Exp\big[ e^{-u Z_j} \big]
  \leq  e^{-\frac{\sqrt{\beta\lambda}}{8}\big( \sqrt{u}+ \frac{\lambda\sqrt{\beta}}{2\sqrt{2}} \big)}  .
\end{equation}
Here we used the bound from Proposition~\ref{prop:transitlaplace} with $\omega=\sqrt{\lambda/8}$.

\begin{proposition}\label{prop:zerobounds}
  Fix $\nu>0$, $\lambda\ge1$, and $\epsilon>0$.  There exist constants $C_{\beta,\epsilon},c>0$ such that, for every $R\ge1$ and all sufficiently large $T$ (depending on $\beta,\epsilon,\nu$), the following holds with probability at least $1-C_{\beta,\epsilon}e^{-c\beta R T^\nu}$:
  \[
\mathbf{N}_{+\infty} \le  (\log T)^{1+\epsilon} \sqrt{R\, T^\nu} 
    \qquad  \text{and}    \qquad
\sum_{k}  \frac{1}{\sqrt{\mathfrak{z}_k}} \le  \sqrt{ R\, T^\nu} 
  \]
\end{proposition}

\begin{proof}
Combining \eqref{eq:genericnbound} and \eqref{Laplace4} with Markov's inequality, we obtain for any $t\ge 0$, $k\in\N$ and $u \ge 0$ sufficiently large,
\[ 
 \P\big[  \mathbf{N}_t  \ge k + 1\big]   \leq  e^{u t - k \frac{\sqrt{\beta\lambda}}{8}\big( \sqrt{u}+ \frac{\lambda\sqrt{\beta}}{2\sqrt{2}} \big)} . 
\]
To minimize the RHS  amounts to choosing 
$ u = \frac{\lambda\beta}{128} \frac{k^2}{t^2}$. 
Hence, if $k\ge  2\sqrt{2} t  (\sqrt{\lambda} \vee \frac{8}{\beta^2})$, this implies that 
\[ 
 \P\big[  \mathbf{N}_t  \ge k + 1\big]   \leq 
 \begin{cases}
   e^{- c\beta \lambda^{3/2} - c\beta k^2 \lambda/t}  & \text{if } k\ge  2\sqrt{2} t  (\sqrt{\lambda} \vee \frac{8}{\beta^2}) \\
  e^{- c\beta \lambda^{3/2}}    & \text{else} 
 \end{cases} , 
\]
 for a small numerical constant $c>0$. 
For $R  \ge 1$, if $T$ is sufficiently large (depending only on $\beta$), we immediately obtain for any $\lambda \ge 1$
\[ 
 \P\big[  \mathbf{N}_{1}  \ge    \sqrt{R\, T^\nu} \big]   \leq  e^{-c \beta R\, T^\nu}
\]
and by the Markov property,
\[ 
 \P\big[  \mathbf{N}_{+\infty}- \mathbf{N}_{(R\, T^\nu)^{2/3}}  \ge 2\big]  
 \le \P\big[  \mathbf{N}_{+\infty}^{(R\, T^\nu)^{2/3}}   \ge 2\big]  
  \leq  e^{-c \beta R\, T^\nu}   .
\] 
By a similar argument, it holds for any $j\in\N_0$, 
 \[ 
 \P\big[  \mathbf{N}_{2^{j+1}} -  \mathbf{N}_{2^j}  \ge \sqrt{R T^\nu}\big] 
\le   \P\big[  \mathbf{N}_{2^j}^{(2^j)}   \ge  \sqrt{R\, T^\nu}\big]
  \leq  e^{-c \beta R\, T^\nu} e^{- c\beta 2^{3j/2} }  .
\] 
 By a union bound, this shows that for any $\lambda\ge 1$  and $\epsilon>0$, if $T$ is sufficiently large (depending only on $\beta, \epsilon$),
 \[
  \P\big[  \mathbf{N}_{+\infty}\ge (\log T)^{1+\epsilon}\sqrt{R T^\nu} \big]
  \le C e^{-c \beta R\, T^\nu}  . 
 \]
Similarly, this also implies that 
  \[
  \P\bigg[  \sum_{k}  \frac{1}{\sqrt{\mathfrak{z}_k}} \ge  \sqrt{ R\, T^\nu} \bigg] \le 
   C e^{-c \beta R\, T^\nu} ,
  \]
 the main contribution coming from the first zeros.
\end{proof}

%\begin{remark} Work of Dumaz and Virag shows that there are in fact no zeros of this diffusion past $C(\log N)^{2/3}$ with overwhelming probability, which we do not capture.
%\end{remark}

\subsection{Control of the nonsingular part: Proof of Proposition~\ref{prop:nsc}} \label{sec:nsc}

Let $h(t) = \rho(t)- \sqrt{t+\lambda}$ for $\lambda\ge 0$ and 
set $\overline{\tau_{k}}= \tau_{k}\wedge T$,  for a given $T>0$.
By definition, $\big(h(t)\big)_{t\ge0}$ is a continuous diffusion which satisfies the SDE:
\[
\d h = -\bigg(h(t)^2+ 2h(t) \sqrt{t+\lambda} - \frac{1/2}{\sqrt{t+\lambda}}\bigg)\d t + \tfrac{2}{\sqrt{\beta}}\,\d B 
\qquad \text{with}\quad h(0) = \omega-\sqrt{\lambda}. 
\]

Our goal is to show that in the long run, with overwhelming probability, the process $h(t)$ does not become too large. 
The idea is to compare $h$ to a simpler diffusion $x$ which is driven by the same Brownian motion and satisfies the SDE:
\begin{equation} \label{xSDE}
\d x = - 2x(t) \sqrt{t+\lambda}  \d t + \tfrac{2}{\sqrt{\beta}}\,\d B +\frac{\d t}{2\sqrt{t+\lambda}}
\qquad \text{with}\quad x(0) = T^{(\nu+1)/4} . 
\end{equation}
In particular, if $\omega \le T^{(\nu+1)/4}$, we have $h(t) \le x(t)$ for all $t \le \mathfrak{z}_1$, that is, up to the first blow-down of the Riccati diffusion. The previous SDE is explicitly solvable using the integrating factor $y_\lambda$ defined in \eqref{def:y}, and its solution is given by
\[
x(t) =  T^{(\nu+1)/4} y_\lambda(0) y_\lambda^{-1}(t)   + y_\lambda^{-1}(t) \Bigg( \int_{0}^t \frac{y_\lambda(u)}{{2}\sqrt{u+\lambda}} \d u
+ \tfrac{2}{\sqrt{\beta}}\int_{0}^t y_\lambda(u)  \d B(u) \Bigg)
\]
Hence,  we obtain the bound 
\[
\int_0^{\overline{\tau_1}}h(t) \d t  \le \int_0^{\overline{\tau_1}}x(t) \d t  \le T^{(\nu+1)/4} y_\lambda(0) \int_0^\infty  y_\lambda^{-1}(t)  \d t
+ \int_0^T \int_{0}^t \frac{y_\lambda^{-1}(t)y_\lambda(u)}{2\sqrt{u+\lambda}} \d u \d t
+  \tfrac{2}{\sqrt{\beta}}\Xi_0   
\]
where $\displaystyle  \Xi_0 = \int_0^{\overline{\tau_1}} y_\lambda^{-1}(t) \int_{0}^t y_\lambda(u)  \d B(u)\d t$ is a  random variable. 
Let us observe that by the (strong) Markov property, we also have $h(\sigma_k+t) \le x_k(t)$ for $t\in [0, \mathfrak{z}_{k+1}-\sigma_k)$  where $x_k$ satisfies the SDE \eqref{xSDE} driven by the Brownian motion $\big(B_k(t) = B(t+\sigma_k)- B(\sigma_k) \Big)_{t\ge 0}$ for $k\ge0$, with $B_0=B$, and with drift $\lambda+\sigma_k$ instead of $\lambda$. 
Hence, this shows that for any $t\in [0,T]$, 
\begin{equation} \label{xest1}
\sum_{k=0}^{\mathbf{N}_t} \int_{\sigma_k}^{\overline{\tau_{k+1}}} h(u)\d u 
\le T^{(\nu+1)/4} \sum_{k=0}^{\mathbf{N}_t} y_{\lambda+\sigma_k}(0)\int_0^\infty  y_{\lambda+\sigma_k}^{-1}(s)  \d s
+ \frac12\sum_{k=0}^{\mathbf N_t}\int_0^T\int_0^s
\frac{y_{\lambda+\sigma_k}^{-1}(s)y_{\lambda+\sigma_k}(u)}
{\sqrt{u+\lambda+\sigma_k}}\,\d u\,\d s
+ \tfrac{2}{\sqrt{\beta}}\sum_{k=0}^{\mathbf{N}_t} \Xi_k
\end{equation}
where
\[
\Xi_k
= \int_0^{(\overline{\tau_{k+1}}-\sigma_k)_+} \hspace{-1em}
  y_{\lambda+\sigma_k}^{-1}(t) \int_{0}^t
  y_{\lambda+\sigma_k}(u)\,\d B_k(u)\,\d t.
\]
Conditionally on $\mathcal F_{\sigma_k}^B$, the integrand in this definition is the Gaussian process from \eqref{def:mathfrakh} with parameter $\lambda+\sigma_k$:
\[
\left(y_{\lambda+\sigma_k}^{-1}(t)\int_0^t
  y_{\lambda+\sigma_k}(u)\,\d B_k(u)\right)_{t\ge0}
\overset{\mathrm{law}\mid\mathcal F_{\sigma_k}^B}{=}
\left(\mathfrak X_{\lambda+\sigma_k,0}(t)\right)_{t\ge0}
.
\]

For the first two terms on the RHS of \eqref{xest1}, we easily verify that there exists a numerical constant $C>0$ such that for any $\lambda\ge0$, 
\begin{equation} \label{xest2}
y_\lambda(0)\int_0^\infty  y_\lambda^{-1}(t)  \d t \le \frac{C}{1+\sqrt{\lambda}}
\qquad\text{and}\qquad  
 \int_0^T \int_{0}^t \frac{y_\lambda^{-1}(t)y_\lambda(u)}{\sqrt{u+\lambda}} \d u\d t \le C\big(1+ \log T \big) . 
\end{equation}

To control the contributions from the random variables $\Xi_k$, we rely on the following lemma which is proved at the end of this section. 

\begin{lemma} \label{lem:X}
Let $\lambda\ge 0$. We can rewrite for any $\tau\ge 0$, 
\[
\int_0^{\tau} 
\mathfrak{X}_{\lambda,0}(t)
\d t
%y_\lambda^{-1}(t) \int_{0}^t y_\lambda(u)  \d B(u)
= \mathfrak{M}_\lambda(\tau)  + \mathfrak{Y}_\lambda(\tau) 
\qquad\text{where}\qquad
\mathfrak{M}_\lambda(\tau) =  \1_{\tau>1} \int_{1}^\tau  \frac{\d B(u)}{2\sqrt{\lambda+u}}
\]
 and $\mathfrak{Y}_\lambda$ is a continuous (adapted) Gaussian process; see formula \eqref{Y0}. 
Moreover, there exists a constant $C_\beta>0$ such that  
\begin{equation} \label{Y1}
\big\VERT  {\textstyle \sup_{\tau\in\R_+}} | \mathfrak{Y}_\lambda(\tau)| \big\VERT_2^2 \le  C_\beta . 
%\P\bigg[ \sup_{\tau\in\R_+}  \big| \mathfrak{N}_\lambda(\tau)   \big| \ge \sqrt{R\, T^\nu} \bigg] \le  CN^{- cR /\beta }  . 
\end{equation}
\end{lemma}

Notice that the martingale  $\mathfrak{M}_\lambda$ is a continuous Gaussian log-correlated process in the sense that its quadratic variation is $\tau \mapsto \tfrac14\log\big(\frac{\lambda+\tau}{\lambda+1}\big)$. Hence, it holds for any $\lambda\ge 0$ 
and $T\ge 2$, 
\[
\big\VERT  {\textstyle \sup_{\tau \le T}}  \big| \mathfrak{M}_\lambda(\tau) \big|  \big\VERT_2^2 \lesssim \log T. 
\]

Taking the supremum over the possible endpoint $\overline\tau_{k+1}-\sigma_k$, this stopping time being less than $T$, the preceding estimate and Lemma~\ref{lem:X} imply, that uniformly over for all $k$ and all $\theta\in\R$,
\begin{equation}\label{eq:XiConditionalMGF}
\E\left[e^{\theta\Xi_k}\mid\mathcal F_{\sigma_k}^B\right]
\le \exp\left(C_\beta\left(|\theta|\sqrt{\log T}+\theta^2\log T\right)\right).
\end{equation}

For a fixed integer $m\ge0$, set $S_m=\sum_{k=0}^m\Xi_k$.  Since $S_{m-1}$ is $\mathcal F_{\sigma_m}^B$-measurable, iterating conditional expectation in \eqref{eq:XiConditionalMGF} gives
\[
\E e^{\theta S_m}
\le \exp\left(C_\beta(m+1)
\left(|\theta|\sqrt{\log T}+\theta^2\log T\right)\right).
\]
Applying this estimate with both signs of $\theta$ and optimizing shows that, for every $x>0$,
\begin{equation}\label{eq:ThetaFixedSum}
\P\left[|S_m|\ge C_\beta(m+1)\sqrt{\log T}+x\right]
\le 2\exp\left(-\frac{c_\beta x^2}{(m+1)\log T}\right).
\end{equation}
Now put
\[
L_T=\left\lceil(\log T)^{1+\epsilon}\sqrt{RT^\nu}\right\rceil,
\qquad
\mathscr A_T=\{\mathbf N_T\le L_T\}.
\]
Proposition~\ref{prop:zerobounds} gives $\P(\mathscr A_T^c)\le C_{\beta,\epsilon}e^{-c\beta RT^\nu}$.  Take $\Lambda=C\sqrt R\,T^{(3\nu+1)/4}$.  Since
\[
(L_T+1)\sqrt{\log T}
\lesssim \sqrt{RT^\nu}(\log T)^{3/2+\epsilon}
=o\left(\sqrt R\,T^{(3\nu+1)/4}\right),
\]
we may increase $C$ and take $T$ sufficiently large so that the deterministic term in \eqref{eq:ThetaFixedSum} is at most $\Lambda/2$ for every $m\le L_T$.  A union bound over these fixed values of $m$ then gives
\[
\begin{aligned}
\P\Bigg[\sup_{t\le T}\bigg|\sum_{k=0}^{\mathbf N_t}\Xi_k\bigg|\ge\Lambda\Bigg]
&\le \P(\mathscr A_T^c)
+2(L_T+1)\exp\left(-\frac{c_\beta\Lambda^2}{(L_T+1)\log T}\right).
\end{aligned}
\]
For $R\le\sqrt T$,
\[
\frac{\Lambda^2}{(L_T+1)\log T}
\gtrsim \frac{C^2\sqrt R\,T^{\nu+1/2}}{(\log T)^{2+\epsilon}}
\ge \frac{C^2T^{1/4}}{(\log T)^{2+\epsilon}}RT^\nu.
\]
The prefactor is therefore absorbed for all sufficiently large $T$, and we conclude that
\[
\P\Bigg[ \sup_{t\le T} \bigg| \sum_{k=0}^{\mathbf{N}_t}  \Xi_k \bigg|
\ge C\sqrt R\,T^{(3\nu+1)/4} \Bigg]
\le C_{\beta, \epsilon}e^{-c\beta R\, T^\nu} .
\]

Then, by \eqref{xest2} and Proposition~\ref{prop:zerobounds}, it also holds with probability at least $1-  C_\beta e^{-c\beta R\, T^\nu} . $
\[
\sup_{t\ge 0} \sum_{k=0}^{\mathbf{N}_t}
y_{\lambda+\sigma_k}(0)\int_0^\infty y_{\lambda+\sigma_k}^{-1}(s)\,\d s
\le C\bigg(1+\sum_{k\ge1}\frac1{\sqrt{\mathfrak z_k}}\bigg)
\le C\sqrt{RT^\nu}.
\]

After multiplication by the Riccati cutoff, this bounds the first term in \eqref{xest1} by $C\sqrt R\,T^{(3\nu+1)/4}$.  The second estimate in \eqref{xest2} and Proposition~\ref{prop:zerobounds} similarly give, on the same event,
\[
\begin{aligned}
&\frac12\sum_{k=0}^{\mathbf N_t}\int_0^T\int_0^s
\frac{y_{\lambda+\sigma_k}^{-1}(s)y_{\lambda+\sigma_k}(u)}
{\sqrt{u+\lambda+\sigma_k}}\,\d u\,\d s\\
&\hspace{2cm}\le C(\mathbf N_t+1)(1+\log T)
\le C_{\epsilon}\sqrt{RT^\nu}(\log T)^{2+\epsilon}
\le C_{\epsilon}\sqrt R\,T^{(3\nu+1)/4}.
\end{aligned}
\]

Hence from the estimate \eqref{xest1}, we conclude that it holds for any fixed $\epsilon>0$ and for $R\in[1,\sqrt T]$,
\begin{equation*} \label{xest3}
\P\Bigg[   \sup_{t\le T} \bigg( \sum_{k=0}^{\mathbf{N}_t} \int_{\sigma_k}^{\overline{\tau_{k+1}}} h(u)\d u  \bigg) \ge C\sqrt R\,T^{(3\nu+1)/4}   \Bigg]
\le C_{\beta, \epsilon}e^{-c\beta R\, T^\nu} .
\end{equation*}

Rescaling $R$ by a fixed constant gives the stated form of Proposition~\ref{prop:nsc}. This ends our proof of Proposition~\ref{prop:nsc}. To complete the argument, it remains to prove Lemma~\ref{lem:X}.

\begin{proof}[Proof of Lemma~\ref{lem:X}] 
First, an integration by parts shows that 
\begin{equation} \label{X1}
  \int_0^{\tau} \tfrac{1}{y_\lambda(t)} \int_{0}^t y_\lambda(u)  \d B(u) \d t
 = -Y_\lambda(\tau)  \int_{0}^\tau y_\lambda(u)\d B(u) +
 \int_{0}^\tau Y_\lambda(u)y_\lambda(u)\d B(u) , 
\end{equation}
where $\displaystyle Y_\lambda(\tau)  = \int_\tau^{\infty} y_\lambda^{-1}(t) \d t$.
Recall  $\displaystyle G_\lambda(\tau)  =\int_0^\tau y_\lambda^{2}(u) \d u.$  By the martingale representation theorem, the first term on the RHS of \eqref{X1} has the same law as the process $\tau \mapsto Y_\lambda(\tau) B_{G_\lambda(\tau)} $.
In particular, we have the following sub-Gaussian norm estimate: for any $T\ge 0$,
\[
\bigg\VERT \sup_{\tau \ge T} \Big|
Y_\lambda(\tau)  \int_{0}^\tau y_\lambda(u)\d B(u)
\Big| \bigg\VERT_2 
=
\bigg\VERT \sup_{t \ge G_\lambda(T)} \Big|
Y_\lambda\big(G_\lambda^{-1}(t)\big)  B(t)
\Big| \bigg\VERT_2 .
\]

To control this norm, we rely on the following facts:
$1)$ There exists an increasing bijection $\zeta : [1,\infty) \to [0,\infty)$ such that 
$\int_0^{\zeta(t)}  \sqrt{v+\lambda}  \d v = \frac14 \log(t)$.
$2)$ We check that $\zeta(t) \ge \big(\tfrac 38\log t \big)^{2/3}$ for all $t\ge 1$. 
$3)$ $G_\lambda$ is an increasing function with
$G_\lambda(\tau) \le  y_\lambda^{2}(\tau)$ (see \ref{eq:Gupper}) so that by \eqref{def:y}, its inverse function satisfies $G_\lambda^{(-1)}(t) \ge \zeta(t) $ for $t\ge1$ (c.f.\ \ref{eq:Gupper}). 
$4)$ $Y_\lambda$ is a decreasing function with 
$Y_\lambda(\tau) \le \frac{1}{2y_\lambda(\tau)\sqrt{\tau+\lambda}} $, 
so that it holds for any $t\ge 1$, 
\[
Y_\lambda\big(G_\lambda^{-1}(t)\big) \le  \frac{1}{2y_\lambda(\zeta(t))\sqrt{\zeta(t)+\lambda}} \le \frac{1}{\sqrt{t}(\log t)^{1/3}} . 
\]
This bound implies that
\[
\bigg\VERT \sup_{\tau >0} \Big(
Y_\lambda(\tau)  \int_{0}^\tau y_\lambda(u)\d B(u)
\Big) \bigg\VERT_2 
\le Y_\lambda(0) \bigg\VERT \max_{t \in[0,2]} \big|B(t)
\big| \bigg\VERT_2 +
\bigg\VERT \sup_{t \ge 2} \frac{|B(t)|}{\sqrt{t} (\log t)^{1/3}}
\bigg\VERT_2 ,
\]
so that  by Lemma~\ref{lem:BMestimate}, we obtain the sub-Gaussian norm estimate:
\begin{equation} \label{X2}
\bigg\VERT \sup_{\tau\in\R_+}  \bigg| Y_\lambda(\tau)  \int_{0}^\tau y_\lambda(u)\d B(u)  \bigg| \bigg\VERT \lesssim 1. 
\end{equation}
where the implied constant is independent of $\lambda>0$ (and may depend on $\beta$). 
We also verify that 
\[
\lim_{T\to\infty} \bigg\VERT \sup_{\tau \ge T} \Big|
Y_\lambda(\tau)  \int_{0}^\tau y_\lambda(u)\d B(u)
\Big| \bigg\VERT_2  =0
\]
which implies that $\displaystyle  \lim_{\tau \to\infty} \bigg( Y_\lambda(\tau)  \int_{0}^\tau y_\lambda(u)\d B(u) \bigg) = 0$ almost surely.

%for some standard Brownian motion $\widehat{B}$ with variance $\E[\widehat{B}(t)^2] =  t$ for $t\ge 0$. 
The function $\tau \in \R_+ \mapsto Y_\lambda^2(\tau)G_\lambda(\tau)$ is smooth, non-negative and uniformly bounded (independently of $\lambda\ge 0$). Indeed, we verify that we have the following asymptotics as $\tau\to\infty$,
\[
Y_\lambda^2(\tau)G_\lambda(\tau) = \frac{1/16}{(\tau+\lambda)^{3/2}} \big(1+\O(\tau^{-1}) \big) . 
\]

For the second term on the RHS of \eqref{X1} which is a martingale, let us observe that
\[
Y_\lambda(t)y_\lambda(t) = \frac{1}{2\sqrt{\lambda+t}} - H_\lambda(t)
\qquad
 H_\lambda(t) : = \frac{1}{4}  \int_t^\infty e^{-\frac43(u^{3/2}-t^{3/2})} \frac{\d u}{(u+\lambda)^{3/2}}  . 
\]
This decomposition shows that indeed we have for any $\tau\ge 0$, 
\[
 \int_0^{\tau} y_\lambda^{-1}(t) \int_{0}^t y_\lambda(u)  \d B(u) \d t
 =  \1_{\tau>1}  \int_{1}^\tau  \frac{\d B(u) }{2\sqrt{\lambda+u}} + \mathfrak{Y}_\lambda(\tau) 
\]
where
\begin{equation} \label{Y0}
\mathfrak{Y}_\lambda(\tau) \coloneqq Y_\lambda(\tau)  \int_{0}^\tau y_\lambda(u)\d B(u) + \int_{0}^{\tau\wedge1} \hspace{-.3cm}Y_\lambda(u)y_\lambda(u)\d B(u)   
+   \1_{\tau>1} \bigg( \int_{1}^\tau  H_\lambda(u) \d B(u) \bigg) 
\end{equation}
is a continuous (adapted) Gaussian process. 
Finally, observe that the function  $H_\lambda \in L^2([1,\infty))$ since it has the following asymptotics
$H_\lambda(\tau) =  \frac{1/8}{(\tau+\lambda)^{2}}\big(1+\O(\tau^{-1}) \big)$ as $\tau\to+\infty$. This shows that the process $\mathfrak{Y}_\lambda$ converges almost surely and 
\begin{equation} \label{limitY}
 \mathfrak{Y}_\lambda(\infty)  \coloneqq \lim_{\tau\to\infty} \mathfrak{Y}_\lambda(\tau) =
\int_0^1 Y_\lambda(u)y_\lambda(u)\d B(u)    +  \int_{1}^\infty  H_\lambda(u) \d B(u) . 
\end{equation}
Moreover, for any $\lambda\ge 0$, we have the following sub-Gaussian norm estimate :
\[
\bigg\VERT \sup_{\tau\ge1}  \bigg| \int_{1}^\tau  H_\lambda(u) \d B(u) \bigg| \bigg\VERT \lesssim 1. 
\]
%where $c = \| H_0 \|_{L^2([1,\infty)}^{-2}$. 
Combining the previous estimate with \eqref{X2} and an analogous bound for the second martingale on the RHS of \eqref{Y0}, we obtain \eqref{Y1}.
\end{proof}

\subsection{Control of the singular part: Proof of Proposition~\ref{prop:sc}} \label{sec:sc}

The proof requires two properties that we need to show hold with overwhelming probability: the intervals $[\tau_k,\sigma_k]$ are small and that both the Dirichlet and Neumann solutions are well-behaved for short times. 
The following proposition summarizes these estimates.

\begin{proposition} \label{prop:singular}
Fix $\nu\in(1,2)$, $\lambda\in[-T,T]$.
For all $T$ sufficiently large and all $k\ge1$, let $\mathcal G_k$ denote the event
\[
\sigma_k \le \tau_k +5T^{-(\nu+1)/4} 
\qquad\text{and}\qquad
\sup_{0\le s\le 5T^{-(\nu+1)/4}}|\v_k(\tau_k+s)|\le C_\nu
\qquad\text{and}\qquad
\sup_{0<s\le 5T^{-(\nu+1)/4}}\frac{|\u_k(\tau_k+s)|}{s}\le C_\nu  . 
\]
There are constants $c_\beta,C_\beta>0$ such that
\[
\P\Big[\{\tau_k\le T\}\cap\mathcal G_k^c\Big]
\le C_\beta e^{-c_\beta T^{(\nu+2)/2}}.
\]
\end{proposition}

Using  Proposition~\ref{prop:zerobounds} and Proposition~\ref{prop:singular}, it is straightforward to control the singular contribution.

\begin{proof}[Proof of Proposition~\ref{prop:sc}]
Recall that
\[
\mathrm{II}_T
=\prod_{k=1}^{\mathbf{N}_{T}}
\left|\frac{\phi_\lambda(\sigma_k)}{\phi_\lambda(\tau_k)}\right|.
\]
For $k\ge1$, set $s_k=\sigma_k-\tau_k$.  Then
\[
\frac{\phi_\lambda(\sigma_k)}{\phi_\lambda(\tau_k)}
=\v_k(\tau_k+s_k)-T^{(\nu+1)/4}\u_k(\tau_k+s_k).
\]
Consequently, Proposition~\ref{prop:singular} implies that, on the event
$\{\mathbf{N}_{T}\le T^{3\nu/4}\}$, with overwhelming probability,
\[
\left|\frac{\phi_\lambda(\sigma_k)}{\phi_\lambda(\tau_k)}\right|\le C_\nu,
\qquad 1\le k\le\mathbf N_T.
\]
Moreover, by Proposition~\ref{prop:zerobounds} with $R = T^{(2-\nu)/2} $, we obtain
\[
\P\big[ \mathbf{N}_{T}\ge T^{3\nu/4} \Big] \le C_\beta e^{- c_\beta  T^{(\nu+2)/2}}
\]
Hence, by a union bound (and adjusting the constants), we conclude that
\[
\P\big[  \mathrm{II}_T  \le C_\nu^{T^{3\nu/4}}\big] \ge 1- C_\beta e^{- c_\beta  T^{(\nu+2)/2}} . \qedhere 
\]
\end{proof}

It remains to prove Proposition~\ref{prop:singular}.
Our analysis is based on the study of the  \emph{inverse Riccati diffusion}  $x = 1/\rho$ that we already encountered in Section~\ref{sec:entrance}. 
First, we show that when started from a small $x_0$, the process $x(t)$ remains approximately linear for a short amount of time with overwhelming probability. These estimates will be instrumental for the rest of the proof. Then, we provide short-time uniform bounds for both Dirichlet and Neumann solutions. Finally, we complete the proof of Proposition~\ref{prop:sc} in the last section. 

\paragraph{Inverse Riccati diffusion: short time estimates.}
We define the \emph{inverse Riccati diffusion}  $x = 1/\rho.$  
Applying It\^o formula to \eqref{Ric0}, we verify that $x$ solves the following SDE:
\begin{equation}\label{eq:inversericatti}
  \d x = (1- (t+\lambda)  x^2 + \tfrac{4}{\beta}x^3) \d t - \tfrac{2}{\sqrt{\beta}}x^2 \d B 
\qquad \text{with}\quad x(0) = x_0 = 1/\omega.
\end{equation}
We are interested in the case where $x_0=0$ (Neumann solution) or $x_0$ is small.
This equation has a strong solution until its first blow-up (which corresponds to the Riccati diffusion hitting 0 if it occurs). 

Note that by definitions \eqref{stoppingtimes}, the (strong) Markov property implies that for any $k\ge 1$,  on the event $\{\tau_k<\infty\}$,  the process  $x_k(t) = 1/\rho(\tau_k+t)$ solves \eqref{eq:inversericatti} with $x_0 =-T^{-(\nu+1)/4}$
and 
\begin{equation} \label{stoppingtimes2}
(\sigma_k-\tau_k) = \inf\big\{ t > 0 : x_k(t) = T^{-(\nu+1)/4} \big\}   
\end{equation}
Hence by showing that if $x_0$ is small, the inverse Riccati diffusion $x(t)$ remains in a linear tube for a short window of time, we can get an estimate for \eqref{stoppingtimes2}.

\begin{lemma}\label{lem:blowdown}
  Let $x$ solve \eqref{eq:inversericatti} with $x_0 \le 0$  and set $\vartheta = \inf \big\{ t \in [0,1] : |x(t) - x_0 - t| \geq \frac{ t^{1+\alpha}-x_0}{2} \big\}$ for $\alpha\in[0,1]$.
  There is a small constant $c_\beta$ so that
  \begin{itemize}
  \item If $x_0 = 0$, it holds for all $\delta \leq \frac{c_\beta}{\sqrt{1+|\lambda|}}$, 
$
    \P[ \vartheta \leq \delta]
    \leq
    2e^{-c_\beta\delta^{-3+2\alpha}}.$
  \item If $\alpha=0$ and $|x_0| \le \frac{c_\beta}{\sqrt{1+|\lambda|}}$, then
    $  \P\big[ \vartheta < 5|x_0| \big]
    \leq
    2e^{-c_\beta|x_0|^{-3}}.$
  \end{itemize}
\end{lemma}
\begin{proof}
  Integrating \eqref{eq:inversericatti}, we have that for $t \leq \vartheta,$
  \[
  t + x_0-  x(t)  = \int_0^t \biggl\{(u + \lambda)x^2(u) - \tfrac{4}{\beta}x^3(u)\,\d u\biggr\} + M(t),
  \]
  where $M(t) = \tfrac{2}{\sqrt{\beta}}\int_0^t x^2(u)\,\d B(u)$ is a martingale. By definition of $\vartheta,$ we verify that we can bound for $t\le \vartheta$,  
  \[
  | x(t) -x_0 -t | \leq \tfrac34|\lambda| (t-x_0)^3 + C_\beta (t-x_0)^4  + |M(t)| .
  \]
  for some constant $C_\beta.$  
If $x_0 = 0$, choosing $\delta>0$ sufficiently small so that $3|\lambda| \delta +  4C_\beta \delta^2 \le 1$, this implies that 
\[
    \P[ \vartheta \leq \delta] \le  \P\big[  |M(t \wedge \vartheta)|  \ge (t \wedge \vartheta)^{1+\alpha}/2 \text{ for a } t\le \delta \big]
  \]

The quadratic variation of the martingale $M(t)$ is bounded  above and below by $t^5$ up to a multiplicative   constants. 
Hence, if we represent $M$ as a time change of a Brownian motion $B,$  we can bound 
  \[
    \sup_{t \leq \delta} \frac{|M(t \wedge \vartheta)|}{(t \wedge \vartheta)^{1+\alpha}}
    \lesssim \sup_{t \leq \delta} \frac{|M(t \wedge \vartheta)|}{ \langle M(t \wedge \vartheta) \rangle^{(1+\alpha)/5}}
    \lesssim \sup_{t \leq C_\beta \delta^5} \frac{ |B(t)| }{ t^{(1+\alpha)/5}} 
  \]
for a constant $C_\beta>0$.  
By  Lemma \ref{lem:BMestimate} with $f(t) = t^{3/10- \alpha/5}$, 
we obtain the  sub-Gaussian norm estimate:
  \(
  \big\VERT \sup_{t \leq \delta} \tfrac{|M(t \wedge \vartheta)|}{(t \wedge \vartheta)^{1+\alpha}} \big\VERT \lesssim_\beta \delta^{3/2-\alpha}.
  \)
This proves the first claim. 
We can use the same argument to obtain the second bound.
 If $27 |\lambda| x_0^2 + 4^3C_\beta |x_0|^3 \le 1/4$, then we verify that
  \[
    \P[ \vartheta \leq 5 |x_0|] \le  \P\big[  |M(t \wedge \vartheta)|  \ge (t\wedge \vartheta-x_0)/2 \text{ for a } t\le 5|x_0|  \big] .
  \]
Taking $\delta=5|x_0|$ in the  previous sub-Gaussian norm estimate also yields the second claim. 
\end{proof}

\paragraph{Control of Dirichlet solutions.}
Our next Lemma deals with the local behavior of the Riccati diffusion started from $0$.  This allows us to get control of the growth of Dirichlet solution for short time. 

\begin{proposition}\label{lem:r-zero}
  Let $\rho$ solve \eqref{Ric0} with $\omega =0$ and $\lambda\in\R$. 
  Given $\eta \ge 1$, let $\vartheta = \inf \{ t>0 :   |\rho(t)| > \eta \}$.
  There is a constant $c_\beta>0$ so that $ \P[ \vartheta \leq \delta] \leq e^{-c \eta^2/\delta}$ for all $\delta \leq  \frac{\eta}{4|\lambda|} \wedge \frac{1}{2\eta}$. 
%Let $\v$ solve \eqref{SA1} with conditions $\v(0)=1$ and $\v'(0)=0$
This implies that  for such small $\delta$, 
  \[
  \P\Big[ \max_{t\in[0,\delta]}| \mathrm{f}_{\lambda,0}(t)| \ge  e^{\delta \eta} \Big]   \le e^{-c_\beta \eta^2/\delta} . 
  \]
\end{proposition}
\begin{proof}
  Integrating the equation \eqref{Ric0}, it holds for $t \in[0, \vartheta]$, 
  \begin{equation*}%\label{eq:rt-i}
    \rho(t) = \int_0^t \big( -\rho^2(u) + u + \lambda \big)\,\d u+ B(t) ,
  \end{equation*}
  since the Riccati diffusion has not yet blown down.
 Thus, we have for $t\le \vartheta \le \delta$,  
  \[
 | \rho(t) | \le  \big| t^2/2+ (\lambda-\eta^2) t \big| + |B(t)| \le 3\eta/4 + |B(t)| ,
  \]
  where we used that $\delta \leq \frac{\eta}{4|\lambda|} \wedge \frac{1}{2\eta}$
 and $\eta\ge 1$ to get the second bound.
  Hence we conclude
    \[
    \P[ \vartheta \leq \delta] \leq
    \P\big[ |B(t)| \ge \eta/4 \text{ for a } t\in[0,\delta] \big]
    \le e^{-c_\beta \eta^2/\delta}. 
  \]
Finally, by \eqref{phipv} applied to $\mathrm{f}_{\lambda,0}$, we obtain the  trivial bound on the event $\{\vartheta\ge \delta\}$, 
$| \mathrm{f}_{\lambda,0}(t)| \le  e^{\delta \eta}$ valid for all $t\le\delta$. 
This proves the second estimate.
\end{proof}

\paragraph{Control of Neumann solutions.}
In this case $\mathrm{g}_{\lambda,0}(0) = 0 $ and the corresponding Riccati diffusion $\rho$ starts from $\omega=\infty$. 
Let $\tau = \inf\{t > 0 :  \mathrm{g}_{\lambda,0}(t) = 0 \}$. 
Applying Lemma~\ref{lem:phipv},  it holds conditionally on the event $\{\epsilon<\tau\}$, for any $t\in [\epsilon,\tau)$, 
\[
\mathrm{g}_{\lambda,0}(t) =  \mathrm{g}_{\lambda,0}(\epsilon) \exp\biggl( \ \int_{\epsilon}^{t}\rho(u)\d u \biggr) .
\]
We can rewrite this formula in terms of the inverse Riccati diffusion $x = 1/\rho$,  then take the limit as $\epsilon\to 0$. 
By continuity of $\mathrm{g}_{\lambda,0}'$, it holds almost surely 
$\displaystyle  \lim_{\epsilon\to0}\mathrm{g}_{\lambda,0}(\epsilon)  \epsilon^{-1}  = \mathrm{g}_{\lambda,0}'(0) =1 $, so that for any $t\le \tau$, 
\begin{equation} \label{Neumannpv}
\mathrm{g}_{\lambda,0}(t) =  t\exp\bigg( \lim_{\epsilon\to 0} \bigg( \int_{\epsilon}^{t} \bigg(\frac{1}{x(v)} - \frac{1}{v} \bigg)\d v \bigg)  \bigg) . 
\end{equation}

By Lemma~\ref{lem:blowdown} with $x_0 =0$ and a small $\alpha>0$, we know that with overwhelming probability $|x(t)-t| \le t^{1+\alpha}/2$ for a short time. Hence for such short time, we expect that $|\mathrm{g}_{\lambda,0}(t) | \lesssim t$. This is the content of our next proposition. 

\begin{proposition} \label{lem:Neumann-growth}
For any $\alpha\in(0,1]$, there exists a constant $C_\alpha>0$ such that it holds for any $\delta \leq \frac{c_\beta}{\sqrt{1+|\lambda|}}$, 
\vspace{-.36cm}
\[
\P\Big[ \max_{t\in[0,\delta]}|\mathrm{g}_{\lambda,0}(t) t^{-1}|  \ge C_\alpha  \Big]      \le 2e^{-c_\beta\delta^{-3+2\alpha}}.
\]
\end{proposition} 

\begin{proof}
Let $\vartheta$ be as in Lemma~\ref{lem:blowdown} with $x_0=0$ --- we plainly have 
$\vartheta \le \tau$.
Observe that conditionally on the event $\{\vartheta \ge t\} $, we have for any $0<\epsilon<t \le 1$, 
\[
\bigg| \int_{\epsilon}^{t} \frac{\d v}{x(v)} - \frac{\d v}{v}  \bigg| \le 
\int_{\epsilon}^{t} \bigg|  \frac{x(v)-v}{x(v)v}   \bigg| \d v \lesssim   \int_\epsilon^1 \frac{\d v}{v^{1-\alpha}} \le c_\alpha
\]
By \eqref{Neumannpv}, we conclude that for any $\delta>0$, on the event $\{\vartheta \ge \delta\} $, 
$\sup_{t\in[0,\delta]}|\mathrm{g}_{\lambda,0}(t) t^{-1}|  \le C_\alpha $. Then, the tail bound follows from the first claim in  Lemma~\ref{lem:blowdown}.
\end{proof}

\paragraph{Proof of Proposition~\ref{prop:singular}.}
Fix $k\ge1$.  On $\{\tau_k\le T\}$, conditionally on
$\mathcal F^B_{\tau_k}$, the shifted inverse Riccati process has parameter
$\lambda+\tau_k$ and satisfies $|\lambda+\tau_k|\le2T$.  Since $\nu>1$,
the condition $T^{-(\nu+1)/4}\le c_\beta/\sqrt{1+2T}$ holds for all
sufficiently large $T$.
Let us first obtain control for the length of the interval $[\tau_k,\sigma_k]$.  
By \eqref{stoppingtimes2} and applying the second estimate from Lemma~\ref{lem:blowdown} with $x_0 = -T^{-(\nu+1)/4}$, we obtain
\[
\P\Big[\{\sigma_k-\tau_k>5T^{-(\nu+1)/4}\}\cap\{\tau_k\le T\}\Big]
\le 2e^{-c_\beta T^{3(\nu+1)/4}}.
\]

We can now get control of the singular contribution. 
Recall that for $k\ge 1$, $\v_k = \mathrm{f}_{\lambda,\tau_k}$ and $\u_k = \mathrm{g}_{\lambda,\tau_k}$
where $\tau_k$ is the stopping time \eqref{stoppingtimes}. Then, conditionally on $\mathcal F^B_{\tau_k}$, the strong Markov property gives
$(\v_k(\tau_k+s),\u_k(\tau_k+s))_{s\ge0} \overset{\rm law}{=} (\mathrm{f}_{\lambda+\tau_k,0}(s), \mathrm{g}_{\lambda+\tau_k,0}(s))_{s\ge0}$.
For the Dirichlet solutions, applying Proposition~\ref{lem:r-zero} with $\delta = 5T^{-(\nu+1)/4}$ and $\eta = T^{(\nu+1)/4}/10$, we obtain 
  \[
  \P\Big[ \sup_{0\le s\le 5T^{-(\nu+1)/4}}|\v_k(\tau_k+s)| \ge  e^{1/2} ,  \tau_k \le T   \Big]   \le e^{-c_\beta T^{3(\nu+1)/4}} .
  \]
Here, we used the assumption $T\ll T^{(\nu+1)/2}$, so that the condition $\delta \leq  \frac{\eta}{8T}$ is satisfied. 
For the Neumann solutions, applying Proposition~\ref{lem:Neumann-growth} with
$\delta = 5T^{-(\nu+1)/4}$ and $\alpha = \frac{\nu-1}{2(\nu+1)}$, we conclude that
\[
\P\Big[ \sup_{0<s\le 5T^{-(\nu+1)/4}}\frac{|\u_k(\tau_k+s)|}{s}  \ge C_\nu  ,  \tau_k \le T \Big]      \le 2e^{-c_\beta T^{(\nu+1)(3-2\alpha)/4}} = 2 e^{- c_\beta  T^{(\nu+2)/2}}.
\]
Combining these three estimates and using
$3(\nu+1)/4\ge(\nu+2)/2$ for $\nu\ge1$ gives
$\P\big[\{\tau_k\le T\}\cap\mathcal G_k^c\big]
\le C_\beta e^{-c_\beta T^{(\nu+2)/2}}$, as claimed.
\qed

% ===== end sections/08-expanding-direction.tex =====

% ===== begin sections/02-parabolic-region.tex =====
\section{Approximation results for the Characteristic Polynomial} \label{sec:quant}

The goal of this section is to go from the three-term recursion for the characteristic polynomials of the tridiagonal matrix model to the stochastic Airy equation. This approximation is valid in a neighborhood of any point $z\in[-1,1]$ in the spectrum, with $z\neq0$, in the parabolic part of the recursion as explained in Section~\ref{sec:para}, and it involves coupling the noise \eqref{def:XY} with a Brownian motion $\boldsymbol{B}^{(z)}$ as we now explain.

To describe the parabolic part of the recursion, we rely on the following conventions: we consider a sequence $z(N)\in [-1,1]$ such that 
\begin{equation}\label{def:hyper}
z(N)\to z \in [-1,1] \quad\text{as $N\to\infty$}, \qquad   N_0 =  \lfloor Nz(N)^2 \rfloor \to \infty \quad\text{as $N\to\infty$}. 
\end{equation}
We abuse notation and write $z = z(N)$.
Here $N_0$ is the \emph{turning point} and we define for $t\in\R$ and $\lambda\in\C$:
\begin{equation} \label{parascale}
z_\lambda = z\bigg(1 +  \frac{\lambda}{2 N_0^{2/3}} \bigg) , \qquad 
N_{t}(z) = \lfloor N z^2 - t \mathfrak{L}   \rfloor , \qquad \mathfrak{L} \coloneqq N_0^{1/3} ,
\end{equation}
where $ \mathfrak{L}$ is the parabolic scale.
Throughout this section, we work with the following conventions: 

\begin{assumption}\label{def:para}
We assume that $z$ satisfies \eqref{def:hyper}. 
As in Theorem~\ref{thm:hypersimple}, we let $\boldsymbol{q} =\boldsymbol{q}_T^{(z)}$ with 
\begin{equation*}
T  \coloneqq (\log N_0)^{1-\epsilon} \qquad\text{for a small $\epsilon>0$}. 
\end{equation*}
There is a standard Brownian motion $\big\{ \boldsymbol{B}^{(z)}_t ; t\in \R \big\}$ such that under the Assumptions~\ref{def:noise}, for any $\alpha>0$, 
\begin{equation} \label{eq:B}
    \Bigg\VERT  
    \sup_{|t| \leq N_0^{\alpha}} 
    \bigg|
    \boldsymbol{B}^{(z)}(t) 
    {+}
    \frac{1}{\sqrt{2}}
    \sum_{k=1}^{\lfloor t \mathfrak{L} \rfloor} 
    \frac{X_{N_0-k}+Y_{N_0-k}}{\sqrt{ \mathfrak{L}}}  
    \bigg|   
    \Bigg\VERT_1 
    \le \frac{ C_{\beta,\alpha} \log \mathfrak{L}}{\sqrt{ \mathfrak{L}}} .
\end{equation}
We consider the stochastic Airy equation driven by the Brownian motion $B$ such that 
\begin{equation}\label{Bcoup}
B_t =
\boldsymbol{B}^{(z)}_{t} 
\quad\text{for $t \le T$} , \qquad 
B_{t+T}- \boldsymbol{B}^{(z)}_T    \quad\text{is independent of $\F_\infty$ for $t\ge 0$}. 
\end{equation}
Finally, let $\mathcal{K} \subset \C$ be a compact set and
$K' = \big\{ \lambda \in \mathcal{K} : |\Im \lambda| \le N_0^\epsilon \big\}$.
\end{assumption}

The bound \eqref{eq:B} is a direct consequence of \cite[Theorem C.1]{LambertPaquette02} and the proof relies on a \emph{strong embedding} of random walk or \emph{KMT} coupling.
Moreover, the construction \eqref{Bcoup} is such that the Brownian motion $B$ is independent of $\F_{N_T}$, that is, independent of the hyperbolic part of the recursion. 

\medskip

In the sequel, we consider the rescaled characteristic polynomial
\begin{equation}\label{eq:Psi}
  \Psi^{(z)}_n(\lambda) \coloneqq  w_{\lfloor n\rfloor}\varphi_{\lfloor n\rfloor} (z_\lambda) , \qquad n\ge 1, \quad \lambda\in\C . 
\end{equation}

The function $(t,\lambda) \in \R\times \C \mapsto \Psi_{N_t}(\lambda)$ is piecewise constant in $t$ and analytic in $\lambda$ and we claim that it is an approximate solution to the stochastic Airy equation backward in time. 
In Section~\ref{sec:approx}, we formulate two different approximations which both hold with \emph{overwhelming probability}:
\begin{itemize}[leftmargin=.5cm]
\item In Theorem~\ref{main:approx}, we compare the rescaled characteristic polynomial and its \emph{discrete derivative} to the stochastic Airy function.
The stochastic Airy function is a universal object, the special solution of the stochastic Airy equation that is stable when the equation runs backward. 
Moreover, the only contribution from the \emph{hyperbolic part} of the recursion is given in terms of a mean-zero Gaussian random variable which is independent of the stochastic Airy function.
\item In Theorem~\ref{thm:para}, we compare the rescaled characteristic polynomial to another solution of the stochastic Airy equation with deterministic initial conditions at some arbitrary time $S\gg 1$ (which may be independent of $N_0$) . We view this solution as an approximation of the stochastic Airy function, that is adapted to the Brownian filtration of $\{\boldsymbol{B}^{(z)}_{t} ; t \le S\}$. The advantage is that the properties of this solution are more directly accessible using the SDE~\eqref{SA1}.
\end{itemize} 

The proofs of these results are given in Section~\ref{sec:fde} and rely on comparing the recursion for the characteristic polynomial in the parabolic regime to the integral equation \eqref{SA2} run backward from time $T\gg 1$ using the stability results for the stochastic Airy equation (e.g.~Proposition~\ref{prop:stability}).

\subsection{Uniform approximations} \label{sec:approx}

In this section, we give our two main approximation results for the rescaled characteristic polynomial  \eqref{eq:Psi} in the parabolic regime. 

\begin{theorem}\label{main:approx}
Under the Assumptions~\ref{def:para}, for  $t \in  [-e^{T},T/2]$ and $\lambda \in K'$,
\[\begin{aligned}
 \mathfrak{L} \big(\Psi_{N_t-1}(\lambda)  - \Psi_{N_t}(\lambda)\big) & =\boldsymbol{\Gamma}_\lambda^{(N,z)}  \Big( \SAi_\lambda'(t) +  \Upsilon_{1,\lambda}^{(N,z)}(t)\Big)
\\
\Psi_{N_{t}}(\lambda) & =  \boldsymbol{\Gamma}_\lambda^{(N,z)}
  \Big( \SAi_\lambda(t) + \Upsilon_{0,\lambda}^{(N,z)}(t)\Big)
\end{aligned}\]
where
\begin{itemize}[leftmargin=*]
\item $\SAi_\lambda$ is the stochastic Airy function, $\mathfrak{X}$  is a Gaussian process as in  \eqref{eq:SAIG} driven by the Brownian motion $B$, and 
\begin{equation} \label{def:Gamma}
\boldsymbol{\Gamma}_\lambda^{(N,z)}  = 
  \frac{\exp\Big(\frac{2}{\sqrt{\beta}}\int_{0}^{T} \mathfrak{X}(u)\d u + \boldsymbol{g} + \varepsilon_{\lambda}^{(N,z)}  \Big)} 
  {\E  \exp\big(\frac{2}{\sqrt{\beta}}\int_{0}^{T} \mathfrak{X}(u)\d u +  \boldsymbol{g} \big)} . 
\end{equation}
\item 
$\boldsymbol{g}$ is a mean-zero Gaussian, independent of $B$ with variance $\tfrac2\beta \log\big(\mathfrak{L}/\sqrt{4T}\big)$ (given explicitly by $\gamma \mathbf{W}_t(z)$ in Proposition~\ref{prop:Gauss}).
In particular, $\E\big( \tfrac{2}{\sqrt{\beta}}\int_{0}^{T} \mathfrak{X}(u)\d u + \boldsymbol{g} \big)^2 = \tfrac{2}{3\beta} \log(N_0) + \O(1)$ as $N \to\infty$. 
\item The error $\lambda\in \C \mapsto \epsilon_{\lambda}^{(N,z)} $ is analytic, measurable with respect $\F_{N_T} \vee \sigma(B)$ and there is a constant $C = C(\beta,\epsilon)$ such that 
\[
\P\bigg[
\max_{\lambda \in\mathcal{K}} | \varepsilon_{\lambda}^{(N,z)}| \ge C T^{2\epsilon-1/6} 
\bigg] \le e^{-(\log N_0)^{1+\epsilon}}.
\]
\item For $i=0,1$, the errors $\Upsilon_{i,\lambda}^{(N,z)}(t)$ are analytic for $\lambda\in \mathcal{K}$ and there are  constants $C = C(\beta,\epsilon,\ell)$ such that for any $\ell \in \N$: 
\begin{equation*}  \label{tail:main}
    \P\bigg[ \max_{t\in [-e^{T},T/2] ,\lambda \in K'}  \big| e^{ \frac23 t_+^{3/2}} \partial_\lambda^{\ell-1} \Upsilon_{i,\lambda}(t) \big|  \ge C  N_0^{\ell\epsilon-1/6}  \bigg] \le e^{-(\log N_0)^{1+\epsilon}}.
  \end{equation*}
\end{itemize} 
\end{theorem}

In particular, taking $z=1$ in Theorem~\ref{main:approx}, we recover the edge asymptotics of Theorem~\ref{main:thm} for the characteristic polynomial. 

\medskip

There are two types of errors in Theorem~\ref{main:approx}: 
\begin{itemize}[leftmargin=*]
\item $\varepsilon_{\lambda}$ which comes from the fact that the initial condition at the start of the parabolic stretch does not exactly match $\SAi$. However, this error is not independent of the Brownian motion $B$. 
\item $\Upsilon_{i,\lambda}(t)$ which comes from the discretization errors. 
\end{itemize}

For future purposes, it will be convenient to compare the characteristic polynomial to a solution of the SDE \eqref{SA1} that is adapted to the Brownian filtration with some explicit initial conditions. The next result provides precisely such a comparison.

\begin{theorem} \label{thm:para}
Under the Assumptions~\ref{def:para}, let $S \in [1,T]$ and $\big\{ \boldsymbol{\phi}_\lambda^{(z)}(t;S) : t\le S, \lambda \in \C \big\}$ be the solutions of the stochastic Airy equation driven by $\boldsymbol{B}$ with initial data at time S:
\begin{equation} \label{IdAiry} 
\boldsymbol{\phi}_\lambda(S;S) = \Ai(\lambda+S) , \qquad  
\boldsymbol{\phi}_\lambda'(S;S)  =  \Ai'(\lambda+S) . 
\end{equation}
Then, for $t \in [-e^S,S/2]$ and $\lambda \in K'$,
\[\begin{aligned}
 \mathfrak{L} \big(\Psi_{N_t-1}(\lambda)  - \Psi_{N_t}(\lambda)\big) & =  \big(  \boldsymbol{\phi}'_\lambda(t;S)+ \Upsilon_{1,\lambda}^{(N,S,z)}(t) \big) \exp\Big(\boldsymbol{q}_S -  \tfrac12 \E \boldsymbol{q}_S^2 + \varepsilon_{\lambda}^{(N,S,z)} \Big) 
\\
\Psi_{N_{t}}(\lambda) & =  \big(  \boldsymbol{\phi}_\lambda(t;S)+ \Upsilon_{0,\lambda}^{(N,S,z)}(t)  \big)  \exp\Big(\boldsymbol{q}_S -  \tfrac12 \E \boldsymbol{q}_S^2 + \varepsilon_{\lambda}^{(N,S,z)} \Big) 
\end{aligned}\]
where $\boldsymbol q_S=\boldsymbol q_S^{(z)}$ denotes the variable from Assumption~\ref{def:para} and \eqref{qrv}.  It is independent of
$\big\{ \boldsymbol{\phi}_\lambda(t;S) : t\le S, \lambda \in \C \big\}$,
approximately Gaussian, and for $1<\nu< \frac32$,  
\begin{itemize}[leftmargin=*] 
\item the error $\lambda\in \C \mapsto \epsilon_{\lambda}$ is analytic and there are constants $C = C(\beta,\nu)$ such that for $S\gg1$, 
\[
\P\bigg[
\max_{\lambda \in\mathcal{K}} | \varepsilon_{\lambda} | \ge C  S^{\frac12-\frac\nu3}
\bigg] \le e^{-S^{\nu}}.
\]
\item for $i=0,1$, the errors $\lambda \in \C \mapsto \Upsilon_{i,\lambda}(t)$ are analytic and, for $0<\alpha<1/6$, there is a  constant $C = C(\beta,\nu,\alpha)$ such that 
\begin{equation*}  
    \P\bigg[ \max_{t\in [-e^{S},S/2] ,\lambda \in K'}  \big|   \Upsilon_{i,\lambda}(t) \big|  \ge C e^{-S^{3/2}/4} \vee  N_0^{-\alpha}  \bigg] \le e^{-S^{\nu}}.
  \end{equation*}
\end{itemize} 
\end{theorem}

\begin{remark} 
One has $\E  \SAi_\lambda(t) = \Ai(\lambda+t)$ for $(t,\lambda) \in \R \times \C$; for any $S\ge 1$, one also has $\E \boldsymbol{\phi}_\lambda(t;S) = \Ai(\lambda+t)$ for $t \le S$ and $\lambda\in \C$. 
\end{remark}

\subsection{Schr\"odinger equation; from discrete to continuous time} \label{sec:fde}
 
The goal of this section is to rewrite the Dumitriu--Edelman recurrence \eqref{eq:recurrence} in the parabolic region as a discrete version of the stochastic Airy equation in order to prove the results from Section~\ref{sec:approx}. 
 Namely, we show that $t\mapsto \Psi_{N_{t}}(\lambda)$ is a solution to a perturbed version of the integral equation  \eqref{SA2}. Then, using the stability Proposition~\ref{prop:stability}, we obtain the following result:

\begin{proposition}\label{prop:para} 
Under the Assumptions~\ref{def:para}, let $\big\{ \phi_\lambda(t) : t\le T, \lambda \in \C \big\}$ be the solution of the stochastic Airy equation \eqref{SA1} driven by the Brownian motion $\big\{B(t) ; t\le T \big\}$  with initial data: 
\begin{equation} \label{id1} 
c_{1,\lambda} = \phi_\lambda(T) = \frac{\Psi_{N_T}(\lambda)}{\Ai(\lambda+T)} e^{-\boldsymbol{q} +\frac12 \E\boldsymbol{q}^2} , \qquad
c_{2,\lambda} = \phi_\lambda'(T) =\frac{\mathfrak{L}\big(\Psi_{N_T-1}(\lambda)  - \Psi_{N_T}(\lambda)\big)}{\Ai(\lambda+T)} e^{-\boldsymbol{q} +\frac12 \E\boldsymbol{q}^2}.
\end{equation}
Then, for any $0<\alpha<1/3$, there are  constants $C=C(\beta,\epsilon,\alpha)$ and $c=c(\beta,\epsilon,\alpha)>0$ so that with probability at least  $1- e^{-c (\log N_0)^{1+\epsilon}}:$ for $t\in [-e^{T},T]$ and $\lambda \in K'$:   
\[
\big| \Psi_{N_t}(\lambda) e^{-\boldsymbol{q} +\frac12 \E\boldsymbol{q}^2} - \phi_\lambda(t) \Ai(\lambda+T)    \big| \le C N_0^{-\alpha} e^{-\frac23 t_+^{3/2} }
\]
and the initial data satisfies uniformly for $\lambda \in \mathcal{K}$: 
\begin{equation} \label{id2} 
c_{1,\lambda} = 1+\O(T^{2\epsilon-1/6}) , \qquad
c_{2,\lambda}/\sqrt{T} = -1+ \O(T^{2\epsilon-1/6}). 
\end{equation}
\end{proposition}

Proposition~\ref{prop:para} compares the rescaled characteristic polynomial to a solution of the stochastic Airy equation with matching initial data (up to a scaling) on an event of overwehlming probability. In particular, under the Assumptions~\ref{def:para}, the initial data \eqref{id1} is independent of $B$. 

Then, using the stability Proposition~\ref{prop:SAiwins} and the asymptotics of $\SAi$, this will be straightforward to deduce Theorem~\ref{main:approx}. 

\begin{proof}[Proof of Proposition~\ref{prop:para}]
In terms of the (normalized) random variables \eqref{def:XY}, we can rewrite the recurrence \eqref{eq:recurrence}:
\[
  \varphi_n = 
  \left( z - \frac{X_{n-1}}{\sqrt{2N \beta}} \right)\varphi_{n-1}
  -\left( \frac{Y_{n-1} \sqrt{2\beta (n-1)} + \beta(n-1)}{4N\beta} \right)\varphi_{n-2}  .
\]
Rewriting this equation in terms of \eqref{parascale}--\eqref{eq:Psi}, we arrive at 
\begin{equation}\label{eq:recurrence2}
\Psi_n- 2\Psi_{n-1} + \Psi_{n-2} = 
\underbrace{2\left(z_\lambda \sqrt{\tfrac Nn}-1 - \tfrac{X_{n-1}}{\sqrt{2n \beta}}\right)}_{=R_{n-1,\lambda}\mathsf{h}}\Psi_{n-1}
  -\underbrace{\left( \tfrac{Y_{n-1} \sqrt{2}}{\sqrt{\beta n}} + \sqrt{\tfrac{n-1}{n}}-1 \right)}_{=-S_{n-1}\mathsf{h}}\Psi_{n-2} 
\end{equation}
where $\mathsf{h} = \mathfrak{L}^{-1}$. 
Observe that the \emph{driving term} satisfies, with $n=N_t$,  
\begin{equation} \label{drive}
2z_\lambda \sqrt{\frac N{N_t}}-2 = 2 \frac{1 + \lambda / 2 \mathfrak{L}^{2}}{\sqrt{1- t/ \mathfrak{L}^{2}}}-2 =  (\lambda+t) \mathsf{h}^{2} + \O\big(\mathsf{h}^{4}\big)
\end{equation}
so that \eqref{eq:recurrence2} is a discretization of the SDE \eqref{SA1} on scale $\mathsf{h}$. 
To make the comparison, we introduce a piecewise $C^1$ function which linearly interpolates   $(\Psi_n)_{n\ge N_T}$: 
\begin{equation}\label{eq:psi}
  \psi_\lambda'(t) \coloneqq \mathfrak{L} \big(\Psi_{N_t-1}(\lambda)  - \Psi_{N_t}(\lambda)\big)
  \quad\text{for $t\le T$} \qquad
  \text{and}
  \qquad
\psi_\lambda(T) =  \Psi_{N_{T}}(\lambda).
\end{equation}
% \begin{equation}\label{eq:Pi}
%   \Pi_\lambda'(t) \coloneqq \mathfrak{L} \big(\Psi_{n-1}(\lambda)  - \Psi_n(\lambda)\big)\one_{n = N_0 -\lfloor t \mathfrak{L}\rfloor }
%   \quad
%   \text{and}
%   \quad
% \Pi_\lambda(T) =  \Psi_{N_{-T}}(\lambda).
% \end{equation}

Without loss of generality, we assume that  $T \in \mathsf{h}\Z$. Then, by construction, the two processes $\Psi_{N_t}(\lambda) = \psi_\lambda(t)$ on the mesh $ t\in (-\infty, T]\cap \mathsf{h}\Z$. 
Moreover, by substituting \eqref{eq:psi} into \eqref{eq:recurrence2} with $n-1 =N_t$, we obtain on the mesh:
\begin{equation}\label{eq:recurrence3}
\psi_\lambda'(t)-\psi_\lambda'(t-\mathsf{h}) = (R_{n,\lambda}+S_n) \psi_\lambda(t) +S_n \psi_\lambda'(t) \mathsf{h}
\end{equation}

We use the following version of Abel's summation formula to turn \eqref{eq:recurrence3} into an integral equation: if $f:\R\to\R$ is an absolutely continuous function, then for any $\mathfrak{L}\ge 1$ and $T>t$:  
\begin{equation*}
  \sum_{\mathfrak{L} t  < k \leq \mathfrak{L} T}  a_k f(\tfrac k{\mathfrak{L}}) =  -  \int_t^T A(u,t) f'(u) \d u  + A(T,t)f(T) , \qquad A(u,t) = \sum_{t\mathfrak{L} < n \le u \mathfrak{L}} a_n . 
\end{equation*} 

Then, summing the equation \eqref{eq:recurrence3}, we obtain on the mesh
\begin{equation}\label{eq:recurrence1}
\psi_\lambda'(T)-\psi_\lambda'(t) = U_\lambda(T,t)\, \psi_\lambda(T) -  \int_t^T U_\lambda(u,t)\, \psi_\lambda'(u)\d u  + \mathsf{h} \sum_{t\mathfrak{L} < n \le T \mathfrak{L}}  S_{N_{n\mathsf{h}}} \, \psi_\lambda'(n \mathsf{h} ) 
\end{equation}
with the  kernel for $t < u$, 
\begin{equation} \label{Udiscrete}
 U_\lambda(u,t)    \coloneqq  \sum_{N_u\le n<N_t}  \big( R_{n,\lambda} + S_{n} \big) .
\end{equation}
In \eqref{eq:recurrence1}, since $t \mapsto \psi_\lambda'(t)$ is piecewise constant, we can also write 
\[
\mathsf{h} \sum_{t\mathfrak{L} < n \le T \mathfrak{L}}  S_{N_{n\mathsf{h}}} \, \psi_\lambda'(n \mathsf{h}) 
= \mathsf{h} S_{N_T}\psi_\lambda'(T)+   \int_t^T  S_{N_u}  \psi_\lambda'(u) \1\{u \ge [t] \} \d u 
\qquad\text{where $ [t] \coloneqq \mathsf{h} \big(\lfloor t \mathfrak{L} \rfloor+1\big)$}.
\]
Moreover, since both sides of \eqref{eq:recurrence1} are piecewise continuous on the mesh, this equation holds for all $t\le T$. 

\medskip

Now, let $\Es_\lambda$ denote the kernel in \eqref{SA2} with the Brownian motion $B=B^{(z)}$, and set
\[
\mathsf{Norm}_\lambda
\coloneqq \Ai(\lambda+T)\exp\bigl(\boldsymbol{q}-\tfrac12\E\boldsymbol{q}^2\bigr).
\]
According to \eqref{id1},
\begin{equation*} %\label{id0} 
%\widehat{c}_1 = 
\psi_\lambda(T) =  \Psi_{N_T}(\lambda) = c_1\mathsf{Norm}_\lambda, \qquad
%\widehat{c}_2 =  
\psi_\lambda'(T)   
= \mathfrak{L}\big(  \Psi_{N_T-1}(\lambda)  - \Psi_{N_T}(\lambda)  \big) 
= c_2\mathsf{Norm}_\lambda .
\end{equation*}
Using these notations, we conclude that
\begin{equation} \label{adiscrete}
\begin{aligned}
\psi_\lambda'(t)  & = \psi_\lambda'(T)   + \int_t^T \big( U_\lambda(u,t)- S_{N_u}  \1\{u \ge [t] \} \big)\, \psi_\lambda'(u)\d u - \psi_\lambda(T) U_\lambda(T,t) - \mathsf{h} \psi_\lambda'(T)  S_{N_T} \\
\widehat{\Phi}_\lambda(t)
& = {c}_2  + \int_t^T \big(\Es_\lambda(u,t)+ \mathcal{D}_{\lambda}^1(t,u)\big) \, \widehat{\Phi}_\lambda(u)\d u - {c}_1 \Es_\lambda(T,t) + \mathcal{D}_{\lambda}^2(t)
\end{aligned}
\end{equation}
where (after rescaling) $\widehat{\Phi}_\lambda(t) = \psi_\lambda'(t)/\mathsf{Norm}_\lambda $ and the \emph{errors}: 
\begin{equation}\label{DPert}
\mathcal{D}_\lambda^1(u,t)   =  U_\lambda(u,t)-\Es_\lambda(u,t) - S_{N_u}\1\{u \ge [t] \}, \qquad 
\mathcal{D}_\lambda^2(t)  =  {c}_1 \big( U_\lambda(T,t)-\Es_\lambda(T,t) \big)  - \mathsf{h} {c}_2  S_{N_T}   . 
\end{equation}

Thus, \eqref{adiscrete} is a perturbation of the integral equation \eqref{SA2} of the form \eqref{eq:SAperturb}.   
Then, the strategy is to compare $\widehat{\Phi}_\lambda(t)$ to a solution of the stochastic Airy equation with initial condition $(c_1,c_2)$ by applying the stability Proposition~\ref{prop:stability} with the following estimates. 

\begin{lemma}\label{lem:Kdiff}
Let $\Es$ be the kernel \eqref{kernelU} driven by  the Brownian motion $\boldsymbol{B}^{(z)}$ from \eqref{eq:B}. 
Let $0<\epsilon <1/6$. 
With probability at least $1- e^{-(\log N_0)^2}$, it holds uniformly for $\lambda \in \C$ with $|\lambda| \le N_0^\epsilon$ and $t, u \in \R$ with $|t|, |u| \le N_0^{\epsilon}$ and $t<u$, 
\[
| \Es_\lambda(u,t) - U_\lambda(u,t)|  \lesssim (\log N_0)^3 N_0^{-1/6}  
\]
and  
\[
\max_{ t \in \mathsf{h} \Z, |t| \le N_0^{\epsilon} }|S_{N_t}| \le (\log N_0)^3 N_0^{-1/6}  
\]
\end{lemma}

\begin{proof}[Proof of Lemma~\ref{lem:Kdiff}]
According to \eqref{eq:recurrence2} and \eqref{Udiscrete}, we can decompose the kernel 
\[
U_\lambda(u,t)
= \mathrm{I}_\lambda +  \mathrm{II} - \mathrm{III} , 
\]
where 
\[
\mathrm{I}_\lambda = 2 \mathsf{h}^{-1}\hspace{-.3cm}\sum_{N_u\le n<N_t}  \left(z_\lambda \sqrt{\tfrac{N}{n+1}}-1 \right) , \qquad 
\mathrm{II} = \mathsf{h}^{-1}\hspace{-.3cm}\sum_{N_u\le n<N_t} 1- \sqrt{\tfrac{n}{n+1}} , \qquad
\mathrm{III} =  \sqrt{\tfrac 2\beta} \mathsf{h}^{-1}\hspace{-.3cm}\sum_{N_u\le n<N_t} \frac{X_n+Y_n}{\sqrt{n+1}} .
\]
The terms $ \mathrm{I}_\lambda $ and $\mathrm{II}$ are deterministic and, by \eqref{drive}, we have uniformly for $|\lambda| \le N_0^\epsilon$ and $|u|, |t| \le N_0^{\epsilon}$ 
\[
\mathrm{I}_\lambda = \int_{t}^{u} (\lambda+s) ds + \O\big(\mathsf{h}^{1-3\epsilon}\big) , \qquad \mathrm{II}  =  \O\big(\mathsf{h}^{1-3\epsilon}\big)
\]
using that $N_0 = \mathsf{h}^{-3}$ and approximating the sum by integrals. 
Here the error is uniform on the stated range: the mesh and floor errors
are $O(N_0^\epsilon\mathsf{h})=O(\mathsf{h}^{1-3\epsilon})$, and the Taylor
remainder contributes $O(\mathsf{h}^{2-9\epsilon})$, which is smaller when
$\epsilon<1/6$.

Similarly, the noise term $\mathrm{III}$ can be compared to the Brownian motion $\boldsymbol{B}$ using~\eqref{eq:B}, with probability at least  $1- e^{-c (\log N_0)^{2}}$:
\[
\big| \mathrm{III} - \tfrac{2}{\sqrt{\beta}}\big(\boldsymbol{B}(t) -  \boldsymbol{B}(u) \big)\big| \lesssim(\log \mathsf{h}^{-1})^3 \mathsf{h}^{1/2}
\]
Indeed, before applying \eqref{eq:B}, replacing $(n+1)^{-1/2}$ by
$N_0^{-1/2}$ in $\mathrm{III}$ costs at most
$O\big((\log\mathsf{h}^{-1})^3\mathsf{h}^{2-9\epsilon/2}\big)$ uniformly on this
range, which is negligible compared with
$(\log\mathsf{h}^{-1})^3\mathsf{h}^{1/2}$ because $\epsilon<1/6$.
Combining these estimates, with $\Es$ denoting the kernel \eqref{kernelU} driven by  $\boldsymbol{B}$, with probability at least  $1- e^{-c (\log N_0)^{1+\epsilon}}$:
\[
|\Es_\lambda(u,t) - U_\lambda(u,t) |  \lesssim (\log N_0)^3 N_0^{-1/6}  
\]
uniformly for  $|\lambda| \le N_0^\epsilon$ and $|u|, |t| \le N_0^{\epsilon}$. 
This proves the first claim. The second claim is a direct consequence of the fact that the random variable $Y_n$ are sub-exponential.
Using that $\mathsf{h}^{-1} N_0^{-1/2} = N_0^{-1/6}$, by a union bound: for any $C>0$, if $N_0 \gg 1$, 
\[
\P \left[
\max_{k\in \Z, |k| \le N_0^{C} } \hspace{-.1cm}  |S_{N_0+k}|  \ge (\log N_0)^3 N_0^{-1/6} 
\right]  \le   e^{-  (\log N_0)^2} . \qedhere
\]
\end{proof}

Then, Lemma~\ref{lem:Kdiff} implies that 
on an event of probability at least
$1- e^{-(\log N_0)^2}$:
\begin{equation*}
\sup_{t\in [-e^{T},T], \lambda\in K}\sup_{u\in[t,T]}|\mathcal{D}_\lambda^1(u,t)| \lesssim (\log N_0)^3 N_0^{-1/6} \, ,\qquad
 \sup_{ t\in [-e^{T},T], \lambda\in K}|\mathcal{D}_\lambda^2(t) | \lesssim  \sqrt{ {c}_1^2+{c}_2^2}\, (\log N_0)^3  N_0^{-1/6} .
\end{equation*}
In particular, on this event, the conditions \eqref{eq:D_stab} are satisfied with $\nu =\frac1{1-\epsilon}$ \big(so that $e^{T^\nu} = N_0$\big) and $0<\alpha<1/6$. 
Hence, by Proposition~\ref{prop:stability}, with probability at least $1 - e^{-c\, T^{1+2\epsilon}}$: for $t\in[-e^{T},T]$ and $\lambda \in K'$, 
\begin{equation*} 
|\widehat{\Phi}_\lambda(t)-  \Phi_\lambda(t)| \le \big(1+\sqrt{{c}_1^2+  {c}_2^2}\big) N_0^{-\alpha} E_\lambda(t,T)
\end{equation*}
where $\Phi_\lambda = \phi'_\lambda$ is the solution of the stochastic Airy equation \eqref{SA2} with initial conditions  $({c}_1,{c}_2)$ at time $T$. 
% In particular, conditionially on $({c}_1,{c}_2)$, the process $t\in (-\infty,T] \mapsto \phi_\lambda$  is adpated to the reverse filtration $\sigma(B_t : t\le T)$. 

Moreover, by Theorem~\ref{thm:hypersimple}, the initial data satisfies \eqref{id1} on an event of probability $1- e^{- c (\log N_0)^{1+\epsilon}}$. 
Going back to  $ \psi_\lambda'(t)=\mathsf{Norm}_\lambda\widehat{\Phi}_\lambda(t) $, we conclude that on an event of probability at least
$1- e^{-c (\log N_0)^{1+\epsilon}}$: for $t\in[-e^{T},T]$ and $\lambda \in K'$,
\begin{equation} \label{para}
|\psi_\lambda'(t) e^{-\boldsymbol{q} +\frac12 \E\boldsymbol{q}^2}  - \Ai(\lambda+T)   \phi_\lambda'(t)| \lesssim N_0^{-\alpha} \underbrace{ E_\lambda(t,T) \Ai(\lambda+T)}_{\lesssim e^{- \frac23 t_+^{3/2}}} .
\end{equation}
Integrating \eqref{para}, we obtain a similar uniform bound for $\psi_\lambda(t) = \Psi_{N_t}(\lambda)$ on the mesh $ t\in [-e^T, T]\cap \mathsf{h}\Z$. 
Finally, observe that by Corollary~\ref{cor:solbound}, the continuous part satisfies for any $\epsilon>0$, $| \Ai(\lambda+T)   \phi_\lambda'(t)| \lesssim N_0^\epsilon e^{- \frac23 t_+^{3/2}} $  for $t\in[-e^{T},T]$  with a similar probability. 
Hence the approximation holds for the piecewise constant interpolation $\Psi_{N_t}(\lambda)$ for all $t\in[-e^{T},T]$. This concludes the proof of Proposition~\ref{prop:para}. 
% \begin{equation} \label{coup3}
% |\psi_\lambda(t) e^{-\boldsymbol{q} +\frac12 \E\boldsymbol{q}^2}  - \Ai(\lambda+T)   \phi_\lambda(t)| \lesssim N_0^{-\alpha} . 
% \end{equation}
\end{proof}

The next step is to replace the solution from Proposition~\ref{prop:para} by the stochastic Airy function by using that the stochastic Airy equation run backward is not sensitive to initial conditions. 

\begin{proof}[Proof of Theorem~\ref{main:approx}]
We apply the stability Proposition~\ref{prop:SAiwins} to relate the solution $\phi_\lambda$ from Proposition~\ref{prop:para}  to the stochastic Airy function:
\[
    \phi_\lambda(t)
    =
    \Theta_\lambda  \SAi_\lambda(t)
    +
    \chi_\lambda(t)
    \qquad
    \text{for all $t \in \R$}
\]
and with probability at least $1 - e^{-T^\nu}$ --- choosing $\nu =\frac{1+\epsilon}{1-\epsilon} <3/2$ --- one has
\[
\Theta_\lambda = \frac{c_{1}-c_{2}/\sqrt{T} + \O(T^{2\epsilon-1/6})}{2\SAi_\lambda(T)\big(1+\O(T^{2\epsilon-1/6})\big)}
= \frac{2+\O(T^{\epsilon/2-1/6})}{2\SAi_\lambda(T)}
=  \frac{1+\O(T^{\epsilon/2-1/6})}{\Ai(T+\lambda)}
\frac{\exp\big(\tfrac{2}{\sqrt{\beta}}\int_0^T \mathfrak{X}(u)\d u\big)}{\E\exp\big(\tfrac{2}{\sqrt{\beta}}\int_0^T \mathfrak{X}(u)\d u\big)}
\]
where we used first the initial data \eqref{id2} and second the asymptotics of Theorem~\ref{thm:asymp} for $\SAi$. This implies that with probability at least  $1- e^{-c (\log N_0)^{1+6\epsilon}}$: for $t\in [-e^{T},T/2]$ and $\lambda\in K'$, 
\[
\Ai(\lambda+T)\phi_\lambda(t)
= \SAi_\lambda(t) \frac{\exp\big(\tfrac{2}{\sqrt{\beta}}\int_0^T \mathfrak{X}(u)\d u +\O(T^{2\epsilon-1/6})\big)}{\E\exp\big(\tfrac{2}{\sqrt{\beta}}\int_0^T \mathfrak{X}(u)\d u\big)}+ \widehat{\chi}_\lambda(t)  
\]
where the error is measurable with respect to $\F_{N_T} \vee \sigma(B)$, independent of $t$,  and $\widehat{\chi}_\lambda(t)=\Ai(\lambda+T){\chi}_\lambda(t)$ satisfies
\begin{equation}\label{errchi}
|\widehat{\chi}_\lambda(t)| \vee |\widehat{\chi}'_\lambda(t)| \lesssim   e^{- T^{3/2}/4}  e^{- \frac23 t_+^{3/2}} . 
\end{equation}
By \eqref{para}, we conclude that with probability at least  $1- e^{-c (\log N_0)^{1+\epsilon}}$: uniformly for $t\in [-e^{T},T/2]$ and $\lambda\in K'$, 
\[
 \psi_\lambda'(t) e^{-\boldsymbol{q} +\frac12 \E\boldsymbol{q}^2}  =  
 \SAi_\lambda'(t) \frac{\exp\big(\tfrac{2}{\sqrt{\beta}}\int_0^T \mathfrak{X}(u)\d u +\O(T^{2\epsilon-1/6})\big)}{\E\exp\big(\tfrac{2}{\sqrt{\beta}}\int_0^T \mathfrak{X}(u)\d u\big)}  + \O\big( N_0^{\epsilon-\frac16} e^{- \frac23 t_+^{3/2}} \big) 
\]
where $\psi_\lambda'(t)$ as in \eqref{eq:psi} and the error \eqref{errchi} is negligible.
Moreover, the error is analytic in $\lambda$, so all its derivatives are controlled uniformly. 

Finally, by Proposition~\ref{prop:Gauss}, we can replace the random variable $\boldsymbol{q}$ by a mean-zero Gaussian $\boldsymbol{g}$  with variance $\tfrac2\beta \log\big(\mathfrak{L}/\sqrt{4T}\big)$, independent of $B$, with an error $\O(N_0^{-\alpha})$ with $0<\alpha<1/6$ with probability at least  $1- e^{-c (\log N_0)^{1+\epsilon}}$.
Moreover, by \eqref{eq:SAIG}, $\int_0^T \mathfrak{X}(u)\d u$ is also Gaussian, mean-zero, with variance $\tfrac14\log T+\O(1)$. So, one has $|\int_0^T \mathfrak{X}(u)\d u |\lesssim (\log N_0)^{\frac 12 + \epsilon}$ with probability at least  $1- e^{-c (\log N_0)^{1+\epsilon}}$. 

By Proposition~\ref{prop:para}, a similar estimate also holds for $\Psi_{N_t}(\lambda)$. 
\end{proof}

Finally, we can also replace $\SAi$ by an adapted solution of the backward stochastic Airy SDE defined at an arbitrary time $S \gg 1$ with $S\le T$. 

\begin{proof}[Proof of Theorem~\ref{thm:para}]
We proceed as in the proof of Theorem~\ref{main:approx}.
Let $1<\nu<3/2$ and $S \le T = (\log N_0)^{1-\epsilon}$ with $\epsilon>0$ sufficiently small.  
Applying Proposition~\ref{prop:SAiwins} to the solution $ \boldsymbol{\phi}_\lambda(\cdot;S)$ with initial condition \eqref{IdAiry}, with probability at least $1 - e^{-cS^{\nu}}$:
\[
\SAi_\lambda(t)=  \Theta_\lambda^{-1}\big( \boldsymbol{\phi}_\lambda(t;S)+ \O(e^{-S^{3/2}/4}) \big)
\]
where the error is a solution controlled uniformly for $\lambda \in K'$ and $t \in [-e^{S},S/2]$ and the coefficient 
\[  
\Theta_\lambda
=  \frac{\Ai(\lambda+S)-\Ai'(\lambda+S)/\sqrt{S} + \O(S^{-\epsilon})}{2\SAi_\lambda(S)\big(1+\O(S^{-\epsilon})\big)} 
= \frac{\exp\big(\tfrac{2}{\sqrt{\beta}}\int_0^S \mathfrak{X}(u)\d u\big)}{\E\exp\big(\tfrac{2}{\sqrt{\beta}}\int_0^S \mathfrak{X}(u)\d u\big)}\big(1+\O(S^{-\alpha})\big) 
\]
where $\alpha = \frac12-\frac\nu3$. 
Here, we used the Airy asymptotics \eqref{Airyasymp} and the asymptotics from Theorem~\ref{thm:asymp}.  
Recall that
$\boldsymbol{\Gamma}_\lambda=\boldsymbol{\Gamma}^{(N,z)}_\lambda$ is defined in
\eqref{def:Gamma}.  Combining this definition with the preceding expression for
$\Theta_\lambda$, one has, with probability at least $1-e^{-cS^\nu}$,
\[
\boldsymbol{\Gamma}_\lambda\Theta_\lambda^{-1}
=\frac{\exp\Big(\tfrac{2}{\sqrt\beta}\int_S^T\mathfrak X(u)\,\d u+\boldsymbol q_T
+\O(S^{-\alpha})\Big)}
{\E\exp\Big(\tfrac{2}{\sqrt\beta}\int_S^T\mathfrak X(u)\,\d u+\boldsymbol q_T\Big)}
=\exp\big(\boldsymbol{q}_S-\tfrac12\E\boldsymbol{q}_S^2
+\O(S^{-\alpha})\big),
\]
where for the second equation we have used Propositions~\ref{prop:qnoise} and~\ref{prop:Gauss}.
We conclude that with probability at least $1 - e^{-cS^{\nu}}$: it holds uniformly for $\lambda \in K'$ and $t \in [-e^{S},S/2]$, 
\[
\Psi_{N_{t}}(\lambda) =  
  \Big( \boldsymbol{\phi}_\lambda(t;S) +\O\big( e^{-S^{3/2}/4} \vee N_0^{\epsilon-1/6}  \big) \Big)
\exp\Big(\boldsymbol{q}_S -  \tfrac12 \E \boldsymbol{q}_S^2 + \O(S^{-\alpha})\Big) 
\]
A similar estimate also holds for the \emph{discrete derivative} $ \mathfrak{L} \big(\Psi_{N_t-1}(\lambda)  - \Psi_{N_t}(\lambda)\big)$ in terms of $\boldsymbol{\phi}_\lambda'(t;S)$.  
\end{proof}

% ===== end sections/09-characteristic-poly-approximation.tex =====

\appendix
% ===== begin sections/A-uniformity-hyperbolic.tex =====
%Uniformity of the error bound in the
\section{Hyperbolic approximation} 
\label{sec:hyperuniform}

In this Appendix, we review the main results from \cite{LambertPaquette02} on the \emph{hyperbolic part of the characteristic polynomial recursion} \eqref{eq:recurrence}. 
In particular, this paper gives an asymptotic approximation of the characteristic polynomial $\varphi_N$ away from the spectral support:

\begin{theorem}[Theorem 1.4 in \cite{LambertPaquette02}]\label{thm:planar}
There is a Gaussian process $\mathbf{W} : \C \setminus [-1,1] \to \C$ such that  $\mathbf{W}(\overline{z}) = \overline{\mathbf{W}(z)}$ with covariance kernel: for $z,w\in \C\setminus[-1,1]$ 
\begin{equation}  \label{eq:Wtz}
  \E[\mathbf{W}(z)\mathbf{W}(w)] = - \log\big(1- J(z)J(w)\big)  , \qquad   \text{where}
  \quad
J(z) \coloneqq z- \sqrt{z^2-1}, 
\end{equation}
such that for any compact set $K \subset \C\setminus[-1,1]$, 
there exist constants $C_{K}, c>0$ so that for all $N\in\N$, 
\[
\P\left[ 
  \sup_{z \in K} \Bigg|
  \tfrac{\varphi_{N}(z) \E\big[\exp\big(\sqrt{\frac{2}{{\beta}}}\mathbf{W}(z)\big)\big]  }{\widetilde\varphi_N(z) \exp\big(\sqrt{\frac{2}{{\beta}}}\mathbf{W}(z)\big)} -1 \Bigg|  \ge C_{K} N^{-\frac{1}{10}}  \right] \le e^{-c N^{1/45}} . 
 \]
\end{theorem}

The mapping $\sqrt{\cdot}$ is chosen so that $J: \C \setminus [-1,1] \to \D$ is a conformal map and 
$\mathbf{W}$ is the analytic extension of a Gaussian log-correlated field on $[-1,1]$ (see \cite[Section 1.4]{LambertPaquette02} for details). 
This Gaussian field completely describes the fluctuations of $z\mapsto \varphi_{N}(z)$ away from the spectral support and it discloses the underlying conformal structure.  
Moreover, the random field $\mathbf{W} = \boldsymbol{W}_1$ is defined using the Brownian embedding from \cite[Appendix~C]{LambertPaquette02} by an explicit formula; see \eqref{gfield} below. 

\medskip

A crucial input to prove our main result, Theorem~\ref{thm:para} (and Theorem~\ref{main:thm} as a consequence), is a generalization of Theorem~\ref{thm:planar} at the boundary of the parabolic region. 

\begin{theorem}\label{thm:hypersimple}
% We consider a sequence $z=z(N)$ such that 
% \[
% z(N) \in [-1,1] , \qquad z(N)\to z \quad\text{as $N\to\infty$}, \qquad \lim_{N\to\infty} N_0(z) =  \lfloor Nz(N)^2 \rfloor =\infty . 
% \]
Let $\mathcal{K} \subset \C$ be a compact set and 
\begin{equation}\label{Time}
T  \coloneqq (\log N_0)^{1-\epsilon} , \qquad\text{for $0<\epsilon<1/6$}. 
\end{equation}
Under the Assumptions~\ref{def:noise}, using the notations \eqref{def:hyper}--\eqref{eq:Psi}, there is a constant $c>0$ such that 
\begin{equation}\label{est1}
\P\Bigg[ \max_{\lambda\in\mathcal{K}} \bigg|  \frac{\Psi_{N_T-1}(\lambda)}{\Ai(\lambda+T)} e^{-\boldsymbol{q} +\frac12 \E\boldsymbol{q}^2} -1 \bigg| \ge (\log N_0)^{\epsilon-1/6} \Bigg] \le e^{- c (\log N_0)^{1+\epsilon/2}}  
\end{equation}
and 
\begin{equation}\label{est2}
\P\Bigg[ \max_{\lambda\in\mathcal{K}} \bigg|  \frac{\mathfrak{L}\big(\Psi_{N_T} - \Psi_{N_T-1}\big)(\lambda)}{\sqrt{T} \, \Ai(\lambda+T)}e^{-\boldsymbol{q} +\frac12 \E\boldsymbol{q}^2} -1 \bigg| \ge (\log N_0)^{\epsilon-1/6} \Bigg] \le e^{- c (\log N_0)^{1+\epsilon/2}} . 
\end{equation}
where the random variable
\begin{equation} \label{qrv}
\boldsymbol{q}= \boldsymbol{q}_T^{(z)} \coloneqq \sum_{k=1}^{N_T-1} \frac1{\sqrt{Nz^2-k}}\frac{X_k +  Y_k J\big(z \sqrt{\tfrac{N}{k}}\big)}{\sqrt{2\beta}} . 
\end{equation}
\end{theorem}

In particular, \eqref{est1} and \eqref{est2} imply \eqref{id2}.

 This statement is a consequence of \cite[Theorem 1.6]{LambertPaquette02} after an appropriate reworking and a suitable meshing argument to control the errors uniformly in $\lambda\in \mathcal{K}$. 

\begin{proof} The proof is divided in several steps:
\begin{itemize}[leftmargin=.5cm]
\item[1] We state the relevant results from \cite{LambertPaquette02} on the product of transfer matrices \eqref{eq:recurrence} in the hyperbolic regime up to a neighborhood of the turning point (in particular, the parameter $T$ is chosen so that the estimates hold with overwhelming probabilities for a fixed $\lambda\in\C$). 
\item[2] We perform some deterministic reductions to express the \emph{mean behavior} of the characteristic polynomial in terms of the Airy function. Using the asymptotics 
\begin{equation} \label{Airyasymp}
   \Ai(\lambda + T) =  \frac{ e^{-\frac23(T+\lambda )^{3/2}}}{2\sqrt{\pi}(T+\lambda)^{1/4}}  \left(1+ \O(T^{-3/2})\right)
   \quad\text{and}\quad
      \Ai'(\lambda + T) = -  \sqrt{T+\lambda}    \Ai(\lambda + T)  \left(1+ \O(T^{-3/2})\right)     , 
\end{equation}
 one can replace $\sqrt{T} \, \Ai(\lambda+T) \leftarrow \Ai'(\lambda+T)$ in \eqref{est2} with a similar error.  
\item[3] We show that the contribution from the noise is given by $e^{\boldsymbol{q} -\frac12 \E\boldsymbol{q}^2}$, up to an error which is controlled uniformly for $\lambda\in\mathcal{K}$. It is also possible to replace $\boldsymbol{q}$ by a Gaussian random variable with a similar error; see Proposition~\ref{prop:Gauss}. 
\item[4] Finally, we obtain uniform control of the approximation errors using meshing argument which relies on the fact that $\lambda\in \C \mapsto \varphi_{n} (z_\lambda) $ is a polynomial of degree $n$. 
\end{itemize}

\paragraph{1. Review of the results from \cite{LambertPaquette02}.}
In terms of the noise \eqref{def:XY}, we work on the \emph{truncation event}
\[
\mathsf{Trunc} \coloneqq \big\{ |X_k|+|Y_k| \le T^{3/2} : k\le N_0/2 \big\} \cap \big\{ |X_k|+|Y_k| \le T^{3/4} : N_0/2 \le k\le N_0\big\}.
\]
These events depend on $z\in [-1,1]$, and one has $\P[\mathsf{Trunc}^{\rm c}] \lesssim N_0 e^{- cT^{3/2}} \lesssim e^{-c(\log N_0)^{3/2-3\epsilon/2}}$ for some $c>0$ (provided that $\epsilon<1/6$).  
On this event, we can use the bounds from \cite[Theorem 3.5]{LambertPaquette02} with formula (3.8), 
we obtain with $n = N_T < N_0$:
\begin{equation}  \label{rec2}
\begin{pmatrix} \varphi_{n} \\ \varphi_{n-1} \end{pmatrix}
= \biggl[ \prod_{k=1}^{n-1}  \lambda_+(\tfrac{k}N)\big(1- \delta_k - \boldsymbol{\eta}_k \big) \biggr] \begin{bmatrix}
     \lambda_+(\tfrac{n}N) & \lambda_-(\tfrac{n}N) \\
      1 & 1 \\
    \end{bmatrix} \left( \begin{bmatrix} 
1 &   0 \\   0 &  0
\end{bmatrix} + \boldsymbol{\varepsilon} \right) {\rm U}
\end{equation}
where the deterministic quantities
\begin{equation}\label{deltaeq}
 \lambda_{\pm}(t;z) = \frac{\sqrt{t} J(z/\sqrt{t})^{\mp 1}}{2} , \qquad 
\delta_k(z) = \frac{1}{4(Nz^2-k)}+O\!\left(\frac1{(Nz^2-k)^2}\right) ,
\end{equation}
the initial condition
\begin{equation}\label{Uini}
{\rm U}(z) =  \begin{pmatrix} 1 \\ 0\end{pmatrix} - \frac{b_1}{\sqrt{Nz^2}} \begin{pmatrix} 1 \\ -1\end{pmatrix} + \O\left(\frac{1}{Nz^2}\right)
\end{equation}
and the random variables $\boldsymbol{\eta_k}$ and $\boldsymbol{\epsilon}$ satisfy (for a small $\epsilon>0$)
\begin{equation} \label{def:eta}
\boldsymbol{\eta}_k(z)=\sqrt{\frac{1/2\beta}{Nz^2-k}}\left( X_k +  Y_k J\big(z \sqrt{\tfrac{N}{k}}\big) \right) , \qquad  
\P\big[ \big\{ \|\boldsymbol{\epsilon}(z)\| \ge (\log N_0)^{\epsilon-1/6} \big\} \cap \mathsf{Trunc} \big] \le e^{- c (\log N_0)^{1+\epsilon/2}}.
\end{equation}
In particular, in \cite[Theorem 3.5]{LambertPaquette02} holds with $N \leftarrow N_T(z)$ , $z \leftarrow z_\lambda= z\big(1 +  \frac{\lambda}{2 N_0^{2/3}} \big)$ in which case the bounds hold with $\Omega = (\log N_0(z))^{1/3-3\epsilon/2}$ \big(i.e.~$T=(\Omega \log N_0)^{2/3}= (\log N_0(z))^{1-\epsilon} $\big).
To check the bounds, with $\mathfrak{L} = N_0(z)^{1/3}$ the parabolic length scale, it suffices that 
\begin{equation} \label{N0exp}
N_0(z_\lambda) =  N_0(z) + \lambda \mathfrak{L} + \O(1), \qquad\text{with }|\lambda| \ll T . 
\end{equation}
Here, we restrict ourselves to $\lambda\in\mathcal{K}$ for a fixed compact set $\mathcal{K}$. 
In particular, by \cite[(3.12)]{LambertPaquette02}, the error $\boldsymbol{\epsilon}(z_\lambda)$ satisfies the pointwise estimate \eqref{def:eta} uniformly for $\lambda\in\mathcal{K}$. 

\paragraph{2. Deterministic part.}
The deterministic product in \eqref{rec2} can be expressed in terms of Hermite polynomials (see \cite[(6.2)]{LambertPaquette02}), with $\widetilde\varphi_n= \E\varphi_n$, 
\begin{equation*}
   \widetilde\varphi_{n-1}  =  \prod_{k=2}^{n}  \lambda_+(\tfrac{k-1}N)\big(1- \delta_{k} \big)  \big( 1+ \O\big(n^{-{\rm c}}\big) \big) . 
\end{equation*}
for some ${\rm c}>0$, where the error is controlled uniformly for $\lambda\in\mathcal{K}$. 
According to \eqref{eq:wn}, we multiply \eqref{rec2} by the matrix
\[
\begin{bmatrix} 
w_n &   -w_{n-1} \\   0 &  w_{n-1}
\end{bmatrix}
=  w_{n-1} \begin{bmatrix} 
\frac{w_n}{w_{n-1}} &   -1 \\   0 &  1
\end{bmatrix}, 
\qquad
\frac{w_{n}}{w_{n-1}} = \frac2{\sqrt{n/N}}
\]
so that by \eqref{deltaeq}: 
\[
\begin{bmatrix} 
w_n &   -w_{n-1} \\   0 &  w_{n-1}
\end{bmatrix}
\begin{bmatrix}
     \lambda_+(\tfrac{n}N) & \lambda_-(\tfrac{n}N) \\
      1 & 1 \\
\end{bmatrix} 
=  w_{n-1} \begin{bmatrix} 
J(z_\lambda/\sqrt{n/N})^{-1} -1 &   J(z_\lambda/\sqrt{n/N}) -1 \\   1 &  1
\end{bmatrix}
\]
where (here $z>0$ and $n= N_T$) 
\[
J\big(z_\lambda/\sqrt{n/N}\big)^{\pm 1} = 1 \mp \frac{\sqrt{T+\lambda}}{\mathfrak L}
+\O\bigg(\frac{T}{\mathfrak L^2}\bigg) . 
\]

Then, using the notation \eqref{eq:Psi} and with $\widetilde\Psi_n= \E\Psi_n$, by \eqref{rec2} evaluated at $z_\lambda$,  this implies that with $n= N_T$:
\begin{equation}  \label{rec3}
\begin{pmatrix}\frac{\mathfrak L}{\sqrt{T+\lambda}}\big(\Psi_{n} - \Psi_{n-1}\big)   \\ \Psi_{n-1} \end{pmatrix}
= \frac{\widetilde\Psi_{n-1}}{1+ \O\big(n^{-{\rm c}}\big) }
\prod_{k=1}^{n-1} \frac{1- \delta_k - \boldsymbol{\eta}_k}{1-\delta_k}
\left( \begin{bmatrix} 
1 &  - 1 \\   1 &  1
\end{bmatrix} + \O\bigg(\frac{\sqrt T}{\mathfrak L}\bigg) \right)
 \left( \begin{bmatrix} 
1 &   0 \\   0 &  0
\end{bmatrix} + \boldsymbol{\varepsilon} \right) {\rm U} . 
\end{equation}
 
The asymptotics of the (normalized) Hermite polynomials $\widetilde\Psi_n$ as $n\to\infty$ are given in terms of the Airy function \eqref{eq:Airy}:

We also require the following quantitative extension: the monic Hermite polynomials have the following uniform asymptotics for $|\lambda| \le T$ and $|n-N_p|\leq  T N_p^{1/3}$: 
\begin{equation}  \label{HermiteAiry}
\widetilde\Psi_{n-1}(z_\lambda) = \Ai(\lambda+T)  \big( 1+ \O\big(T/\mathfrak{L}\big) \big) 
\end{equation}
we refer e.g. to the proof of \cite[Lemma 7.2]{AM05} (these asymptotics are uniform for $\lambda\in\mathcal{K}$). 

Finally, by \eqref{Uini}, the initial condition satisfies ($\VERT b_1\VERT_1 \le c$) for $\gamma>0$
\[
\max_{\lambda \in \mathcal{K}}\, \P\big[ \| {\rm U}(z_\lambda) - \left(\begin{smallmatrix}1\\ 0 \end{smallmatrix}\right)\| \ge  \gamma \big] \lesssim e^{-c \gamma \sqrt{ N_0} } .
\]
We conclude that with $n=N_T$, 
\begin{equation}  \label{rec4}
\begin{pmatrix}\frac{\mathfrak L}{\sqrt{T+\lambda}}\big(\Psi_{n} - \Psi_{n-1}\big)   \\ \Psi_{n-1} \end{pmatrix}
= \Ai(\lambda+T)
\prod_{k=1}^{n-1} \frac{1- \delta_k - \boldsymbol{\eta}_k}{1-\delta_k}
\left(   \begin{pmatrix} 1 \\ 1\end{pmatrix} + \boldsymbol{\Delta} \right)  
\end{equation}
where, like \eqref{def:eta}, the error satisfies (for $\epsilon>0$ sufficiently small):
\begin{equation}  \label{deltaer}
\max_{\lambda \in \mathcal{K}}\,\P\big[ \big\{\|\boldsymbol{\Delta}(\lambda)\| \ge (\log N_0)^{\epsilon-1/6}\big\}  \cap \mathsf{Trunc}  \big] \le e^{- c (\log N_0)^{1+\epsilon/2}}.
\end{equation}
Note that the deterministic errors, the error coming from the initial condition and $\P[\mathsf{Trunc}^{\rm c}]$ are negligible. 

\paragraph{3. Random part.}
Let $\boldsymbol{q} = \boldsymbol{q}_T^{(z)}$ be as in \eqref{qrv}, $\mathcal{K} \subset \C$ be a compact set, and $0<\epsilon<1/6$.  
The goal of this section is to prove the following estimates: 
\begin{equation} \label{prodrandest}
\P\Bigg[ \bigg\{\max_{\lambda \in \mathcal{K}} \bigg| 
e^{-\boldsymbol{q} +\frac12 \E\boldsymbol{q}^2} \,  \prod_{k=1}^{n-1} \frac{1- \delta_k - \boldsymbol{\eta}_k}{1-\delta_k}(z_\lambda) -1 \bigg|  \ge (\log N_0)^{\epsilon-1/6}\bigg\}  \cap \mathsf{Trunc}   \Bigg]
\lesssim e^{-(\log N_0)^{1+\epsilon}} . 
\end{equation}

\begin{proof}
Using \eqref{deltaeq} and \eqref{def:eta}, by a linearization, 
\[
\prod_{k=1}^{n-1} \frac{1- \delta_k - \boldsymbol{\eta}_k}{1-\delta_k}
= \exp\left( \sum_{k<n} -\boldsymbol{\eta}_k - \frac{\boldsymbol{\eta}_k^2}{2}+ \O(\delta_k\boldsymbol{\eta}_k)+ \O(\boldsymbol{\eta}_k^3)  \right)
\]
where $\boldsymbol{\eta}_k$ are independent random variables and, according to \eqref{qrv}, $ \boldsymbol{q} = \sum_{k=1}^{n-1} \boldsymbol{\eta}_k(z)$.
On the event $\mathsf{Trunc}$, one has $|\boldsymbol{\eta}_k(z_\lambda)| \lesssim \frac{T^{3/2}}{\sqrt{N_0-k +\O(\mathfrak{L})}} \le \frac{T}{\sqrt{\mathfrak{L}}}$ so that
\begin{equation} \label{prodrand}
\max_{\lambda \in \mathcal{K}} \bigg| 
e^{-\boldsymbol{q} +\frac12 \E\boldsymbol{q}^2} \,  \prod_{k=1}^{n-1} \frac{1- \delta_k - \boldsymbol{\eta}_k}{1-\delta_k}(z_\lambda) -1 \bigg| 
\lesssim \left(   \boldsymbol{\xi} + \boldsymbol{\zeta} + \O\left(\frac{T^C}{\sqrt{\mathfrak{L}}}\right)   \right)
\exp\big(\boldsymbol{\xi} + \boldsymbol{\zeta}\big)
\end{equation}
with a deterministic error (for a constant $C\ge 1$), where 
\[
\boldsymbol{\xi}\coloneqq   \sum_{k<N_T} \max_{\lambda \in \mathcal{K}}  \big|\boldsymbol{\eta}_k(z_\lambda)- \boldsymbol{\eta}_k(z)\big|
, \qquad \boldsymbol{\zeta}\coloneqq \max_{\lambda \in \mathcal{K}} \bigg| \E\boldsymbol{q}^2 -  \sum_{k<N_T} \boldsymbol{\eta}_k^2(z_\lambda) \bigg|.
\]

Observe that $z_\lambda/z= 1+ \O(N_0^{-2/3})$ so that by \eqref{N0exp}, one has
\[
\boldsymbol{\xi}_k \coloneqq 
\max_{\lambda\in\mathcal{K}}\big|\boldsymbol{\eta}_k(z_\lambda) -  \boldsymbol{\eta}_k(z)\big| \lesssim \frac{\big(|X_k|+|Y_k|\big) \mathfrak{L}}{\big(N_0-k +\O(\mathfrak L)\big)^{3/2}} + \frac{|Y_k| N_0^{-1/6}}{N_0-k +\O(\mathfrak L)} 
\]
using the simple bound $|J'(z)| \le \frac1{\sqrt{|z^2-1|}}$. This implies that for $k<N_0$:
\[
\E[\boldsymbol{\xi}_k]\lesssim \frac{\mathfrak{L}}{\big(N_0-k +\O(\mathfrak L)\big)^{3/2}} + \frac{N_0^{-1/6}}{N_0-k +\O(\mathfrak L)} , 
\]
\[
\VERT \boldsymbol{\xi}_k\VERT_{2,\mathsf{Trunc}}^2 \lesssim   
\frac{T^3 N_0^{-1/3}}{\big(N_0-k +\O(\mathfrak L)\big)^2} +
\frac{\mathfrak{L}^2}{\big(N_0-k +\O(\mathfrak L)\big)^3}\1\big\{k > N_0/2\big\} . 
\]
Since $\boldsymbol{\xi}_k$ are independent random variables, by \cite[Proposition 2.6.1]{Vershynin},
\[
\VERT \boldsymbol{\xi} - \E \boldsymbol{\xi} \VERT_{2,\mathsf{Trunc}}^2 \lesssim \sum_{k<N_T} \VERT \boldsymbol{\xi}_k\VERT_{2,\mathsf{Trunc}}^2 \lesssim T^{-2}, \qquad
\E \boldsymbol{\xi} \lesssim \sum_{k<N_T} \E[\boldsymbol{\xi}_k] \lesssim T^{-1/2} . 
\]
Similarly, using that $\E\boldsymbol{q}^2 = \sum_{k=1}^{n-1} \E \boldsymbol{\eta}_k^2(z)$ and 
\[
\boldsymbol{\zeta} \lesssim  \frac{T^{3/2}  \boldsymbol{\xi}}{\sqrt{\mathfrak{L}}}  + |\boldsymbol{\zeta}'| , \qquad \boldsymbol{\zeta}' = \sum_{k<n} \big(\boldsymbol{\eta}_k^2(z) - \E \boldsymbol{\eta}_k^2(z) \big)  ,
\]
where $\boldsymbol{\zeta}'$ is a sum of independent random variables with variance
$\sum_{k=1}^{n-1} \E \big(\boldsymbol{\eta}_k^2 - \E \boldsymbol{\eta}_k^2 \big)^2 \lesssim \sum_{k=1}^{n-1}\frac1{(Nz^2-k)^2} \lesssim \frac{1}{\mathfrak{L}}$ and 
$\big|\boldsymbol{\eta}_k^2 - \E \boldsymbol{\eta}_k^2 \big| \lesssim \frac{T^2}{\mathfrak{L}} $ on $\mathsf{Trunc}$, so that by Freedman's inequality:
\[
\VERT \boldsymbol{\zeta}'\VERT_{1,\mathsf{Trunc}}  \lesssim \frac{T^2}{\mathfrak{L}}
\]

Going back to \eqref{prodrand}, since $\E\boldsymbol{\xi} \lesssim T^{-1/2} \ll (\log N_0)^{-1/3}$ (since $\epsilon<1/6$) this implies that for $ (\log N_0)^{-1/3}<\delta <{\rm c}$
\[\begin{aligned}
\P\Bigg[
\begin{gathered}
\bigg\{ \max_{\lambda \in \mathcal{K}} \bigg|
e^{-\boldsymbol{q} +\frac12 \E\boldsymbol{q}^2} \,  \prod_{k=1}^{n-1} \frac{1- \delta_k - \boldsymbol{\eta}_k}{1-\delta_k}(z_\lambda) -1 \bigg|  \ge 2\delta \bigg\}\\[-0.2em]
{}\cap \mathsf{Trunc}
\end{gathered}
\Bigg]
& \lesssim \P\big[ \{\boldsymbol{\xi} \ge \delta\}  \cap \mathsf{Trunc} \big] \\
&\quad +\P\big[ \{ \boldsymbol{\zeta} \ge \delta \}   \cap \{ \boldsymbol{\xi} \le \delta \}  \cap \mathsf{Trunc} \big] \\
&\lesssim \P\big[ \{\boldsymbol{\xi} - \E\boldsymbol{\xi}  \ge \delta/2 \}  \cap \mathsf{Trunc} \big] + \P\big[ \{  \boldsymbol{\zeta}' \ge \delta/2 \}  \cap \mathsf{Trunc} \big] .
\end{aligned}\]
The second term (involving $\boldsymbol{\zeta}'$) is negligible and we conclude that there is a constant $c>0$ so that 
\[
\P\Bigg[ \bigg\{ \max_{\lambda \in \mathcal{K}} \bigg| 
e^{-\boldsymbol{q} +\frac12 \E\boldsymbol{q}^2} \,  \prod_{k=1}^{n-1} \frac{1- \delta_k - \boldsymbol{\eta}_k}{1-\delta_k}(z_\lambda) -1 \bigg|  \ge 2\delta \bigg\} \cap \mathsf{Trunc} \Bigg]
\lesssim \exp\big(- c T^2  \delta^2\big)
\]
Choosing $2 \delta = (\log N_0)^{\epsilon-1/6}$ this proves the estimate \eqref{prodrandest} \big(recall that $\P[\mathsf{Trunc}^{\rm c}] \lesssim e^{-(\log N_0)^{1+3\epsilon/2}}$\big). 
\end{proof}

\paragraph{4. Uniformity.}
Let $\mathcal{K} \subset \C$ be a compact set. We define $\boldsymbol{\Upsilon}, \boldsymbol{\Upsilon}'$ implicitly by, with $n=N_T(z)$, 
\[
\begin{pmatrix}\frac{\mathfrak L}{\sqrt{T+\lambda}}\big(\Psi_{n} - \Psi_{n-1}\big)(\lambda)  \\ \Psi_{n-1}(\lambda) \end{pmatrix}
= \Ai(\lambda+T) e^{\boldsymbol{q} -\frac12 \E\boldsymbol{q}^2} \left( \begin{pmatrix} 1 \\ 1\end{pmatrix} +
\begin{pmatrix} {\boldsymbol{\Upsilon}'(\lambda)} \\  {\boldsymbol{\Upsilon}(\lambda)} \end{pmatrix} \right) .
\]
Since $T \gg 1$, if $ \mathcal{S}$ is an open neighborhood of $\mathcal{K}$,  $\lambda \in \mathcal{S} \mapsto \boldsymbol{\Upsilon}(\lambda)$ and $\lambda \in \mathcal{S} \mapsto \boldsymbol{\Upsilon}'(\lambda)$ are analytic functions. 
Moreover, by \eqref{rec4}, combining the estimates \eqref{deltaer} and \eqref{prodrandest}, we have 
\begin{equation}  \label{Upser}
\max_{\lambda \in \mathcal{K}}\,\P\big[ \big\{\|\boldsymbol{\Upsilon}(\lambda)\| \ge (\log N_0)^{\epsilon-1/6}\big\}  \cap \mathsf{Trunc}  \big] \le e^{- c (\log N_0)^{1+\epsilon/2}}
\end{equation}
and similarly for $\boldsymbol{\Upsilon}'$. 
The goal of this section is to turn \eqref{Upser} into a uniform bound.

Let $\varphi : \lambda \in \mathcal{S} \mapsto \frac{w_{n-1}(z_\lambda)}{\Ai(\lambda+T)}$ with $w_{n-1}$ as in \eqref{eq:wn}.
We record that 
\[\begin{cases}
\boldsymbol{\Upsilon}(\lambda) = \varphi(\lambda) Q(\lambda) -1 \\
\boldsymbol{\Upsilon}'(\lambda) = \varphi(\lambda) Q'(\lambda)-1 
\end{cases}\qquad
\begin{cases}
Q(\lambda) = \varphi_{n-1}(z_\lambda) e^{-\boldsymbol{q} +\frac12 \E\boldsymbol{q}^2} \\
Q'(\lambda) = \frac{\mathfrak{L}e^{-\boldsymbol{q} +\frac12 \E\boldsymbol{q}^2}}{ \sqrt{\lambda+T}} \big(\tfrac2{\sqrt{n/N}}\varphi_{n}(z_\lambda)- \varphi_{n-1}(z_\lambda)\big)
\end{cases}
\]
where $\lambda\mapsto Q(\lambda)$  and $\lambda\mapsto Q'(\lambda)$  are polynomials of degree $\le n$. 
Then, to show that the errors can be controlled uniformly for $\lambda\in\mathcal{K}$, we rely on the following polynomial interpolation result.

\begin{lemma} \label{lem:polycontrol}
Let $\mathcal{K} \subset \C$ be a compact set and $\widehat{\mathcal{K}} = \mathcal{S} \cap n^{-R} (\Z+\i\Z)$ for $R>1$ where $ \mathcal{S}$ is a neighborhood of $\mathcal{K}$. 
Suppose that $\varphi$ is a function such that for a constant ${\rm C} \ge 1$ and $\sigma\ge 0$, it holds for all $z, w\in  \mathcal{S}$: 
\begin{equation} \label{ratioasymp}
|\varphi(z)| \le {\rm C} |\varphi(w)| \qquad\text{if }|z-w| \le n^{-\sigma} . 
 \end{equation}
 Then, for any polynomial $Q$ of degree $\le n$, if $n$ is sufficiently large:  
\[
\max_{\substack{ u,v \in \mathcal{K}\\ |u-v| \le 2 n^{-R}}} \Big|\varphi(u)\big(Q(u)-Q(v) \big) \Big|  \lesssim {\rm C} (\log n)   n^{1+\sigma-R}  \max_{v \in \widehat{\mathcal{K}}} |\varphi(v)Q(v)| .
\]
\end{lemma}

\begin{proof}
We rely on basic interpolation estimates for polynomials. By Bernstein's inequality, for any $\eta>0$,
\[
  \max_{|z| \leq \eta} |Q'(z)| \leq  \eta^{-1} n\,  \max_{|z| \leq \eta} |Q(z)|.
\]
By the max-modulus principle, we may take the maximum over the square of side-length $2\eta$. 
Then, using the condition \eqref{ratioasymp}, this implies that if $\eta \le n^{-\sigma}/2$,
\[
  \max_{\substack{|u|, |v| \leq \eta \\ |u-v| \le 2n^{-R}}} \Big|\varphi(u)\big(Q(u)-Q(v) \big) \Big|  \leq  2{\rm C} \eta^{-1} n^{1-R}  \max_{z\in \square_\eta} | \varphi(z)Q(z)|.
\]
  Using \cite[Theorem 1]{Rakhmanov}, we may bound the maximum of a polynomial of degree $n$ on a segment by $\O(\log n)$ (for some universal constant) times the maximum of this polynomial over any equally spaced mesh of cardinality $n$, so that with $ \eta = {\rm c} n^{-\sigma}$ for an appropriate constant ${\rm c}$:
 \[
  \max_{\substack{|u|, |v| \leq \eta \\ |u-v| \le 2n^{-R}}} \Big|\varphi(u)\big(Q(u)-Q(v) \big) \Big|  \lesssim {\rm C} (\log n)n^{1+\sigma-R}  \max_{\substack{z\in \square_\eta \\ z\in  n^{-R} (\Z+\i\Z)}} | \varphi(z)Q(z)|.
\]
Hence, by covering $\mathcal{K}$ with balls of radius $\eta$ (here $\eta \ll 1$), we conclude that  
  \[
 \max_{\substack{ u,v \in \mathcal{K}\\ |u-v| \le 2n^{-R}}} \Big|\varphi(u)\big(Q(u)-Q(v) \big) \Big|  \lesssim {\rm C} (\log n) n^{1+\sigma-R}  \max_{z\in \widehat{\mathcal{K}}} |\varphi(z)Q(z)| .  \qedhere
  \]
\end{proof}

We also need the following deterministic bounds.

\begin{lemma}  \label{lem:varphicont}
Let $\varphi : \lambda \in \mathcal{S} \mapsto \frac{w_{n-1}(z_\lambda)}{\Ai(\lambda+T)}$ and $\delta \le \mathfrak L^{-1}$. 
One has for $\lambda, \mu \in \mathcal{S}$: 
\[
\max_{|\lambda-\mu| \le \delta} \bigg| \frac{\varphi(\lambda)}{\varphi(\mu)}-1\bigg| \lesssim \delta \mathfrak L . 
\]
\end{lemma}

\begin{proof}
Using the Airy function asymptotics \eqref{Airyasymp}: for $|\lambda-\mu| \le \delta$, 
\[
\bigg| \frac{\Ai(\lambda+T)}{\Ai(\mu+T)}-1 \bigg| \lesssim \delta \sqrt{T}
\]
According to \eqref{eq:wn} and using \eqref{N0exp}, uniformly for $\lambda,\mu \in \mathcal{S}$ with $|\lambda-\mu| \le \delta \le \mathfrak L^{-1}$,
\[
\frac{w_n(z_\lambda)}{w_n(z_\mu)}
=
\exp\!\big(-\mathfrak L (\lambda-\mu)\big)
\left(1+O(\mathfrak L^{-2})\right)
\]
Then 
\[
\bigg| \frac{w_n(z_\lambda)}{w_n(z_\mu)}  -1 \bigg| \lesssim \delta \mathfrak L . 
\]
This proves the claim. 
\end{proof}

Lemma~\ref{lem:varphicont} yields the condition \eqref{ratioasymp} if $\sigma \ge 1/3$ (since $n<N_0$ and $\mathfrak{L} = N_0^{1/3}$) and it also implies that, with $R>4/3$, for $v\in\mathcal{K}$: 
\begin{equation} \label{ratiocont}
\max_{|\lambda-v| \le 2 n^{-R}} \bigg| \frac{\varphi(\lambda)}{\varphi(v)}-1\bigg| \lesssim N_0^{-1} . 
\end{equation}

Let $\widehat{\mathcal{K}}$ be as in Lemma \ref{lem:polycontrol} with $n \le N_0$. 
By \eqref{Upser} and a union bound (adjusting the constant $c>0$ depending on $R>1$ and $\epsilon>0$): 
\begin{equation}  \label{Upseruni}
\P\bigg[  \bigg\{\max_{v \in \widehat{\mathcal{K}}} \|\boldsymbol{\Upsilon}(v)\| \ge (\log N_0)^{\epsilon-1/6}\bigg\}  \cap \mathsf{Trunc}  \bigg] \le e^{- c (\log N_0)^{1+\epsilon/2}} . 
\end{equation}

Now, for $v \in \widehat{\mathcal{K}}$, one has with $\eta = 2n^{-R}$, 
\[\begin{aligned}
\max_{|\lambda-v| \le \eta} \big|  \varphi(\lambda) Q(\lambda)-1 \big|
&\le  \max_{|\lambda-v| \le \eta} \big|  \varphi(\lambda) Q(v)-1 \big| +\max_{|\lambda-v| \le \eta} \big| \varphi(\lambda) \big(Q(\lambda)-Q(v)\big) \big|\\
&\le  \big|  \varphi(v) Q(v)-1 \big|
+ \max_{|\lambda-v| \le \eta} \bigg| \frac{\varphi(\lambda)}{\varphi(v)}-1\bigg| \cdot \big|  \varphi(v) Q(v) \big|
+   \max_{|\lambda-v| \le \eta} \big| \varphi(\lambda) \big(Q(\lambda)-Q(v)\big) \big| 
\end{aligned}\]
Then, by Lemma~\ref{lem:polycontrol} with $\sigma=1/3$, on the event $\big\{ \max_{v \in \widehat{\mathcal{K}}} |\boldsymbol{\Upsilon}(v)|  \le \delta \big\}$ with $0<\delta\le 1$, we obtain (by taking the max over $v \in \widehat{\mathcal{K}}$) for $R>4/3$:
\[
\max_{\lambda\in\mathcal{K}} \big|\underbrace{\varphi(\lambda) Q(\lambda)-1}_{\boldsymbol{\Upsilon}(\lambda)} \big|
\le \delta + 2 \max_{\substack{|\lambda-v| \le 2n^{-R} \\ v\in \widehat{\mathcal{K}}}} \bigg| \frac{\varphi(\lambda)}{\varphi(v)}-1\bigg| +\O\left( (\log n)   n^{4/3-R} \right)
\]
By \eqref{ratiocont} and \eqref{Upseruni}, choosing $\delta = (\log N_0)^{\epsilon-1/6}$,  we conclude that if $N_0 \gg 1$ (here $n =N_T$)
\[
\P\bigg[ \bigg\{\max_{\lambda\in\mathcal{K}} \big| \boldsymbol{\Upsilon}(\lambda)\big| \ge 2\delta \bigg\}  \cap \mathsf{Trunc} \bigg] \le  
\P\bigg[ \bigg\{\max_{v\in\widehat{\mathcal{K}}} \big| \boldsymbol{\Upsilon}(v)\big| \ge \delta\bigg\} \cap \mathsf{Trunc}  \bigg] \le e^{- c (\log N_0)^{1+\epsilon/2}} . 
\]
This concludes the proof. 
\end{proof}

Using the strong embedding from \cite[Appendix~C]{LambertPaquette02}, let $(\widehat{\mathbf{X}},\widehat{\mathbf{Y}})$ denote the independent Brownian motions coupled to the noise sequences $(X_k,Y_k)$.  We rescale these Brownian motions by defining
\[
({\mathbf{X}}_t : t \in [0,1]) = (N^{-1/2}\widehat{\mathbf{X}}_{tN} : t \in [0,1]) , \qquad 
({\mathbf{Y}}_t : t\in[0,1]) = (N^{-1/2}\widehat{\mathbf{Y}}_{tN} : t \in [0,1]) .
\]
Then, $(\mathbf{X},\mathbf{Y})$ remain standard independent Brownian motions and, as in \cite[(1.9)]{LambertPaquette02} we define a Gaussian process for $t \in [0,1)$ and $z \in \C{\setminus[-1,1]}$: 
\begin{equation} \label{gfield}
  \boldsymbol{W}_t(z) \coloneqq -\frac{1}{2}\int_0^t \frac{\d \mathbf{X}_u+ \mathrm{J}(z/\sqrt{u}) \d\mathbf{Y}_u}{\sqrt{u}\sqrt{z^2/u-1}},
  \qquad
  \qquad
  \text{where $\mathrm{J}(z) \coloneqq z- \sqrt{z^2-1}$ as in \eqref{eq:Wtz}}
\end{equation}
and $z \in \C\setminus [-1,1] \mapsto \sqrt{z^2/u-1}$ is the same $\sqrt{\cdot}$ as in $\mathrm{J}$.

\begin{proposition} \label{prop:Gauss}
Let $\gamma = \sqrt{2/\beta}$.
Let $z\in[-1,1]$ with $N_0(z) \gg 1$ and $T\ge 1$. 
According to \eqref{qrv} with  $t = \tfrac{N_T-1}{N_0}$: 
\[
\VERT \boldsymbol{q}_T^{(z)} -  \gamma\boldsymbol{W}_t(z)  \VERT_1 \lesssim  \frac{\log N_0}{\sqrt{T\mathfrak{L}}}. 
\]
In particular, the random variable $\boldsymbol{q}_T^{(z)}$ is approximately Gaussian with mean-zero and variance $\frac2\beta \log\big(\mathfrak{L}/\sqrt{4T}\big)$. 
\end{proposition} 

\begin{proof}
Observe that with  $t = \tfrac{N_T-1}{N_0}$, we can rewrite
\[
\boldsymbol{q} = 
\frac{1}{\sqrt{2\beta N_0}} \sum_{k=1}^{N_0t} \big( f_1(\tfrac{k}{N_0}) X_k +  f_2(\tfrac{k}{N_0}) Y_k \big)
\]
where $f_1 : u \in [0,1] \mapsto \frac{1}{\sqrt{1-u}}$, $f_2 : u \in [0,1] \mapsto  \frac{J(1/\sqrt{u})}{\sqrt{1-u}}$. Then, by \cite[Theorem C.1]{LambertPaquette02} and \cite[Proposition C.2]{LambertPaquette02}, 
\[
\VERT \boldsymbol{q}_T - \gamma \boldsymbol{W}_t  \VERT_1 \lesssim \frac{\log N_0}{\sqrt{N_0}}
\max_{j=1,2}\|f_j\|_{\operatorname{TV},t} , \qquad \|f\|_{\operatorname{TV},t}
    =  |f(t)| + \int_{0}^t |f'(s)|\,ds . 
\]
Since both $f_j$ are increasing on $[0,1]$, $\|f\|_{\operatorname{TV},t} \le 2 f_j(t) \le 
\frac{2}{\sqrt{1-t}}$ for $t\in[0,1)$, we conclude that 
\[
\VERT \boldsymbol{q}_T -  \gamma\boldsymbol{W}_t  \VERT_1 \lesssim \frac{\log N_0}{\sqrt{N_0(1-t)}} \le \frac{\log N_0}{\sqrt{N_0-N_T}} = \frac{\log N_0}{\sqrt{T\mathfrak{L}}}.
\]
This proves the first claim. Then by \cite[(1.11)]{LambertPaquette02}:
\[
\E[ \boldsymbol{W}_t^2] = - \log\big(1 - J(1/\sqrt t)^2 \big) = \log\big(\mathfrak{L}/\sqrt{4T}\big) + \O\big(\sqrt{T}/\mathfrak{L}+(T\mathfrak{L})^{-1}\big)
\]
using that \(1-J(1/\sqrt{u}) \simeq \sqrt{1-u}\) as $u\to1$ and $1-t \simeq T \mathfrak{L}^{-2}$.  
\end{proof}

\begin{remark} \label{rk:sym}
Let $z = \pm 1$ and $t = \tfrac{N_T-1}{N}$ with $T = (\log N)^{1-\epsilon}$. 
We claim that the random variables $\boldsymbol{W}_t(\pm 1)$ from \eqref{gfield} are not independent. These variables are mean-zero Gaussian, and their covariance is computed below.
Observe that $J(-x)= -J(x)$ for $x \in \R\setminus[-1,1]$. Then, as stochastic integrals driven by two independent Brownian motions $(\widehat{\mathbf{X}},\widehat{\mathbf{Y}})$, we have 
\[
\E\big[\boldsymbol{W}_t(+1) \boldsymbol{W}_t(-1)\big] 
= \frac14 \bigg( \int_0^t \frac{- \d u}{1-u} + \int_0^t \frac{ \mathrm{J}(1/\sqrt{u})^2 \d u}{1-u} \bigg) = \frac12 \int_0^t \bigg(1-\frac{1}{\sqrt{1-u}} \bigg) \frac{\d u}{u}
\] 
using that $\frac{\mathrm{J}(z)^2 -1}{2} = z^2-1- z\sqrt{z^2-1}$ so that $\frac{\mathrm{J}(1/\sqrt{u})^2 -1}2 = \frac{1-u - \sqrt{1-u}}{u}$ for $u\in[0,1]$. Using that $t\simeq 1$ and substituting $v= \sqrt{1-u}$ in the above integral, we obtain 
\[
\E\big[\boldsymbol{W}_t(+1) \boldsymbol{W}_t(-1)\big] 
= - \int_{\sqrt{1-t}}^1 \frac{1-v}{1-v^2} \d v
\simeq - \big[ \log(1+v) \big]_{0}^1  = -\log 2. 
\] 
\end{remark}

\begin{proposition} \label{prop:qnoise}
Consider the random variable $\boldsymbol{q}_S= \boldsymbol{q}_S^{(z)}$ for $S\ge 1$; \eqref{qrv}.
Let $\mathfrak{X}(u) =\int_0^u e^{\frac{4}{3}( t^{3/2} - u^{3/2})} \d \boldsymbol{B}(t)$ for $u\ge 0$ as in  \eqref{eq:SAIG} where $\boldsymbol{B} = 
\boldsymbol{B}^{(z)}$ is the Brownian motion \eqref{Bcoup}. 
For $0<\alpha<1/6$, there is a constant $c= c(\beta,\alpha)$  such that for any $S\in [1,T]$ with $T \le (\log N_0)^2$ and $0<\eta< N_0^\alpha$:  
\[
\P\bigg[
\bigg| \boldsymbol{q}_T + \tfrac{2}{\sqrt{\beta}}\int_S^T \mathfrak{X}(u)\d u - \boldsymbol{q}_S \bigg| \ge \eta \bigg]
\lesssim \exp\big( - c S^{3/2} \eta^2 \big) . 
\]

\end{proposition}
\begin{proof}
On the one hand, using that 
$J\big(z \sqrt{\tfrac{N}{N_0- k}}\big) = 1- \O\big(\sqrt{k}\mathfrak{L}^{-1}\big)$ for $k\ge 1$, 
 by \eqref{qrv}, one has for $0<S<T$:
\[
\boldsymbol{q}_S = \boldsymbol{q}_T + \sum_{k=N_T}^{N_S-1} \frac1{\sqrt{N_0-k}}\frac{X_k +  Y_k J\big(z \sqrt{\tfrac{N}{k}}\big)}{\sqrt{2\beta}} 
=  \boldsymbol{q}_T +  \sum_{k = \lfloor S\mathfrak{L} \rfloor+1}^{\lfloor T \mathfrak{L} \rfloor} \frac{X_{N_0-k}+ Y_{N_0-k}}{\sqrt{2\beta k }}
+ \O\bigg(\sqrt{\tfrac{T}{\mathfrak L}}\bigg)
\]
where the error is controlled in sub-exponential norm (under the Assumptions~\ref{def:noise}). 

On the other hand, applying stochastic Fubini to \eqref{eq:SAIG}, one has
\begin{equation}
\int_S^T \mathfrak{X}(u) \d u =
\int_0^S \d B(t) \bigg(
  \int_S^T e^{\frac{4}{3}( t^{3/2} - u^{3/2})} \d u \bigg) +{}
\int_S^T \d B(t) \bigg(
  \int_t^T e^{\frac{4}{3}( t^{3/2} - u^{3/2})} \d u \bigg) = \int_S^T \frac{\d B(t)}{2\sqrt t}
+ \O(S^{-3/4}) 
\end{equation}
where the error is controlled in sub-Gaussian norm. Indeed, the missing first integral equals $a_{S,T}\mathfrak X(S)$, where
\begin{align*}
a_{S,T}
&:=\int_S^T e^{-\frac43(u^{3/2}-S^{3/2})}\,\d u
\lesssim S^{-1/2},\qquad
\|\mathfrak X(S)\|_{\psi_2}
\lesssim
\bigg(\int_0^S e^{\frac83(t^{3/2}-S^{3/2})}\,\d t\bigg)^{1/2}
\lesssim S^{-1/4}.
\end{align*}
Thus $\|a_{S,T}\mathfrak X(S)\|_{\psi_2}\lesssim S^{-3/4}$.
Then, by the coupling \eqref{eq:B}, this implies that 
\[
\frac{2}{\sqrt{\beta}}\int_S^T \frac{\d B(t)}{2\sqrt t}
= -\sqrt{\frac2{\beta\mathfrak L}}
\sum_{k = \lfloor S\mathfrak{L} \rfloor+1}^{\lfloor T \mathfrak{L} \rfloor} \frac{X_{N_0-k}+ Y_{N_0-k}}{2\sqrt{k/\mathfrak{L}}} + \O\bigg( \frac{\log \mathfrak{L}}{\sqrt{ \mathfrak{L}}} \bigg) 
\]
where the error is controlled in  sub-exponential norm.  

Combining these sub-exponential / sub-Gaussian estimates, we conclude that  for $1 \le S \le T \le (\log N_0)^2$, it holds for any $\eta>0$: 
\[
\P\bigg[
\bigg| \boldsymbol{q}_T + \tfrac{2}{\sqrt{\beta}}\int_S^T \mathfrak{X}(u)\d u - \boldsymbol{q}_S \bigg| \ge \eta \bigg]
\lesssim \exp\Bigg( - \frac{c \eta^2}{ S^{-3/2} + \eta  (\log N_0) N_0^{-1/6}} \bigg) . 
\]
In the range of the proposition, $\eta<N_0^\alpha$ with $\alpha<1/6$ and $S\leq (\log N_0)^2$, the second term in the denominator is at most a constant times $S^{-3/2}$ for all sufficiently large $N_0$, which gives the stated tail bound.
This concludes the proof. 
\end{proof}

% ===== end sections/A-uniformity-hyperbolic.tex =====

% ===== begin sections/B-brownian-motion.tex =====
\section{Estimates for Brownian motion} \label{sec:bm}

Let us recall that $B$ is a standard Brownian motion, so that
$\E[B(t)^2] =  t$ for $t\ge 0$.

Let us consider the following functions
\[ g(u) =  \begin{cases} 
f(u)  &\text{if}  f(u)=u^\alpha \text{ for a }\alpha>0 \\
\frac{f(u)}{(\log(1+\log(1+ u)))^{3/2}} &\text{if}  f(u)= \log(1+u)^\alpha \text{ for a }\alpha>0 
\end{cases} . \]

\begin{lemma}\label{lem:BMestimate}
    There is a constant $C=C_{\alpha}$ so that for any $T> 0$, 
    \[
 \Big\VERT \sup_{t \geq T} \tfrac{B(t)}{\sqrt{t} f(t)} \Big\VERT_2 \leq \tfrac{C}{g(T)} \quad
      \text{and}
      \quad
      \Big\VERT \sup_{t \leq 1/T} \tfrac{B(t)}{\sqrt{t}f(1/t)} \Big\VERT_2 \leq \tfrac{C}{g(T)} .
    \]
\end{lemma}
%\noindent For our purposes, the most important examples of functions $f$ are $u \mapsto \log^{\alpha}(u)$ for $\alpha \geq \frac{1}{3}$ and $u \mapsto u^{\epsilon}$ for $\epsilon >0.$  
\begin{proof}
  The second statement follows from the first by time inversion symmetry of Brownian motion (i.e.\ $t\mapsto tB(1/t)$ is again a Brownian motion with the same variance).  Hence we just need to control the sub-Gaussian norm of the random variable $X_T =  \sup_{t \geq T} \big| \tfrac{B(t)}{\sqrt{t} f(t)} \big|$.  
 From the Borell--TIS inequality,
\begin{equation}\label{BorellTIs}
	     \VERT X_T -\E X_T \VERT_2 \lesssim   \sqrt{  \sup_{t \geq T}  \tfrac{\E B(t)^2}{t f(t)^2} } \lesssim \frac{1}{f(T)}
 \end{equation}
 since $f$ is an increasing function. 
 Hence, it suffices to bound $\E X_T$. 
 The proof relies on the fact that $W(t) := B(Te^t)/\sqrt{Te^t}$ is an Ornstein--Uhlenbeck process and its maximum on an interval of length $\ell$ has mean at most $C\sqrt{\log(1+\ell)}$.\footnote{ To see this: cut the time interval into $\ell$ intervals of length $1.$  On an interval of length $1,$ \eqref{eq:Dudley2} gives the sub-Gaussian norm is bounded by a constant.  Now the max of $\ell$ variables of uniformly bounded sub-Gaussian norm has expectation bounded by $\sqrt{\log \ell}$ with at most constant order sub-Gaussian fluctuations.}
 Partitioning $[0,\infty)$ into intervals of length $2^{k},$ we obtain
  \[
    \Exp \sup_{t \geq T} \Big| \frac{B(t)}{\sqrt{t}f(t)} \Big|
    =\Exp \sup_{t \geq 0} \Big| \frac{W(t)}{f(Te^t)} \Big|
    \lesssim \sum_{k=0}^\infty \frac{\sqrt{k+1}}{f(Te^{2^k})}
  \]
where we used that the function $f$ is increasing. 
In case $f(u) = u^\alpha$, this immediately implies that 
$ \E X_T \lesssim_\alpha \frac{1}{f(T)}$. 
In case $f(u) = \log^\alpha(u)$, the main contribution to the sum comes from $k: 2^k \leq 1+\log(1+ T),$ and we obtain 
$ \E X_T \lesssim_\alpha \frac{(\log(1+\log(1+T)))^{3/2}}{f(T)}$. 
By \eqref{BorellTIs} and the triangle inequality, this completes the proof. 
 \end{proof}

 We also need a bound for certain quadratic functionals of Brownian motion:
 \begin{lemma} \label{lem:S}
Consider $\mathbf{S} = \int_T^\infty L(u)^2 B(u)^2 \frac{\d u}{u^2}$ where $L$ is a positive, eventually non-increasing slowly varying function and $T$ is sufficiently large that $L$ is non-increasing on $[T,\infty)$, so that
\[
m_T =  \E  \mathbf{S} = \int_T^\infty \frac{L(u)^2}{u}\d u < \infty .
\]
Then there exists a small constant $c>0$ depending on $L$ so that 
\[
\P[\mathbf{S} \ge 2 m_T] \le \exp\bigg(- c\frac{m_T}{L(T)^2} \bigg) 
\] 
\end{lemma}

\begin{proof}
It holds for any $0<\lambda < \|\mathcal{K}\|^{-1}$,
\begin{equation}\label{eq:legaultfredholm}
  \Exp( e^{\lambda \mathbf{S}/2}) = \frac{1}{\sqrt{ \det( \Id - \lambda \mathcal{K})}}
  =   \exp\biggl(
  \sum_{\ell\in\N}\frac{\lambda^\ell}{2\ell}
  \tr [\mathcal{K}^{\ell}]
  \biggr),
\end{equation}
where  $\mathcal{K}$ is an integral operator acting on $L^2([T,\infty))$ with kernel $\mathcal{K}(u,v) = \frac{L(u)}{u} \min\{u,v\}\frac{L(v)}{v}$. Observe that if $q(v) = v^{-1/2}$, it holds for $u \ge  T $, 
\[
\mathcal{K} q (u) =  \int_T^\infty \mathcal{K}(u,v) v^{-1/2} \d v
=  \frac{L(u)}{u} \bigg(\int_T^u L(v) v^{-1/2} \d v + u \int_{u}^\infty  L(v) \frac{\d v}{v^{3/2}}  \bigg) 
\le C  L(u)^2 q(u)
\]
The inequality uses the standard integral estimates for slowly varying functions. Since $L(u)\le L(T)$ for $u\ge T$, by Schur's test, this implies that $\|\mathcal{K}\| \lesssim L(T)^2 $. 
Hence, using that $\E  \mathbf{S}  = \tr \mathcal{K}$ and $  \tr [\mathcal{K}^{\ell+1}] \le \|\mathcal{K}\|^{\ell}\tr \mathcal{K}  $ for any $\ell \in\N$, this implies that if $0<\lambda \le c L(T)^{-2} $ for a sufficiently small constant $c$, then
\[
 \Exp\big[ e^{\lambda( \mathbf{S}-\E  \mathbf{S})/2}\big] \le \exp\bigg(\frac{\lambda \E  \mathbf{S}}{4}\sum_{\ell\in\N}\lambda^\ell \|\mathcal{K}\|^{\ell} \bigg)
 \le  \exp\bigg(C \lambda^2 m_T L(T)^2\bigg) 
\]
Hence, by Markov's inequality choosing $\lambda = 4 c L(T)^{-2} $ for a possibly smaller $c$, 
we conclude that
\[
\P[\mathbf{S} \ge 2 m_T] 
\le  \Exp\big[ e^{\lambda( \mathbf{S}-\E  \mathbf{S})/2}\big] e^{-m_T \lambda/2}
\le \exp\bigg(- c m_T L(T)^{-2} \bigg) . \qedhere
\]
\end{proof}

% ===== end sections/B-brownian-motion.tex =====

% ===== begin sections/C-log-derivative.tex =====
\section{Logarithmic derivative} \label{sec:ld}

In order to relate a solution of the stochastic Airy equation to its logarithmic derivative, we rely on the following deterministic lemma. 

\begin{lemma}\label{lem:pv}
Suppose that $\phi \in C^1(\R_+\to\R)$ has only simple isolated zeros $\{\mathfrak{z}_k\}$ and that  $\phi(0) >0$. 
Let $\rho = \phi'/\phi$  and $\mathbf{N}_{t} = \#\{  k\in\N:  \mathfrak{z}_k \in (0,t) \}$ for $t\in\R_+$.  
Then the function $\rho: \R_+\to \overline{\R}$, the one point compactification of $\R$, is well defined, continuous and 
\begin{equation} \label{exppv}
\phi(t) = \phi(0) \exp\biggl( \operatorname{pv} \int_{0}^t \rho(s)\d s+  \i \pi \mathbf{N}_{t} \biggr)   
\end{equation}
where 
$\displaystyle
\operatorname{pv} \int_{0}^t \rho(s)\d s 
:=  \lim_{\epsilon\to0} \int_{S_{t,\epsilon}}   \rho(s)\d s $ and $S_{t,\epsilon} : =  [0,t] \setminus \bigcup_{k\le \mathbf{N}_{t}}  [\mathfrak{z}_k-\epsilon , \mathfrak{z}_k+\epsilon]$. 
\end{lemma}
\begin{proof}
By assumption, the zero set $ \mathfrak{N} = \{\mathfrak{z}_k\}$  is locally finite (hence countable and ordered), and we verify that 
$\lim_{\epsilon \to 0 } \rho(\mathfrak{z} \mp \epsilon ) = \mp \infty$ for all $\mathfrak{z} \in \mathfrak{N} $. This shows that we can turn  $\rho: \R_+\to \overline{\R}$ into a continuous function such that  $\mathfrak{N} = \big\{  u \in \R_+: \rho(u)=\infty  \big\}$. 
For any $\eta>0$, we have
\[
\int_0^t \frac{\phi'(u)}{\phi(u) +\i \eta} \d u = \log\big(\phi(t) + \i \eta\big) - \log\big(\phi(0) + \i \eta\big)  
\]
for the principal branch of $\log(\cdot)$. 
Upon taking the limit as $\eta\to0$ and using the continuity of $\exp(\cdot)$, this shows that 
\[
\phi(t) = \phi(0)  \exp\Bigg( \lim_{\eta\to0} \int_0^t \frac{\phi'(u)}{\phi(u) +\i \eta} \d u  \Bigg) . 
\]
In particular, the above limit exists. 
Suppose that $t\notin \mathfrak{N}$ and choose $\epsilon>0$ sufficiently small so that
\[
 \int_0^t \frac{\phi'(u)}{\phi(u) +\i \eta} \d u
 = \int_{S_{t,\epsilon}} \frac{\phi'(u)}{\phi(u) +\i \eta} \d u
 +  \sum_{k\le \mathbf{N}_{t}} \int_{\mathfrak{z}_k-\epsilon}^{\mathfrak{z}_{k}+\epsilon}
  \frac{\phi'(u)}{\phi(u) +\i \eta} \d u .
 \]
As above, we verify that 
\[
\lim_{\eta\to0}\int_{S_{t,\epsilon}} \frac{\phi'(u)}{\phi(u) +\i \eta} \d u
=\int_{S_{t,\epsilon}} \rho(u) \d u
\qquad\text{and}\qquad
\lim_{\eta\to0} 
\int_{\mathfrak{z}_k-\epsilon}^{\mathfrak{z}_{k}+\epsilon}  \frac{\phi'(u)}{\phi(u) +\i \eta} \d u
  = \log\bigg|\frac{\phi(\mathfrak{z}_{k}+\epsilon)}{\phi(\mathfrak{z}_{k}-\epsilon)} \bigg| \mp \i \pi
\] 
where $\pm = \operatorname{sgn}\big(\phi'(\mathfrak{z}_{k})\big)$. 
These elementary computations show that for all $\epsilon>0$ sufficiently small, 
\[
 \exp\Bigg( \lim_{\eta\to0} \int_0^t \frac{\phi'(u)}{\phi(u) +\i \eta} \d u  \Bigg) 
=   \exp\Bigg(\int_{S_{t,\epsilon}} \rho(u) \d u +  \sum_{k\le \mathbf{N}_{t}}   \log\bigg|\frac{\phi(\mathfrak{z}_{k}+\epsilon)}{\phi(\mathfrak{z}_{k}-\epsilon)}\bigg| +  \i \pi \mathbf{N}_{t}   \Bigg)  
\]
Note that by a Taylor expansion and using that $\phi'(\mathfrak{z}_{k}) \neq 0$, we have for any $k$, 
$\displaystyle   \log\bigg|\frac{\phi(\mathfrak{z}_{k}+\epsilon)}{\phi(\mathfrak{z}_{k}-\epsilon)}\bigg| = \log|1+ \o(1)| $ as $\epsilon\to0$. 
This implies that $\displaystyle
\operatorname{pv} \int_{0}^t \rho(s)\d s 
:=  \lim_{\epsilon\to0} \int_{S_{t,\epsilon}}   \rho(s)\d s $
exists in $\R$ and that
\[
 \exp\Bigg( \lim_{\eta\to0} \int_0^t \frac{\phi'(u)}{\phi(u) +\i \eta} \d u  \Bigg) =
 \exp\biggl( \operatorname{pv} \int_{0}^t \rho(s)\d s+  \i \pi \mathbf{N}_{t} \biggr)  .
 \]
This completes the proof of \eqref{exppv} in case  $t\notin \mathfrak{N}$. On the other hand, if $t\in\mathfrak{N}$, then we have seen that $\displaystyle \lim_{u\to t_-} \rho(u) =-\infty$ so that by definition $\displaystyle \operatorname{pv} \int_{0}^t \rho(s)\d s =-\infty$. 
Hence, formula  \eqref{exppv} holds in this case as well. 
\end{proof}

% ===== end sections/C-log-derivative.tex =====

\printbibliography[heading=bibliography]%,title={References}]

\end{document}